**Oscar Sheynin**

**Theory of Probability. A Historical Essay**

**Revised and Enlarged Edition**

Berlin
2017

**Contents**

**Annotation**

This book covers the history of probability up to Kolmogorov with essential additional coverage of statistics up to Fisher. Based on my work of ca. 50 years, it is the only suchlike book. Gorrochurn (2016) is similar but his study of events preceding Laplace is absolutely unsatisfactory. Hald (1990; 1998) are worthy indeed but the Continental direction of statistics (Russian and German statisticians) is omitted, it is impossible to find out what was contained in any particular memoir of Laplace and the explanation does not always explain the path from, say, Poisson to a modern interpretation of his results. Finally, the reader ought to master modern math. statistics.

I included many barely known facts and conclusions, e. g., Gauss' justification of least squares (yes!), the merits of Bayes (again, yes!), the unforgivable mistake of Laplace, the work of Chebyshev and his students (merits and failures) etc., etc.

The book covers an extremely wide field, and is targeted at the same readers as any other book on history of science. Mathematical treatment is not as difficult as it is for readers of Hald.

# Preface

> To disregard bygone years and only
> kow-tow to the present is mean,
> barbarous and ignorant
> A. S. Pushkin

> I do feel how wrongful it was to work for so
> many years at statistics and neglect its history
> K. Pearson (1978, p. 1)

I have attempted to include everything essentially interesting. A historical essay such as this one can (and hopefully will) be also used for methodological purposes, so that the narrative should not be dry, but it is not for me to determine whether my book is now suitable enough for a comparatively broad circle of readers.

The book is intended for those interested in the history of mathematics or statistics and more or less acquainted with the latter. It will also be useful for statisticians. My exposition is based, in the first place, on my own investigations published over some 50 years. True, I am not satisfied with a few of them anymore. Note also, that I was unable to check the proofs of some of my papers which are therefore corrupted by misprints. I bear in mind a manuscript that was smuggled out of the Soviet Union (Sheynin 1989a) as well as my Russian articles in the *Istoriko-Matematicheskie Issledovania* from 1993 onward whose proofs I never saw. A good few years ago this journal, the only Russian outlet for papers on history of mathematics, was killed by the transformers of the Academy of Sciences.

I describe the origin of the notions of randomness and subjective or logical probability in antiquity, discuss how laymen comprehended the main notions of the theory of probability, dwell on the birth of political arithmetic and study the history of the theory of probability proper. I also trace the development of statistics and its penetration into natural sciences as well as the history of the mathematical treatment of observations (**Ptolemy**, **Al-Biruni**, **Kepler**, the classical error theory). I stop at the axiomatization of probability and at the birth of the real mathematical statistics, i.e., at **Kolmogorov** and **Fisher**.

From adjoining general sources[1] written from a modern point of view, I mention Stigler (1986), **Hald** (1990; 1998) and Farebrother (1999). The first of these, in spite of its title, only dwells on separate chapters of the history of statistics and is utterly corrupted by slandering the memory of **Euler** and **Gauss**. And since it was universally hailed with **no one** defending those giants, it showed that the scientific community was (and is) seriously ill.

The next two books are specimens of an exposition of a mathematical subject, but they are intended for really qualified readers; then, some topics in **Hald** (1998), especially the description of the work of Russian mathematicians, are omitted and the exposition is not always connected with specific contributions; thus, it is difficult to see what had Laplace to say in any of his early memoirs. Finally,

Farebrother's book dwells on the treatment of observations. My own booklet (2006/2009) is a collection of 832 short dictums, – of pronouncements made by great many scholars over the centuries on probability, statistics, theory of errors and randomness. I see it as a very important scientific and methodological supplement to "ordinary" treatises.

During the last years, quite a few worthless or mediocre contributions to my subject have appeared which was apparently made possible by unsatisfactory preliminary reviewing (and then justified by subsequent superficial abstracting). I do not mention such literature and I also note that in 1915 the Petersburg Academy of Sciences awarded a gold medal to **Chuprov** for a review written on its request (Sheynin 1990c/2011, p. 50). Then, I quote **Truesdell** (1984, p. 292):

*By definition, now, there is no learning, because truth is dismissed as an old-fashioned superstition. Instead* […] *there is perpetual 'research' on anything and everything. In virtue of the Parkinson's law, the professional historian must keep on publishing. Whiteside's monument to Newton, like Wren's masterpiece for St. Paul, will soon be hidden by towering concrete hives of new bureaus and new slums.*

The general situation is bad. For example, it is difficult to publish an honest critical review of a new book. Periodicals receive free copies of them for reviewing from the publishers, and Editors are therefore obliged only to approve sweet nothings. As an example, I advise readers to compare my reviews of two books (Sheynin 2006c; 2006d) with any other of their published reviews[2]. See also § 0.5.

With sincere gratitude I recall the late Professors **Youshkevitch**, who was always favourably disposed towards me, and **Truesdell**, the Editor of the *Archive for History of Exact Sciences*, who had to busy himself with my English and compelled me to pay due attention to style. In 1991, after moving to Germany, I became able to continue my work largely because of Professor **Pfanzagl's** warm support. In particular, he secured a grant for me (which regrettably dried up long ago) from Axel-Springer Verlag. In my papers, I had acknowledged the help of many colleagues including the late Doctors Chirikov (an able mathematician whose bad health thwarted his scientific career) and Eisenhart.

The reader should bear in mind that even **Markov** did not always distinguish between strict and non-strict inequalities. A second similar fact is that the distinction between a magnitude sought (for example, an unknown constant) and its statistical estimate had not been explicitly indicated until perhaps the end of the 19$^{th}$ century (and still later in the Biometric school). Then, the expression such as $P(x = m)$, used for example by **Laplace**, should be understood as $P(m < x < m + dm)$ with one or both inequalities being possibly non-strict.

I am using the following abbreviations: CLT – central limit theorem; LLN – law of large numbers; and MLSq – method of least squares. W-i means Gauss, *Werke*, Bd. i (reprint: Hildesheim, 1973 – 1981) and W/Erg-i is Gauss, *Werke*, Ergänzungsreihe, Bd. i, the reprint of the correspondence of Gauss (Hildesheim, 1975 – 1987).

Bände 1, 3, 4 and 5 are his correspondences with Bessel, Gerling, Olbers and Schumacher respectively. Finally notation **S, G**, i means that the source in question is available as a downloadable file i on my website [www.sheynin.de](www.sheynin.de) which is being copied by Google (Google, Oscar Sheynin, Home). I attach this notation to a source if I have provided its English translation from other languages or if that source is rare.

When describing the contributions of previous years and centuries I sometimes use modern terms but indicate them in square brackets. Thus, [probability] implies that the appropriate author had not applied that expression.

I have gathered the Notes at the end of the chapter in question. I am mentioning many Russian sources, some of them, translated by myself into English. A double date provided in a reference indicates the dates of both the original and the later edition to which I refer.

I have managed to publish abroad a substantial number of papers while still living under a dictatorial regime in Moscow, and the difficulties which I had to overcome have been unknown to the external world. Without Youshkevitch they would have been insurmountable and in any case in addition to English texts (even of the reviews for *Zentralblatt MATH*) I had to submit their Russian versions for the censors. Add to this all the humiliation meted out to a Jew, and you will begin to understand my former situation. I described it in my Russian autobiography, see my site.

*Acknowledgement*. Dr. Valerii Salov had encouraged me to publish this contribution and essentially helped me to submit it.

**Notes**

**1.** Since I also dwell on population statistics, I ought to mention J. & M. Dupâquier (1985). Among other issues, they describe the history of national and international societies and institutions.

**2.** Here is a typical case. In 1998 Desrosières stated that **Poisson** had formulated the strong LLN, and **Gauss** had derived the normal distribution as the limit of the binomial (see my review in *Isis*, vol. 92, 2001, pp. 184 – 185) whereas Stigler (1999, p. 52) called him a first-class scholar.

# 0. Introduction
## 0.1. The Stages

**Kolmogorov** (1947, p. 54) "tentatively" separated the history of probability into four stages: the creation of its "elements" (from **Pascal** and **Fermat** to **Jakob Bernoulli**); the 18th, and the commencement of the 19th century (from **De Moivre** to **Poisson**); the second half of the 19th century (**Chebyshev, Markov; Liapunov** and the origin of mathematical statistics); and the beginning of the 20th century. **Gnedenko** (1958) and **Prokhorov** & **Sevastianov** (1999) offered, roughly speaking, the same pattern and connected the fourth period with the introduction of the ideas and methods of the set theory and the theory of functions of a real variable.

I stress two points. First, I think that there existed an initial version of the theory of probability whose acme were the LLN, the **De Moivre – Laplace** theorem (in essence proved by the former), and the inverse **Bayes** theorem (§ 5.2). Second, the modern stage of the theory, considered up to **Kolmogorov**, began with **Chebyshev**, and this fact should be here more clearly reflected. And so, my pattern of the history of probability is as follows.

1. Its antenatal period (from **Aristotle** to the mid-17th century).
2. The creation of its initial version (finally achieved by **Jakob Bernoulli, De Moivre** and **Bayes**).
3. Its development as an applied mathematical discipline (from **Bayes** to **Laplace** and **Poisson** to **Chebyshev**).
4. A rigorous proof of its limit theorems (**Chebyshev, Markov, Liapunov**) and its gradual transition to the realm of pure mathematics.
5. Axiomatization.

Laplace reasonably transferred the theory of probability from pure (as understood by his predecessors) to applied mathematics and Poisson and Poincaré followed suit. During the third, and partly the fourth stage, mathematicians barely regarded the theory of probability as a serious scientific discipline. And even now they barely recognize the theory of errors.

## 0.2. Mathematical Statistics

Its separation from probability or from statistics in general is difficult. It originated in the early years of the 20th century as the result of the work of the Biometric school and the Continental direction of statistics. Its aim is the systematization, processing and utilizing statistical data (**Kolmogorov** & **Prokhorov** 1988/1990, p. 138). They added a definition of statistical data:

"Information on the number of objects in some more or less extensive collection that have some specific properties".

They apparently excluded the theory of errors and it is unclear whether they meant raw or corrected (by exploratory data analysis) information.

Theoretical statistics is wider since it additionally studies the collection and that same analysis of data.

For **Pearson**, statistics remained an applied mathematical discipline whereas **Fisher** managed to create its theory, mathematical statistics.

### 0.3. The Theory of Errors

From its origin in the mid-18th century and until the 1920s the stochastic theory of errors had been a most important chapter of probability theory. Not without reason had **P. Lévy** (1925, p. vii) maintained that without it his main work on stable laws of distribution would have no *raison d'être*[1]. Actually, for the theory of errors that book was meaningless (Sheynin 1995c) and it is incomprehensible why he had not noticed it.

In turn, mathematical statistics borrowed its principles of maximum likelihood and minimal variance from the error theory. Today, the stochastic theory of errors is the application of the statistical method to the treatment of observations[2].

The history of the theory of errors has its own stages. In ancient times, astronomers were dealing with observations as they saw fit. At the second stage, beginning perhaps with **Tycho Brahe**, observations ceased to be "private property", but their treatment was not yet corroborated by quantitative considerations. This happened during the third stage (**Simpson, Lambert**), and the final, fourth stage was the completion of the classical theory of errors (**Laplace** and especially **Gauss**) although later **Helmert** fruitfully continued the relevant investigations.

The main peculiarity of the error theory is the usage of the notion of *real* (*true*) *value* of the constant sought, see § 6.3.3 and Sheynin (2007a), and I emphasize that it is really needed in experimental science rather than being outdated and on its way out. **Fourier** (1826/1890, p. 534) defined it as the limit of the arithmetic mean, which, incidentally, provides a new dimension to the **Gaussian** postulate of the mean [an expression due to **Bertrand** (1888a, p. 176)], see § 9A.2-2, and to the attempts to justify its usage by the notion of consistency, cf. §§ 9A.4-7, 11.2-8, 13.2-7 and 14.4-2. See also § 6.3.3. It is a peculiar fact that mathematicians neglected (and still neglect) the reasonable basis of Gauss' attitude and preferred (prefer) the hardly useful Laplacean approach.

### 0.4. The Statistical Method

It might be thought that statistics and statistical method are equivalent notions; it is usual, nevertheless, to apply the former term when studying population and to use the latter in all other instances and especially when applying statistics to natural sciences. However, there also exist such expressions as *medical* and *stellar statistics*, and, to recall, theory of errors (§ 0.3).

I understand that, respectively, they are applications of the statistical method to medicine, stellar astronomy and treatment of observations. This explanation is in line with Pearson's statement (1892, p. 15): *The unity of all science consists alone in its method*. And to my mind the statistical method is mathematical (or theoretical) statistics. It is tempting to suggest that mathematics is the application of the mathematical method, i. e., of the introduction and study of systems which possibly bear no relation to reality; a simple example: complex numbers.

Three stages might be distinguished in the history of the statistical method. At first, conclusions were being based on (statistically) noticed qualitative regularities, a practice which conformed to the qualitative essence of ancient science. Here, for example, is the statement of the Roman scholar **Celsus** (1935, p. 19):

*Careful men noted what generally answered the better, and then began to prescribe the same for their patients. Thus sprang up the Art of medicine.*

The second stage (**Tycho** in astronomy, **Graunt** in demography and medical statistics) was distinguished by the availability of statistical data. Scientists had then been arriving at important conclusions either by means of simple stochastic ideas and methods or even directly, as before. During the present stage, which dates back to the end of the 19$^{th}$ century, inferences are being checked by quantitative stochastic rules.

### 0.5. Some Nasty Facts

I briefly discuss some facts concerning my subject, but I will hardly err when connecting them with history of mathematics or perhaps with history of science as well.

1. Disdainful attitude towards information. This is contrary to the proper statement of Shaw & Austin (§ 10.8.3) and mostly due to the influence from beyond. The scientologists (no connection with religion) who vainly attempt to estimate the worth of scientists by numerical measurements apparently do not consider reviewing or abstracting as real scientific work. The scientific rat race (publish or perish!) *therefore* prevents careful and honest reviewing. The situation with respect to reviewing is horrible and the more so since for most researchers the main abstracting journals are now out of reach.

2. The same rat race hinders honest work in general.

3. Standardization, or the sledgehammer law engenders robots and is the curse on science. No one ever compelled Chekhov to write just like Tolstoy did or required Agatha Christie to imitate Edgar Poe. Science is however standardized by a sledgehammer. The format of the presentation of papers must be standardized however different are their aims, styles or lengths. And why? To justify the very existence of those scientologists? A most simple example: In a manuscript, I referred both to K. and E. S. Pearson, but a worker of the editorial staff deleted the initials! Spelling of names is required to disregard the rights of authors or their publications abroad. An eminent Russian author, Bernstein published many notes in France, signed his name as I mentioned. At the very least, it had thus become his pen name, but present authors are required to mention its ugly version, Bernshtein.

And now other aspects of the sledgehammer law: periodicals do not publish translations, period! In the West, Russian literature is badly known. If a journal has a thousand readers (which is more than generous) only one or two will be able to come across a reference to an important Russian source, get hold of it and understand something. One or two out of a thousand! Crass robotic stupidity.

Finally, no one will ever know how many worthy contributions have not been published since their authors were unable to overcome the sledgehammer law! See also my Preface.

The late professor Truesdell edited 49 volumes of his prestigious periodical, *Archive for History of Exact Science* without standardizing anything, and nothing bad had happened, but after his death that journal was promptly driven into the general fold.

I have described all this in detail (Sheynin 2017) and insistently recommend the readers to see it. My epigraph was a phrase from Shakespeare: *Something is rotten in the state of Denmark*.

### Notes

**1.** In 1887 **Chebyshev** (§ 13.1-4) indicated that the CLT can substantiate the MLSq although his statement had nothing to do with the Gaussian approach and **Poincaré** (§ 11.2-8), in the last years of his life, stated that the theory of errors had *naturellement* been his main aim in probability. However, his contribution to that theory is barely significant.

**2.** I especially note its application to metrology (Ku 1969). There also existed a determinate branch of the error theory, unconnected with stochastic considerations. Here is one of its problems: Compile such a program for measuring angles in the field that the unavoidable errors, both systematic and random, will least corrupt the final results (the coordinates of the stations of the trigonometric network). Another problem: bearing in mind the same goal, determine the optimal form of a triangle of triangulation. It is opportune to mention **Cotes** (1722) who solved 28 pertinent problems concerning plane and spherical triangles with various sets of measured elements. I leave such problems aside although, in principle, they can be included in the province of the design of experiments. The determinate error theory is also related to the exploratory data analysis that aims at uncovering the underlying structures (e.g., the systematic errors). Or, rather, that branch of the error theory was swallowed by both of those new disciplines.

# 1. The Antenatal Stage
## 1.1. Randomness, Probability, Expectation

Is an infinite (much more difficult: finite) sequence random or not? This is a fundamental problem which I leave aside. Chaitin (1975) is a popular discussion of the efforts to answer this question. The role of randomness (for example, in studying the evolution of species or in the kinetic theory of gases) is obvious. In statistics, a random variable ought to be stable (best of all, possess an invariable law of distribution) which is not necessary in natural science[1]. Thus, Kolmogorov (1983/1992, p. 515) stated:

*We should distinguish randomness in the wider sense* (*absence of any regularity*) *and stochastic ransom events* (*which constitute the subject of probability theory*).

**1.1.1. Aristotle.** Ancient scholars repeatedly mentioned randomness and (logical or subjective) probability[2]. The first notion implied lack of aim or law, such as a discovery of a buried treasure (Aristotle, *Metaphys.* 1025a) or a sudden meeting of two persons known to each other (Aristotle, *Phys*. 196b 30). In the second example randomness might be interpreted as an intersection of two chains of determinate events[3] and in both cases a small change in the action(s) of those involved (in these chains) would have led to an essential change of the result. Thus, the treasure would have remained hidden, the meeting would not have taken place, cf. § 11.2-9 for a link with modernity; true, randomness is not anymore connected with lack of intention.

In each of these illustrations the sudden event could have been (although was not) aimed at; Aristotle would not have called random a meeting with a stranger, or a discovery of a rusty nail. I would have asked: Suppose that a venomous snake bites a man in his own kitchen, should have Aristotle called such an event random?

The examples above also mean that randomness is a possibility, and Aristotle (*Methaphys.* 1064b – 1065a) indeed said so. It was **Hegel** (1812/1978, pp. 383 – 384) who formulated the converse statement. Suppose that a discrete random variable takes values $x_i$, $i = 1, 2, …, n$, with certain probabilities. Then, according to Hegel, any $x_i$ itself is random.

Aristotle's special example (*Phys*. 199b 1; also see *De generatione animalium* 767b5) considered deviations from law, monstrosities. The first departure of nature from the type

*Is that the offspring should become female instead of male; … as it is possible for the male sometimes not to prevail over the female …*

His was the first, not really convincing example of a dialectical conflict between randomness and necessity[4].

Aristotle (*De Caelo* 283b 1, also see *Phys*. 196b and *Rhetorica* 1369a) stated that "The products of chance and fortune are opposed to what is, or comes to be, always or usually". Similar assertions are made about accidents (*Metaphys.* 1025a, 1065a), chance conjunctions (*Anal. Post*. 87b) and coincidences (*Parva naturalia* 463b). Finally,

Aristotle (*Phys*. 197b 0, 197b 14 and 197a 5) distinguishes between chance (which is "in the sphere of moral actions") and spontaneity ("in lower animals and in many inanimate objects").

Junkersfeld (1945) was and possibly remains the author of the best account of Aristotle's concept of chance, and she (p. 22) stated that by chance he meant something that takes place

*Occasionally; has the character of an end; is such that it might have been the object of a natural or of a rational appetite; was not in fact the object of any appetite but came into being by accident*.

The circularity of this definition is actually due to Aristotle.

That, according to Aristotle, chance takes place rarely, is also indirectly seen in his statement (*Metaphys*. 1064b 15, 1026b, 1027a; *Ethica Eudemia* 1247b) that "none of the traditional sciences busies itself about the accidental […] but only sophistic […]". He (*Ethica Eud*. 1247a) notes, however, that "chance largely enters" navigation and strategy, and even that navigation "is as in throwing dice". Regrettably, he did not say anything about sophistic. I also note that neither does the theory of probability consider the accidental, but rather studies the laws of randomness. In the words of Schiller (*Spaziergang*, 1795),

"Der Weise […] sucht das vertraute Gesetz in des Zufalls grausenden Wundern, sucht den ruhenden Pol in der Erscheinungen Flucht".

**Plato** (Cioffari 1935, p. 30) has a better opinion about navigation: although

*Chance is almost everything in the arts of the* […] *pilot, and the physician, and the general*, […] *yet in a storm there must surely be a great advantage in having the aid of the pilot's art*.

And here is another example connected with astronomy (*De caelo* 292a 30 and 289b 22) that shows, incidentally, that games of chance even then provided examples of stochastic considerations:

*Ten thousand Coan throws* [whatever that meant] *in succession with the dice are impossible* and it is therefore *difficult to conceive that the pace of each star should be exactly proportioned to the size of its circle*[5].

**Cicero** (Franklin 2001, p. 164) came to a similar conclusion unconnected, however, with natural sciences.

Aristotle (*Anal. Priora* 70a 0) also stated that probability

*Is a generally approved proposition: what men know to happen or not to happen, to be or not to be, for the most part thus and thus* […], *e. g., "the envious hate"*[6].

He (*Rhetorica* 1402a 5; *De Poetica* 1461b) believed that "it is possible that improbable things" will happen. Yes, but his was of course only a qualitative formula. And, while discussing poetry, he (*De Poetica* 1460a 25) introduced a rudimentary scale of subjective probabilities: "a likely impossibility is always preferable to an unconvincing possibility". Understandably, he (*Rhetorica* 1376a 19) recommended the use of probabilities in law courts. He (*Problemata* 951b 0) also thought that it was better to acquit a "wrong-doer" than to condemn an innocent person[7] so that the statistical idea about errors of the two kinds is seen here.

Aristotle several times mentions luck (*Metaphys.* 1065a; *Rhetorica* 1361b) and fortune (*Magna Moralia* 1206b, 1270a) as notions expressing deviations from reasonable qualitative expectation.

Finally, Aristotle (for example, *Ethica Nicomachea* 1104a 24) believed that mean behaviour, moderation possessed optimal properties. Analogous statements had appeared even earlier in ancient China; the doctrine of means is attributed to a student of **Confucius** (Burov et al 1973, pp. 119 – 140). Again, a similar teaching existed in the Pythagorean school (Makovelsky 1914, p. 63), and **Nicomachus of Gerasa** (1952, p. 820) stated that a perfect number was in the realm of equality, was a mean between numbers the sum of whose divisors was less, and greater that the number itself; was between excess and deficiency. In medicine (§ 6.2) the mean was considered as the ideal state (of health), while in games of chance (§ 5) the (arithmetic) mean was believed to possess certain stochastic properties.

In the new time, the arithmetic mean became the main estimator of the constants sought (§ 1.2.4) and has been applied in civil life (§ 2.1.2). In addition, it is obviously connected with the appropriate expectation.

Both **Plato** and **Aristotle**, as witnessed by **Simplicius** (Sambursky 1956, p. 37; he only provided an exact reference on p. 3 of the reprint of this paper), called natural sciences "the science of the probable [eikotologia]". A later scholar, **Levi ben Gerson**, thought that the determinism of natural laws was only approximate and probable (Rabinovitch 1973, p. 77, with a reference to Levi's work), and, similarly (Ibidem, p. 166), **Maimonides** held that natural philosophy only offered probable theories.

**1.1.2. The Bible and the Talmud** (Sheynin 1998b). The earlier part of the Talmud is an interpretation of the first five books of the Old Testament. Called *Mishna*, it is subdivided into more than 60 treatises. The other part of the Talmud is made up of later commentaries on the Mishna. The Jerusalem Talmud was essentially completed in the fourth century and precedes the more influential Babylonian Talmud by about a century. I have seen the English edition of the Babylonian version (six volumes; London, 1951 – 1955) and I also refer to the German edition of the Talmud (12 vols, Berlin, 1930 – 1936) and the French edition of the Jerusalem Talmud (six volumes; Paris, 1960). When referring to the Talmud, I abbreviate it; thus, T/Avoth means treatise Avoth from the Talmud.

My reasoning is not linked with religious belief, but it might be emphasized that the sacred texts represented the general feelings of ancient communities. Therefore, any statements there, which complement and sometimes precede Aristotle's teaching, characterize the contemporaneous knowledge stored by mankind.

Ancient commentators of the Mishna hardly set high store by natural science or mathematics. Thus, Rabbi Elieser ben Chisma (T/Avoth3[18]) argued that the laws concerning

*Bird offerings and the onset of menstruation are essential traditional ordinances, but the calculation of the seasons and geometry are but the after-course of wisdom.*

**Maimonides** (1977, pp. 118 – 129), however, was a later exception:

*Our sages confirmed [that mathematical astronomy was] the true wisdom in the sight of the people. But the theories of the astrologists are devoid of any value.*

It is generally known that theology had hindered the advancement of science as well as the introduction of preventive measures against deadly small-pox (§ 6.2.3). On the other hand, the most celebrated scholars including **Kepler** and **Newton** had been inspired by wishing to discover Divine laws of nature. Scientists believed that those laws directly showed the acts of God whereas the Bible was, after all, written by mortals.

In any case theology demanded logical thought. Thus, the Old Testament (Proverbs 14:28) contains a direct and a contrary proposition: "In a multitude of people is the glory of a king, but without people a prince is ruined". See also Matthew 13:34, but the example above is interesting also in that it could have well been formulated in the mid-17th century by **John Graunt** or **William Petty**, the fathers of political arithmetic, or by **Johann Süssmilch** in the next century. Consider also the questions listed by Moses (Numbers 13: 17 – 20) which also provide a link with political arithmetic. He sent spies to the land of Canaan to find out

*Whether the people who dwell in it are strong or weak, whether they are few or many*, […] *whether the land is rich or poor* […].

Several authors preceded me. Hasofer (1967) studied the use of lots as described in the Talmud and Rabinovitch (1973), drawing on his earlier papers, dwelt on stochastic considerations in the same source. He aptly refuted an author of the *Encyclpaedia Hebraica* (Tel Aviv, 1962, XIV, pp. 920 – 921) who had maintained that

*Ancient Jewish thought did not know* [did not make use of] *the concept of probability.*

I am greatly indebted to Rabinovich, but am not satisfied either with his selection of cases or explanations. He discusses marginal problems, finds non-existing subjects (axioms of probability theory) or makes too much of some other topics (law of large numbers, sampling).

I am adding other important facts and, unless stated otherwise, offering my own comments.

**Randomness** is mentioned several times in the Old Testament:

*By chance I happened to be on mount Gilbo'a* (2 Samuel 1:6)
*Now there happened to be there a worthless fellow* (2 Samuel 20:1)
*A certain man drew his bow at a venture and struck the King of Israel* (1 Kings 22:34, 2 Chronicles 18:33)

Here, and in the other cases, randomness implicitly meant lack of purpose, cf. § 1.1.1. Ecclesiastes 9:11, however, states much more: Time and chance determine the fate of man.

The Talmud makes indirect use of **probability** (Rabinovitch 1973, Chapter 4). Thus, in certain cases prohibited fruit could not have exceeded 1/101 of its general quantity. Under less rigid demands it was apparently held that only the sign of the deviation from 1/2 was essential (Ibidem, p. 45).

Rabinovitch (1973, p. 44) noted that in the absence of a "clear majority" the Talmud regards doubts as half and half which is certainly only valid for subjective probabilities[8]. He refers to two tracts but in one of them (Makhshirin $2^{3-11}$) the statement is rather that the majority is equivalent to the whole. And at least in one case that principle seems to have been misapplied (Rabinovitch p. 45):

*If nine shops sell ritually slaughtered meat, and one sells meat that is not ritually slaughtered* [obviously, meat of cattle, twice] *and he* [and someone] *bought in one of them and does not know which one – it is prohibited because of the doubt; but if meat was found* [in the street] *one goes after the majority*.

Another example (Rabinovitch p. 40) obviously shows that subjective probability can lead to sophisms, see § 8.1 and Sheynin (2002b). I bear in mind the opinion of Rabbi Shlomo ben Adret (1235 – 1310), whom Rabinovitch mentions time and time again. There are several pieces of meat, all of them kosher except one. Eating the first one is allowed, since it is likely kosher; the same with the second one etc., and when only two pieces of meat are left, the forbidden piece was likely already eaten and they are also allowed.

Both the Bible and the Talmud provide examples of attempts to **distinguish between randomness and causality** and to act accordingly. Genesis 41:1 – 6 discusses cows and ears of corn as seen by the Pharaoh in two consecutive dreams. The dreams, essentially the same, differed in form. Both described an event with an extremely low probability (more precisely, a miracle), and they were thus divine rather than random. Then, Job (9:24 and 21:17 – 18) decided that the world was "given over to the wicked" [this being the cause] since the alternative had a low pseudo-statistical probability: "their lamp was put out rarely". This statement possibly brought forward a commentary (T/Avoth $4^{15}$) to the effect that we are unable to explain his conclusion.

A random event of sorts with a rather high probability seems nevertheless to have been presented in Exodus 21:29: an attack by an ox will be likely if, and only if, he "has been accustomed to gore".

And here are several examples from the Talmud. If in three consecutive days (not all at once, or in four days) three (nine) persons died in a town "bringing forth" 500 (1500) soldiers, the deaths should be attributed to a plague, and a state of emergence must be declared (Taanid $3^4$)**[9]**. The probability of death of an inhabitant during three days was apparently considered equal to 1/2, see above, so that a random and disregarded death of three people (in the smaller town) had probability 1/8. The case of all three dying at once was for some reason left aside. An early commentator, Rabbi Meir, lamely explained the situation by mentioning the goring ox (the German edition of the Talmud, Bd. 3, p. 707)**[10]**. A similar and simpler example concerned an amulet (Sabbath $6^2$): for being approved, it should have healed three patients consecutively. No follow-up checks had been demanded.

The second example (Leviticus 6:3 – 10) concerned the annual Day of Atonement when the High Priest brought up two lots, for the Lord and "for Azazel", one in each hand, from an urn. During 40 years, the first lot invariably came up in the Priest's right hand and that was regarded as a miracle and ascribed to his special merit [to a cause]. A special example concerned the redemption of the first born by lot (Jerus. Talmud/Sangedrin $1^4$). **Moses** wrote "Levite" on 22, 273 ballots and added 273 more demanding five shekels each. The interesting point here is that only 22,000 "Levite" ballots were needed so that Moses ran the risk of losing some of the required money. Nevertheless, the losing ballots turned up at regular intervals, which was regarded as a miracle. The existence of the superfluous ballots was not explained; I believe that the Israelites were afraid that the last 273 of them to draw the lots will be the losers. Such misgivings are, however, unfounded, see § 8.4 for an explanation based on subjective probabilities and Tutubalin (1972, p. 12) for a real demonstration.

Rabinovitch (1977, p. 335) provided a similar example, again from the Talmud. I cite finally the ruling about abandoned infants (Makhshirin $2^7$). A child, found in a town whose population was mostly gentile, was supposed to be gentile, and an Israelite otherwise (also when the two groups were equally numerous).

**1.1.3. Medicine** (Sheynin 1974, pp. 117 – 121). **Hippocrates** described a large number of [case histories] containing qualitative stochastic considerations in the spirit of the later **Aristotle** (and ancient science in general). Thus (1952a, pp. 54 – 55): "It is probable that, by means of […], this patient was cured"; or (1952b, p. 90),

*To speak in general terms, all cases of fractured bones are less dangerous than those in which …*

He (see below) also understood that healing depended on causes randomly changing from one person to another. Then, Hippocrates formulated qualitative [correlational] considerations, as for example (1952c, No. 44):

*Persons who are naturally very fat are apt to die earlier than those who are slender.*

He (1952c, § 69 p. 116; § 71 p. 117) also states that much depends on the constitution and general condition of the patient. Thus, in the second case:

*Men's constitutions differ much from one another as to the facility or difficulty with which dislocations are reduced* […].

Aristotle left similar reasoning, for example (*Problemata* 892a 0): "Why is it that fair men and white horses usually have grey eyes?" Or, noting some dependence between the climate or weather and human health (Ibidem 859b 5, 860a 5).

**Galen** also made use of stochastic reasoning. Most interesting is his remark (1951, p. 202) which I shall recall in § 11.2-9:

*… in those who are healthy* […] *the body does not alter even from extreme causes; but in old men even the smallest causes produce the greatest change.*

One of his pronouncements (Ibidem, p. 11) might be interpreted as stating that randomness was irregular: The body has "two sources of deterioration, one intrinsic and spontaneous, the other [which is] extrinsic and accidental", affects the body occasional[ly], irregular[ly] and not inevitabl[y]. In principle, one of his conclusions (1946, p. 113) is connected with clinical trials:

*What is to prevent the medicine which is being tested from having a given effect on two* [on three] *hundred people and the reverse effect on twenty others, and that of the first six people who were seen at first and on whom the remedy took effect, three belong to the three hundred and three to the twenty without your being able to know which three belong to the three hundred, and which to the twenty* […] *you must needs wait until you see the seventh and the eighth, or, to put it shortly, very many people in succession.*

Galen (1951, Book 1, Chapter 5, p. 13) also thought that a mean state, a mean constitution were best, cf. § 1.1.1:

*Health is a sort of harmony* […] *all harmony is accomplished and manifested in a two-fold fashion; first, in coming to perfection* […] *and second, in deviating slightly from this absolute perfection* […].

There also (pp. 20 – 21 of Book 1, Chapter 6) he says, for example, that "A good constitution [is] a mean between extremes".

**1.1.4. Astronomy.** Astronomers understood that their observations were imperfect; accordingly, they attempted to **determine some**

**bounds for the constants sought.** Thus (Toomer 1974, p. 139), the establishment of bounds "became a well-known technique […] practised for instance by **Aristarchus, Archimedes** and **Erathosthenes**". Here, for example, is Aristarchus(1959, p. 403):

*The diameter of the sun has to the diameter of the earth a ratio greater than [… 19:3] but less than [… 43:6], –*

greater than 6.33 and less than 7.17.

This did not however exclude the need to assign point estimates which was done by taking into account previous data (including bounds), qualitative considerations and convenience of further calculations. Modern mathematical statistics teaches that in case of large errors any observation can be chosen as the final result. Then, concerning the last-mentioned circumstance Neugebauer (1950, p. 252) remarks:

*The 'doctoring' of numbers for the sake of easier computation is evident in innumerable examples of Greek and Babylonian astronomy. Rounding-off in partial results as well as in important parameters can be observed frequently often depriving us of any hope of reconstructing the original data accurately.*

And, again (Neugebauer 1975, p. 107):

*In all ancient astronomy, direct measurements and theoretical considerations are […] inextricably intertwined […] ever present numerical inaccuracies and arbitrary rounding […] repeatedly have the same order of magnitude as the effects under consideration*[11].

**Babbage** (1874) seems to have been the only author in the new times who paid special attention to various methods of *doctoring*. He distinguished *hoaxing* or *forging* (downright deceiving), *trimming* (leaving the mean intact but "improving" the precision of measurement) and *cooking* (selecting observations at will). I return to errors of observation below.

Theoretical considerations partly replaced measurements in the Chinese meridian arc of 723 – 726 (Beer et al 1961, p. 26; Needham 1962, p. 42). At the end of their paper the authors suggested that it was considered more suitable to present harmonious results.

Special attention was being paid to the **selection of the optimal conditions for observation**; for example, to the determination of time intervals during which an unavoidable error least influenced the final result[12]. Thus (**Ptolemy** 1984, IX, 2, p. 421, this being the *Almagest*, his main contribution): on certain occasions (during *stations*) the motion of planets is too small to be observable. No wonder that he (Ibidem, III 1, p. 137) "abandoned" observations "conducted rather crudely". **Al-Biruni** (1967, pp. 46 – 51), the only Arab scholar to surpass Ptolemy and to be a worthy forerunner of **Galileo** and **Kepler**, to whom I return below, rejected four indirect observations of the latitude of a certain town in favour of its single and simple direct measurement.

Studying this aspect, Aaboe & De Solla Price (1964, pp. 2 and 3) argued that

*In the pre-telescopic era there is* […] *a curious paradox that even a well-graduated device* [their estimate: its error was 5 ] *for measuring celestial angles* […] *is hardly a match for the naked and unaided eye judiciously applied* […].

The function of smaller antique instruments, as they believe, was to serve as a means for avoiding calculations while the

*Characteristic type of measurement depended not on instrumental perfection but on the correct choice of crucial phenomena.*

They even, and I would say, mistakenly, called antique observations qualitative.

Hipparchus was aware that, under favourable conditions, a given error of observation can comparatively little influence the unknown sought (Toomer 1974, p. 131). And (Neugebauer 1950, p. 250) Babylonian astronomers of the Seleucid period possessed some understanding of such phenomena: their lunar and planetary computations "were based on an exceedingly small number of observations" and a "very high accuracy" of these observations was not needed. He continued:

*It seems to me one of the most admirable features of ancient astronomy that all efforts were concentrated upon reducing to a minimum the influence of the inaccuracy of individual observations with crude instruments by developing to the farthest possible limits the mathematical consequences of a very few basic elements.*

Elsewhere he (1948, p. 101) argued that in antiquity

*Observations were more qualitative than quantitative; <u>when angles are equal</u> may be decided fairly well on an instrument but not <u>how large are the angles</u> says Ptolemy with respect to the lunar and solar diameter*[s].

His reference is not exactly correct; it should have been *Almagest* (V 14, p. 252; H 417). The wording in the edition of 1984 is somewhat different.

**Regular observations** constituted the third main feature of ancient astronomy. Neugebauer (1975, p. 659) credited Archimedes and Hipparchus with systematic observation of the apparent diameters of the sun and the moon. Moreover, otherwise Hipparchus could have hardly been able to compile his star catalogue. **Ptolemy** (1984, III, 1, pp. 132 and 136; IV, 9, p. 206) apparently made regular observations as well. For example, in the second instance Ptolemy mentions his "series of observations" of the sun.

**Al-Biruni** (1967) repeatedly tells us about his own regular observations, in particular (p. 32) for predicting dangerous landslides (which was hardly possible; even latitude was determined too crude). Then, **Levi ben Gerson** (Goldstein 1985, pp. 29, 93 and 109) indirectly but strongly recommended them. In the first two cases he maintained that his regular observations proved that the declination of the stars and the lunar parallax respectively were poorly known. Thus, already in those times some contrast between the principle of regular (and therefore numerous) observations and the selection of the best of these began to take shape, also see § 1.2.2.

Astronomers undoubtedly knew that some errors, for example those caused by refraction, acted one-sidedly (e. g., Ptolemy 1984, IX 2); I also refer to Lloyd (1982, p. 158, n 66) who argues that Ptolemy "implicitly recognized what we should describe as systematic errors" and notes that he "has a [special] term for *significant* or *noteworthy* differences". Nevertheless, the separation into random and systematic errors occurred only at the end of the $18^{th}$ century (§ 6.3.1). But, for example, the end of one of Ptolemy's statements (1956, III, 2, p. 231) hints at such a separation:

*Practically all other horoscopic instruments [...] are frequently capable of error, the solar instruments by the occasional shifting of their positions or of their gnomons, and the water clocks by stoppages and irregularities in the flow of the water from different causes and by mere chance.*

**Al-Biruni** (1967, pp. 155 – 156) formulated a similar statement about water clocks.

Here are my conclusions derived from the *Almagest*:

Ptolemy knew that errors of observation were unavoidable or almost so.
He mentioned a number of sources of error.
He knew how to minimize the influence of some of them.
He gave thought to the choice of methods of observation.
He knew that some errors acted systematically.

It is possible that, when selecting a point estimate for the constants sought, ancient astronomers were reasonably choosing almost any number within the appropriate bounds (see above). Indeed, modern notions about treating observations, whose errors possess a "bad" distribution, justify such an attitude, which, moreover, corresponds with the qualitative nature of ancient science. Ptolemy's cartographic work corroborates my conclusion: he was mainly concerned with semblance of truth [I would say: with general correctness] rather than with mathematical consistency (Berggren 1991). A related fact pertains even to the Middle Ages (De Solla Price 1955, p. 6):

*Many medieval maps may well have been made from general knowledge of the countryside without any sort of measurement or estimation of the land by the 'surveyor'.*

Many authors maintained that **Ptolemy** had borrowed observations from **Hipparchus**, and, in general, doctored them while R. R. Newton (1977, p. 379) called him "the most successful fraud in the history of science". Yes, he likely borrowed from Hipparchus, but in good faith, in accordance with the day's custom. No, he had not doctored any observations, but rejected, adjusted or incorporated them "as he saw fit" (Gingerich 1983, p. 151; also see Gingerich 2002), he was an opportunist ready "to simplify and to fudge" (Wilson 1984, p. 43).

I adduce now three noteworthy statements which partly shaped my opinion above.

1) **Kepler** (1609/1992, p. 642/324):

*We have hardly anything from Ptolemy that we could not with good reason call into question prior to its becoming of use to us in arriving at the requisite degree of accuracy.*

2) **Laplace** (1796/1884, p. 413):

*Hipparchus among all ancient astronomers deserves the gratitude of astronomy for the large number and precision of his observations, for important conclusions which he had ben able to make by comparing them with each other and with earlier observations, and for the witty methods by which he guided himself in his research. Ptolemy, to whom we are mostly indebted for acquainting us with his work, had invariably based himself on Hipparchus' observations and theories. He justly appraised his predecessor…*

Indeed, Ptolemy (1984, IX 2, p. 421; H 210) also see (1984 III 1, p. 136; H 200) called him "a great lover of truth". And Laplace on the next page:

*His Tables of the Sun, in spite of their imperfection, are a durable monument to his genius and Ptolemy respected them so much that he subordinated his own observations to them.*

3) **Newcomb** (1878, p. 20): "... all of Ptolemy's *Almagest* seems to me to breathe an air of perfect sincerity"[13].

Al-Biruni (1967, p. 152) was the first to consider, although only qualitatively, the propagation of computational errors and the combined effect of observational and computational errors:

*The use of sines engenders errors which become appreciable if they are added to errors caused by the use of small instruments, and errors made by human observers.*

One of his statements (Ibidem, p. 155) on the observation of lunar eclipses for determining the longitudinal difference between two cities testified to his attempt to exclude systematic influences from final results: Observers of an eclipse should *Obtain all its times* [phases] *so that every one of these, in one of the two towns, can be related to the corresponding time in the other. Also, from every pair of opposite times, that of the middle of the eclipse must be obtained.*

Such a procedure would have ensured some understanding of the systematic influences involved, cf. **Boscovich**' calculation of a latitudinal difference in § 6.3.2.

For **Al-Biruni**, see **Al-Khazini** (1983, pp. 60 – 62), the arithmetic mean was not yet the universal estimator of the constants sought; when measuring the density of metals he made use of the [mode], the [midrange] as well as of some values situated within the extreme measurements[14], also see Sheynin (1992; 1996a, pp. 21 – 23).

Levi ben Gerson (Goldstein 1985, p. 28) stated that (at least) "all Ptolemy's predecessors preferred" to use no instruments at all. He also distrusted the precision of Ptolemy's invention, "the instrument of the rings" as well as astrolabes and quadrants in general. He knew that small errors of observation can produce large errors in the resulting position of the stars. Regrettably he (Ibidem, p. 29) attributed this fact to the "third kind" of errors thus spoiling their classification. And he even claimed to have constructed and tested an instrument which allowed to make error-free observations.

Levi (Ibidem, pp. 29, 93 and 109) several times mentions his regular observations and maintains (Rabinovitch 1974, p. 358) that repeat observations ought to be made "as often as is required". This is interesting but not definite enough. Moreover, the same author (1973, p. 77) quoted Levi's belief in the universal existence of uncertainty, and thus in the impossibility of such an instrument:

*Perfect knowledge of a thing is to know it as it is* […], *that is, to know that aspect that is determined and bounded and to know also that indeterminacy which is in it*.

**1.1.5. Maimonides and Thomas Aquinas.** In accordance with the Talmud, the consumption of some foods was allowed only for priests and in many other cases the part of the forbidden food should not have exceeded certain limits, cf. §1.1.2. In this connection **Maimonides** (Rabinovitch 1973, p. 41) listed seven relevant ratios, i.e., seven different probabilities of eating the forbidden food.

The Talmud also qualitatively discussed the estimation of prices for quantities depending on chance (Franklin 2001, p. 261). It is opportune to recall here that the Roman lawyer **Ulpianus** (170 – 228) compiled a table of expectations of life, common for men and women (Sheynin 1977b, pp. 209 – 210), although neither the method of its compilation, nor his understanding of expectation are known. His table was being used for determining the duration of some allowances (Kohli & **van der Waerden** 1975, p. 515). At least methodologically his table constituted the highest achievement of demographic statistics until the 17th century.

**Maimonides** (Rabinovitch 1973, p. 164) mentioned expectation on a layman's level. He noted that a marriage settlement (providing for a widow or a divorced wife) of 1000 *zuz* "can be sold for 100 [of such monetary units], but a settlement of 100 can be sold only for less than 10". It follows that there existed a more or less fixed expected value of a future possible gain[15].

A marriage settlement is a particular case of insurance; the latter possibly existed in an elementary form even in the 20th century BC (Raikher 1947, p. 40). Another statement of Maimonides (Rabinovitch 1973, p. 138) can also be linked with jurisprudence and might be considered as an embryo of **Jakob Bernoulli's** thoughts about arguments (cf. §3.2.1):

*One should not take into account the number of doubts, but rather consider how great is their incongruity and what is their disagreement with what exists. Sometimes a single doubt is more powerful than a thousand other doubts.*

Incongruity and disagreement, however, rather have to do with opinions.

We also find an embryo of a random variable in Maimonides' works (Rabinovitch 1973, p. 74):

*Among contingent things some are very likely, other possibilities* [of other things?] *are very remote, and yet others are intermediate.*

In the new time one of the first to follow suit in natural science (not in lotteries) was **Maupertuis** (1745, vol. 2, pp. 120 – 121) who effectively explained instances when a child resembled one of his remote ancestors, as well as [mutations] by *non-uniform* randomness. At the same time, however, while discussing the origin of eyes and ears in animals, he (1751, p. 146) only compared "une attraction uniforme & aveugle" with some "principe d'intelligence" and came out in favour of design.

An important digression from statistics is not out of order. While considering a case of uncertain paternity, the Talmud (Makkot 3[15–16]; Rabinovitch 1973, p. 120) declared: "What we see we may presume, but we presume not what we see not", cf. the examples of direct and contrary propositions at the beginning of this subsection. Compare **Newton** (1729/1960, p. 547):

*I have not been able to discover the cause of these properties of gravity from phenomena, and I frame no hypotheses* [*hypotheses non fingo*].

Maimonides (1975, p. 123) advised physicians and judges to test, so to say, simpler hypotheses, and then, only when necessary, to go on to more complex assumptions. Cure ought to be attempted "through food", then by "gentle" medicine whereas "strong drugs" should be considered as the last means. Just the same, a judge should try to "effect a settlement" between litigants, then judge "in a pleasant manner" and only then become "more firm". Compare Newton's Rule No. 1 of reasoning in philosophy (Ibidem, p. 398):

*We are to admit no more causes of natural things than such as are both true and sufficient to explain their appearances*[16].

Moreover, *Occam's razor* (William Occam, 1287 – 1347) comes to mind: when two explanations are possible, the simpler is usually better.

**Thomas** (1952) was the main commentator of **Aristotle** and he strove to adapt the pagan Philosopher to Christianity. Just as his hero, he believed that random events occurred in the minority of cases and were due to some hindering causes (Sheynin 1974, p. 103):

*Casual and chance events* are such as *proceed from their causes in the minority of cases and are quite unknown*:

*Some causes are so ordered to their effects as to produce them not of necessity but in the majority of cases, and in the minority to fail in producing them …* [which] *is due to some hindering cause.*

Again following Aristotle, Thomas illustrated this statement by the "production of woman" which was nevertheless "included in nature's intention"[17]. He (Ibidem, p. 108) also maintained that law courts should guide themselves by stochastic considerations.[18]

*In the business affairs of men* […] *we must be content with a certain conjectural probability*.

On the introduction of moral certainty see §§ 2.1.2, 2.2.2 and 3.2.2.

Finally, Thomas (Ibidem p. 107) attributed ranks and degrees to miracles; on this point see also Kruskal (1988). At a stretch, this meant an introduction of qualitative probabilities. Byrne (1968) studied the work of Thomas and argued (pp. 202 – 208) that he had used a rudimentary frequency theory of probability but I fail to see this. Anyway, Thomas provided a link between medieval and modern science.

### 1.2. Mathematical Treatment of Observations
#### 1.2.1. Theory of Errors: General Information.
I introduce some notions and definitions which will be needed in § 1.4. Denote the observations of a constant sought by

$$x_1, x_2, …, x_n, x_1 \quad x_2 \quad … \quad x_n. \tag{1}$$

It is required to determine its value, optimal in some sense, and estimate the residual error. The classical theory of errors considers independent observations (see § 9A.4-4) and, without loss of generality, they might be also regarded as of equal weight. This problem is called *adjustment of direct observations.*

Suppose now that $k$ unknown magnitudes $x, y, z, …$ are connected by a redundant system of $n$ physically independent equations ($k < n$)

$$a_i x + b_i y + c_i z + \ldots + s_i = 0 \qquad (2)$$

whose coefficients are given by the appropriate theory and the free terms are measured. The approximate values of $x, y, z, \ldots$ were usually known, or they could have been calculated by solving any subsystem of $k$ equations from (1) hence the linearity of (2). The equations are linearly independent (a later notion), so that the system is inconsistent (which was perfectly well understood). Nevertheless, a solution had to be chosen which should have left small enough residual free terms (call them $v_i$). More precisely, those $v_i$'s should have obeyed the properties of *usual* random errors (approximately equal number of positive and negative errors and more errors small in absolute value than large errors).

The values of the unknowns thus obtained are their estimates ($x_0, y_0, \ldots$) and this problem, again with an evaluation of their errors, is called *adjustment of indirect measurements*.

Since the early 19$^{th}$ century the usual condition for solving (2) was that of least squares

$$W = \Sigma v_i^2 = [vv] = v_1^2 + v_2^2 + \ldots + v_n^2 = \min^{19}, \qquad (3)$$

so that

$$\partial W/\partial x = \partial W/\partial y = \ldots = 0. \qquad (4)$$

Conditions (4) easily lead to a system of *normal equations*

$$[aa]x_0 + [ab]y_0 + \ldots + [as] = 0,\ [ab]x_0 + [bb]y_0 + \ldots + [bs] = 0, \ldots, \quad (5)$$

having a positive definite and symmetric matrix. For direct measurements the same condition (3) leads to the arithmetic mean. Another no less important and barely known to statisticians pattern of adjusting indirect observations is described in § 9A.4-9.

**1.2.2. Regular Observations.** I have mentioned them in § 1.1.4. They are necessary for excluding systematic, and compensating the action of random errors. **Kepler** is known to have derived the laws of planetary motion by issuing from **Tycho's** regular observations, who had thought that they provided a means for averaging out "random, instrumental and human error" (Wesley 1978, pp. 51 – 52). Note, however, that both instrumental and human errors can well be partly random which testifies to Wesley poor knowledge of the theory of errors. Wesley also states that Tycho (somehow) combined measurements made by different instruments and that exactly that approach might have produced most favourable results[20].

Nevertheless, it seems that, when compiling his star catalogues, **Flamsteed**, the founder of the Greenwich observatory, made use of only a part of his observations (Baily 1835, p. 376):

*Where more than one observation of a star has been reduced, he has generally assumed that result which seemed to him most*

*satisfactory at the time, without any regard to the rest. Neither […] did he reduce the whole (or anything like the whole) of his observations. […] And, moreover, many of the results, which have been actually computed,* […] *have not been inserted in any of his MS catalogues.*

Reduction, however, was a tiresome procedure and, anyway, Flamsteed likely considered his results as preliminary which is indirectly testified by the last lines of the passage just above and by his own pronouncements (Sheynin 1973c, pp. 109 – 110), for example by his letter of 1672, see Rigaud (1841, pp. 129 – 131):

*I* […] *give you the sun's diameters, of which I esteem the first, third and fourth too large* […] *the rest I esteem very accurate, yet will not build upon them till I have made some further trials with an exacter micrometer.*

**Bradley's** principle of treating observations remains somewhat unexplained (Sheynin 1973c, p. 110; Rigaud 1832). In one case (Rigaud, p. 78) he derived the arithmetic mean of 120 observations, and he (Ibidem, p. 17) supplemented his discovery of nutation of the Earth's axis by stating that

*This points out to us the great advantage of cultivating* [astronomy] *as well as every other branch of natural knowledge by a regular series of observations and experiments.*

He also discovered the aberration of light. At the same time he (Rigaud 1832, p. 29) reported that

*When several observations have been taken of the same star within a few days of each other, I have either set down the mean result, or that observation which best agrees with it.*

And Boyle, the cofounder of scientific chemistry and co-author of the **Boyle – Mariotte** law, kept to his own rule (Boyle 1772, p. 376):

*Experiments ought to be estimated by their value, not their number*; […] *a single experiment* […] *may as well deserve an entire treatise* […]. *As one of those large and orient pearls may outvalue a very great number of those little* […] *pearls, that are to be bought by the ounce* […].

So are series of observations needed? All depends on the order of the random errors, their law of distribution, on the magnitude of systematic influences, the precision and accuracy required (the first term concerns random errors, the second one describes systematic corruption) and on the cost of observation. In any case, it is hardly advisable to dissolve a sound observation in a multitude of worse measurements. In geodesy, observations ought to be carried out under differing (but good enough) meteorological conditions to lessen the influence of systematic errors. Mendeleev (§ 10.9.3) likely held the same opinion. In general, that same aim requires a prior knowledge of the number of observations which prevents the application of sequential analysis.

**1.2.3. Galileo.** When treating discordant observations of the parallax of the New Star of 1572 made by several astronomers, Galileo (1632, Day Third) formulated some propositions of the not yet existing error theory[21] and, first of all, indicated the properties of *usual* random errors (also known to **Kepler**, see § 1.2.4). The method of observation was of course worthless: in those days, even annual star parallaxes remained unyielding to measurement. Astronomers were only interested in placing the New Star either "beneath" the Moon or "among" the fixed stars (but even that problem was really insoluble). In essence, Galileo compared two natural scientific hypotheses with each other and chose the latter. His test, later applied by **Boscovich** (§ 6.3.2), was the minimal sum of absolute values of the parallaxes. Because of computational difficulties, however, Galileo only took into account some of the pairs of observations.

Buniakovsky (1846, Chapter on history of probability) mentioned his investigation in a few lines but did not provide a reference; Maistrov (1967/1974, pp. 30 – 34) described Galileo's reasoning, but see **Hald** (1990, pp. 149 – 160) for a detailed and rigorous discussion.

Apparently during 1613 – 1623 Galileo wrote a note about the game of dice first published in 1718 (F. N. David 1962, pp. 65 – 66; English translation on pp. 192 – 195). He calculated the number of all the possible outcomes (and therefore, indirectly, the appropriate probabilities) and testified that gamblers were believing that 10 or 11 points turned out more often than 9 or 12. If only these events are considered (call them A and B respectively), then the difference between their probabilities

$P(A) = 27/52$, $P(B) = 25/52$, $\Delta P = 1/26 = 0.0385$

can apparently be revealed. On determination of such small differences see also Note 10 to Chapter 2.

In 1610 – 1612 several astronomers discovered sunspots. Daxecker (2004; 2006) singled out the contribution of Christopher Scheiner and discussed his book *Rosa ursina sive Sol* of 1626 – 1630. Much more known is Galileo's achievement (1613) who managed to separate the regular rotation of the spots with the Sun itself from their random proper movement relative to the Sun's disk and estimated the period of the Sun's rotation as one lunar month; the present-day estimate is 24.5 – 26.5 days.

**Humboldt** (1845 – 1862, 1858, p. 64n) thought that sunspots could have well been observed much earlier, for example "On the coast of Peru, during the *garua* [?] […] even with the naked eye", but he was unable to substantiate his opinion.

Here, however, is the great traveller Marco Polo (Jennings 1985, p. 648). He described his conversation with the astronomer Jamal-ud-Din, a Persian, and his team of Chinese astronomers. They discovered the sunspots (and apparently observed them repeatedly) when "the desert dust veiled the Sun". Marco Polo narrative appeared in 1319 and his conversation with astronomers took place in the last quarter of the 13[th] century somewhere near the present-day Chinese city Tianjin. Marco Polo offered no comment, and it seems that he had not become interested at all.

The subsequent history of sunspots is connected, first and foremost, with their diligent observer, Schwabe, but it was Humboldt's description of his work that turned the attention of astronomers to them (Clerke 1885/1893, p. 156).

While denying any possibility of an irregular motion of a heavenly body, Galileo (1623/1960, p. 197) apparently denied randomness in general:

*Those lines are called regular which, having a fixed and definite description, have been susceptible of definition and of having their qualities and properties demonstrated.* […] *But irregular lines are those which have no determinacy whatever* […] *and hence indefinable.* […] *The introduction of such lines is in no way superior to the sympathy, antipathy, occult properties, influences and other terms employed by some philosophers as a cloak for the correct reply, which would be* <u>I do not know</u>.

This was possibly directed against some of Kepler's utterances and then, even if to a small extent, explains Galileo's still mysterious failure to recognize the Keplerian laws of planetary motion.

**1.2.4. Kepler. Randomness and the Treatment of Observations.** Randomness played a certain part in Kepler's astronomical constructions. True, he (1606/2006, p. 163) denied it:

*But what is randomness? Nothing but an idol, and the most detestable of idols; nothing but contempt of God sovereign and almighty as well as of the most perfect world that came out of His hands.*

Nevertheless, his laws of planetary motion were unable to justify the values of the eccentricity of their orbits. He finally had to consider them random, caused by disturbances, which was quite in the **Aristotelian** spirit (§ 1.1.1)[22], occasioned by deviation from (Divine) laws of nature. In this connection I quote **Poincaré** (1896/1912, p. 1) who most clearly formulated the dialectical link between randomness and necessity in natural sciences:

*There exists no domain where precise laws decide everything, they only outline the boundaries within which randomness may move. In accordance with this understanding, the word randomness has a precise and objective sense.*

Kepler (1596) first encountered those eccentricities when attempting to construct a model of the solar system by inserting the five regular solids between the spheres of the then six known planets: they, the eccentricities, and, for that matter, unequal one to another, much worried him. At that time, Kepler understood eccentricity as the preordained eccentric position of the Sun as measured from the centre of the circular orbit of a given planet.

He then changed his (actually, ancient) definition of eccentricity:

*In the motion of the planets, the variables, if any (such as, in the motion of the planets, the varying distance from the sun, or the eccentricity), arise from the concurrence of extrinsic causes* (Kepler 1609/1992, Chapter 38, p. 405).

On the same page he illustrated his opinion by obstacles which prevent rivers from descending "towards the centre of the earth", and finally, on the next page, he concluded that "other causes are conjoined with the motive power from the sun" [affect their motion].

Kepler (1620 – 1621/1952, p. 932) voiced his main statement in a later contribution:

*If the celestial movements were the work of mind, as the ancients believed, then the conclusion that the routes of the planets are perfectly circular would be plausible.* […]
*But the celestial movements are* […] *the work of* […] *nature* […] *and this is not proved by anything more validly than by observation of the astronomers, who* […] *find that the elliptical figure of revolution is left in the real and very true movement of the planet.* […]
*Because in addition to mind there was then need of natural and animal faculties* [which] *followed their own bent* […].

Later, however, Kepler returned to his still cherished explanation of the construction of the Solar system by means of the regular solids and argued that

*The origin of the eccentricities of the individual planets* [is] *in the arranging of the harmonies between their motions* (Kepler 1619/1997, title of Chapter 9 of Book 5 on p. 451).

On that very page he explained that God had combined the planetary motions with the five regular solids and thus created the only most perfect prototype of the heaven.

Again in the same chapter, in Proposition 5, on p. 454, he indirectly mentioned in this connection his second law of planetary motion; for that matter, he could have referred to it in his previous source, see above. So it was this law according to which the preordained eccentricities regulated the planetary motions. Had Kepler realized that his contrived theory of the construction of the Solar system was patently wrong, he would have probably returned to his previous opinion about the "faculties" which "followed their own bent". True, it is difficult to see how, even acknowledging his theory, could have eccentricities regulated planetary motions.

On the contrary, Kepler repeated his statement about the just mentioned preordained eccentricities in the second edition of 1621 of his theory (Kepler 1596/1963, pp. 97 and 104). Note, however, that he had not abandoned his opinion about the possibility of external forces corrupting Divine intention.

Without referring to anyone, **Kant** (1755/1910, p. 337) repeated Kepler's pronouncement about external forces:

*The multitude of circumstances that participate in creating each natural situation, does not allow the preordained regularity to occur. Why are their* [the planets'] *paths not perfectly circular? Is it not seen clearly enough, that the cause that established the paths of celestial bodies […] had been unable to achieve completely its goal? […] Do we not perceive here the usual method of nature, the invariable deflection of events from the preordained aim by various additional causes?*

In § 3.2.4 I mention Kepler in connection with his ideas about the end of the world.

**Kepler** (1604/1977, p. 337) also decided that a possible (that is, an aimless) appearance of the New star in a definite place and at a definite moment (both the place and the moment he, in addition, considered remarkable) was so unlikely that it should have been called forth by a cause [it had an aim], cf. § 1.1.1.

Kepler (Sheynin 1974, § 7) considered himself the founder of scientific astrology, of a science of [correlational] rather than strict influence of heaven on men and states. Thus (Kepler 1619/1952, book 4, pp. 377 – 378), his heavenly bodies were not Mercury, but **Copernicus** and **Tycho Brahe**, and the constellations at his birth had only awakened rather than hardened his spirit and the abilities of his soul. And (1610/1941, p. 200), "heaven and earth are not coupled as cog-wheels in a clockwork". Before him Tycho likely held the same view (Hellman 1970, p. 410). As an astrologer, **Ptolemy** (1956, I, 2 and I, 3) also believed that the influence of the heaven was a tendency rather than a fatal drive, that astrology was to a large extent a science of qualitative correlation, and **Al-Biruni** (1934, p. 232) likely thought the same way: "The influence of Venus is towards …", "The moon tends …"[23]

They both thus anticipated **Tycho** and **Kepler**. For the latter, the main goal of astrology was not the compilation of horoscopes concerning individuals, but the determination of tendencies in the development of states for which such circumstances as geographical position, climate, etc., although not statistical data, should also be taken into account. Cf. the approach of political arithmeticians (§ 2.1.4).

A few decades ago physicists and mechanicians began to recognize randomness as an essentially more important agent, a fact which I am leaving aside, as well as chaotic phenomena, that is, exponentially deviating paths of, for example, a ball on a non-elliptical billiard table. Owing to unavoidable small uncertainty of its initial conditions, the path becomes a cloud which fills a certain region. This chaos defies simple quantitative definition[24].

I qualitatively explain the difference between *ordinary* randomness and chaotic motion: however long and complicated is the fall of a tossed coin, the outcomes of this toss do not change. However, the instability of chaotic motion rapidly increases with time and an uncountable infinity of possible paths (a cloud) inevitably appears.

Kepler had to carry out enormous calculations and, in particular, to adjust both direct and indirect measurements. The most interesting example in the first case (Kepler 1609/1992, p. 200/63) was the adjustment of the four following observations (I omit the degrees)

$x_1 = 23\ 39$ , $x_2 = 27\ 37$ , $x_3 = 23\ 18$ , $x_4 = 29\ 48$ .

Without any explanation, Kepler selected $x = 24\ 33$ as the "medium ex aequo et bono" (in fairness and justice). A plausible reconstruction (Filliben, see Eisenhart 1976) assumes that $x$ was a generalized arithmetic mean with weights of observations being 2, 1, 1, and 0 (the fourth observation was rejected). But the most important circumstance here is that the Latin expression above occurred in **Cicero's** *Pro A. Caecina oratio* and carried an implication *Rather than according to the letter of the law*. Rosental & Sokolov (1956, p. 126), in their Latin textbook intended for students of law, included that expression in a list of legal phrases and adduced Cicero's text (p. 113; German translation, see Sheynin 1993c, p. 186). In other words, Kepler, who likely read Cicero, called the ordinary arithmetic mean *the letter of the law*, i.e., the universal estimator [of the parameter of location].

It might be thought that such a promotion was caused by increased precision of observations; their subjective treatment (§ 1.1.4) became anachronistic. In addition, astronomers possibly began perceiving the mean as an optimal estimator by analogy with the ancient idea of the expediency of "mean" behaviour (§ 1.1.1)[25].

Overcoming agonizing difficulties, Kepler repeatedly adjusted indirect measurements. I dwell here on two points only. And, first of all: How had he convinced himself that **Tycho's** observations were in conflict with the **Ptolemaic** system of the world? I believe that Kepler applied the minimax principle (§ 6.3.2) demanding that the residual free term of the given system of equations, maximal in absolute value, be the least from among all of its possible "solutions". His equations were not even algebraic, but, after some necessary successive approximations, they might have been considered linear. He apparently determined such a minimum, although only from among some possibilities, and found out that that residual was equal to 8 which was inadmissible, see his appropriate statement (1609/1992, p. 286/113):

*The divine benevolence had vouchsafed us Tycho Brahe, a most diligent observer, from whose observations the* 8 *error in this Ptolemaic computation is shown* […] [and, after a few lines] *because they could not have been ignored, these eight minutes alone will have led the way to the reformation of all of astronomy, and have constituted the material for a great part of the present work.*

Then, when actually adjusting observations, he (Ibidem, p. 334/143) corrupted them by small arbitrary corrections. He likely applied elements of what is now called statistical simulation, but in any case he must have taken into account the properties of "usual" random errors, i.e., must have chosen a larger number of small positive and negative corrections and about the same number of the corrections of each sign. Otherwise, Kepler would have hardly achieved success.

I stress finally that Kepler (1609/1992, p. 523/256), although not from the very beginning, came to understand that each observation had to be taken into consideration:

*Since the first and the third position […] agree rather closely, some less thoughtful person will think that it* [the constant sought] *should be established using these, the others being somehow reconciled. And I myself tried to do this for rather a long time.*

Here is an overview (Pannekoek 1961/1989, pp. 339 – 340):

*In former centuries the astronomer selected from among his observations those that seemed the best; this made him liable to bias or inclined to select such data as showed a possibly unreal agreement. […] A new attitude was brought into being, typical of the nineteenth-century scientist towards his material; it was no longer a mass of data from which he selected what he wanted, but it was the protocol of an examination of nature, a document of facts to which he had to defer.*

## Notes

**1.** See Note 22 (examples of randomness in the general sense).

**2.** I leave aside the views of **Democritus**, **Epicurus** and **Lucretius** since I think that their works are not sufficiently understandable. **Russell** (1962, p. 83) considered them "strict determinists", but many other scholars were of the opposite opinion. Thus, **Kant** (1755/1910, p. 344) remarked that the random mutual movement of the Lucretius' atoms had not created the world. Many ancient scientists reasoned on randomness, see Note 3. And here is a strange statement of **Strabo** (1969, 2.3.7), a geographer and historian:

*Such a distribution of animals, plants and climates as exists, is not the result of design – just as the difference of race, or of language, is not, either – but rather of accident and chance.*

**Chrysippus** (Sambursky 1956/1977, p. 6) held that chance was only the result of ignorance and **St. Augustinus**, and, much later, **Spinoza** and **D'Alembert** expressed similar thoughts (**M. G. Kendall** 1956/1970, p. 31, without an exact reference). Kendall also maintains that **Thomas** (Sheynin 1974, p. 104) stated something similar: a thing is fortuitous with respect to a certain, but not to a universal cause but I think that that proposition is rather vague.

**3.** Bru (**Cournot** 1843/1984, p. 306) noted that a number of ancient scholars expressly formulated such an explanation of chance. An example taken from ancient Indian philosophy (Belvalkar & Ranade 1927, p. 458) admits of the same interpretation:

*The crow had no idea that its perch would cause the palm-branch to break, and the palm-branch had no idea that it would be broken by the crow's perch; but it all happened by pure Chance.*

Lack of aim or intersection of chains of events may be seen in **Hobbes'** remark (1646/1840, p. 259):

*When a traveller meets with a shower, the journey had a cause, and the rain had a cause [...]; but because the journey caused not the rain, nor the rain the journey, we say they were contingent one to another.*

He added, however, that the rain was a random event since its cause was unknown, cf. above.

**4.** More interesting is a by-law pronounced in ancient India, between the 2$^{nd}$ century BC and 2$^{nd}$ century of our era (Bühler 1886/1967, p. 267):

*The witness* [in law-suits pertaining to loans]*, to whom, within seven days after he has given evidence, happens* [a misfortune through] *sickness, a fire or the death of a relative, shall be made to pay the debt and a fine.*

This was an attempt to isolate necessity (a speedy divine punishment) from chance. Another example (Hoyrup 1983) describes the death of 20 murderers in an accident with only one person (the one who had unsuccessfully tried to prevent the murders) surviving. This story, concerning the year 590 or thereabouts, was possibly invented, but in any case it illustrated the same attempt.

**5.** Apparently: an invariable mutual arrangement of the stars cannot be random. **Levi ben Gerson** (1999, p. 48) left a similar but less direct statement, but, strictly speaking, such arguments are not convincing. It is impossible to say beforehand which outcomes of ten (say) throws of a coin exhibit regularity, and which are a result of chance: all of them are equally probable. It is opportune to recall the **D'Alembert – Laplace** problem: The word *Constantinople* is composed of printer's letters; was the composition random? D'Alembert (1768a, pp. 254 – 255) stated that all the arrangements of the letters were equally probable only in the mathematical sense, but not "physically" so that the word was not random. Laplace (1776/1891, p. 152 and 1814/1995, p. 9) more correctly decided that, since the word had a definite meaning [had an aim], its random [aimless] composition was unlikely. He thus reasonably refused to solve this problem strictly. **Poisson** (1837a, p. 114) provided an equivalent example and made a similar conclusion. Matthiesen (1867), however, reported an extremely rare event: in a game of whist it occurred that each of the four gamblers received cards of only one suit. It is of course impossible to check this story and, anyway, it is reasonable to follow Laplace and Poisson.

**6. Cicero** understood the probable just as **Aristotle** did (Franklin 2001, p. 116). Much later, in the *Digest* (the Roman code of civil laws, 533), the same interpretation was indirectly repeated (Ibidem, p. 8). I supplement this information about Cicero by another of his statements (1997, Book 1, § 12, p. 7):

*Many things are probable and [...] though these are not demonstrably true, they guide the life of the wise man because they are so significant and clear-cut.*

**7.** Later authors repeatedly expressed the same idea; I name **Maimonides** (Rabinovitch 1973, p. 111), **Thomas Aquinas** (Byrne 1968, pp. 223 and 226) and even **Peter the Great** (in 1716, his *Kriegs-Reglement*, see Sheynin 1978b, p. 286, Note 39). I do not think, however, that practice followed such statements.

**8.** Half-proofs were mentioned in law courts apparently in the 1190s (Franklin 2001, p. 18).

**9.** These numbers indirectly indicated the population of the towns. Deaths of infants hardly counted here.

**10.** And here is **Markov's** opinion of 1915 about graduates of Russian Orthodox seminaries (Sheynin 1993a, p. 200):

*They are getting accustomed by their schooling to a special kind of reasoning. They must subordinate their minds to the indications of the Holy fathers and replace their minds by the texts from the Scripture.*

**11.** To compare all this with the practice of the 19$^{th}$ and partly 20$^{th}$ century: Gauss calculated his measured angles to within 0.″001 (actual precision about 700 hundred times lower); **Karl Pearson** kept to the same habit, even if not to the same extent, and at least one similar example pertained to Fisher, see a discussion of this subject in *Science*, vol. 84, 1936, pp. 289 – 290, 437, 483 – 484 and 574 – 575.

**12.** A good example about the determinate error theory (Chapter 0, Note 2).

**13.** And here is a general estimate of **Ptolemy** written by a Soviet dinosaur during Russia's brighter (!) years (Chebotarev 1958, p. 579): his system "held mankind in spiritual bondage for fourteen centuries".

**14.** In accordance with the Talmud (Kelim 17$^6$), the volume of a "standard" hen's egg, which served as the unit of volume, was defined as the mean of the "largest" and the "smallest" eggs [from a large batch]. The Talmud also provided elementary considerations about linear measurements and some stipulations regarding their admissible errors (Sheynin 1998b, p. 196).

**15.** Large payments were thus valued comparatively higher and this subjective attitude can also be traced in later lotteries up to our days: although large winnings are unrealistic, gamblers are apt to hope for them. The expectations of the various winnings in the Genoese lottery, which had been carried out from the mid-15$^{th}$ century, confirm the conclusion made above: they decreased with the increase in the theoretically possible gain (**N. Bernoulli** 1709/1975, p. 321; Biermann 1957; Bellhouse 1981). Recent experiments (Cohen et al 1970, 1971) seem to suggest that psychological subjective probabilities differ from objective statistical probabilities.

An embryo of the notion of expectation might be seen in the administration of justice in 11$^{th}$-century India (**Al-Biruni** 1887, vol. 2, pp. 158 – 160). Thus,

*If the suitor is not able to prove his claim, the defendant must swear. […] There are many kinds of oath, in accordance with the value of the object o the claim.*

The oath apparently became ever more earnest as that value increased; the probability of lying with impunity multiplied by the value in question was the expectation of fraudulent gain. However, when

*the object of claim was of some importance, the accused man was invited to drink some kind of a liquid which in case he spoke the truth would do him no harm.*

**16.** At least once Ptolemy (1984, III 4, p. 153; H 232) stated that a simpler hypothesis "would seem more reasonable".

**17.** Hardly anyone later recalled that doubtful example; **Lamarck** (1815, p. 133), however, believed that there existed deviations from divine design in the tree of animal life and explained them by "une … cause accidentelle et par consequent variable". In 1629, William Ames, a theologian, stated that random events might occur even with probability $p$   1/2, see Bellhouse (1988, p. 71) who does not elaborate and provides no exact reference.

**18.** The *Laws of Manu* (Ibidem) and the ancient Chinese literature (Burov et al 1972, p. 108) contain examples of decisions based on elementary stochastic considerations, e.g., accept as true the statement of the majority.

**19.** I am using the **Gauss** apt notation

$[ab] = a_1b_1 + a_2b_2 + … + a_nb_n.$

**20.** But how did Tycho treat his observations when one or even more of his instruments had to be temporarily taken out of service? Its (their) removal from the pool could have led to a systematic shift in the mean measurement.

**21.** Possibly somewhat exaggerating, Rabinovitch (1974, p. 355), who described the legal problems and rituals of the Judaic religion, concluded that these propositions (not formulated by **Ptolemy**) were known even in antiquity.

**22.** Chance also began to be recognized in biology (and even **Aristotle** thought that monstrosities were random, see § 1.1.1). **Harvey** (1651/1952, p. 338) stated that spontaneous generation (then generally believed in) occurred accidentally, as though aimlessly, again see § 1.1.1:

*Creatures that arise spontaneously are called automatic […] because they have their origin from accident, the spontaneous act of nature.*

I would say that Harvey considered randomness an intrinsic feature of nature.
**Lamarck** (1809, p. 62), also see Sheynin (1980, p. 338), kept to the same opinion. Harvey (1651/1952, p. 462) also believed that the form of hen's eggs was "a mere accident" and thus indicated an example of intraspecific variation and Lamarck (1817, p. 450) explained such variations by accidental causes.

**23. Ptolemy's** source was devoted to astrology. Thorndike (1923, p. 115) testifies that in the Middle Ages the *Tetrabiblos* exerted "great influence".

**24.** Ekeland (2006) provides pictures of such clouds. Chaos is one of the main subjects of his book, but its description is not good enough, and he (p. 125) even compares it with a game of chance. Randomness, which is another of his subjects, is treated hardly better (p. 86): reality "lies somewhere between" order and dependence of everything on everything else. There are other defects as well, and one of them is an extremely bad documentation.

**25.** Astronomers certainly applied the arithmetic mean even before **Kepler** did. **Tycho** (Plackett 1958/1970, pp. 122 – 123) combined 24 of his observations into (12) pairs and calculated the (generalized arithmetic) mean of these pairs, and of three separate observations assigning equal weight to each of the 15 values thus obtained. He chose the pairs in a manner allowing the elimination of the main systematic errors and, apparently, so as to estimate, even if qualitatively, the influence of random errors in 12 cases out of the 15. The separate observations could have been somehow corrected previous to the adjustment. I shall describe a similar case in § 6.3.2.

The mean occurs in early approximate calculations of areas of figures and volumes of bodies so as to compensate the errors of approximate formulas and/or the deviations of the real figures and bodies from their accepted models (Colebrooke 1817, p. 97).

# 2. The Early History
## 2.1. Stochastic Ideas in Science and Society

**2.1.1. Games of Chance.** They fostered the understanding of the part of chance in life and, even in antiquity, illustrated practically impossible events (§ 1.1.1) whereas mathematicians discovered that such games provided formulations of essentially new problems. Furthermore, although **Pascal** did not apply his relevant studies to any other domain, he (1654b/1998, p. 172) had time to suggest a remarkable term for the nascent theory, – *Aleae geometria, La Géométrie du hazard,* – and to indicate his desire to compile a pertinent tract. Later **Huygens** (1657/1920) introduced his study of games of chance by prophetically remarking that it was not a simple "jeu d'esprit" and that it laid the foundation "d'une spéculation fort intéressante et profonde". **Leibniz** (1704/1996, p. 506) noted that he had repeatedly advocated the creation of a "new type of logic" so as to study "the degrees of probability" and recommended, in this connection, to examine all kinds of games of chance[1]. In 1703 he wrote to **Jakob Bernoulli** (Kohli 1975b, p. 509):

*I would like that someone mathematically studies different games (in which excellent examples of* [the doctrine of estimating probabilities] *occur). This would be both pleasant and useful and not unworthy either of you or of another respected mathematician.*

Indeed, games of chance, and possibly only they, could have at the time provided models for posing natural and properly formulated stochastic problems. And, in addition, they were in the social order of the day.

**Rényi** (1969) attempted to conjecture the essence of **Pascal's** proposed tract. He could have been right in suggesting its subject-matter but not with regard to the year, – 1654, – when Pascal or rather Rényi described it. Another shortcoming of Rényi's attempt is that in spite of treating philosophical issues, Pascal (again, Rényi) had not mentioned **Aristotle**.

The theory of probability had originated in the mid-17[th] century rather than earlier; indeed, exactly then influential scientific societies came into being and scientific correspondence became usual. In addition, during many centuries games of chance had not been sufficiently conducive to the development of stochastic ideas (**M. G. Kendall** 1956/1970, p. 30). The main obstacles were the absence of "combinatorial ideas" and of the notion of chance events, superstition, and moral or religious "barriers" to the development of stochastic ideas, see also § 6.2.3. In essence, combinatorial analysis dates back to the 16[th] century although already **Levi ben Gerson** (Rabinovitch 1973, pp. 147 – 148) had created its elements.

**Montmort** (1713/1980, p. 6) had testified to the superstition of gamblers; **Laplace** (1814/1995, pp. 92 – 93) and **Poisson** (1837a, pp. 69 – 70) repeated his statement (and adduced new examples). When a number has not been drawn for a long time in the French lottery, **Laplace** says, "the mob is eager to bet on it" and he adds that an opposite trend is also noticeable. The same illusions exist in our time although **Bertrand** (1888a, p. XXII) had convincingly remarked that the roulette wheel had "ni conscience ni mémoire". Even a just game (with a zero expectation of loss) is ruinous (§ 6.1.1) and is therefore based on superstition while lotteries are much more harmful. Already **Petty** (1662/1899, vol. 1, p. 64) stated that they were "properly a Tax upon unfortunate self-conceited fools" and **Arnauld** & **Nicole** (1662/1992, p. 332) indicated that an expectation of a large winning in a lottery was illusory. In essence, they came out against hoping for unlikely favourable events, cf. § 1.1.5, Note 15.

**2.1.2. Jurisprudence.** I mentioned it in §§ 1.1.1 and 1.1.5 and, in particular, I noted that one of the first tests for separating chance from necessity was provided for the administration of justice. It seems, however, that the importance of civil suits and stochastic ideas in law courts increased exactly in the mid-17$^{th}$ century[2]. **Descartes** (1644/1978, pt. 4, § 205, p. 323) put moral certainty into scientific circulation, above all apparently bearing in mind jurisprudence[3]. **Arnauld** & **Nicole** (1662/1992) mentioned it; **Jakob Bernoulli** (1713a/1999, pt. 4, Chapter 3) later repeated one of their examples (from their pt. 4, Chapter 15). True, **Leibniz** (§3.1.2) doubted that observations might lead to it.

In the beginning of the 18$^{th}$ century, **Niklaus Bernoulli** (§ 3.3.2) devoted his dissertation to the application of the "art of conjecturing" to jurisprudence. Elements of numerical measures in that science are known even in the Roman canon law: it "had an elaborate system of full proofs, half proofs, and quarter proofs" (Garber & Zabell 1979, p. 51, Note 23). On the Roman code of civil law see also Note 6 in Chapter 1 and Franklin (2001, p. 211).

**Leibniz** (1704/1996, pp. 504 – 505) mentioned degrees of proofs and doubts in law and in medicine and indicated that

*Our peasants have since long ago been assuming that the value of a plot is the arithmetic mean of its estimates made by three groups of appraisers*[4].

**2.1.3. Insurance of Property and Life Insurance.** Marine insurance was the first essential type of insurance of property but it lacked stochastic ideas or methods. In particular, there existed an immoral and repeatedly prohibited practice of betting on the safe arrivals of ships. Anyway, marine insurance had been apparently based on rude and subjective estimates. Chaufton (1884, p. 349) maintained that in the Middle Ages definite values were assigned to risks in marine operations but he possibly meant just such estimates.

And here is a quote from the first English Statute on assurance (Publicke Acte No. 12, 1601; *Statutes of the Realm*, vol. 4, pt. 2, pp. 978 – 979):

*And whereas it hathe bene tyme out of mynde an usage amongste merchants, both of this realme and of forraine nacyons, when they make any great adventure,* […] *to give some consideracion of money to other persons* […] *to have from them assurance made of their goodes, merchandizes, ships, and things adventured,* […] *whiche course of dealinge is commonly termed a policie of assurance* […].

Now, a lucid (regrettably, unsubstantiated) statement about the indirect importance of marine insurance (O'Donnell (1936, p. 78):

*Wide research concedes that Life Insurance came into its own not by a front-door entrance, but by marine insurance porthole.*

Life insurance exists in two main forms. Either the insurer pays the policy-holder or his heirs the stipulated sum on the occurrence of an event dependent on human life; or, the latter enjoys a life annuity. Life insurance in the form of annuities was known in Europe from the $13^{th}$ century onward although later it was prohibited for about a century until 1423 when a Papal bull officially allowed it (Du Pasquier 1910, pp. 484 – 485). The annuitant's age was not usually taken into consideration either in the mid-$17^{th}$ century (Hendriks 1853, p. 112), or even, in England, during the reign of William III [1689 – 1702] (**K. Pearson** 1978, p. 134). Otherwise, as it seems, the ages had been allowed for only in a very generalized way (Sheynin 1977b, pp. 206 – 212; Kohli & **van der Waerden** 1975, pp. 515 – 517; **Hald** 1990, p. 119). It is therefore hardly appropriate to mention expectation here, but at the end of the $17^{th}$ century the situation began to change.

It is important to note that in the $18^{th}$, and even in the mid-$19^{th}$ century, life insurance therefore hardly essentially depended on stochastic considerations[5]; moreover, the statistical data collected by the insurance societies as well as their methods of calculations remained secret. A special point is that more or less honest business based on statistics of mortality hardly superseded downright cheating before the second half of the $19^{th}$ century. Nevertheless, beginning at least from the $18^{th}$ century, the institute of life insurance strongly influenced the theory of probability, see §§ 4.2 and 6.1.1c.

I single out the *Memorandum* of **De Witt** (1671). He separated four age groups (5 – 53; 53 – 63; 63 – 73; and 73 – 80 years) and assumed that the chances of death increased in a definite way from one group to the next one but remained constant within each of them. According to his calculations, the cost of an annuity for "young" men should have been 16 times higher than the yearly premium (not 14, as it was thought).

Eneström (1896/1897, p. 66) called De Witt's exposition unclear. Thus, the risk of dying always concerned a three year old infant, which was not only not explained, but expressed in a misleading way; and the proposed chances of death at various ages were contrary to what was actually calculated.

An appendix to the main text (Hendriks 1853, pp. 117 – 118) contains an interesting observation. Examining "considerably more than a hundred different classes, each class consisting of about one hundred persons", De Witt found that

*For young lives each of these classes always produced to the annuitants […] a value of more than sixteen florins of capital arising from one florin of annual rent. […] Thus […] in practice, when the purchaser of several life annuities comes to divide his capital […] upon several young lives – upon ten, twenty, or more – this annuitant*

*may be assured, without hazard or risk of the enjoyment of an equivalent, in more than sixteen times the rent which he purchases*.

This statement belongs to the prehistory of the LLN and it also shows a peculiar aspect of business in those times.

In the same year, in a letter to another mathematician, **Hudde**, De Witt (Hendriks 1853, p. 109) in an elementary way calculated the cost of annuity on several lives (an annuity that should be paid out until the death of the last person of the group; usually, of a married couple). In the process, he determined the distribution of the maximal term of a series of observations [obeying a uniform law]. Kohli & **van der Waerden** (1975) described in detail the history of the institution of life insurance including the work of De Witt and **Huygens** (§ 2.2.2), and I only note that the former had not justified his assumed law of mortality. A likely corollary of De Witt's work was that the price of annuities sold in Holland in 1672 – 1673 depended on the age of the annuitants (Commelin 1693, p. 1205).

The first estimation of the present worth of life annuities, based on a table of expectations of life, was made by the Praetorian Prefect **Ulpianus**, see § 1.1.5. **Leibniz** (1986, pp. 421 – 432), also see Leibniz (2000), in a manuscript written in 1680, described his considerations about state insurance, see Sofonea (1957a). He had not studied insurance as such, but maintained that the "princes" should care about the poor, remarked that the society ought to be anxious for each individual etc. Much later **Süssmilch** (§ 6.2.2) formulated similar ideas.

Tontines constituted a special form of mutual insurance. Named after the Italian banker Laurens Tonti (Hendriks 1863), they, acting as a single body, distributed the total sums of annuities among their members still alive, so that those, who lived longer, received considerable moneys. Tontines were neither socially accepted nor widespread "on the assumed rationale that they are too selfish and speculative" (Hendriks 1853, p. 116). Nevertheless, they did exist in the 17th century, see Sheynin (1977b, pp. 210 – 211).

**2.1.4. Population Statistics.** The Old Testament (Numbers, Chapter 1) reports on a general census, or, more precisely, on a census of those able to bear arms. To recall (§ 1.1.2), the Talmud estimated the population of towns only by the number of soldiers "brought forth" [when needed]. In China, in 2238 BC or thereabouts, an estimation of the population was attempted and the first census of the warrior caste in Egypt occurred not later than in the 16th century BC (Fedorovitch 1894, pp. 7 – 21). In Europe, even in 15th century Italy, for all its achievements in accountancy and mathematics (**M. G. Kendall** 1960),

*counting was by complete enumeration and still tended to be a record of a situation rather than a basis for estimation or prediction in an expanding economy.*

Only **Graunt** (1662) and, to a lesser extent, **Petty** (1690) can be called the fathers of population statistics. They studied population, economics, and commerce and discussed the appropriate causes and connections by means of elementary stochastic considerations, also see Urlanis (1963) and **K. Pearson** (1978, pp. 30 – 49). It was Petty who called the new discipline *political arithmetic*. He (Petty 1690/1899, vol. 2, p. 244) plainly formulated his denial of "comparative and superlative Words" and attempted to express himself in "Terms of Number, Weight, or Measure…"; Graunt undoubtedly did, if not said the same.

Petty (1927, vol. 1, pp. 171 – 172) even proposed to establish a "register generall of people, plantations & trade of England". In particular, he thought of collecting the accounts of all the

*Births, Mariages, Burialls* […] *of the Herths, and Houses* […] *as also of the People, by their Age, Sex, Trade, Titles, and Office.*

As noticed by Greenwood (1941 – 1943/1970, p. 61), the scope of this Register was to be "wider than that of our existing General Register Office".

Strictly speaking, neither Petty, nor, as it seems, his followers ever introduced any definition of political arithmetic. However, without violating Petty's thoughts quoted above, it is possible to say that the aims of this new scientific discipline were to study from a socio-economic point of view states and separate cities (or regions) by means of (rather unreliable) statistical data on population, industry, agriculture, commerce etc.

At least 30 from among Petty's manuscripts (1927) pertained to political arithmetic. This source shows him as a philosopher of science congenial in some respects with **Leibniz**, his younger contemporary. I adduce one passage (Ibidem, pp. 39 – 40); also see Sheynin (1977b, pp. 218 – 220):

*What is a common measure of Time, Space, Weight, & motion? What number of Elementall sounds or letters, will […] make a speech or language? How to give names to names, and how to adde and subtract sensata, & to ballance the weight and power of words; which is Logick & reason.*

**Graunt** (1662) studied the weekly bills of mortality in London which began to appear in the 16[th] century and had been regularly published since the beginning of the 17[th] century. For a long time his contribution had been attributed to **Petty**. However, according to Hull (Petty 1899, vol. 1, p. lii), Petty

*Perhaps suggested the subject of inquiry*, […] *probably assisted with comments upon medical and other questions here and there* […], *procured* [some] *figures* […] *and may have revised, or even written the Conclusion…*

If so, Petty still perhaps qualifies as co-author, but I shall not mention him anymore. And I also quote his *Discourse* (1674):

> ***I have also (like the author of those Observations** [like Graunt]*) ***Dedicated this Discourse to [...] the Duke of Newcastle***.

And so, Graunt was able to use the existing fragmentary statistical data and estimated the population of London and England as well as the influence of various diseases on mortality. His main merit consisted in that he attempted to find definite regularities in the movement of the population. I only indicate that he established that both sexes were approximately equally numerous (which contradicted the then established views) and that out of 27 newly born about 14 were boys. When dealing with large numbers, Graunt did not doubt that his conclusions reflected objective reality which might be seen as a fact belonging to the prehistory of the LLN. Thus, the statistical ratio 14:13 was, in his opinion, an estimate of the ratio of the respective [probabilities].

Nevertheless, he uncritically made conclusions based on a small number of observations as well and thought that the population increased in an arithmetical progression. Many authors had since replaced that law by the geometrical progression definitely introduced by Süssmilch and Euler (§ 6.2.2).

In spite of the meagre and sometimes wrong information provided in the bills about the age of those dying, Graunt was able to compile the first mortality table (common for both sexes). He calculated the relative number of people dying within the first six years and within the next decades up to age 86. According to his table, only one person out of a hundred survived until that age. Now, how exactly had Graunt calculated his table? Opinions vary, but, in any case, the very invention of the mortality table was the main point here. The indicated causes of death were also incomplete and doubtful, but Graunt formulated some important relevant conclusions as well (although not without serious errors)[6]. His general methodological (but not factual) mistake consisted in that he assumed, without due justification, that statistical ratios during usual years (for example, the per cent of yearly deaths) were stable.

It is generally known that Graunt had essentially influenced later scholars and here are a few pertinent opinions (**Huygens**, letter of 1662/1888 – 1950, 1891, p. 149; **Süssmilch**, see **K. Pearson** 1978, pp. 317 – 318; Willcox, see Graunt 1662/1939, p. x; **Hald** 1990, p. 86).

1. *Grant's* [!] *discourse really deserves to be considered and I like it very much. He reasons sensibly and clearly and I admire how he was able to elicit all his conclusions from these simple observations which formerly seemed useless.*

2. The discovery of statistical regularities *was even as feasible as that as America, it wanted only its Columbus…*

3. *Graunt is memorable mainly because he discovered [...] the uniformity and predictability of many biological phenomena taken in the mass [...] thus he, more than any* [other] *man, was the founder of statistics.*

4. Graunt, in his own words, reduced the data from "several great confused Volumes into a few perspicuous Tables" and analysed them in " a few succinct Paragraphs" which is exactly the aim of statistics. Hald could have referred here to **Kolmogorov** & **Prokhorov**, see § 0.2.

**Halley** (1693a; 1693b), a versatile scholar and an astronomer in the first place, compiled the next mortality table. He made use of statistical data collected in Breslau[7], a city with a closed population. Halley applied his table for elementary stochastic calculations connected with life insurance and he was also able to find out the general relative population of the city. Thus, for each thousand infants aged less than a year, there were 855 children from one to two years of age, …, and, finally, 107 persons aged 84 – 100. After summing up all these numbers, Halley obtained 34 thousand (exactly) so that the ratio of the population to the newly born occurred to be 34. Until 1750 his table remained the best one (**K. Pearson** 1978, p. 206).

For a long time, the explanation of the calculations in his main (first) paper had not been properly understood. Then, the yearly rate of mortality in Breslau was 1/30, the same as in London, and yet Halley considered that city as a statistical standard. If only such a notion is appropriate, standards of several levels ought to be introduced. Finally, Halley thought that the irregularities in his data

*Would rectify themselves, were the number of years* [of observation] *much more considerable*.

Such irregularities could have well been produced by systematic influences, but, anyway, Halley's opinion shows the apparently widespread belief in an embryo of the LLN.

Halley's second note is indeed interesting as a reasoning on the welfare of the population.

In spite of the deficiencies mentioned, great success came immediately. K. Pearson (1978, p. 78) indicated that Halley had made "all the use that a modern actuary could" of his data and that he had computed his life-table "as we should do it today". Sofonea (1957b, p. 31*) called Halley's contribution "the beginning of the entire development of modern methods of life insurance", and Hald (1990, p. 141) stated that it "became of great importance to actuarial science".

Later scholars applied Halley's contribution; the most important example is De Moivre's introduction (1725/1756/1967) of the uniform law of mortality for ages beginning at 12 years.

In 1701 Halley (Chapman 1941, p. 5) compiled a chart of Northern Atlantic showing the lines of equal magnetic declinations so that he (and of course **Graunt**) might be called the founders of exploratory data analysis, a most important, even if elementary stage of statistical investigations, cf. § 10.8.3.

In 1680 – 1683 **Leibniz** wrote several manuscripts mostly pertaining to the so-called statecraft (§ 6.2.1) and first published in 1866 (Leibniz 1986, pp. 340 – 349, 370 – 381, and 487 – 491), see also Sheynin (1977b, pp. 224 – 227). He recommended the publication of "state tables" (numerical or not?) of remarkable facts and their comparison, year with year, or state with state. Their compilation, as he suggested, should have been the duty of special recording offices and, as it seems, for such offices Leibniz (disorderly) listed 56 questions from which I mention the number of inhabitants of a state and the comparison of the birth rate and mortality. Then, he thought it advisable to collect information about scientific achievements, "clever ideas" and medical and meteorological observations, and to establish "sanitary boards" for compiling data on a wide range of subjects (meteorology, medicine, agriculture).

One of Leibniz' manuscripts (Leibniz 1986, pp. 456 – 467, or, with a German translation, 2000, pp. 428 – 445) was devoted to political arithmetic. There, he introduced the *moyenne longueur de la vie humaine*[8], necessary, as he remarked, for calculating the cost of annuities; assumed without substantiation several regularities, for example, that the ratio of mortality to population was equal to 1:40; and wrongly stated that the mortality law for each age group including infants was uniform. Following **Arnauld** & **Nicole** (1662/1992, pp. 331 and 332), he discussed *apparence* or *degré de la probabilité* and *apparence moyenne* [expectation]. See §6.2.1 about the history of university statistics.

Leibniz began by stating that in a throw of two dice the occurrence of seven points is three times (actually six times) more probable than of twelve points. Todhunter (1865, p. 48) indicated his second mistake of the same kind. And Leibniz ended his essay by arguing that, since hardly 1/9 or 1/10 of fertile women conceive yearly, the birth-rate can be nine or ten times higher than it actually is. I simply do not understand this argument.

Population statistics owed its later development to the general problem of isolating randomness from Divine design. **Kepler** and **Newton** achieved this aim with regard to inanimate nature, and scientists were quick to begin searching for the laws governing the movement of population, cf. **K. Pearson's** appropriate remark in § 2.2.3.

## 2.2. Mathematical Investigations

**2.2.1. Pascal and Fermat.** In 1654 Pascal and Fermat exchanged several letters (Pascal 1654a) which heralded the beginning of the formal history of probability. They discussed several problems; the most important of them was known even at the end of the 14$^{th}$ century. Here it is: Two or three gamblers agree to continue playing until one of them scores $n$ points; for some reason the game is, however, ínterrupted on score $a:b$ or $a:b:c$ ($a$, $b$, $c < n$) and it is required to divide the stakes in a reasonable way[9]. Both scholars solved this problem (the problem of points; see Takácz 1994) by issuing from one and the same rule: the winnings of the gamblers should be in the same ratio(s) as existed between the expectations of their scoring the $n$ points, see for example Sheynin (1977b, pp. 231 – 239). The actual introduction of that notion, expectation, was their main achievement. They also effectively applied the addition and the multiplication theorems[10].

The methods used by Pascal and Fermat differed. In particular, Pascal solved the above problem by means of the arithmetic triangle composed, as is well known, of binomial coefficients of the development $(1 + 1)^n$ for increasing values of *n*. Pascal's relevant contribution (1665) was published posthumously, but Fermat was at least partly familiar with it. Both there, and in his letters to Fermat, Pascal in actual fact made use of partial difference equations (**Hald** 1990, pp. 49 and 57).

The celebrated **Pascal** wager (1669/2000, pp. 676 – 681), also published posthumously, was in essence a discussion about choosing a hypothesis. Does God exist, rhetorically asked the devoutly religious author and answered: you should bet. If He does not exist, you may live calmly [and sin]; otherwise, however, you can lose eternity. In the mathematical sense, Pascal's reasoning[11] is vague; perhaps he had no time to edit his fragment. Its meaning is, however, clear: if God exists with a fixed and however low probability, the expectation of the benefit accrued by believing in Him is infinite.

**2.2.2. Huygens.** Huygens was the author of the first treatise on probability (1657). Being acquainted only with the general contents of the **Pascal – Fermat** correspondence, he independently introduced the notion of expected random winning and, like those scholars, selected it as the test for solving stochastic problems. Note that he went on to prove that the "value of expectation", as he called it, of a gambler who gets *a* in *p* cases and *b* in *q* cases was

$$\frac{pa + qb}{p + q}. \qquad (1)$$

**Jakob Bernoulli** (1713a/1999, p. 9) justified the expression (1) much simpler than Huygens: if each of the *p* gamblers gets *a*, and each of the *q* others receives *b*, and the gains of all of them are the same, then the expectation of each is equal to (1). After Bernoulli, however, expectation began to be introduced formally: expressions of the type of (1) followed by definition.

**Huygens** solved the problem of points under various initial conditions and listed five additional problems two of which were due to **Fermat**, and one, to **Pascal**. He solved them later, either in his correspondence, or in manuscripts published posthumously. In essence, they demanded the use of the addition and multiplication theorems, the actual introduction of conditional probabilities and the formula (in modern notation)

$P(B) = P(A_i)P(B/A_i), i = 1, 2, …, n.$

I describe two of the five additional problems. Problem No. 4 was about sampling without replacement. An urn contained 8 black balls and 4 white ones and it was required to determine the ratio of chances that in a sample of 7 balls 3 were, or were not white. Huygens determined the expectation of the former event by means of a partial difference equation (**Hald** 1990, p. 76), cf. Korteweg's remark about Huygens' analytical methods below. Nowadays such problems leading to the hypergeometric distribution (**J. Bernoulli** 1713a/1999, pp. 167 – 168; **De Moivre** 1712/1984, Problem 14 and 1718/1756, reprint 1967, Problem 20) appear in connection with statistical inspection of mass production.

    **Pascal's** elementary Problem No. 5 was the first to discuss the gambler's ruin. Gamblers *A* and *B* undertake to score 14 and 11 points respectively in a throw of 3 dice. They have 12 counters each and it is required to determine the ratio of the chances that they be ruined. The stipulated numbers of points occur in 15 and 27 cases and the ratio sought is therefore $(5/9)^{12}$.

    In 1669, in a correspondence with his brother Lodewijk, **Huygens** (1895), see Kohli & **van der Waerden** (1975), discussed stochastic problems connected with mortality and, to be sure, life insurance. So it happened that the not yet formed theory of probability spread over new grounds. Issuing from **Graunt's** mortality table (§ 2.1.4), Huygens (pp. 531 – 532) introduced the probable duration of life (but not the term itself) and explained that it did not coincide with expected life. On p. 537 he specified that the latter ought to be used in calculations of annuities and the former for betting on human lives. Indeed, both he (pp. 524 – 526) and Lodewijk (pp. 484 – 485) mentioned such betting. Christiaan also showed that the probable duration of life could be determined by means of the graph (a continuous curve passing through empirical points given by **Graunt's** table of mortality; plate between pp. 530 and 531) of the function

$$y = 1 - F(x),$$

where, in modern notation, $F(x)$ was a remaining unknown integral distribution function with admissible values of the argument being $0 \leq x \leq 100$.

    Also in the same correspondence Huygens (p. 528) examined the expected period of time during which 40 persons aged 46 will die out; and 2 persons aged 16 will both die. The first problem proved too difficult, but Huygens might have remarked that the period sought was 40 years (according to **Graunt**, 86 years was the highest possible age). True, he solved a similar problem but made a mistake. He assumed there that the law of mortality was uniform and that the number of deaths will decrease with time, but for a distribution, continuous and uniform in some interval, *n* order statistics will divide it into ($n$ + 1) approximately equal parts and the annual deaths will remain about constant. In the second problem Huygens applied conditional expectation when assuming that one of the two persons will die first. Huygens never mentioned **De Witt** (§ 2.1.3) whose work (an official and classified document) had possibly been remaining unknown[12].

    When solving problems on games of chance, Huygens issued from expectations which varied from set to set rather than from constant probabilities. He was therefore compelled to compose and solve difference equations (Korteweg, see Huygens 1888 – 1950, 1920, p. 135) and he (like **Pascal**, see § 2.2.1) should be recalled in connection with their history. See also Shoesmith (1986).

While developing the ideas of **Descartes** and other scholars about moral certainty (§ 2.1.2), **Huygens** maintained that proofs in physics were only probable and should be checked by appropriate corollaries and that common sense should determine the required degree of certainty of judgements in civil life. In a letter of 1691 Huygens (1888 – 1950, t. 10, p. 739) had indeed mentioned Descartes and, without justification, dismissed probabilities of the order of $p = 10^{-11}$ although he hardly applied this, or any other number as a criterion. Note that **Borel** (1943/1962, p. 27) proposed $p = 10^{-6}$ and $10^{-15}$ as insignificant *on the human* and the *terrestrial scale* respectively. However, a single value is hardly suitable for any arbitrary event. Also see Sheynin (1977b, pp. 251 – 252).

**2.2.3. Newton.** Newton left interesting ideas and findings pertaining to probability (Sheynin 1971a), but much more important were his philosophical views. Here is the opinion of **K. Pearson** (1926):

*Newton's idea of an omnipresent activating deity, who maintains mean statistical values, formed the foundation of statistical development through Derham* [a religious philosopher], *Süssmilch* [§ 6.2.2], *Niewentyt* [a statistician], *Price* [Chapt. 5] *to Quetelet* [§ 10.5] *and Florence Nightingale*[13]…

Newton had not stated such an idea (although he thought that God regularly delivered the system of the world from accumulating corruptions, see below). In 1971 **E. S. Pearson** answered my question on this point:

*From reading* [the manuscript of K. Pearson (1978)] *I think I understand what K. P. meant* […]. *He has stepped ahead of where Newton had got to, by stating that the laws which give evidence of Design, appear in the stability of the mean values of observations …*

I have since found that **K. Pearson** (1978, pp. 161 and 653) had attributed to **De Moivre** (1733/1756, reprint 1967, pp. 251 – 252) the Divine *stability of statistical ratios*, *that is, the original determination or original design* and referred to **Laplace** who (1814/1995, p. 37) had indeed formulated a related idea:

*In an indefinitely continued sequence of events, the action of regular and constant causes ought, in the long run, to outweigh that of irregular causes.*

However, as I also note in § 7.1-3, Laplace never mentioned Divine design.

K. Pearson (1926) then went over to **De Moivre** (§ 4.4) and **Bayes** (Chapter 5) and maintained that their work was motivated by theological and sociological causes rather than by mathematics itself.

And here is **Newton's** most interesting pronouncement (1704/1782, Query 31):

*Blind fate could never make all the planets move one and the same way in orbs concentrick, some inconsiderable irregularities excepted, which may have risen from the mutual actions of comets and planets upon one another, and which will be apt to increase, till this system wants a reformation. Such a wonderful uniformity in the planetary system must be allowed the effect of choice. And so must the uniformity in the bodies of animals.*

I have indicated that such considerations are logically imperfect but practically certain (Note 5 in Chapter 1). The idea of a divine reformation of the system of the world was later abandoned, but Newton's recognition of the existence and role of its random disturbances is very important. Random, I specify, in the same sense as the outcome of coin-tossing is. But at the same time Newton (1958, pp. 316 – 318), just like **Kepler** (§ 1.2.4), denied randomness and explained it by ignorance of causes. It was the future theologian **Bentley**, who, in 1693, expressed his thoughts after discussing them with Newton. The texts of two of his sermons, of Newton's letters to him, and an article on Newton and Bentley are in Newton (1958). Finally, Newton's remark (Schell 1960), that the notion of chance might be applied to a single trial, has a philosophical side[14].

When studying the chronology of ancient kingdoms, **Newton** (1728, p. 52) left an interesting statement:

*The Greek Chronologers […] have made the kings of their several Cities […] to reign about 35 or 40 years a-piece, one with another; which is a length so much beyond the course of nature, as is not to be credited. For by the ordinary course of nature Kings Reign, one with another, about 18 or 20 years a-piece; and if in some instances they Reign, one with another, five or six years longer, in others they reign as much shorter: 18 or 20 years is a medium.*

Newton derived his own estimate from other chronological data and his rejection of the twice longer period was reasonable. Nevertheless, a formalized reconstruction of his decision is difficult: within one and the same dynasty the period of reign of a given king directly depends on that of his predecessor. Furthermore, it is impossible to determine the probability of a large deviation of the value of a random variable from its expectation without knowing the appropriate variance (which Newton estimated only indirectly and in a generalized way). **K. Pearson** (1928a) described Newton's later indication of the sources of his estimate and dwelt on **Voltaire's** adjoining remarks, and, especially, on the relevant work of **Condorcet**.

I am now mentioning Newton's manuscript (1967, pp. 58 – 61) written sometime between 1664 and 1666. This is what he wrote:

*If the Proportion of the chances […] bee irrational, the interest [expectation] may bee found after ye same manner.*

He thought of a ball falling upon the centre of a circle divided into sectors whose areas were in "such proportion as 2 to √5". If the ball

"tumbles" into the first sector, a person gets *a*, otherwise he receives *b*, and his "hopes is worth"

(2*a* + *b* 5) ÷ (2 + 5).

This was a generalization of expectation as defined by **Huygens** (§ 2.2.2) and the first occurrence of geometric probability (§ 6.1.6). Newton's second example was a throw of an irregular die. He remarked that [nevertheless] "it may bee found how much one cast is more easily got than another". He hardly thought about analytic calculations, he likely bore in mind statistical probabilities. I can only add that Newton may well have seen **Graunt's** contribution (§ 2.1.4).

In 1693, when answering a question, Newton (Gani 1982) determined the [probability] of throwing not less than one, two, and three sixes with six, 12 and 18 dice respectively (cf. the De Méré problem in Note 10). In the last-mentioned case, for example, his calculations can be described by the formula

$$P = 1 - C_{18}^2 (5/6)^{16} (1/6)^2 - C_{18}^1 (5/6)^{17} (1/6) - (5/6)^{18}.$$

**2.2.4. Arbuthnot.** Arbuthnot (1712) collected the data on births (more precisely, on baptisms) in London during 1629 – 1710. He noted that during those 82 years more boys (*m*) were invariably born than girls (*f*) and declared that that fact was *not the Effect of Chance but Divine Providence, working for a good End*. Indeed, as he added, boys and men were subject to greater dangers and their mortality was higher than that of the females, "as Experience convinces us". Even disregarding such [hardly exhibited] regularities as the "constant Proportion" *m:f* and "fix'd limits" of the difference (*m* – *f*), the "Value of Expectation" of a random occurrence of the observed inequality was less than $(1/2)^{82}$, he stated, see however Note 5 in Chapter 1.

Arbuthnot could have concluded that the births of both sexes obeyed [the binomial distribution], which, rather than the inequality *m > f*, manifested Divine design; and could have attempted to estimate its parameter (approximately equal to 14:13, see § 2.1.4). Then, he had not remarked that baptisms were not identical with births; that Christians perhaps somehow differed from other people and, again, that London was perhaps an exception. **Graunt** (1662, end of Chapter 3) had indeed stated that during 1650 – 1660 less than half of the general population had believed that baptism was necessary.

And Arbuthnot had not known the comparative mortality of the sexes. Nevertheless, later authors took note of his paper, continued to study the same ratio *m:f* and, by following this subject, made important stochastic findings (see especially § 4.4). **Freudenthal** (1961, p. xi) even called Arbuthnot the author of the first publication on mathematical statistics. From among many other recent commentators I name Shoesmith (1987) and H. A. David & Edwards (2001, pp. 9 – 11) and I note that Arbuthnot was the first to publish a trick equivalent to the application of a generating function of the binomial distribution although only for its particular case. **Jakob**

**Bernoulli** (§ 3.1.2) applied a generating function before Arbuthnot, but his book only appeared in 1713.

In 1715, 'sGravesande (**K. Pearson** 1978, pp. 301 – 303; **Hald** 1990, pp. 279 – 280) improved on Arbuthnot's reasoning and discussed it with **Niklaus Bernoulli**, cf. § 3.3.4.

Bellhouse (1989) described Arbuthnot's manuscript written likely in 1694. There, the author examined the game of dice, attempted to study chronology (two examples, cf. § 2.2.3) and to a certain extent anticipated his published note of 1712.

The study of the sex ratio at birth by most eminent scholars was based on Arbuthnot's data.

### Notes

**1.** He himself began studying games of chance in 1675 (Biermann 1955). See Gini (1946), Kohli (1975b) and Sylla (1998) for the main body of his correspondence with **Jakob Bernoulli** (with comments). Also published (in the original Latin) is the entire correspondence (**Leibniz** 1971, pp. 10 – 110; Weil et al 1993).

**2.** The volte-face of the public mood can be perceived when considering the problem of absentees about whom nothing is known. So as not to violate God's commandment [which one?], **Kepler** (1610/1941, p. 238) as an astrologer refused to state whether a particular absentee was alive or not. **Jakob Bernoulli** (1713/1999, p. 235), however, suggested to study, in such cases, the pertinent stochastic considerations as also did **Niklaus Bernoulli** (§ 3.3.2).

**3.** This notion was introduced about 1400 for the solution of ethical problems (Franklin 2001, p. 69). Recall also Thomas' "conjectural probability" in § 1.1.5.

**4.** In 1705 he repeated his statement about the appraisal of lots (end of § 2.1.2) in a letter to **Jakob Bernoulli** (Kohli 1975b, p. 512). Much earlier, he included it in a Latin manuscript of 1680 – 1683 (Leibniz 2000, pp. 428 – 445, with a German translation).

**5.** Several authors mentioned the practice of insuring a number of healthy infants possibly substantiated by intuitive stochastic reasoning, see § 3.2.3.

**6.** It might be thought that **Graunt** attempted to allow for systematic corruptions of the data. Thus, he reasonably supposed that the number of deaths from syphilis was essentially understated [not only because of the difficulty in diagnosing but also] out of ethical considerations.

**7.** He (as well as **Leibniz**) obtained them from Caspar Neumann, a *Magister der Philosophie* and Member of the *Societät der Wissenschaften* in Berlin. In a letter of 1692 Leibniz (1970, p. 279) stated that the data were interesting. On **Halley** see Böckh (1893), Chapman (1941) and Cook (1998).

**8.** Perhaps he was not acquainted with the correspondence of **Huygens** (§ 2.2.2).

**9.** About 1400 an anonymous Italian author (Franklin 2001, pp. 294 – 296) correctly solved this problem for the case of two gamblers, but had not sufficiently justified his solution and made a mistake in a more complicated instance. It is much more important to note, however, that he had not introduced the notion of expectation (cf. below).

**10.** The term *probability* does not appear in the extant part of the correspondence, only *chance* is applied. De Méré, a man of the world, had unintentionally initiated that correspondence by asking **Pascal** why the chances of two apparently equivalent events were different. An elementary calculation shows that either the gamblers were able to reveal a difference of probabilities equal to 0.0264, cf. § 1.2.3, or, as **Ore** (1960, pp. 411 – 412) and **van der Waerden** (1976) believed, De Méré was able to calculate the appropriate probabilities, – but still thought that achieving a six in four throws of a die and two sixes in 24 throws of two dice should have been equally probable since 24/36 = 4/6. Actually,

$$P_1 = 1 - (5/6)^4 \approx 0.5177, \quad P_2 = 1 - (35/36)^{24} \approx 0.4913.$$

Later **Jakob Bernoulli** (1713a/1999, p. 32) considered the same problem.

A queer episode concerning De Méré occurred in the 19<sup>th</sup> century. **Georg Cantor** mistakenly thought that by his conclusion the man of the world had wished to destroy science. Accordingly, he privately called **Kronecker** (who had denied the emerging set theory) "Herr von Méré" (Fraenkel 1930, p. 199).

**11.** For a similar reasoning see **Arnauld** & **Nicole** (1662/1992, p. 334).

**12.** During the last years of his life **Jakob Bernoulli** vainly asked **Leibniz** to procure for him a copy of **De Witt's** work.

**13.** I would have excluded **Niewentyt** from the **Pearsonian** chain; **Derham**, however, had been very influential.

**14.** Here is a similar statement formulated in the 4<sup>th</sup> century BC (Burov et al 1972, p. 203):

*Who even before battle gains victory by military estimations has many chances* […] *who has many chances gains victory*; *who has few chances does not gain victory*; *all the less he who has no chances at all.*

### 3. Jakob Bernoulli and the Law of Large Numbers

I consider Bernoulli's main work, the *Ars conjectandi* (AC)[1], published posthumously in Latin and touch on his Diary (*Meditationes*) for 1684 – 1690. Only the stochastic part of the latter was published together with the AC (both in their original Latin), with other materials and comments (J. Bernoulli 1975). I also discuss related topics and dwell on Bernoulli's contemporaries. The AC was translated into German and its separate parts have also appeared in other living languages. In 1913 **Uspensky** translated Part 4 into Russian (reprint: J. Bernoulli 1986). My own translation, or rather rendition of the same part has also appeared (J. Bernoulli 2005). An English translation of the entire *Ars* has been published; I reviewed it (Sheynin 2006a) and concluded that the *Cobbler should stick to his last*. In the sequel, I refer to the German translation of Bernoulli (1899, reprinted in 1999).

#### 3.1. Bernoulli's Works

**3.1.1. The Diary.** There, Bernoulli studied games of chance and the stochastic side of civil law. He (1975, p. 47) noted that the probability[2] of a visitation of the plague in a given year was equal to the ratio of the number of these visitations during a long period of time to the number of years in that period. I stress that Bernoulli thus applied the definition of probability of an event (of statistical probability!) rather than making use of chances. An interesting point in this connection is that he (p. 46, marginal note) wrote out the imprint of a review published in 1666 of **Graunt's** book (§ 2.1.4) which Bernoulli possibly had not seen; he had not referred to it either in the *Meditationes* itself or in the AC. But the most important part of his Diary is a (fragmentary) proof of the LLN. This fact means that Bernoulli proved it not later than in 1690.

**3.1.2. The *Art of Conjecturing* (1713a, b). Its Contents.** The last part of this book is entitled *The use and application of the previous doctrine to civil, moral and economic affairs* (J. Bernoulli 1713a/1999, p. 229) but nothing of the sort had appeared there[3]. Interesting problems are solved in parts 1 and 3 (the study of random sums for the uniform and the binomial distributions, a similar investigation of the sum of a random number of terms for a particular discrete distribution, a derivation of the distribution of the first order statistic for the discrete uniform distribution and the calculation of probabilities appearing in sampling without replacement). The author's analytical methods included combinatorial analysis and calculation of expectations of winning in each set of finite and infinite games and their subsequent summing.

Part 1 is a reprint of **Huygens'** tract (§ 2.2.2) including the solution of his five additional problems, one of which Bernoulli (1713a/1999, p. 167) had, however, carried over to Part 3, together with vast and valuable commentaries. Nevertheless, this form again testifies that he was unable to complete his contribution. Also in Part 1 Bernoulli (pp. 22 – 28), while considering a game of dice, compiled a table which enabled him to calculate the coefficients of $x^m$ in the development of $(x + x^2 + \ldots + x^6)^n$ for small natural values of $n$.

Part 2 did not bear on probability. It dealt with combinatorial analysis and it was there that the author introduced the *Bernoulli numbers.*

Part 4 contained the LLN. There also we find a not really formal "classical" definition of probability (a notion which he had not applied when formulating that law), a reasoning, in Chapter 2, on the aims of the art of conjecturing (determination, as precise as possible, of probabilities for choosing the best solutions of problems, apparently in civil life) and elements of stochastic logic[4].

Bernoulli likely considered the art of conjecturing as a mathematical discipline based on probability as a measure of certainty and on expectation and including (the not yet formally introduced) addition and multiplication theorems and crowned by the LLN.

Bernoulli informed **Leibniz** about the progress in his work in a letter of 3 Oct. 1703 (Kohli 1975b, p. 509). He was compiling it for many years with repeated interruptions caused by his "innate laziness" and worsening of health; the book still lacked its "most important part", the application of the art of conjecturing to civil life; nevertheless, he, Bernoulli, had already shown his brother [**Johann**] the solution of a "difficult problem, special in its own way" [§ 3.2.3], that justified the applications of the art of conjecturing.

Most important both in that letter and in the following correspondence of 1703 – 1705[5] (Ibidem, pp. 510 – 512) was the subject of statistical probabilities, also see §§ 3.2.2 – 3.2.3. Leibniz never agreed that observations could secure moral certainty but his arguments were hardly convincing. Thus, he in essence repeated the statement of **Arnauld** & **Nicole** (1662/1992, pp. 304 and 317) that the finite (the mind; therefore, observations) cannot always grasp the infinite (for example, God, but also, as Leibniz stated, any phenomenon depending on innumerable circumstances).

These views were possibly caused by his understanding of randomness as something "whose complete proof exceeds any human mind" (Leibniz 1686/1960, p. 288). His heuristic statement does not contradict a modern approach to randomness founded on complexity and he was also right in the sense that statistical determinations cannot definitively corroborate a hypothesis.

In his letter of 3 Dec. 1703 **Leibniz** (Gini 1946, p. 405) had also maintained that the allowance for all the circumstances was more important than subtle calculations, and **Bortkiewicz** (1923, p. 12) put on record **Keynes'** favourable attitude towards this point of view and indicated the appropriate opinion of **Mill** (1843/1886, p. 353) who had sharply contrasted the consideration of circumstances with "elaborate application" of probability. Mill could have mentioned applied mathematics in general; the *circumstances* there are the extents to which its models are compatible with reality. Calculations (not necessarily stochastic) are hardly more important than such circumstances, cf. **Gauss'** opinion in §§ 9A and 9A.5-2. But circumstances should not be contrasted with calculations, the less so because they can initially remain insufficiently known.

Bernoulli paid due attention to Leibniz' criticism: more than a half of Chapter 4 of Part 4 of the AC in essence coincided with the respective passages from his letters to Leibniz; in that chapter, Bernoulli (1713a/1999, p. 250), in particular, discussed the objections made by "scientists", that is, by Leibniz[6].

### 3.2. The *Art of Conjecturing*, Part 4: Its Main Propositions

**3.2.1. Stochastic Assumptions and Arguments.** Bernoulli examined these in Chapters 2 and 3 but did not return to them anymore; he possibly thought of applying them in the unwritten pages of his book. The mathematical aspect of his considerations consisted in the use of the addition and the multiplication theorems for combining various arguments.

Unusual was the non-additivity of the [probabilities][7]. Here is one of his examples (p. 244): "something" possesses 2/3 of certainty but its opposite has 3/4 of certainty; both possibilities are probable and their probabilities are as 8:9. Koopman (1940) resumed, in our time, the study of non-additive probabilities whose sources can be found in the medieval doctrine of probabilism that considered the opinion of each theologian as probable. Franklin (2001, p. 74) attributed the origin of probabilism to the year 1577, or, in any case (p. 83), to 1611. Nevertheless, similar pronouncements on probabilities of opinion go back to **John of Salisbury** (the 12$^{th}$ century) and even to **Cicero** (Garber & Zabell 1979, p. 46).

I note a "general rule or axiom" concerning the application of arguments (pp. 234 and 236): out of two possibilities, the safer, the more reliable, or at least the more probable should be chosen[8]; gamblers, however, always acted that same way (§ 1.2.3) if only they did not follow superstitious beliefs (§ 2.1.1).

**3.2.2. Statistical Probability and Moral Certainty.** Before going on to prove his LLN, Bernoulli (p. 246) explained that the theoretical "number of cases" was often unknown, but what was impossible to obtain beforehand, might at least be determined afterwards, i.e., by numerous observations. The application of statistical [probabilities], he maintained, was not new at all and referred to the "celebrated" **Arnauld**, the co-author of Arnauld & Nicole (1662/1992)[9]. In his Diary, Bernoulli indirectly mentioned **Graunt** (§ 3.1.1) and, furthermore, quite to the point, in connection with the impossibility of determining, without applying statistical data, how much more probable was the case of a youth outliving an old man than the opposite instance[10].

He informed **Leibniz** about his opinion (cf. § 3.1.2) and added that exactly that consideration led him to the idea of replacing, when necessary, prior knowledge by posterior. Recall also Bernoulli's reasoning on the statistical probability of a plague epidemic (§ 3.1.1).

I discussed moral certainty in §§ 2.1.2 and 2.2.2. Bernoulli (p. 238) maintained that it ought to be admitted on a par with absolute certainty and that judges must have firm instructions about what exactly (for example, 0.99 or 0.999 of certainty) constituted moral certainty. The latter idea was hardly ever put into practice; furthermore, the probability of a just sentence must be the higher the more important it is. On p. 249 Bernoulli mentioned moral certainty once more. His theorem will show, he declared, that statistical [probability] was a morally certain [a consistent, in modern terms] estimator of the theoretical [probability][11].

**3.2.3. The Law of Large Numbers.** Bernoulli proved a proposition that, beginning with **Poisson**, is called the LLN (§ 8.7). Let $r$ and $s$ be natural numbers, $t = r + s$, $n$, a large natural number, $\nu = nt$, the number of independent[12] trials in each of which the studied event occurs with [probability] $r/t$, $\mu$ – the number of the occurrences of the event (of the successes). Then Bernoulli proved without applying mathematical analysis that

$$P\left(|\frac{\mu}{\nu} - \frac{r}{t}| \leq \frac{1}{t}\right) \geq 1 - \frac{1}{1+c} \qquad (1)$$

and estimated the value of $\nu$ necessary for achieving a given $c > 0$. In a weaker form Bernoulli's finding meant that

$$\lim P = \left(|\frac{\mu}{\nu} - \frac{r}{t}| < \varepsilon\right) = 1, \quad \nu \to \infty, \qquad (2)$$

where, as also in (1), $r/t$ was the theoretical, and $\mu/\nu$, the statistical probability.

**Markov** (Treatise, 1924, pp. 44 – 52) improved Bernoulli's estimate mainly by specifying his intermediate inequalities and **K. Pearson** (1925), by applying the **Stirling** formula, achieved a practically complete coincidence of the Bernoulli result with the estimate that makes use of the normal distribution as the limiting case of the binomial law[13]. In addition, Pearson (p. 202) considered Bernoulli's estimate of the necessary number of trials in formula (1) "crude" and leading to the ruin of those who would have applied it. He also inadmissibly compared Bernoulli's law with the wrong **Ptolemaic** system of the world (and **De Moivre** with **Kepler** and **Newton**):

*Bernoulli saw the importance of a certain problem; so did Ptolemy, but it would be rather absurd to call Kepler's or Newton's solution of planetary motion by Ptolemy's name*!

The very fact described by formulas (1) and (2) was, however, extremely important for the development of probability and statistics[14]; and, anyway, should we deny the importance of existence theorems?

And so, the LLN established a correspondence between the two probabilities[15]. Bernoulli (p. 249) had indeed attempted to ascertain whether or not the statistical probability had its "asymptote"; whether there existed such a degree of certainty, which observations, no matter how numerous, were unable to ensure. Or, in my own words, whether there existed such positive numbers $\varepsilon$ and $\alpha < 1$, that

$$\lim P\left(|\frac{\mu}{\nu} - \frac{r}{t}| < \varepsilon\right) \leq 1 - \alpha, \quad \nu \to \infty.$$

He answered his question in the negative: no, such numbers did not exist. He thus established, within the boundaries of stochastic knowledge, a relation between deductive and inductive methods and combined statistics with the art of conjecturing. Strangely enough, statisticians for a long time had not recognized this fact. Haushofer (1872, pp. 107 – 108) declared that statistics, since it was based on induction, had no "intrinsic connections" with mathematics based on deduction (consequently, neither with probability). A most noted German statistician, **Knapp** (1872, pp. 116 – 117), expressed a strange idea: the LLN was hardly useful since statisticians always made only one observation, as when counting the inhabitants of a city. And even later, Maciejewski (1911, p. 96) introduced a "statistical law of large numbers" instead of the Bernoulli proposition that allegedly impeded the development of statistics. His own law qualitatively asserted that statistical indicators exhibited ever lesser fluctuations as the number of observations increased.

All such statements definitely concerned the **Poisson** law as well (European statisticians then hardly knew about the **Chebyshev** form of the LLN) and Maciejewski's opinion likely represented the prevailing attitude of statisticians. Here, indeed, is what **Bortkiewicz** (1917, pp. 56 – 57) thought: the expression *law of large numbers* ought to be used only for denoting a "quite general" fact, unconnected with any definite stochastic pattern, of a higher or lower degree of stability of statistical indicators under constant or slightly changing conditions and given a large number of trials. Even **Romanovsky** (1912, p. 22; 1924, pt. 1, p. 15; 1961, p. 127) kept to a similar view. Thus, in the last-mentioned contribution he stressed the natural-scientific essence of the law and called it physical.

The LLN has its prehistory. It was thought, long before **Bernoulli**, that the number of successes in $n$ "Bernoulli" trials with probability $p$ was approximately equal to

$$\mu = np. \qquad (3)$$

**Cardano** (**Ore** 1963, pp. 152 – 154 and 196), for example, applied this formula in calculations connected with games of dice. When compiling his mortality table, **Halley** (§ 2.1.4) assumed that "irregularities" in his data would have disappeared had he much more observations at his disposal. His idea can be interpreted as a statement on the increase in precision of formula (3) with $n$, see however § 2.1.4. Also see **Graunt's** reasoning (same subsection).

A second approach to the LLN took shape in astronomy when the arithmetic mean became the universal estimator of the constant sought (§ 1.2.4), cf. § 6.3.3. If the expectation of each of the magnitudes (1.1) is equal to that constant, i.e., if systematic errors are absent, and if (as always was the case) their variances are bounded, it could be thought that $a$ was approximately equal to the arithmetic mean of the observations.

Similar but less justified statements concerning sums of magnitudes corrupted by random errors had also appeared. Thus, **Kepler** (Sheynin

1973c, p. 120) remarked that the total weight of a large number of metal money (the mean weight!) of the same coinage did not depend on the inaccuracy in the weight of the separate coins. Then (§ 2.1.3), when buying annuities "upon several young lives", the expectation of a gain $Ex_i$ from each such transaction was obviously positive; if constant, the buyer could expect a total gain of $n\ Ex_i$. Much later **Condorcet** (1785a, p. 226) testified that those engaged in such "commerce" (and apparently ignorant of the LLN) had regarded it as "sûre".

There also likely existed a practice of an indirect participation of (petty?) punters in many games at once. At any rate, both **De Moivre** (1718/1756, reprint 1967, Problem 70) and **Montmort** (1708/1980, p. 169) mentioned in passing that some persons bet on the outcomes of games[16]. The LLN has then been known but that practice could have existed from much earlier times. And, finally, Gower (1993, p. 272) noted that **Boscovich** (1758, § 481) had [somewhat vaguely] maintained that the sum of random magnitudes decreased with an increase in the number of terms.

The attitude of Kepler and Boscovich is understandable: a similar wrong opinion about random sums had persisted at least until the 19th century, so that **Helmert** (1905, p. 604) thought advisable to refute it in his own field.

**3.2.4. Randomness and Necessity.** Apparently not wishing to encroach upon theology, Bernoulli (beginning of Chapter 1) refused to discuss the notion of randomness. Then, in the same chapter, he offered a subjective description of the "contingent" but corrected himself at the beginning of Chapter 4 where he explained randomness by the action of numerous complicated causes. Finally, the last lines of his book contained a statement to the effect that some kind of necessity was present even in random things. He referred to **Plato** who had taught that after a countless number of centuries everything returned to its initial state. Bernoulli likely thought about the archaic notion of the Great Year whose end will cause the end of the world with the planets and stars returning to their positions at the moment of creation. Without justification, he widened the boundaries of applicability of his law and his example was, furthermore, too complicated. It is noteworthy that **Kepler** (1596) believed that the end of the world was unlikely. In this, the first edition of this book, his reasoning was difficult to understand but later (1621) he substantiated his conclusion by stating, in essence, like **Oresme** (1966, p. 247) did before him, that two [randomly chosen] numbers were "probably" incommensurable[17]. Bernoulli (end of Chapter 1) also borrowed **Aristotle's** example (§ 1.1.1) of finding a buried treasure, but, unlike him, had not connected it with randomness.

### 3.3. Bernoulli's Contemporaries

I dwell somewhat on the ideas and findings of some of Bernoulli's contemporaries but I postpone the discussion of **De Moivre**, whose first publication had appeared before the AC did, until Chapter 4.

**3.3.1. Arnauld.** Arnauld & Nicole anonymously put out their book, *Art of Reasoning* (1662/1992)[18]. Arnauld was its main author and I mentioned him in Note 11 of Chapter 2 (in connection with the **Pascal** wager), in § 2.1.2 (moral certainty) and § 2.1.1 (an advice to neglect unlikely favourable events). Then (§§ 3.2.2 and 3.2.1), I noted that Bernoulli mentioned him when justifying the use of statistical probabilities and borrowed his principle of behaviour. Finally, Arnauld repeatedly, although without a formal definition, applied the term *probabilité* (for example, on pp. 331 and 332), and *degrez de probabilité*. Recall that **Leibniz** (§§ 3.1.2 and 2.1.4), in turn, borrowed from him a reasoning and a term.

**3.3.2. Niklaus Bernoulli.** He published a dissertation on the application of the art of conjecturing to jurisprudence (1709/1975). An utterly unsuitable English translation is available in the Internet: Latin phrases up to 40 lines long are rendered in even longer sentences and in general very much is incomprehensible.

The dissertation contained

a) The calculation of the mean duration of life for persons of different ages.

b) A recommendation of its use for ascertaining the value of annuities and estimating the probability of death of absentees about whom nothing is known. He proposed to consider them dead when the probability of their death is twice higher than that of the alternative.

c) Methodical calculations of expected losses in marine insurance.

d) The calculation of expected gains (more precisely, of expected losses) in the Genoese lottery.

e) Calculation of the probability of truth of testimonies.

f) The determination of the life expectancy of the last survivor of a group of men (pp. 296 – 297; **Todhunter** 1865, pp. 195 – 196). Assuming a continuous uniform law[19] of mortality, he calculated the expectation of the appropriate [order statistic]. He was the first to use, in a published work, both this distribution and an order statistic.

g) A comment on the introduction of expectation by **Huygens** (p. 291; Kohli 1975c, p. 542), see expression (2.1). Bernoulli interpreted it as a generalized arithmetic mean and the centre of gravity "of all probabilities" (this is rather loose).

Apparently in accordance with his subject he had not discussed the treatment of observations, cf. § 2.1.4. Bernoulli's work undoubtedly fostered the spread of stochastic notions in society (cf. § 2.1.2), but I ought to add that not only did he pick up some hints included in the manuscript of the *Ars conjectandi*, he borrowed separate passages both from it and even from the *Meditationes* (Kohli 1975c, p. 541), never intended for publication. His numerous general references to **Jakob** do not excuse his plagiarism.

**3.3.3. Montmort.** He is the author of an anonymous book (1708/1980), important in itself and because of its obvious influence upon **De Moivre** as well as on **Niklaus Bernoulli**, the correspondence with whom **Montmort** included in 1713 in the second edition of his work. In the Introduction (1708/1713, p. iii) he indicated that in practical activities and considerations it was desirable to be guided by "geometry" rather than by superstition, cf. § 2.1.1. However, he (p. xii) added that, since he was unable to formulate appropriate "hypotheses", he was not studying the applications of [stochastic] methods to civil life[20].

Henny (1975) and **Hald** (1990) examined Montmort's findings. The latter, on his p. 290, listed Montmort's main methods: combinatorial analysis, recurrent formulas and infinite series; and on p. 209 Hald added the method (the formula) of inclusion and exclusion

$$P(\cup A_i) = \sum P(A_i) - \sum P(A_i \cdot A_j) + \sum P(A_i A_j A_k) - \ldots, \qquad (4)$$

where $A_1, A_2, \ldots, A_n$ were events and $i < j < k < \ldots$ This formula is an obvious stochastic corollary of a general proposition about arbitrarily arranged sets.

Here are some problems solved by Montmort (1708/1980, pp. 244 – 246/NNo. 46 – 50 and 203 – 205/NNo. 200 – 202; 130 – 143), see Hald (1990, pp. 196 – 198; 206 – 213; 292 – 297; and 328 – 336) respectively:

1) The problem of points. Montmort arrived at the negative binomial distribution. He returned to this problem in his correspondence with Niklaus Bernoulli (Hald 1990, pp. 312 – 314).

2) A study of throwing $s$ points with $n$ dice, each having $f$ faces. Montmort applied the combinatorial method and formula (4).

3) A study of arrangements and, again, of a game of dice. Montmort arrived at the multivariate hypergeometric, and the multinomial distributions.

4) A study of occupancies. Tickets numbered $1, 2, \ldots, n$, are extracted from an urn one by one without replacement. Determine the probability that at least one ticket with number $k$, $1 \leq k \leq n$, will occur at the $k$-th extraction. Montmort derived the appropriate formulas

$$P_n = 1 - 1/2! + 1/3! - \ldots + (-1)^{n-1}/n!, \lim P_n = 1 - 1/e, n \to \infty.$$

Niklaus Bernoulli and De Moivre returned to this problem, see H. A. David & Edwards (2001, pp. 19 – 29).

**3.3.4. Montmort and Niklaus Bernoulli: Their Correspondence.** I outline their correspondence of 1710 – 1713 (Montmort 1708/1980, pp. 283 – 414).

1) The strategic game *Her* (**Hald** 1990, pp. 314 – 322). The modern theory of games studies it by means of the minimax principle. Nevertheless, already Bernoulli indicated that the gamblers ought to keep to [mixed strategies]. See also Bellhouse et al (2015).

2) The gambler's ruin. Montmort wrote out the results of his calculations for some definite initial conditions whereas Bernoulli indicated, without derivation, the appropriate formula (an infinite series). Hald believes that he obtained it by means of the method of inclusion and exclusion. On this point and on the appropriate findings of Montmort and **De Moivre** see also Thatcher (1957), Takácz (1969) and Kohli (1975a).

3) The sex ratio at birth (Montmort 1708/1713, pp. 280 – 285; Shoesmith 1985a). I only dwell on Bernoulli's indirect derivation of the normal distribution (Sheynin 1968[21]; 1970a, pp. 201 – 203). Let the sex ratio be $m/f$, $n$, the total yearly number of births, and μ and $(n - μ)$, the numbers of male and female births in a year. Denote

$$n/(m + f) = r, m/(m + f) = p, f/(m + f) = q, p + q = 1,$$

and let $s = 0(\sqrt{n})$. Then Bernoulli's derivation (Montmort 1708/1980, pp. 388 – 394) can be presented as follows:

$$P(|\mu - rm| \geq s) \leq (t - 1)/t,$$
$$t \leq [1 + s(m + f)/mfr]^{s/2} \approx \exp[s^2(m + f)^2/2mfn],$$
$$P(|\mu - rm| \leq s) \geq 1 - \exp(s^2/2pqn),$$

$$P[|\mu - np|/\sqrt{npq} \leq s] \geq 1 - \exp(-s^2/2).$$

This result does not however lead to an integral theorem since $s$ is restricted (see above) and neither is it a local theorem; for one thing, it lacks the factor $\sqrt{2/\pi}$ [22].

4) The Petersburg game. In a letter to Montmort, Bernoulli (Ibidem, p. 402) described his invented game. *B* throws a die; if a six arrives at once, he receives an *écu* from *A*, and he obtains 2, 4, 8, … *écus* if a six only occurs at the second, the third, the fourth, … throw. Determine the expectation of *B*'s gain. Gabriel Cramer insignificantly changed the conditions of the game; a coin appeared instead of the die, and the occurrence of heads (or tails) has been discussed ever since. The expectation of gain became

$$E = 1·1/2 + 2·1/4 + 4·1/8 + … = \infty, \quad (5)$$

although a reasonable man will never pay any considerable sum in exchange for it.

This paradox is still being examined. Additional conditions were being introduced; for example, suggestions were made to neglect unlikely gains, i.e., to truncate the series (5); to restrict beforehand the possible payoff; and, the most interesting, to replace expectation by *moral expectation*[23]. In addition, **Condorcet** (1784, p. 714) noted that the possibly infinite game nevertheless provided only one trial and that only some mean indicators describing many such games could lead to an expedient solution. Actually issuing from the same idea, **Freudenthal** (1951) proposed to consider a number of games with the role of the gamblers in each of them to be decided by lot. Finally, the Petersburg game caused **Buffon** (1777, § 18) to carry out the apparently first statistical experiment. He conducted a series of 2048 games; the mean payoff was 4.9 units, and the longest duration of play (in six cases), nine throws[24]. From a theoretical point of view, the game was interesting because it introduced a random variable with an infinite expectation.

### Notes

1. **Bernoulli** (1713a/1999, p. 233) had additionally explained this expression by the Greek word *stochastice* which **Bortkiewicz** (1917, p. X), with a reference to him, put into scientific circulation. Already **Wallis** (1685, p. 254) had applied the expression *stochastic* (iterative) process and Prevost & Lhuilier (1799, p. 3) mentioned stochastics, or "l'art de conjecturer avec rigueur sur l'utilité et l'étendue [and the extensiveness] du principe par lequel on estime la probabilité des causes". Hagstroem (1940) indicated that **Plato** and **Socrates** had applied the term *stochastics* and that the *Oxford English Dictionary* had included it with a reference to a source published in 1662.

   **2. Bernoulli** had not invariably applied this term, see § 3.1.2.

   **3.** Bearing in mind the published subject of that part, it would have been expedient to isolate the mentioned applications of the art of conjecturing.

   **4.** Concerning the use of probable inferences in civil life see also the opinion of **Cicero** in Note 6 to Chapter 1. And, in connection with **Leibniz'** understanding of randomness (below), I recall Cicero (1991, Buch 2, § 17, p. 149) once more:

   "Nothing is more opposed to calculation and regularity than chance".

The publishers appended to the AC the author's French contribution *Lettre à un Amy sur les Parties du Jeu de Paume* (**J. Bernoulli** 1975) in which he calculated the expectations of winnings in a time-honoured variety of tennis.

   **5.** I mentioned it in §§ 2.1.1 and 2.1.2.

   **6.** In 1714, in a letter to one of his correspondents, **Leibniz** (Kohli 1975b, p. 512) softened his doubts about the application of statistical probabilities and for some reason added that the late **Jakob Bernoulli** had "cultivated" [the theory of probability] in accordance with his, Leibniz', "exhortations".

   **7.** On **Bernoulli's** non-additive probabilities see Shafer (1978) and Halperin (1988).

   **8.** See **Arnauld** & **Nicole** (1662/1992, p. 327): we should choose the more probable.

   **9.** I have found there only one appropriate but not really convincing example, see § 2.1.2. On their p. 281 these authors mention the possibility of posterior reasoning.

   **10.** Cf. his Example 4 (**Bernoulli** 1713a/1999, p. 236).

   **11. Leibniz** (1668 – 1669 (?)/1971, pp. 494 – 500 reasoned on the application of moral certainty in theology. One of its chapters should have included the expression "infinite probability or moral certainty". In a later manuscript of 1693 he (**Couturat** 1901, p. 232), unfortunately, as it seems, isolated logical certainty, physical certainty or "logical probability", and physical probability. His example of the last-mentioned term, the southern wind is rainy, apparently described a positive correlative dependence.

   **12. De Moivre** (§ 4.1) was the first to mention independence.

   **13. Markov** had not applied that formula apparently because **Bernoulli** did not yet know it.

   **14.** Strengthened by the prolonged oblivion of **De Moivre's** finding (§ 4.3).

   **15.** Throughout Part 4, **Bernoulli** considered the derivation of the statistical probability of an event given its theoretical probability and this is most definitely seen in the formulation of his *Main Proposition* (the LLN) in Chapter 5. However, both in the last lines of that chapter, and in Chapter 4 he mentioned the inverse problem and actually alleged that he solved it as well. I return to this point in my Chapter 5.

   **16. Cournot** (1843, §11) also mentions them vaguely.

   **17.** For us, **Oresme's** understanding of incommensurability is unusual, but I do not dwell on this point. Before him, **Levi ben Gerson** (1999, p. 166) stated that the heavenly bodies will be unable to return to their initial position if their velocities were incommensurable. He had not, however, mentioned the end of the world.

   **18. Bernoulli** possibly thought about that expression when choosing a title for his book (and for the new discipline of the same name, the predecessor of the theory of probability).

   **19. Huygens'** appropriate reasoning (§ 2.2.2) appeared in print much later.

   **20.** Cf. **Daniel Bernoulli's** moral expectation (§ 6.1.1).

   **21.** Only in its reprint of 1970 (p. 232).

   **22.** Nevertheless, **A. P. Youshkevich** (1986) reported that at his request three mathematicians, issuing from the description offered by **Hald**, had concluded that **Bernoulli** had come close to the local theorem. Neither had Hald (1998, p. 17) mentioned that lacking factor.

   **23.** See **Daniel Bernoulli's** memoir of 1738 in § 6.1.1. He published it in Petersburg, hence the name of the game.

   **24.** O. Spieß (1975) dwelt on the early history of the Petersburg game and Jorland (1987) and Dutka (1988) described later developments. Dutka also adduced the results of its examination by means of statistical simulation.

## 4. De Moivre and the De Moivre – Laplace Limit Theorem
### 4.1. *The Measurement of Chance* (1712)

In his first probability-theoretic work De Moivre (1712/1984) justified the notion of expected random gain by common sense rather than defining it formally as has been done later, cf. § 2.2.2; introduced the multiplication theorem for chances (mentioning independence of the events) and applied the addition theorem, again for chances; and, in solving one of his problems (No. 26), applied the formula (3.4) of inclusion and exclusion. I describe some of his problems; I have mentioned Problem 14 (repeated in De Moivre's *Doctrine of chances*) in § 2.2.2.

1) Problem No. 2. Determine the chances of winning in a series of games for two gamblers if the number of remaining games is not larger than *n*, and the odds of winning each game are *a/b*. De Moivre notes that the chances of winning are as the sums of the respective terms of the development of $(a + b)^n$.

2) Problem No. 5. The occurrence of an event has *a* chances out of $(a + b)$. Determine the number of trials (*x*) after which it will happen, or not happen, with equal probabilities[1]. After determining *x* from the equation

$$(a + b)^x - b^x = b^x,$$

De Moivre assumes that $a/b = 1/q$, $q$ , and obtains

$$1 + x/q + x^2/2q^2 + x^3/6q^3 + \ldots = 2, \quad x = q\ln 2, \qquad (1)$$

whose left part resembles the **Poisson** distribution.

3) A lemma. Determine the number of chances for the occurrence of *k* points in a throw of *f* dice each having *n* faces. Later De Moivre (1730, pp. 191 – 197; 1718 and 1756, Problem No. 3, Lemma) solved this problem by means of a generating function of a sequence of possible outcomes of a throw of one die.

4) Problem No. 9 (cf. **Pascal's** problem from § 2.2.2). Gamblers *A* and *B* have *p* and *q* counters, and their chances of winning each game are *a* and *b*, respectively. Determine the odds of their ruining. By a clever trick that can be connected with the notion of martingale (Seneta 1983, pp. 78 – 79) De Moivre obtained the sought ratio:

$$P_A/P_B = a^q(a^p - b^p) \div b^p(a^q - b^q). \qquad (2)$$

He left aside the elementary case of $a = b$.

5) Problem No. 25. Ruining of a gambler during a finite number of games played against a person with an infinite capital. De Moivre described the solution in a generalized way; its reconstruction is due to **Hald** (1990, pp. 358 – 360). The conditions of the problem had not been formulated explicitly.

### 4.2. Life Insurance

De Moivre first examined life insurance in the beginning of the 1720s and became the most influential author of his time in that field.

Issuing from **Halley's** table (§ 2.1.4), he (1725/1756, pp. 262 – 263) assumed a continuous uniform law of mortality for all ages beginning with 12 years and a maximal duration of life equal to 86 years. I describe some of those of his numerous findings which demanded the application of the integral calculus.

1) Determine the expected duration of life for a man of a given age if the maximal duration, or complement of life is $n$ ($n = 86 -$ age). The answer is $n/2$ (p. 288). Reconstruction:

$$\int_0^n xdx/n = n/2.$$

2) Determine the probability of one person outliving another one if the complements of their lives are $n$ and $p$, $n > p$ (p. 324). Here, in essence, is De Moivre's solution. Let the random durations of the lives of $A$ and $B$ be ξ and η. Then, since at some moment $x$ the complement of $A$'s life is $(n - x)$,

$$P(\eta \leq x, \xi = x) = [(n - z)/n]dz/p, \ P(\xi > \eta) =$$
$$\int_0^p [(n - z)/n]dz/p = 1 - p/2n.$$

3) Determine the expected time Eτ during which both men with the same complements of life as in the previous problem do not die (p. 288). De Moivre only provided the answer; a reconstruction (**Czuber**, Note 22 to the German translation of 1906 of De Moivre's work) is as follows.

$$P(x \leq \xi \leq x + dx \text{ or } x \leq \eta \leq x + dx) =$$
$$[(n - x)/n]dx/p + [(p - x)p]dx/n,$$
$$E\tau = \int_0^p \{[(n - x)/np] + [p - x)p/n]\}dx = p/2 - p^2/6n.$$

Note that probabilities of the type of $P (\xi \leq x)$ easily lead to integral distribution functions.

**Hald** (1990, pp. 515 – 546) described in detail the work of **De Moivre** and of his main rival, **Simpson** (1775), in life insurance. Simpson improved on, and in a few cases corrected the former's findings. After discussing one of the versions of mutual insurance, Hald (p. 546) concluded that Simpson's relevant results represented "an essential step forward". However, the relations between De Moivre and Simpson became horrible and K. Pearson (1978) justly called Simpson *a most disreputable character* (p. 145), *an unblushing liar* and *a thorough knave at heart* (p. 184).

### 4.3. The *Doctrine of Chances* (1718, 1738, 1756)

This work published in three editions, in 1718, 1738, and, posthumously, in 1756 (reprinted in 1967), was De Moivre's main achievement. He developed it from his previous memoir (§ 4.1) and intended it for gamblers so that many results were provided there without proof. This fact together with other circumstances[2] caused this extremely important book, whose translation into French contemplated both **Lagrange** and **Laplace**[3], to remain barely known for many decades. I refer to the reprint of its last edition.

In his Introduction, De Moivre listed his main methods: combinatorial analysis, recurrent sequences (whose theory he himself developed) and infinite series; in particular, he applied appropriately truncated divergent series. Also in the Introduction, on pp. 1 – 2 he provided the "classical" definition of probability, usually attributed to Laplace, but kept to the previous reasoning on expectation (§ 4.1) and even introduced the *value of expectation* (p. 3); formulated the multiplication theorem for probabilities (not for chances, as previously) and, in this connection, once more mentioned independence. Two events, *A* and *B*, were independent, if, as he stated,

$P(B) = P(B/A)$, $P(A) = P(A/B)$

(modern notation here and below). For dependent events (p. 6), three in number (say),

$P(A·B·C) = P(A) P(B/A) P(C/A·B)$.                 (3)

I list now some of the problems from the *Doctrine* mentioned by **Hald** (1990, pp. 409 – 413) without repeating those described in § 4.1 and, for the time being, leaving aside the normal distribution.

1) The **Huygens** additional Problem No. 4 (§ 2.2.2) including the multivariate case. The appearance of the hypergeometric distribution: Problems NNo. 20 and 26.

2) Runs of successes in *n* **Bernoulli** trials including the case of $n \to \infty$: Problems NNo. 34 and 74.

3) Coincidences. A generalization of **Montmort's** findings (§ 3.3.3) by the method of inclusion and exclusion: Problems 35 and 36.

4) The gambler's ruin: Problems 58 – 71.

5) Duration of game: Problems 58 – 64, 68 – 71.

For the general reader the main merit of the *Doctrine* was the study of many widely known games whereas De Moivre himself, in dedicating its first edition to **Newton** (reprinted in 1756 on p. 329), perceived his main goal in working out

*A Method of calculating the Effects of Chance* […] *and thereby fixing certain rules, for estimating how far some sort of Events may rather be owing to Design than Chance* […] [so as to learn] *from your Philosophy how to collect, by a just Calculation, the Evidences of exquisite Wisdom and Design, which appear in the Phenomena of Nature throughout the Universe*.

I stress that De Moivre wrote this dedication before proving his limit theorem (§ 4.4). See **Pearson's** statement on Newton's influence in § 2.2.3. The aim of De Moivre's theory of probability was therefore the separation of necessity and randomness of any kind.

### 4.4. The De Moivre – Laplace Theorem

In 1730 De Moivre published his *Miscellanea analytica*. Later he appended two supplements; I am interested in the second one (1733)[4], which he printed in a small number of copies and sent out to his colleagues. In 1738 De Moivre translated it into English and included in the second, and then, in an extended form, in the third edition of the *Doctrine* (pp. 243 – 254 in 1756). Its title includes the words *binomial* $(a + b)^n$ which means that, although studying the particular case of the symmetric binomial, De Moivre thought about the general case. He (p. 250) also expressly and justly stated that the transition to the general case was not difficult. Strangely enough, even recently some authors (Schneider 1988, p. 118) maintained that De Moivre had only considered the particular case.

The date of the compilation of this supplement is known: on the first page of its Latin original De Moivre stated that he had concluded (at least its mathematical part) about 12 years earlier, *i.e.*, much earlier than he published his *Misc. anal.*, but we can only follow that contribution (Sheynin 1970a).

1) In Book 5 of the *Misc. anal.* De Moivre determined the ratio of the middle term of the symmetric binomial to the sum of all of its terms, and in the first supplement to that work he derived, independently from, and simultaneously with **Stirling**, the so-called Stirling formula. Only the value of the constant, $\sqrt{2}$, the latter communicated to him[5].

2) In the same Book, De Moivre calculated the logarithm of the ratio of the middle term of the binomial $(1 + 1)^n$ to the term removed by $l$ from it:

$(m + l – 1/2)\ln(m + l – 1) + (m – l + 1/2)\ln(m – l + 1)$
$– 2m\ln m + \ln[(m + l)/m], m = n/2.$

However, only in the second supplement De Moivre transformed this expression, obtaining, as $n \to \infty$, $– 2l^2/n$. The ratio itself thus became equivalent to

$1 – 2l^2/n + 4l^4/2n^2 – ...$           (4)

Actually, as corroborated by his further calculations, De Moivre thought about the inverse ratio.

3) Also in the same supplement, after integrating (4), De Moivre calculated the ratio of the sum of the terms between the middlemost and the one removed from it by $l$ to the sum of all the terms. It was equal to

$[2/\sqrt{2 n}]\, (l – 2l^3/1\cdot 3n + 4l^5/2\cdot 5n^2 – …).$       (5)

He then calculated this sum either by numerical integration, or, for $l < n/2$, by leaving only a few of its first terms. For $n \to \infty$ his main result can be written in modern notation as

$$\lim P[a \le \frac{\mu - np}{\sqrt{npq}} \le b] = \frac{1}{\sqrt{2\pi}} \int_a^b \exp(-z^2/2)\, dz. \qquad (6)$$

This is the integral De Moivre – **Laplace** theorem (see § 7.1-3), as **Markov** (1924, p. 53) called it, – a particular case of the CLT, a term introduced by **Polya** (1920). Note that neither De Moivre, nor Laplace knew about uniform convergence that takes place here.

**Todhunter** (1865, pp. 192 – 193) inadequately described the essence of De Moivre's finding. He failed to note that De Moivre had actually considered the general case of $p \ne q$ and only stated that "by the aid of **Stirling's** Theorem the value of **Bernoulli's** Theorem is [was] largely increased"[6]. De Morgan (1864) and then Eggenberger (1894) were the first to note that De Moivre had arrived at the [normal distribution]. Concerning De Morgan see Note 1 to Chapter 8.

De Moivre (1718/1756, p. 252) mentioned the study of the sex ratio at birth (§ 2.2.4) and illustrated it by imagined throws of 14 thousand dice each having 35 faces and painted in two colours, 18 faces white and 17 black[7]. His reasoning (and his general considerations on p. 251) meant that, for him, the binomial distribution was a divine law of nature, stochastic only because of possible deviations from it. De Moivre thus actually recognized the mutual action of necessity and randomness, cf. § 3.2.4.

### Notes

**1. De Moivre** thus made use both of chances and probability.

**2. De Moivre's** symbolism soon became dated; the English language had been little known on the Continent; **Todhunter**, the most influential historian of probability of the 19[th] century, inadequately described De Moivre's main finding (§ 4.4); and, last but not least, **Laplace** (1814/1995, p. 119) did not sufficiently explain it.

**3. Lagrange's** letter to **Laplace** of 30 Dec. 1776 in t. 14 of his *Oeuvres* (1892), p. 66.

**4.** I call this memoir a supplement only for the sake of tradition: its extant copies in large libraries were bound to the *Misc. anal.*

**5.** In the same supplement **De Moivre** included a table of lg $n$! for $n$ = 10 (10) 900 with 14 decimals; reprinted: (1718/1756, p. 333). Eleven or twelve decimals were correct; a misprint occurred in the value of lg 380!.

**6.** In 1740, **Simpson** directly dwelt on the general case (**Hald** 1990, pp. 21 – 23).

**7.** Regular 35-hedrons do not exist, but this is not here important. **De Moivre** thought about 14 thousand yearly births with $m:f$ = 18:17.

## 5. Bayes
### 5.1. The Bayes Formula and Induction

I dwell on the posthumous memoir (Bayes 1764 – 1765) complete with the commentaries by **Price**. In its first part Bayes introduced his main definitions and proved a few theorems; note that he defined probability through expectation. There was no hint of the so-called Bayes theorem

$$P(A_i/B) = \frac{P(B/A_i)P(A_i)}{\sum_{j=1}^{n} P(B/A_j)P(A_j)}. \qquad (1)$$

Lubbock et al (1830, p. 48) first applied the term itself, as noted by H. A. David, see H. A. David & Edwards (2001, p. 215). He referred to an oral communication from Hald. I return to formula (1) in §§ 7.1-1 and 9.2-2. Here I indicate that Bayes had in essence introduced induction into probability and that his approach that assumed the existence of prior probabilities or distributions (see below) greatly influenced the development of mathematical statistics[1].

Bayes then studied an imaginary experiment, a ball falling on point *r* situated in a unit square *ABCD*, "to the left" or "to the right" of some straight line *MN* parallel to, and situated between *AB* and *CD*. If, after (*p* + *q*) trials, the point *r* occurred *p* times to the right of *MN* and *q* times, to the left of it, then

$$P(b \leq r \leq c) = \int_b^c u^p(1-u)^q\,du \div \int_0^1 v^p(1-v)^q\,dv \qquad (2)$$

where *bc* is a segment within *AD*. Bayes derived the denominator of (2) obtaining the value of the [beta-function] $B(p + 1; q + 1)$ and spared no effort in estimating its numerator. The right side of (2) is now known to be equal to the difference of two values of the incomplete beta-function

$$I_c(p + 1; q + 1) - I_b(p + 1; q + 1).$$

Thus, given the results of the experiment, and assuming a uniform prior distribution[2] of the location of *MN* and *r*, he determined the appropriate theoretical probability considered as a random variable.

In his covering letter, **Price** provided purely methodical illustrations; the most interesting of them required the [probability] of the next sunrise observed $10^6$ times in succession. Formula (2) indirectly answers his question if $b = 1/2$ and $c = 1$ are chosen; it also provides the probability of the contrary event if $b = 0$ and $c = 1/2$. Price (Bayes 1764/1970, pp. 149 and 150 – 151) also solved the same question for $p = 1$ and $q = 0$ and obtained $P = 3/4$ which is doubtful: knowing nothing about the essence of a phenomenon we should have got $P = 1/2$ (cf. **Poisson's** reasoning in § 8.1). In this case, formula (2) is wrong. Note also that **Chebyshev** (1879 – 1880/1936, p. 158) formulated the same problem on an everyday level: To determine the probability of a student's successful answer to the next question after his previous successes.

I dwell somewhat on the above. The probability (this time, the actual probability) of the next sunrise is

$$P = \int_0^1 x^{p+1} dx \div \int_0^1 x^p dx = \frac{p+1}{p+2}$$

(cf. §7.1-5). Now, $P(p = 1) = 2/3$, not good enough either.

**Polyá** (1954, p. 135) remarked that each consecutive success (sunrise) provided ever less justification for the next one.

**Cournot** (1843, § 93) considered a similar problem: A woman gave birth to a boy; determine the probability that her next child will also be a boy. Without justification, he stated that "perhaps" the odds were 2:1 but that it was impossible to solve that problem. See the opinions of **Laplace** (§§ 7.1-1 and 7.1-5)**, Gauss** (§ 9A.2-2) and **Chebyshev** (§ 13.2-7) about the Bayesian approach. Another point concerned the Bayesian treatment of an unknown constant $r$ in formula (2) as a random variable (Neyman 1938/1967, p. 337).

Beginning with the 1930s and perhaps for three decades English and American statisticians had been denying Bayes. I am, however, leaving aside that period and I only note that the first and the main critic of the Bayes "theorem" or formula was **Fisher** (1922, pp. 311 and 326) but that he had not specified what exactly did he refuse to comply with. It seems that he disagreed with the introduction of hardly known prior probabilities and/or with the assumption that they were equal to one another, cf., however, **Laplace's** general statement about rectifying hypotheses (§ 7.2-1). The papers of Cornfield (1967) and Barnard (1967) were also insufficiently definite; the former figuratively remarked, on his p. 41, that Bayes had returned from the cemetery.

It is methodologically important to note that the *inverse probability* defined by formula (1) is tantamount to conditional probability given that the stipulated condition has indeed been fulfilled.

### 5.2. The Limit Theorem

I dwell now on the case of $n = (p + q)$ which Bayes had not expressly discussed. **Price**, however, remarked that, for a finite $n$, **De Moivre's** results were not precise. I indicate in addition that another posthumous note written by Bayes was published at the same time, 1764, as his main memoir, and in the same periodical, and there he warned mathematicians about the danger of applying divergent series. He had not named De Moivre, but he apparently had in mind the derivation of the De Moivre – **Laplace** theorem (4.6) as well.

Timerding, the Editor of the German translation of the Bayes memoir, nevertheless went on to consider the limiting case. He issued from Bayes' calculations made for large but finite values of $p$ and $q$. Applying a clever trick, he proved that, as $n \to \infty$, the probability of the ball falling to the right of *MN* obeyed the proposition

$$\lim P\left\{\frac{|\ -a|}{\sqrt{pq/n^3}}\ z\right\} = \frac{1}{\sqrt{2\pi}} \int_0^z \exp(-w^2/2)dw, \qquad (3)$$

where (not indicated by Timerding) $a = p/n = $ E , $pq/n^3 = $ var . The functions in the left sides of formulas (4.6) and (3) are random variables, centred and normed in the same way, and it is remarkable that Bayes, without knowing the notion of variance, apparently understood that (4.6) was not sufficiently precise for describing the problem inverse to that studied by **De Moivre**. Anyway, **Price** (Bayes 1764/1970, p. 135) stated that he knew

*of no person who has shewn how to deduce the solution of the converse problem* […]. *What Mr De Moivre has done therefore cannot be thought sufficient* …

**Jakob Bernoulli** (Note 15 in my Chapter 3) maintained that his formulas were also fit for solving the inverse problem – but how precisely? True, **De Moivre** (1718/1756, p. 251) also mentioned the inverse problem:

*Conversely, if from numberless Observations we find the Ratio of the Events to converge to a determinate quantity* […], *then we conclude* [how precisely?] *that this Ratio expresses the determinate Law according to which the Event is to happen.*

In my opinion, the Bayes' insufficiently known proposition is very important. Together with the integral **De Moivre – Laplace** theorem it completed the creation of the first version of the theory of probability. **Bayes was actually the main predecessor of Mises** who had not, however, mentioned him.

### 5.3. Additional Remark

Stigler (1983) quoted a curious statement (Hartley 1749, pp. 338 – 339) and interpreted it as a testimony against Bayes' priority. After referring to **De Moivre**, Hartley wrote, in part:

*An ingenious friend has communicated to me a solution of the inverse problem of determining the probability of an event given the number of times it happened and failed.*

Later Stigler (1986, pp. 98, 132) recalled Hartley and his own earlier paper of 1983, but did not definitively repeat his previous inference. Then, however, he (1999, pp. 291 – 301) reprinted that paper and added a tiny footnote brushing aside all the criticism published by that time.

Stigler inferred that the author of the Bayes' theorem was Saunderson (1682 – 1739), and by applying formula (1), he even found that his conclusion was three times more probable than the former opinion. However, he based himself on subjective estimates and in particular assumed that the prior probabilities of the authorship of Bayes and Saunderson were the same. This means that the extra-mathematical arguments (for example, the evidence of Price, a close friend of Bayes) are not considered at all. And it is opportune to recall the opinion of Gauss (§ 3.1.2): applications of the theory of probability can be greatly mistaken if the essence of the studied object is disregarded.

In addition, not only a honest personality like Saunderson, but almost any pretender will be able to claim equal prior rights with an established author (or a politician of the past). For my part, I have argued that it was Bayes himself who communicated to Hartley the solution "of the inverse problem". Referring to later sources, Zabell (1989, p. 316) concluded that Stigler's opinion was not serious.

### Notes

**1.** A modern encyclopaedia (**Prokhorov** 1999b) contains 14 items mentioning him; for example, Bayesian estimator, Bayesian approach, etc. There also, on p. 37, the author of the appropriate entry mistakenly attributes formula (1) to **Bayes**.

**2. Bayes** himself had not stated that this distribution was uniform, but it is nevertheless necessary to make this assumption (**K. Pearson** 1978, p. 364). Without any explanation provided, **Mises** (1919, § 9.2) remarked that Bayes had considered the general case as well. Following **Czuber**, to whom he referred, Mises proved that the influence of non-uniformity of the prior distribution weakens with the increase in the number of observations.

## 6. Other Investigations before Laplace
### 6.1. Stochastic Investigations

**6.1.1. Daniel Bernoulli.** He published a number of memoirs pertaining to probability and statistics, and, before that, he (1735) provided a stochastic reasoning on the structure of the Solar system. The inclinations of the orbits of the five (excepting the Earth) then known planets with respect to the Earth (considered as random variables with a continuous uniform distribution) were small, and the probability of a "random" origin of that circumstance, as he concluded, was negligible. I dwelt on the logic of such considerations in § 1.1.1; here, however, enters a new dimension (see § 10.8.4): it was possible to study, instead of the inclinations, the arrangement of the poles of the orbits (**Todhunter** 1865, p. 223).

In this subsection, I consider only some of Bernoulli's memoirs and I postpone the study of his other work until §§ 6.2.3 and 6.3, but my general conclusion is that he, together with **De Moivre**, was the main predecessor of **Laplace**.

1) Moral expectation. While attempting to explain the paradoxical nature of the Petersburg game (§ 3.3.4), **Bernoulli** (1738) suggested that the gain $y$ of a gambler was determined by his winnings $x$ in accord with the differential equation (the first such equation in probability)

$dy = cdx/x$, $c > 0$, so that $y = f(x) = c\ln(x/a)$

where $a$ was the gambler's initial capital.

Bernoulli also proposed that the expected winnings $(p_1 x_1 + p_2 x_2 + \ldots + p_n x_n) / \sum p_i$ where $p_i$ were the appropriate probabilities be replaced by their "moral expectation"

$\sum p_i f(x_i) / \sum p_i$.

He indicated but had not proved (see § 7.1-9) that even a "just" game with a zero expected loss for each participant became disadvantageous because the moral expectation of winnings, again for each, was negative, and that the infinite expected gain in the Petersburg game (3.5) can be replaced by a finite moral expectation. Then, applying his innovation to a study of marine shipping of freight, he maintained (again, without proof, see same subsection below) that the freight should be evenly distributed among several vessels.

Bernoulli appended the text of a letter of 1732 from Gabriel Cramer to **Nikolaus Bernoulli** which contained his (not Daniel's) term *moral expectation*. Cramer also indirectly suggested the choice of

$f(x) = \min(x; 2^{24})$ or $f(x) = \sqrt{x}$.

In a letter of 1742 Bernoulli left a curious statement concerning politics (P. N. Fuss 1843/1968, t. 2, p. 496): "I believe that mathematics can also be rightfully applied in politics". After having referred to a positive opinion of **Maupertuis**, he continued:

*An entirely new science will emerge provided that as many observations are made in politics as in physics*.

With all respect, this is at least unclear.

Moral expectation had become popular and **Laplace** (1812/1886, p. 189) therefore proposed a new term for the previous "usual" expectation calling it *mathematical*; his expression persists to no purpose at least in the French and Russian literature and I refuse to apply it. Modern statistical damage functions are akin to moral expectation and moreover, at the end of the 19th century, issuing from Bernoulli's idea, economists began to develop the theory of marginal utility thus refuting **Bertrand's** opinion (1888a, p. 66) that moral expectation was useless:

*The theory of moral expectation became classical, and never was a word used more exactly. It was studied and taught, it was developed in books really celebrated. With that, the success came to a stop; no application was made, or could be made, of it*.

2) A limit theorem. While studying the same problem concerning the sex ratio at birth (§§ 2.2.4, 3.3.4, 4.4), Bernoulli (1770 – 1771), in the first part of his memoir, assumed that male and female births were equally probable. It followed that the probability that the former constituted a half of 2$N$ births will be

$P = [1\ 3\ 5\ …\ (2N – 1)] \div [2\ 4\ 6\ …\ 2N] = q(N)$.

He calculated this fraction not by the **Wallis** formula but by means of differential equations. After deriving $q(N – 1)$ and $q(N + 1)$ and the two appropriate values of $q$, he arrived at

$dq/dN = – q/(2N + 2), dq/dN = – q/(2N – 1)$

and, "in the mean", $dq/dN = – q/(2N + 1/2)$. Assuming that the solution of this equation passed through point $N = 12$ and $q(12)$ as defined above, he obtained

$q = 1.12826/\sqrt{4N + 1}$.

Application of differential equations was Bernoulli's usual method in probability, also see § 6.1.1-1.

Bernoulli also determined the probability of the birth of approximately $m$ boys (see below):

$P(m = N \pm \mu) = q\exp(– \mu^2/N)$ with $\mu = 0(\sqrt{N})$.          (1)

In the second part of his memoir Bernoulli assumed that the probabilities of the birth of both sexes were in the ratio of $a:b$. Equating the probabilities of $m$ and $(m + 1)$ boys being born, once more for 2*N* births, he thus obtained the [expected] number of male births

$$Em = M = \frac{2Na - b}{a + b} \approx \frac{2Na}{a + b}$$

which was of course evident. More interesting was Bernoulli's subsequent reasoning for determining the probability of an arbitrary *m* (for µ of the order of $\sqrt{N}$):

$$P(m = M + \mu + 1) - P(m = M + \mu) \quad d \quad =$$
$$-(a/b)\frac{2N - M - \mu}{M + \mu + 1}d\mu,$$
$$-d \ / \quad = \frac{\mu + 1 + \mu a/b}{m + \mu + 1}d\mu.$$

The subsequent transformations included the expansion of ln[(*M* + 1 + µ)/(*M* + 1)] into a power series. Bernoulli's answer was

$$P(m = M \pm \mu) = \quad = P(m = M) \exp[-\frac{(a+b)\mu^2}{2bM}],$$

hence (1). Note that Bernoulli had not applied the local **De Moivre** (– **Laplace**) theorem.

Issuing from some statistical data, Bernoulli compared two possible ratios *a/b* but had not made a final choice in favour of either of them. He also determined such a value of µ that the sum of probabilities (1), beginning from µ = 0, equalled one half. Applying summation rather than integration, he had not therefore arrived at an integral limit theorem and (also see above) he did not refer to, and apparently had not known about **De Moivre's** findings. This shows, once again (cf. § 4.4), that they had for a long time been forgotten.

3) Urn problems. I consider two of these. An urn contains *n* pairs of white and black stripes. Determine the number (here and below, actually, the expected number) of paired stripes left after (2*n* – *r*) extractions without replacement. By the combinatorial method Bernoulli (1768a) obtained

$$x = r(r-1)/(4n-2); \text{ and } x = r^2/4n \text{ if } n = \quad .$$

He derived the same result otherwise: when *r* decreases by *dr* the corresponding *dx* is either zero [(*r* – 2*x*) cases] or *dr* (2*x* cases) so that

$$dx = [(r - 2x) \cdot 0 + 2x \cdot dr]/r, \ x = r^2/4n \text{ since } r = 2n \text{ if } x = n.$$

Bernoulli then generalized his problem by considering unequal probabilities of extracting the stripes of different colours and he (1768b) applied his findings to study the duration of marriages, a subject which was directly linked with insurance of joint lives.

Suppose now that each of two urns contains an equal number $n$ of balls, white and black, respectively. Determine the number of white balls in the first urn after $r$ cyclic interchanges of one ball. Bernoulli (1770) solved this problem by the same two methods. Thus, he issued from the differential equation

$$dx = - xdr/n + [(n - x)/n]dr \text{ so that } x \quad (1/2)n[1 + e^{-2r/n}].$$

Bernoulli then combinatorially considered the case of three urns with balls of three different colours. He noted that the number of white balls in the first urn was equal to the sum of the first, the fourth, the seventh, … terms of the development of $[(n - 1) + 1]^r$ divided by $n^{r-1}$. For the other urns he calculated, respectively, the sums of the second, the fifth, the eighth, …, and the third, the sixth, the ninth, … terms. For the first urn he obtained

$$A = \frac{1}{n^{r-1}}[(n - 1)^r + C_r^3 (n - 1)^{r-3} + C_r^6 (n - 1)^{r-6} + ...] \quad ne^{-r/n} S. \text{ (2)}$$

The expression designated by $S$ obeyed the differential equation

$$Sdr^3/n^3 = d^3S$$

and was therefore equal to

$$S = ae^{r/n} + be^{-r/2n} \sin(r \ 3/2n) + ce^{-r/2n} \cos(r \ 3/2n)$$

where, on the strength of the initial conditions, $a = 1/3$, $b = 0$, $c = 2/3$.

Bernoulli derived similar expressions for the other urns, calculated the number of extractions leading to the maximal number of white balls in the first urn, and, what is extremely interesting, he also noted the existence of a limiting state, of an equal number of balls of each colour in each urn. This can be easily verified by referring to the theorem on the limiting transition matrix in homogeneous **Markov** chains. Physicists could have noticed here a stochastic model of the previously received notion of the thermal death of (a finite) universe.

Bernoulli obtained formula (2) by issuing from differential equations

$$dx = - xdr/n + [n - (x + y)]dr/n, dy = - ydr/n + xdr/n$$

where $x$, $y$, and $[n - (x + y)]$ were the numbers of white balls in the urns after $r$ interchanges[1]. I return to this problem in §7.1-3; here, I note that **Todhunter** (1865, pp. 231 – 234) simplified Bernoulli's solution and made it more elegant. He wrote the differential equations as

$$dx = (dr/n) (z - x), dy = (dr/n) (x - y), dz = (dr/n) (y - z)$$

and noted that the sum $S$ was equal to

$$S = (1/3)[e^{r/n} + e^{r/n} + e^{r/n}]$$

with , , being the values of $\sqrt[3]{1}$.

Bernoulli's *x* in his first problem, and his *S* and *A* from (2) depend on discrete "time" *r/n*, which is characteristic of stochastic processes with non-homogeneous time. The same is true about the ratio *s*/ξ in § 6.2.3 below.

**6.1.2. D'Alembert.** In the theory of probability, he is mostly known as the author of patently wrong statements[2]. Thus, he (1754) maintained that the probability of heads appearing twice in succession was equal to 1/3 rather than to 1/4. Then, he (1768a) reasoned on the difference between "mathematical" and "physical" probabilities[3], stating without justification that, for example, after one of two contrary events had occurred several times in succession, the appearance of the other one becomes physically more probable. He was thus ridden by prejudices which **Montmort** had already mentioned and which **Bertrand** later refuted by a few words (§ 2.1.1). At the same time, D'Alembert recommended to determine probabilities experimentally but had not followed his own advice (which saved him from revealing his mistakes). Finally, he (1768b) denied the difference (perfectly well understood by **Huygens**, § 2.2.2) between the mean, and the probable durations of life and even considered its existence as an (additional) argument against the theory of probability itself.

It is opportune to recall **Euler's** opinion as formulated in one of his private letters of 1763 (**Juskevic** et al 1959, p. 221): D'Alembert tries "most shamelessly to defend all his mistakes". Anyway, D'Alembert (1768d, pp. 309 – 310) did not ascribe the theory of probability to "a precise and true calculus with respect either to its principles or results"[4].

On the other hand, D'Alembert thought that, in a single trial, rare events should be considered unrealizable (**Todhunter** 1865, § 473) and that absolute certainty was qualitatively different from "the highest probability". It followed from the latter statement that, given a large number of observations, an unlikely event might happen (cf. the strong law of large numbers), and, taken together, his considerations meant that the theory of probability ought to be applied cautiously. D'Alembert (1768c) also reasonably objected to **Daniel Bernoulli's** work on prevention of smallpox and formulated his own pertinent ideas (§ 6.2.3). I ought to add that D'Alembert was indeed praiseworthy for his work in other branches of mathematics (and in mechanics); note also that Euler had not elaborated his likely correct remark. On D'Alembert's work see also Yamazaki (1971).

**6.1.3. Lambert.** He was the first follower of **Leibniz** in attempting to create a doctrine of probability as a component of a general teaching of logic. Like **D'Alembert** (Note 3 to this Chapter), Lambert explained randomness by ignorance of causes, but he also stated that all digits in infinite decimal developments of irrational numbers were equally probable, which was an heuristic approach to the notion of normal numbers, and he formulated a modern-sounding idea about the connection of randomness and disorder (Lambert 1771, § 324; 1772 –

1775), also see Sheynin (1971a, pp. 238 – 239; 1971b, p. 246; 1974, pp. 136 – 137).

His considerations were forgotten, only Cournot (1851/1975, § 33 note) and Chuprov (1909/1959, p. 188) mentioned them.

Lambert did not go out of the confines of *uniform randomness.* To put his ideas in perspective, I ought to add that the philosophical treatises of the 18th century testify to the great difficulties experienced in generalizing the notion of randomness, also see § 2.2.4. One example: even in the 19th century, many scientists, imagining that randomness was only uniform, refused to recognize the evolution of species, and two authors (**Baer** 1873, p. 6; **Danilevsky** 1885, pt. 1, p. 194) independently mentioned the philosopher depicted in *Gulliver's Travels* (but borrowed by **Swift** from **Raymond Lully**, 13th – 14th centuries). That "inventor", hoping to get to know all the truths, was putting on record each sensible chain of words that appeared from among their uniformly random arrangements.

**6.1.4. Buffon.** He is mostly remembered for his definitive introduction of geometric probabilities (§ 6.1.6). Then, he reasonably suggested that the value of winnings in a game of chance diminished with the increase of the gambler's capital (cf. § 6.1.1) and experimentally studied the Petersburg game (§ 3.3.4), proposed the value 1/10,000 as a universally negligible probability, attempted to determine the probability of the next sunrise (see Chapter 5)[5], cf. § 7.1-5, and compiled tables of mortality which became popular.

Negligible, as he thought, was the probability of death of a healthy man aged 56 during the next 24 hours, but his figure was apparently too low; **K. Pearson** (1978, p. 193) thought that 1/1,000 would have been more appropriate. In addition, negligibility ought to be only chosen for a particular event rather than assigned universally. All the above is contained in Buffon's main work (1777). There also (§ 8, Note) he published the text of his letter of 1762 to **Daniel Bernoulli** which contained an embryo of **Quelelet's** celebrated Average man (see my § 10.5):

*Mortality tables are always concerned with the average man, that is, with people in general, feeling themselves quite well or ill, healthy or infirm, robust or feeble.*

**6.1.5. Condorcet.** He attempted to apply the theory of probability to jurisprudence in the ideal and tacitly assumed case of independent judgements made by jurors or judges. He also estimated the trustworthiness of testimonies and critically considered electoral problems. His main method was the application of difference equations. **Todhunter** (1865, pp. 351 – 410) described the work of Condorcet in detail and concluded (p. 352) that in many cases it was "almost impossible to discover" what he had meant[6]: "The obscurity and self contradiction are without any parallel …" He, Todhunter, will provide some illustrations, "but no amount of examples can convey an adequate impression of the extent of the evils". At the very least, however, **Laplace** and **Poisson** continued to apply probability to jurisprudence and certainly profited to some extent from the work of Condorcet. Poisson (1837a, p. 2) mentioned his ideas quite favourably.

I note however that, while discussing games of chance, Condorcet (1785b/1847, p. 561) expressed himself rather unfortunately, and stated on the next page without any justification that **Daniel Bernoulli** had not removed all the objections to the "rule" of expectation and that this was allegedly achieved by **D'Alembert**. In 1772, in a letter to Turgot, he (Henry 1883/1970, pp. 97 – 98) told his correspondent that he was "amusing himself" by calculating probabilities, had compiled "a booklet on that subject" and was keeping to the opinions of D'Alembert. That booklet is unknown. On Condorcet see also Zabell (1988b) and Yamazaki (1971).

**6.1.6. Geometric Probabilities.** These were decisively introduced in the 18$^{th}$ century although the definition of the notion itself, and, for that matter, only on a heuristic level, occurred in the mid-19$^{th}$ century (§ 10.3). **Newton** (§ 2.2.3) was the first to think about geometric probability; **Daniel Bernoulli** (§ 6.1.1) tacitly applied it in 1735 as did somewhat later (possibly even earlier) **De Moivre** (1725/1756, p. 323), and then **T. Simpson** (1757) (§ 6.3.1) and **Bayes** (§ 5.1). Dealing with the continuous uniform distribution, De Moivre assumed, for example, that if $0 < \xi < b$ and $0 < a < b$, then

$P(0 < \xi < a) = [0; a] \div [0; b]$.

Simpson noted that in his case (a continuous triangular distribution) probabilities were proportional to the areas of the appropriate figures. Bayes assumed that, for a continuous uniform distribution, the probabilities of a ball falling on different equal segments were equal to one another.

The **Michell** (1767) problem became classical: Determine the probability that two stars from all of them, uniformly distributed over the celestial sphere, were situated not farther than 1° from each other. Choose an arbitrary point (*A*) on a sphere with centre *O* and imagine a circle perpendicular to *OA* having distance 1° from *A*. The probability sought is the ratio of the surface of the spherical segment thus obtained to that of the sphere. **Newcomb** and **Fisher** calculated the expected number of closely situated stars (§ 10.8-4) and general issues were also debated by others. Thus, **Proctor** (1874, p. 99) wished to determine

*What peculiarities of distribution might be expected to appear among a number of points spread over a plane surface perfectly at random.*

His was a question which concerned the deviations of an empirical density curve from its theoretical counterpart and now belonged to mathematical statistics. And **Bertrand** (1888a, pp. 170 – 171) remarked that without studying other features of the sidereal system it was impossible to decide whether stars were arranged randomly.

**Buffon** (§ 6.1.4) expressly studied geometric probability; the first report on his work (Anonymous 1735) had appeared long before his contribution. Here is his main problem: A needle of length 2*r* falls

"randomly" on a set of parallel lines. Determine the probability *P* that it intersects one of them. It is easily seen that

$$P = 4r/a \qquad (3)$$

where $a > 2r$ is the distance between adjacent lines. Buffon himself had, however, only determined the ratio *r/a* for *P* = 1/2. His main aim was (Buffon 1777/1954, p. 471) to "put geometry in possession of its rights in the science of the accidental [du hasard]". Many commentators described and generalized the problem above. The first of them was **Laplace** (§7.1-4) who noted that formula (3) enabled to determine [with a low precision] the number . I treat the further history of geometric probability in my Chapter 12.

### 6.2. Statistical Investigations

**6.2.1. Staatswissenschaft (Statecraft, University Statistics).** In mid-18th century **Achenwall** (Sheynin 1997b) created the Göttingen school of *Staatswissenschaft* which described the climate, geographical situation, political structure and economics of separate states and estimated their population by issuing from data on births and mortality but did not study relations between quantitative variables. Achenwall referred to **Süssmilch** (§ 6.2.2), advised state measures fostering the multiplication of the population and recommended censuses without which (1763/1779, p. 187) a "probable estimate" of the population could be still got, see above. He (1752/1756, Intro.) also left an indirect definition of statistics:

*In any case, statistics is not a subject that can be understood at once by an empty pate. It belongs to a well digested philosophy, it demands a thorough knowledge of European state and natural history taken together with a multitude of concepts and principles, and an ability to comprehend fairly well very different articles of the constitutions of present-day kingdoms* [*Reiche*].

Achenwall's student **Schlözer** (1804, p. 86) figuratively stated that "History is statistics flowing, and statistics is history standing still". In spite of its elegance, this phrase is meaningless: Schlözer himself (§ 14-3) stated that it was necessary to compare one state with another and the changes in a state in time. And even Leibniz recommended such comparisons in his manuscripts of the 1680s (Sheynin 1977b, p. 224).

However, for those keeping to *Staatswissenschaft* this pithy saying became the definition of statistics which was thus not compelled to study causal connections in society or discuss possible consequences of innovations; which thus failed to adhere to the goals of political arithmetic (§ 2.1.4). The second distinction between the two disciplines consisted in that only political arithmetic was seriously interested in studying population. Finally, the methods of investigation were also different: not numbers, but wordy descriptions of *remarkable features* lay at the heart of the works of the Göttingen school.

Knies (1850, p. 24) quoted unnamed German authors who had believed, in 1806 and 1807, that the issues of statistics ought to be the national spirit, love of freedom, the talent and the characteristics of great personalities and ordinary people of a given state. This criticism, of course, had to do with the limitations of mathematics in general. Recall (§ 1.1.2) that Moses wished to find out not only whether the population of Canaan was numerous or not, but also whether the people there were "strong or weak".

Tabular statistics which had originated with **Anchersen** (1741) could have served as an intermediate link between words and numbers, but Achenwall was apparently opposed to it. Anyway, he (1752, Intro.) stated that he had "experienced a public attack" against the first edition of that book (published in 1749 under a previous title) by Anchersen. "Tabular" statisticians continued to be scorned, they were even called *Tabellenfabrikanten* and *Tabellenknechte* (slaves of tables) (Knies 1850, p. 23). When describing Russia in 1734, I. K. Kirillov used many tables, but his work was only published in 1831 (Ploshko & Eliseeva 1990, pp. 65 – 66). Another contribution which I have not seen is Golitzin (1807).

By the end of the 19th century, owing to the heterogeneity of its subject, *Staatswissenschaft* disintegrated. **K. Pearson** (1978, p. 125) remarked that political economy (**Adam Smith**) was the first discipline to break off from it and that the "evolution of Political Philosophers" had further curtailed the *Staatswissenschaft*. All this means that statistics, in its modern sense, owes its origin to political arithmetic. Consequently, I dwell below on contributions which had not belonged to the former subject, but, to the contrary, were mathematical or, in any case, issued from statistical data.

Apparently in the mid-19th century the Staatswissenschaft had regenerated and partly transferred to political arithmetic. It exists nowadays, at least in Germany, it is taught there in some universities and it can be considered as the application of the statistical method to the various aspects of the life of states.

**K. Pearson** ( p. 29) named **Edward Chamberlayne** (1616 – 1703) as the "English Achenwall" but he also noted that Chamberlayne had "copied" his book from a French work of 1661 (which he did not see).

**6.2.2. Population Statistics. Süssmilch** (1741) adhered to the tradition of political arithmetic. He collected vast statistical data on the movement of population and attempted (as **Arbuthnot** did, see § 2.2.4) to reveal in it divine providence but he treated his materials rather loosely. Thus, when taking the mean of the data pertaining to towns and rural districts, he tacitly assumed that their populations were equally numerous; in his studies of mortality, he had not attempted to allow for the differences in the age structure of the populations of the various regions etc. Nevertheless, it is possible to believe that his works paved the way for **Quetelet** (§ 10.5); in particular, he studied issues which later came under the province of moral statistics (e.g., illegitimate births, crime, suicides). And his tables of mortality had been in use even in the beginning of the 19th century. Like **Graunt**, Süssmilch discussed pertinent causes and offered conclusions. Thus, he (1758) thought of examining the dependence of mortality on climate and geographical position and he knew that poverty and ignorance were conducive to the spread of epidemics.

Süssmilch's main contribution, the *Göttliche Ordnung*, marked the origin of demography. Its second edition of 1765, included a chapter "On the rate of increase and the period of doubling [of the population]" written jointly with **Euler**. Partly reprinted in the latter's *Opera omnia* (in t. 7 of ser. 1, 1923), it served as the basis of one of Euler's memoirs (Euler 1767). Süssmilch naturally thought that the multiplication of mankind was a divine commandment and that, therefore, rulers must take care of their subjects. Quite consistently, he condemned wars and luxury and indicated that the welfare of the poor was to the advantage of both the state, and the rich. His collaboration with Euler and frequent references to the latter in his *Göttliche Ordnung* certainly mean that Euler had shared his general social views.

Both as a statistician and as a deeply religious person Süssmilch opposed polygamy. In this, he followed Graunt and perhaps other statisticians. And it seems that **Daniel Bernoulli**, the bachelor, was the only partisan of that practice. He (1768c, p. 103) stated, without any substantiation, that polygamy "certainly" fostered the multiplication of population. His was a letter to Euler, another deeply religious person, who likely did not comment on that point. The biographers of Bernoulli were not at all sure that he was really religious. On Süssmilch see Birg (1986) and **Pfanzagl** & Sheynin (1997).

**Malthus** (1798) picked up one of the conclusions made in the *Göttliche Ordnung*, viz., that the population increased in a geometric progression (which is still a more or less received statement).

Statisticians have been hard put to reconcile the Biblical command (Genesis 1:28) *Be fruitful and increase in number* … with the registered growth of the population. Pearson (1978, p. 337) commented on this difficulty experienced by an early British statistician:

*Apparently his view is that while the Creator would not approve of starvation for thinning humanity, he would have no objections to plague or war! It reminds me of the farmer's wife.* She saved an earwig from her hot iron and threw it in the kitchen fire saying *I nearly killed ye with the hot iron*!

Euler is known to have left hardly any serious contribution to the theory of probability (see also § 6.3.1 devoted to the theory of errors), but he published a few memoirs on population statistics collected in the same volume 7 of his works. When treating statistical data, he did not introduce any stochastic laws (for example, laws of mortality), but such concepts as increase in population and the period of its doubling are due to him, and his reasoning was always elegant and methodically interesting, in particular for life insurance (Paevsky 1935). See Sheynin (2007b).

**Lambert** published a mainly methodical study in population statistics (1772). Without due justification, he proposed there several laws of mortality (§ 9). One of them was represented by a sum of two terms and he explained that they described physical processes; now, we also see that they belonged to types IX and X of the **Pearson** curves. Then, he formulated the problem concerning the duration of marriages, statistically studied children's mortality from smallpox and the number of children in families (§ 108). See Sheynin (1971b) and Daw (1980) who also appended a translation of Lambert's discussion of the smallpox issue.

When considering the last-mentioned subject, Lambert issued from data on 612 families having up to 14 children and, once more without substantiation, somehow adjusted his materials. It is remarkable that he arbitrarily increased the total number of children by one half and that the new data, as he maintained, were "smoother". It might be thought that Lambert attempted to allow for stillbirths and the death of children. Elsewhere in his work he (§ 68) indicated that statistical investigations should reveal [and explain] irregularities.

**6.2.3. Medical Statistics**. Population statistics was not closely linked with any embryo of medical statistics, Graunt and Süssmilch (§ 6.2.2) notwithstanding. Schlözer passed medical problems over in silence. Statisticians have been ignorant of probability theory, recognized only *Bernoulli trials* (§ 3.2.3). Interesting enough was Poisson's introduction of the (soon forgotten) term *social arithmetic* (§ 8), i. e., population and medical statistics.

Medical statistics originated in the 19$^{th}$ century, partly because of the need to combat the devastating visitations of cholera. I note that the expression *medical probability* appeared not later than in the mid-18$^{th}$ century (Mendelsohn 1761, p. 204). At the end of that century **Condorcet** (1795/1988, p. 542) advocated collection of medical observations[7] and Black (1788, p. 65) even compiled a possibly forgotten "Medical catalogue of all the principle diseases and casualties by which the Human Species are destroyed or annoyed" that reminded of **Leibniz'** thoughts (§ 2.1.4). Descriptions belonging to other branches of natural sciences as well have actively been compiled (mostly later) and such work certainly demanded preliminary statistical efforts[8]. Some authors mistakenly stated that their compilations ruled out the need for theories (cf. **D'Alembert's** opinion in Note 7). Until the beginning of the 20$^{th}$ century, the partisans of complete descriptions continued to deny sampling in statistics proper (§ 10.8-2).

Especially important was the study of inoculation, of the method of preventing smallpox by communicating a mild form of smallpox from one person to another. The history of smallpox epidemics and inoculation is described in various sources (**Condamine** 1759, 1763, 1773; Karn 1931). In his first memoir, Condamine listed the objections against inoculation, both medical and religious. He ended his contribution thus:

*Since the inoculation of the Royal Family in England, this practice became generally accepted in France.*

His statement did not however square with **Daniel Bernoulli's** appeals in its favour. In his second memoir Condamine (p. 464) mentioned the Bernoulli family:

*In Basel, Messrs Bernoulli, whose name can by itself authorize doubtful opinions on many points, have not rested content to declare themselves publicly in favour of inoculation. They obtained approval for the first experimentation from the faculties of medicine and theology in Basel. The younger of the two brothers* [**Johann II**, 1710 – 1790] *and the only one of them married, decided to participate, and in 1756 he inoculated his two younger sons, and, a year ago, their elder brother*.

Johann II had indeed many sons, not all of them famous, whereas his brother was Daniel Bernoulli. An approval from theologians was really needed. White (1896/1898) described the "warfare of science with theology" including, in vol. 2, pp. 55 – 59, examples of fierce opposition to inoculation (and, up to 1803, to vaccination of smallpox). Many thousands of Canadians perished in the mid-19[th] century only because, stating their religious belief, they had refused to be inoculated. White clearly distinguished between theology, the opposing force, and "practical" religion.

Condamine (1773) includes the text of a third memoir on the same subject (pp. 221 – 282) read in 1764 and the author's correspondence, in particular with Daniel Bernoulli. He had given Bernoulli data on smallpox epidemics which the latter used in his research (below).

Karn stated at the very beginning of her article that

*The method used in this paper for determining the influence of the death-rates from some particular diseases on the duration of life is based on suggestions which were made in the first place by D. Bernoulli.*

Bernoulli (1766) justified inoculation. That procedure, however, spread infection, was therefore somewhat dangerous for the neighbourhood and prohibited for some time, first in England, then in France. Referring to statistical data, but not publishing it, Bernoulli suggested that $1/n$ was the yearly rate of the occurrence of smallpox in those who have not had it before; that $1/m$ was the corresponding mortality; that $m = n = 8$ and that the inoculation itself proved fatal in 0.5% of cases.

He formed the appropriate differential equation whose solution

$$s = \frac{m}{1+(m-1)e^{x/n}}$$

showed the relation between age $x$ (in years) and the number of people of the same age, , of which $s$ had not contacted smallpox. Also by means of a differential equation he derived a similar formula for a population undergoing inoculation, that is, for its 99.5% which safely endured it and were not anymore susceptible to the disease. It occurred that inoculation lengthened the mean duration of life by 3 years and 2 months and that it was therefore, in his opinion, extremely useful. The **Jennerian** vaccination, –

*The inestimable discovery by Jenner, who has thereby become one of the greatest benefactors of mankind* (**Laplace** 1814/1995, p. 83), –

was introduced at the end of the 18th century. Its magnificent success had not however ruled out statistical studies. Thus, Simon (1887, vol. 1, p. 230) formulated a question about the impermanence of protection against post-vaccinal smallpox and concluded that only comprehensive national statistics can provide an answer.

**D'Alembert** (1761b; 1768c) criticized **Daniel Bernoulli**[9]. Not everyone will agree, he argued, to lengthen his mean duration of life at the expense of even a low risk of dying at once of inoculation; then, moral considerations were also involved, as when inoculating children. Without denying the benefits of that procedure, D'Alembert concluded that statistical data on smallpox should be collected, additional studies made and that the families of those dying of inoculation should be indemnified or given memorial medals.

He also expressed his own thoughts, methodologically less evident but applicable to studies of even unpreventable diseases. Dietz & Heesterbeek (2002) described Bernoulli's and D'Alembert's investigations on the level of modern mathematical epidemiology and mentioned sources on the history of inoculation. For his part, **K. Pearson** (1978, p. 543) stated that inoculation was

*Said to have been a custom in Greece in the 17th century and was advocated* […] *in the* <u>Phil. Trans. of the Royal Society</u> *in 1713.*

Also see Sheynin (1972a/1977, pp. 114 – 116; 1982, pp. 270 – 272).

**6.2.4. Meteorology.** In § 2.1.4 I noted that **Leibniz** recommended regular meteorological observations**.** And, indeed (Wolf 1935, p. 312),

*Observations of barometric pressure and weather conditions were made at Hanover, in 1678, and at Kiel, from 1679 to 1714, at the instigation of Leibniz.*

The *Societas meteorologica Palatina* in Pfalz (a principality in Germany) was established in 1780, and, for the first time in the history of experimental science, it organized cooperation on an international scale (Sheynin 1984b, § 3.1). At about the same time the *Société Royale de Médecine* (Paris) conducted observations in several European countries (Kington 1974). And even in the 1730s – 1740s they were carried out in several towns in Siberia in accordance with directions drawn up by **Daniel Bernoulli** in 1733 (Tikhomirov 1932). In the second half of the 18th century several scholars (the meteorologist **Cotte**, **Lambert** and **Condorcet**) proposed plans for comprehensive international meteorological studies.

**Lambert** (1773) studied the influence of the Moon on the air pressure. **Daniel Bernoulli**, in a letter of 1759 (Radelet de Grave et al 1979, p. 62), encouraged him and noted that since the Moon's distance was not constant, its influence on the air ought to vary in time, cf. § 7.2-8. Then, **Toaldo** (1775; 1777) statistically studied the connections between phenomena concerning meteorology and stated that the weather depended on the configurations of the Moon. His opinion was not abandoned until the mid-19$^{th}$ century (Muncke 1837, pp. 2052 – 2076), but either then, or later, in the second half of that century, for example when the connection between cyclones and solar activity had been studied (Sheynin 1984a, § 4.2), no embryo of correlation theory was established, see § 10.7.

**Lamarck**, the most eminent biologist of his time, seriously occupied himself with physics, chemistry and meteorology. In meteorology, his merits had for a long time been ignored (Muncke 1837), but he is now remembered for his "pioneer work in the study of weather" (Shaw & Austin 1942, p. 130) and I (Sheynin 1984b, § 6) quoted several of his important pronouncements. He repeatedly applied the term *météorologie statistique* (e.g., 1800 – 1811, t. 4, p. 1) whose aim (Ibidem, t. 11, p. 9 – 10) was the study of climate, or, as he (Ibidem, t. 4, pp. 153 – 154) maintained elsewhere, the study of the climate, of regularities in the changes of the weather and of the influence of various meteorological phenomena on animals, plants and soil.

### 6.3. Mathematical Treatment of Observations

In modernity, mathematical treatment of observations became necessary after regular astronomical observations (**Tycho Brahe**, § 1.2.2) had begun. A new problem of natural sciences, the determination of the figure and the size of the Earth (of the Earth's ellipsoid of rotation), presented itself in the second half of the 18$^{th}$ century. By means of meridian arc measurements the lengths of those arcs were calculated (indirectly, by triangulation). After determining the length of one degree of the meridian in two different and observed latitudes it becomes possible to calculate both parameters of the ellipsoid whereas redundant measurements lead to equations of the type of (1.2) in these unknowns which can then be derived more precisely[10].

The term "Theory of errors" (*Theorie der Fehler*) is due to **Lambert** (1765a, Vorberichte and § 321) who defined it as the study of the relations between errors, their consequences, circumstances of observation and the quality of the instruments. He isolated the aim of the "Theory of consequences" as the study of functions of observed (and error-ridden) quantities. In other words, he introduced the determinate error theory (Note 2 to § 0.3) and devoted to it §§ 340 – 426 of his contribution. Neither **Gauss**, nor **Laplace** ever used the new terminology, but **Bessel** (1820, p. 166; 1838b, § 9) applied the expression "theory of errors" without mentioning anyone and by the mid-19$^{th}$ century it became generally known. As far as that theory is concerned, **Lambert** was Gauss' main predecessor (see § 6.3.1).

I shall separately consider the adjustment of direct and indirect measurements; note, however, that scientists of the 18$^{th}$ century recognized the common character of these problems. Thus, in both cases the unknowns were called by the same term, "Mittel" (Lambert 1765b, § 6) or "milieu" (**Maire** & **Boscovich** 1770, pp. 484 and 501), also see the method of averages (§ 6.3.2).

**6.3.1. Direct Measurements**. The first to touch on this case was **Cotes** (1722, see Gowing 1983, p. 107). Without any justification he advised to regard the weighted arithmetic mean, which he compared with the centre of gravity of the system of points, – of the observations,– as the "most probable" estimator of the constant sought:

*Let p be the place of some object defined by observation, q, r, s the places of the same object from subsequent observations. Let there also be weights P, Q, R, S reciprocally proportional to the displacements arising from the errors in the single observations, and which are given by the limits of the given errors; and the weights P, Q, R, S are conceived as being placed at p, q, r, s, and their centre of gravity Z is found; I say the point Z is the most probable place of the object.*

Cotes appended a figure (perhaps representing a three-dimensional picture) showing nothing except these four points. He had not explained what he meant by *most probable*, nor did he exemplify his rule. Nevertheless, his authority apparently gave support to the existing common feeling (§ 1.2.4). Without mentioning Cotes and putting forth qualitative considerations, **Condamine** (1751, p. 223) recommended to apply it. Then, **Laplace** (1814/1995, p. 121) stated that "all calculators" followed the Cotes rule. Elsewhere Laplace (1812/1886, pp. 351 – 353) remarked that astronomers had begun to follow Cotes after **Euler** (1749), but even before Cotes **Picard** (1693/1729, pp. 330, 335, 343) called the arithmetic mean *véritable*.

**T. Simpson** (1756) applied, for the first time ever, stochastic considerations to the adjustment of measurements; for that matter, he made use of generating functions. He aimed to refute some authors (left unnamed) who had maintained that one good observation was as plausible as the mean of many of them, cf. § 1.2.2. Simpson assumed that the chances of observational errors

$$-v, -v+1, \ldots, -2, -1, 0, 1, 2, \ldots, v-1, v$$

were equal [proportional] either to

$$r^{-v}, r^{-v+1}, \ldots, r^{-2}, r^{-1}, 1, r, r^2, \ldots, r^{v-1}, r^v$$

or to

$$r^{-v}, 2r^{-v+1}, \ldots, (v-1)r^{-2}, vr^{-1}, (v+1), vr, (v-1)r^2, \ldots, 2r^{v-1}, r^v.$$

He assumed that the observational errors obeyed some density law (taking $r = 1$ he thus introduced the uniform and the triangular discrete distributions) and his was the first actual introduction of random errors.

Denote the observational errors by $\varepsilon_i$, and by $N$, the number of some chances. Then, as Simpson noted,

$N(\varepsilon_1 + \varepsilon_2 + \ldots + \varepsilon_n = m)$ was the coefficient of $r^m$ in the expansions of

$$(r^{-v} + \ldots + r^0 + \ldots + r^v)^n = r^{-vn}(1-r)^{-n}(1-r^{2v+1})^n,$$

$$(r^{-v} + 2r^{-v+1} + \ldots + (v+1)r^0 + \ldots + 2r^{v-1} + r^v)^n =$$
$$r^{-vn}(1-r)^{-2n}(1-r^{v+1})^{2n}.$$

The left sides of these two equalities were generating functions with unit coefficients in the first case, and coefficients

1, 2, …, $v + 1$, … 2, 1

in the second instance.

For both these cases Simpson determined the probability that the absolute value of the error of the arithmetic mean of $n$ observations was less than some magnitude, or equal to it[11]. Consequently, he decided that the mean was always [stochastically] preferable to a separate observation and he thus arbitrarily and wrongly generalized his proof. Simpson also indicated that his first case was identical with the determination of the probability of throwing a given number of points with $n$ dice each having $(v + 1)$ faces. He himself (1740, Problem No. 22), and earlier **Montmort** (§ 3.3.3), although without introducing generating functions, and **De Moivre** (1730, pp. 191 – 197) had studied the game of dice.

Soon Simpson (1757) reprinted his memoir adding to it an investigation of the continuous triangular distribution. He passed over to the continuous case by assuming that $|v| \to \infty$ leaving the magnitude $(m/n)/v$ constant. Here, the fraction in the numerator was the admissible error of the mean and $n$, as before, the number of observations. Simpson's graph however represented a finite $v$ and a continuous argument (the observational errors) and the curve of the error of the mean did not possess the distinctive form of the normal distribution.

Simpson naturally had no knowledge of the variance and the calculation of the probability that the error of the mean exceeded the error of a single observation occurred to be difficult (Shoesmith 1985b).

Without mentioning Simpson, **Lagrange** (1776a) studied the error of the mean for several other and purely academic distributions, also by applying generating functions (even for continuous laws, thus anticipating the introduction of characteristic functions). A possible though inadequate reason for leaving out Simpson was the heated dispute over priority between **De Moivre** and him. Lagrange apparently had not wanted to be even indirectly involved in it. De Moivre was a scholar of a much higher calibre (a fact clearly recognized by Simpson) and 43 years the senior. At least on several important occasions Simpson did not refer to De Moivre and, after being accused by the latter (1725; only in edition of 1743, p. xii) of *Mak*[ing] *a Shew of new Rules, and works of mine, appeal*[ed] *to all mankind, whether in his treatment of me* [of Simpson], *he has* [not] *discovered an air of self-sufficiency, ill-nature, and inveteracy, unbecoming a gentleman* (Simpson, posth. publ. 1775, p. 144).

See Pearson's opinion about him at the end of § 4.2.

Lagrange's memoir contained other findings of general mathematical interest. He was the first to use integral transformations, and, in Problem 6, he derived the equation of the multivariate normal distribution (**K. Pearson** 1978, p. 599). In his § 18 he introduced the term *courbe de la facilité des erreurs.* Also see Sheynin (1973a, § 2).

**Lambert** (1760, §§ 271 – 306) described the properties of "usual" random errors, classified them in accordance with their origin (§ 282), unconvincingly proved that deviating observations should be rejected (§§ 287 – 291) and estimated the precision of observations (§ 294), again lamely but for the first time ever. He then formulated an indefinite problem of determining a [statistic] that with maximal probability least deviated from the real value of the constant sought (§ 295) and introduced the principle of maximal likelihood, but not the term itself, for a continuous density (§ 303), maintaining, however (§ 306), that in most cases it will provide estimates little deviating from the arithmetic mean. The translator of Lambert's contribution into German left out all this material claiming that it was dated[12].

Lambert introduced the principle of maximum likelihood for an unspecified, more or less symmetric and [unimodal] curve, as shown on his figure, denote this curve by $(x - x_o)$ where $x_o$ was the sought [parameter of location]. Let the observations be $x_1, x_2, \ldots, x_n$, and, somewhat simplifying his reasoning, write his [likelihood function] as

$$(x_1 - x_o)\ (x_2 - x_o) \ldots (x_n - x_o).$$

When differentiating this function, Lambert had not indicated that the argument here was the parameter $x_o$, etc.

In a few years Lambert (1765a) returned to the treatment of observations. He attempted to estimate the precision of the arithmetic mean, but did not introduce any density and was unable to formulate a definite conclusion. He also partly repeated his previous considerations and offered a derivation of a density law of errors occurring in pointing an instrument (§§ 429 – 430) in accordance with the principle of insufficient reason: it was a semi-circumference (with an unknown radius) simply because there were no reasons for its "angularity".

**Johann III Bernoulli** (1785) published a passage from a manuscript of **Daniel Bernoulli** (1769/1997) which he had received in 1769 but which was written, as its author had told him, much earlier. There, Daniel assumed the density law of observational errors as a semi-ellipse or semi-circumference of some radius *r* ascertained by assigning a reasonable maximal error of observation and the [location parameter] equal to the weighted arithmetic mean with posterior weights

$$p_i = r^2 - (\bar{x} - x_i)^2. \tag{4}$$

Here, $x_i$ were the observations and $\bar{x}$, the usual mean. If required, successive approximations could have been made.

The first to apply weighted or generalized arithmetic means was Short (1763). This estimator demanded a subjective selection of weights and, moreover, it only provided a correction to the ordinary arithmetic mean which tended to vanish for even density functions.

In his published memoir Daniel Bernoulli (1778) objected to the application of the arithmetic mean which (§ 5) only conformed with an equal probability of all possible errors and was tantamount to shooting blindly[13]. Instead, he suggested [the maximum likelihood estimator of the location parameter] and supported his idea (§ 9) by indicating that, when one out of several possible and incompatible events had occurred, it should be thought that it was the event that possessed the highest probability.

Listing a few reasonable restrictions for the density curve (but adding to these the condition of its cutting the abscissa axis almost perpendicularly), he selected a semi-circumference with radius equal to the greatest possible, for the given observer, error. He then (§ 11) wrote out the [likelihood function] as

$$\{[r^2 - (x - x_1)^2]\,[r^2 - (x - x_2)^2]\,[r^2 - (x - x_3)^2]\,\ldots\}^{1/2},$$

where, in somewhat different notation, $x$ was the unknown abscissa of the centre of the semi-circumference, and $x_1, x_2, x_3, \ldots$, were the observations. Preferring, however, to calculate the maximum of the square of that function, Bernoulli thus left the semi-circumference for an arc of a parabola. He certainly had not known that the variance of the result obtained will therefore change.

For three observations his [likelihood equation], as it occurred, was of the fifth degree. Bernoulli numerically solved it in a few particular instances with some values of $x_1, x_2$ and $x_3$ chosen arbitrarily (which was admissible for such a small number of them). In turn, I present his equation as

$$\frac{x - x_1}{r^2 - (x - x_1)^2} + \frac{x - x_2}{r^2 - (x - x_2)^2} + \ldots = 0$$

so that the maximum likelihood estimate is

$$x_0 = \frac{[px]}{\sum p_i},\; p_i = \frac{1}{r^2 - (x_0 - x_i)^2} \tag{5; 6}$$

with unavoidable use of successive approximations. These formulas are lacking in Bernoulli's memoir although the posterior weights (6) were the inverse of the weights (4) from his manuscript. This fact heuristically contradicted his own preliminary statement about shooting skilfully. Neither would have astronomers of his time approved weights increasing towards the tails of a distribution. It is now known, however, that these are expedient in case of some densities. I also note that, according to Bernoulli, the properly normed density was

$$y = (3/4r^3) [r^2 - (x - x_o)^2], \quad x_o - r \leq x \leq x_o + r$$

and that the weights (6) should be corrected accordingly.

Daniel Bernoulli is invariably described as a mathematician and natural scientist who was able to achieve great success with comparatively little mathematics, and I fail to understand why he did not indicate formulas (5) and (6).

**Euler** (1778) commented on Bernoulli's memoir. He (§ 6) objected to the [principle of maximum likelihood] but his reasoning was not good enough. He argued that, on the one hand, the result of an adjustment should barely change whether or not a deviating observation was adopted, but, on the other hand, that the value of the [likelihood function] essentially depended on that decision.

Euler then (§ 7) remarked that, in general, there was no need

*To have recourse to the principle of the maximum, since the undoubted precepts of the theory of probability are quite sufficient to resolve all questions of this kind.*

He had not noticed that his requirement led to the choice of the median.

In the positive part of his commentary, Euler recommended, instead of the arithmetic mean, the estimate (5) with posterior weights (4) and he mistakenly assumed that Bernoulli had actually chosen these same weights. While developing his thoughts, and denoting the $n$ observations by $+a$, $+b$, $+c$, …, where

$$a + b + c + ... = 0, \qquad (7)$$

he formed the equation

$$nx^3 - nr^2x + 3Bx - C = 0,$$
$$B = a^2 + b^2 + c^2 + ..., \quad C = a^3 + b^3 + c^3 + ...,$$

from which the estimate $+ x$ should have been calculated with $x$ equal to its root least in absolute value. Condition (7) meant that the estimate sought was the closest possible to the arithmetic mean; Euler himself (§ 9) justified his choice of the root by noting that $x = 0$ as $r \to \infty$, that is, as $n \to \infty$, also see below.

Euler (§ 11) also remarked that estimate (5) with weights (4) can be obtained from the condition

$$[r^2 - (x_o - a)^2]^2 + [r^2 - (x_o - b)^2]^2 + [r^2 - (x_o - c)^2]^2 + ... = \max. \quad (8)$$

The magnitudes in parentheses are the deviations of observations from the estimate sought and their fourth powers are negligible so that condition (8) is equivalent to the requirement

$$(x_o - a)^2 + (x_o - b)^2 + (x_o - c)^2 + \ldots = \min, \qquad (9)$$

whence, in accordance with condition (7), follows the arithmetic mean. Condition (9) is heuristically similar to the principle of least squares (which in case of one unknown indeed leads to the arithmetic mean) and condition (8) with weights (4) resembles the **Gaussian** principle of maximum weight (of least variance). True, if the density of the observational errors is known (which was the unrealistic assumption of both Bernoulli and Euler), then other estimates can be better than the arithmetic mean, cf. the opinion of Poincaré (§ 11.2).

A small deviation from condition (9) does exist and it is easy to see that it is occasioned by inevitable deviations of the observations from the proposed (or tacitly assumed) symmetrical law. Bernoulli himself noted this fact when (see above) numerically adjusting several sets of three observations. So, I say once more that in actual fact Bernoulli proposed the general arithmetical mean.

In his last memoir **Daniel Bernoulli** (1780) separated, for the first time ever, observational errors into random (*momentanearum*) and systematic (*chronicarum*), although not for observations in general. Thus, he assumed an equal number of *slow* and *quick* swings of the pendulum and said nothing about the possible interdependence of consecutive swings and his systematic error was the simplest possible. But Daniel was the first to apply the normal distribution for calculations in a field belonging to the theory of errors, to use the probable error as a measure of precision in the same field and to isolate the two kinds of errors.

As I remarked in § 1.1.4, even ancient astronomers undoubtedly knew that some errors were *systematic*. And here is the opinion of **D. T. Whiteside** (private communication, 1972):

*Newton in fact* (*but not in explicit statement*) *had a precise understanding of the difference between random and structurally 'inbuilt' errors. He was certainly, himself, absorbed by the second type of 'inbuilt' error, and many theoretical models of differing types of physical, optical and astronomical phenomena were all consciously contrived so that these structural errors should be minimized. At the same time, he did, in his astronomical practice, also make suitable adjustment for 'random' errors in observation …*

I return to Bernoulli. Since he considered pendulums[14], he indicated that these errors acted proportionally to the square root of, and to the time itself respectively. Making use of his previous findings (§ 6.1.1, formula (1)), Bernoulli justified his investigation by the [normal distribution] which thus first occurred in the theory of errors, although only as a limiting law.

The number of vibrations of a seconds pendulum during a day is $2N = 86,400$; Bernoulli assumed that $(N + \mu)$ of them were slower, and $(N - \mu)$ faster than stipulated, with periods of $(1 + )$ and $(1 - )$ respectively. His simple pattern meant that the number of positive (say) errors possessed a symmetric binomial distribution and that the error of the pendulum accumulated after a large number of vibrations had a normal distribution.

In his previous work Bernoulli (1770 – 1771) noted that, for $N = 10,000$,

$$[2/\sqrt{N}]\int_0^{\tilde{}} \exp(-x^2 N)dx = \frac{1}{2}$$

if $\mu = 47.25$. Now, having $N = 43,200$, he obtained, for the same probability of 1/2,

$$\mu = 47.25 \sqrt{4.32} \quad 100.$$

It was this calculation that caused his conclusion (above) about the behaviour of random errors. Already in the 19$^{th}$ century, however, it became known that such errors can possess other laws of distribution (e.g., § 10.9.4).

Note also that Bernoulli came close to introducing the probable error; to recall (§ 2.2.2), **Huygens** discussed the probable duration of life. Bernoulli was also the first to introduce elementary errors. I do not however set high store by this fact; indeed, this notion is not necessary for proving the CLT. I conclude by remarking that Bernoulli had not investigated the more general pattern of an unequal number of the slower and the faster vibrations although it corresponded to the case of unequal probabilities of male and female births, also studied by him. Neither had he said anything about the possible dependence between the periods of successive vibrations.

**6.3.2. Indirect measurements.** Here, I consider the adjustment of redundant systems

$$a_i x + b_i y + \ldots + s_i = v_i, \; i = 1, 2, \ldots, n \tag{10}$$

in $k$ unknowns ($k < n$) and residual free terms $v_i$ (see § 1.2.1).

1) In case of two unknowns (cf. beginning of § 6.3) astronomers usually separated systems (10) into all possible groups of two equations each and averaged the solutions of these groups. In other words, if $(x_{ij}; y_{ij})$ is the solution of group $(i; j)$, $i, j = 1, 2, \ldots, n, i < j$, then, in accordance with this method of *combinations*, the final estimates of the unknowns were

$$x_o = (1/C_n^2) \; x_{ij}, \; y_o = (1/C_n^2) \; y_{ij}.$$

The residual free terms were thus neglected.

In 1757 and later **Boscovich** (Cubranic 1961, pp. 90 – 91; Maire & Boscovich 1770, pp. 483 – 484) applied this method but it did not satisfy him, see below. Interestingly enough, in the first case he (Cubranic 1961, p. 46) derived the arithmetic mean of four latitudinal differences in an unusual way: he first calculated the half-sums of all six pairwise differences and then took their mean. He apparently attempted to exclude, without changing the final result, the unavoidable systematic errors and thus to ensure a (qualitative) estimation of the order of random errors[15].

In the 19th century, it was discovered that the MLSq led to the same result as the method of combinations, although only if the particular solutions were appropriately weighted (**Whittaker** & **Robinson** 1949, p. 251).

2) For the case of three unknowns the method of combinations becomes unwieldy. In an astronomical context, **Mayer** (1750) had to solve 27 equations in three unknowns. He separated them into three groups of nine equations each, calculated three particular solutions (see below), and, finally, averaged them. The plausibility of the results thus obtained depended on the expediency of the separation and it seems (Stigler 1986, pp. 21 – 25) that Mayer had indeed made a reasonable choice. Being mostly interested in only one unknown, he included the equations with their greatest and smallest in absolute value coefficients in the first, and the second group respectively. Note also that Mayer believed that the precision of results increased as the number of observations, but in his time this mistake was understandable.

Mayer solved each group of equations under an additional condition

$$v_i = 0, \qquad (11)$$

where *i* indicates the number of an equation; if the first group includes the first nine of them, then $i = 1, 2, …, 9$. **Biot** (1811, pp. 202 – 203) testified that before the advent of the MLSq astronomers had always applied the method of Mayer and **Laplace** (1812/1886, pp. 352 – 353) stated almost the same.

In a letter of 1850 **Gauss** (W/Erg-5, No. 3/6, p. 90) remarked that **Mayer** had only calculated by means of primitive combinations. He referred to Mayer's manuscripts, but it is likely that Mayer's trick was almost the same in both cases. And Gauss himself, in an earlier letter of the same year (Ibidem, pp. 66 – 67), recommended a similar procedure for calibrating an aneroid.

Condition (11) determines the method of *averages* and **Lambert's** recommendation (1765b, § 20) about fitting an empirical straight line might be interpreted as its application. Lambert separated the points (the observations) into two groups, with smaller and larger abscissas, and drew the line through their centres of gravity. He employed a similar procedure when fitting curves by separating the points into several groups.

The method of averages was intuitively considered as conforming to the equal possibility of errors of each sign (**Maire** & **Boscovich** 1770, p. 501), and, apparently, as leading in case of one unknown to the arithmetic mean. See § 10.1 for its further history.

3) The Boscovich method. He (Maire & Boscovich, Ibidem) adjusted systems (10) under additional conditions

$$v_1 + v_2 + \ldots + v_n = 0, \quad |v_1| + |v_2| + \ldots + |v_n| = \min, \qquad (12; 13)$$

the first of which determined the method of averages. It can be allowed for by summing all the equations and eliminating one of the unknowns from the expression thus obtained. The mean (*milieu*), as Boscovich remarked, should be connected

*Par une certaine loi aux règles des combinaisons fortuites et du calcul des probabilités*[16].

He was unable, however, to explain how his conditions conformed to his aim.

Boscovich's second condition (13)[17] linked his method with the [median]. Indeed, his geometric adjustment of systems (10) consisted in constructing a straight line whose slope occurred to be equal to the median of some fractions. In other words: for meridian arc measurements systems (10) are

$$a_i x + y + s_i = v_i. \qquad (14)$$

After allowing for condition (12), we have

$$[a_i - (1/n)\ a_i]x + [s_i - (1/n)\ s_i] = 0.$$

Calculate the *n* values of *x* and choose as the estimate their median.

**Laplace** (§ 7.2-6) also made use of the Boscovich method.

4) The minimax method. According to it, systems (10) are solved under the additional condition

$$|v_{\max}| = \min,$$

with the minimum being determined from all possible and expedient solutions[18]. In § 2.1.4 I indicated that **Kepler** had apparently made use of some elements of this method (true, not even for algebraic equations). It does not ensure optimal, in any sense, results, but allows to check whether the theory, underlying the given system (10), is correct. Indeed, any other method of its solution will lead to a greater value of $|v_{\max}|$, the gap between theory and observation will be wider, and the correctness of the former might mistakenly be questioned.

Gusak (1961) described the history of the minimax method from 1778, when **Euler** had applied it to an important but, in my context, irrelevant study, to **Chebyshev**. However, Euler (1749) made use of the rudiments of that method much earlier. When solving systems of the type of (10), he compared only a few "solutions" with each other[19]. Then, **Lambert** (1765a, § 420) recommended the same method but owned that he did not know how to apply it "in a general manner and without many roundabout ways". **Laplace** (1789/1895, pp. 493, 496 and 506 and elsewhere) applied the minimax method for preliminary investigations, – for checking whether or not the results of meridian arc measurements and pendulum observations contradicted the theory according to which the Earth was an oblate ellipsoid of rotation. Since the method of minimax has no stochastic underpinning, I am not describing the appropriate algorithms introduced by Laplace; I note, however, that it is applied in the theory of statistical decision-making (Lehmann 1959, Chapter 9).

5) Euler (1749, 1755, 1770) had to treat indirect measurements as well. At least in the first two instances his goal was much more difficult than that outlined in § 1.2.1 where the underlying theory was supposed to be known for sure. Concerning the first of his contributions, Wilson (1980, p. 262n) remarked that Euler was

*Stymied by the finding that, for certain of the variables, the equations led to <u>wildly</u> different values, so that averages* [values adjusted in accordance with some definite rule] *would seem meaningless*.

Euler did not attempt to build a more or less general theory; he rather restricted his efforts to achieving practical results. In accordance with the considerations above, he turned in some cases to the minimax principle. On another important occasion Euler's attitude (1770) is difficult to explain: he did not keep to any definite method and, moreover, combined equations in a doubtful manner. So as to eliminate one unknown, he subtracted each equation from (say) the first one, thus tacitly regarding it as possessing much more weight. See Sheynin (2007b, § 3.5). This procedure seems to have been necessary. Ivory (§ 10.9.1) acted the same way.

**6.3.3. The True Value of a Measured Constant** (Sheynin 2007a). Astronomers, geodesists, metrologists and other specialists making measurements have always been using this expression. Mathematical statistics has done away with true values and introduced instead parameters of densities (or distribution functions), and this was a step in the right direction: the more abstract was mathematics becoming, the more useful it proved to be.

**Fisher** was mainly responsible for that change; indeed, he (1922, pp. 309 – 310) defined the notions of consistency, efficiency and sufficiency of statistical estimators without any reference to true values. But then, on p. 311, he accused the Biometric school of applying the same names to

*The true value which we should like to know* […] *and to the particular value at which we happen to arrive*.

So the true value was then still alive and even applied, as in the lines above, to objects having no existence in the real world. Incidentally, the same can be said about **Gauss** (1816, §§ 3 and 4) who repeatedly considered the true value of a measure of precision of observations. And **Hald** (1998) mentioned the *true value* many times in Chapters 5 and 6; on p. 91 he says: "the estimation of the true value, the location parameter…"

So what is a true value? **Markov** (1924, p. 323) was the only mathematician who cautiously, as was his wont, remarked:

*It is necessary in the first place to presume the existence of the numbers whose approximate values are provided by observations*.

This phrase first appeared in the 1908 edition of his *Treatise* (and perhaps in its first edition of 1900). He did not attempt to define *true value*, but this is exactly what **Fourier** (1826/1890, p. 534) had done more than a century before him. He determined the *véritable objet de la recherche* (the constant sought, or its "true" value) as the limit of the arithmetic mean of *n* appropriate observations as $n \to \infty$.

Many authors, beginning perhaps with Timerding (1915, p. 83) [and including **Mises** (1919/1964a, pp. 40 and 46)], without mentioning Fourier and independently from each other, introduced the same definition. One of them (Eisenhart 1963/1969, p. 31) formulated the unavoidable corollary: the mean residual systematic error had to be included in that "true" value:

*The mass of a mass standard is* […] *specified* […] *to be the mass of the metallic substance of the standard plus the mass of the average volume of air adsorbed upon its surface under standard conditions.*

However, even leaving systematic influences aside, the precision of observations is always restricted (§ 11.2-8) so that the term "limit" in the Fourier definition (which is in harmony with the Mises definition of probability) must not be understood literally. I indicate also that **Gauss** (W-9, pp. 278 – 281), see also Schreiber (1879, p. 141) measured each angle in the field until becoming convinced that further work was meaningless.

The concept of true value is not, however, universally recognized. Chatterjee (2003, p. 264) mentioned the *true value syndrome* which was *ultimately left behind*.

### Notes

**1. Lagrange** (1777) solved a similar problem for a finite number of urns and balls of two colours as well as some other stochastic problems by means of partial difference equations.

**2.** He published many memoirs and papers on the theory of probability and its applications (§ 6.2.3) and it is difficult to organize them bibliographically; on this point see Paty (1988). **Todhunter** (1865) devoted an entire chapter to **D'Alembert.**

**3.** Cf. the **D'Alembert – Laplace** problem (Note 4 in Chapter 1). In 1750 D'Alembert declared that randomness was only caused by ignorance (Note 2 in Chapter 1). The denial of randomness, also upheld by **Kepler** (§ 1.2.4) and **Laplace** (§ 7.3), although only by mere words, proved fruitless.

**4.** Regarding his really strange attitude towards medicine see Note 7.

**5.** His unsubstantiated conclusion was absolutely wrong. Loveland (2001) attempted to reconstruct **Buffon's** reasoning.

**6.** Recall however (§ 3.3.4) that **Condorcet** reasonably remarked on the Petersburg game.

**7. D'Alembert** (1759/1821, p. 163) should also be mentioned. The first edition of this contribution published in 1759 apparently had not contained any such statement. Note, however, that he died in 1783 so that he formulated his similar desire in the 18[th] century. D'Alembert even stated that a physician was a blind man who can strike either the disease or the patient by his club and added, on p. 167, that the best doctor was the one who least believed in medicine. See also § 10.8.

**8.** The same author, Black, appended a "Chart of all the fatal diseases and casualties in London during […] 1701 – 1776" to his book. It will warn us *to make* *the best disposition and preparation for defence.* However, in his previous book Black (1782) stated contradictory ideas.

**9.** In the first case he discussed **Bernoulli's** report; I stress that the latter's memoir appeared only in 1766. Later **D'Alembert** rewrote his memoirs. See **Todhunter** (1865, pp. 265 – 271, 277 – 278 and 282 – 286) for a detailed description of his proposals.

**10.** Not the semiaxes of the ellipsoid, $a$ and $b$ ($a > b$), were determined, but rather $a$ and the flattening $(a - b)/a$. The flattening had also been derived from pendulum observations; see § 10.9.1 where I describe the pertinent work of **Ivory**.

**11.** The distributions introduced by **Simpson**, if considered continuous, can be directly compared with each other in the sense that the respective variances are $v^2/3$ and $v^2/6$.

**12**. In a letter of 1971 **E. S. Pearson** informed me that "curiously" his father's *Lectures* (1978), – then not yet published, – omitted **Lambert**. He explained:

*It was not because* [Lambert's] *writings were in German of which my father was an excellent scholar. I suppose […] that he selected the names of the personalities he would study from a limited number of sources, e.g., Todhunter, and that these did not include Lambert's name.* [**Todhunter** did refer to Lambert but had not described his work.] *Of course, K. P. was over 70 by the time his history lectures passed the year 1750, and no doubt his exploration was limiting itself to the four Frenchmen, Condorcet D'Alembert, La Grange and Laplace.*

**13.** Here, however, is **K. Pearson's** reasonable qualitative statement (1978, p. 268): small errors are more frequent and have their due weight in the mean.

**14.** For this reason his memoir was attributed to practical mechanics and until my publication (Sheynin 1972a) its stochastic nature had not been noticed.

**15. Tycho's** example (Note 25 in Chapter 1) is more convincing.

**16.** The last term deserves attention: it was hardly used before **Boscovich**.

**17. Galileo** (§ 1.2.3) and **Daniel Bernoulli** (1735/1987, pp. 321 – 322) applied this condition in the case in which the magnitudes such as $v_i$ were positive by definition. The latter derived the plane of the solar equator in such a way that the sum of the inclinations of the planetary orbits, considered positive, relative to the equator, was minimal. **W. Herschel** (1805) determined the movement of the Sun by issuing from the apparent motion of the stars. The sum of these motions depends on the former and its minimal value, as he assumed, provided an expedient estimation of that movement, cf. Kepler's similar method at the end of § 1.2.4.

Note that in those times the motion of a star could have been discovered only in the plane perpendicular to the line of vision. Here is W. Herschel's earlier reasoning (1783/1912, vol. 1, p. 120):

*We ought […] to resolve that which is common to all the stars […] into a single real motion of the Solar system, as far as that will answer the known facts, and only to attribute to the proper motions of each particular star the deviations from the general law the stars seem to follow …*

Such, he added, were "the rules of philosophizing". Compare now **Newton's** Rule No. 1 of reasoning in philosophy (§ 1.1.5).

When treating direct measurements W. Herschel (1806) preferred the [median] rather than the arithmetic mean (Sheynin 1984a, pp. 172 – 173).

**18.** It is remarkable that the minimax method corresponds, as **Gauss** (1809b, § 186) noted, to the condition

$$\lim(v_1^{2k} + v_2^{2k} + ... + v_n^{2k}) = \min, k$$

**19.** Stigler (1986, pp. 27 – 28) called **Euler's** memoir (1749) a "statistical failure" and, in his opinion, Euler was a mathematician who "distrusted" the combination of equations. Without perceiving the main goal of the method of minimax, and mentioning a classic in a free and easy manner, Stigler got into a mess. See Wilson's statement at the end of § 6.3.2. In his second book Stigler (1999, pp. 317 – 318) unblushingly called Euler a great statistician but did not notice his lame argument against the [method of maximum likelihood] (§ 6.3.1). Moreover, as stated there, in the particular instance under discussion Euler should have opted for the not yet named median (which he did not mention) rather than the mean (which he preferred).

For that matter, in the 18[th] century practitioners experienced difficulties when deciding how to adjust their observations (Bru 1988, pp. 225 – 226); and at the turn of that century **Laplace** and **Legendre** simply refused to adjust a triangulation chain laid out between two baselines. Instead, likely fearing the propagation of large errors, they decided to calculate each half of the chain starting from its own baseline (Sheynin 1993b, p. 50). Much later Laplace (ca. 1819/1886, pp. 590 – 591) defended their decision by the previous ignorance of the "vraie théorie" of adjustment and added that his justification of the MLSq had changed the situation.

I supplement Bru's description by indicating that **Maupertuis** (1738, p. 160; 1756b, pp. 311 – 319) calculated his triangulation twelve times (each time taking into account differing sets of measured angles), selected two of his results and adopted their mean value.

It is instructive to note that, before the adjustment proper of the Soviet primary triangulation, each of its chains situated between baselines and astronomically determined azimuths was replaced by an appropriate geodetic line (cf. beginning of § 9B). Only these lines were then adjusted after which each chain was finally dealt with independently from one another. One of the benefits of this procedure was that it prevented the systematic errors from "freely walking" over the entire network, as **Izotov**, the leading assistant of **Krasovsky**, the calculator of the *Krasovsky spheroid*, explained ca. 1950 in one of his lectures at the Moscow Geodetic Institute which I attended.

# 7. Laplace
## 7.1. Theory of probability

Laplace devoted a number of memoirs to the theory of probability and later combined them in his *Théorie analytique des probabilités* (abbreviation: TAP) (1812). When referring to it, I often indicate only the page numbers. I describe its second *Livre*; in the first one he studied the calculus of generating functions with application to the solution of ordinary and partial difference equations and the approximate calculation of integrals.

1) In **Chapter 1** Laplace provided the "classical" definition of probability (introduced by **De Moivre**, see § 4.3), formulated the addition and multiplication theorems for independent events as well as theorems concerning conditional probabilities. He described the same material in his *Essai philosophique …*[1] where he (1814/1995, p. 10), in addition, included the so-called **Bayes** theorem, see formula (5.1), calling it a principle. Much earlier he (1774/1891, p. 29) introduced a "fundamental principle", – the same theorem for the case of constant prior probabilities $P(A_i)$:

$$P(A_i/B)/P(A_j/B) = P(B/A_i)/P(B/A_j).$$

2) In **Chapter 2** Laplace solved a number of problems by means of difference, and partial difference equations. I consider three other problems.

a) In an astronomical context Laplace studied sampling with replacement. Tickets numbered from 0 to $n$ are extracted from an urn. Determine the probability that the sum of $k$ numbers thus extracted will be equal to $s$ (p. 257). Let these numbers be $t_1, t_2, …, t_k$, then

$$t_1 + t_2 + … + t_k = s. \qquad (1)$$

Laplace calculated the number of combinations leading to equality (1) allowing for the condition $t_i \quad n, i = 1, 2, ..., k$ by assigning to these $t_i$ probabilities

$$(1 - l^{n+1})/(n + 1) \qquad (2)$$

with $l = 0$ for $t_i \quad n$ and $l = 1$ otherwise. Earlier, I (Sheynin 1973a, pp. 291 – 298) discussed Laplace's use of discontinuity factors in somewhat more detail. I also described his similar method which he applied in 1810 and which dates back to **De Moivre** and **Simpson** (Ibidem, pp. 278 – 279). If, for example, two (three) of the *t*'s exceed *n*, that factor in (2) is raised to the second (to the third) power etc.

Laplace calculated the probability sought and considered the case of $s, n \quad$ and his formula on p. 260 for the distribution of the sum of independent, continuous variables obeying the uniform law on interval [0; 1] corresponds with modern literature (**Wilks** 1962, § 8.3.1) which does not, however, demand large values of *s* and *n*.

Also in an astronomical context, already in 1776, Laplace solved a problem concerning such distributions by very complicated recursion relations (Sheynin 1973a, pp. 287 – 290). Note, however, that even **Simpson** and **Lagrange** (§ 6.3.1) obtained similar findings in the theory of errors.

Laplace treated the two other problems alluded to above in the same way as he did earlier in 1781.

b) Non-negative [random variables] $t_1, t_2, …, t_k$ with differing laws of distribution $\varphi_i(t)$ are mutually independent and their sum is $s$. Determine the integral

$$\int \psi(t_1; t_2; …; t_k) \varphi_1(t) \varphi_2(t) … \varphi(t) dt_1 dt_2 … dt_k$$

over all possible values of the variables; $\psi$ is not yet chosen. Laplace then generalizes his very general problem still more by assuming that each function $\varphi_i(t)$ can be determined by different formulas on different intervals of its domain.

When solving this problem, Laplace made use of the same discontinuity factor as above and derived the **Dirichlet** formula (an expression for the $n$-tuple integral of the product of $n$ power functions over the area in which the sum of the arguments, i. e., of the bases of those functions, is restricted by the interval [0, 1]), even in a more general version. The case of $\psi \equiv 1$ enabled him to determine the probability of equality (1) (which interested Laplace here also). He then once more specified his problem by assuming that

$$\varphi_i(t) = a + bt + ct^2.$$

When solving that problem, Laplace derived a multiple integral of $u, u_1, u_2, …$ over the area

$$0 \le u + u_1 + u_2 + … \le s$$

and differentiated it with respect to $s$, – with respect to that area! He had not mentioned that he calculated that derivative rather than the integral and only provided the final answer. I note that a simple transformation $u = sx$, $u_1 = sx_1$, $u_2 = sx_2$, … saves us from that unusual differentiation (Sheynin 1973a, p. 292).

c) An interval $OA$ is divided into equal or unequal parts and perpendiculars are erected to the interval at their ends. The number of perpendiculars is $n$, their lengths (moving from $O$ to $A$) form a non-increasing sequence and the sum of these lengths is given. Suppose now that the sequence is chosen repeatedly; what, Laplace asks, would be the mean broken line connecting the ends of the perpendiculars? The mean value of a current perpendicular? Or, in the continuous case, the mean curve? Each curve might be considered as a realization of a stochastic process and the mean curve sought, its expectation. Laplace was able to determine this mean curve (Sheynin 1973a, p. 297) by issuing from his previous problem[2] and, in 1781, he attempted to apply this finding in the theory of errors (§ 7.2) and for studying expert opinions. Suppose that some event can occur because of $n$ mutually exclusive causes. Each expert arranges these in an increasing (or decreasing) order of their [subjective] probabilities, which, as it occurs, depend only on *n* and the number of the cause, *r*, and are proportional to

$$\frac{1}{n} + \frac{1}{n-1} + ... + \frac{1}{n-r+1}.$$

The comparison of the sums of these probabilities for each cause allows to show the mean opinion about its importance. To be sure, different experts will attribute differing perpendiculars to one and the same cause.

3) The **third Chapter** is devoted to the integral "**De Moivre – Laplace**" theorem and to several interesting problems connected with the transition to the limit. In proving that theorem (§ 4.4) Laplace applied the **Euler – MacLaurin** summation formula, and, a second innovation, calculated the remainder term to allow for the case of large but finite number of trials. His formula was ( $= \sqrt{n/2xx'}$ ):

$$P(|\mu - np - z| \leq l) =$$

$$(2/\sqrt{\pi}) \int_0^{\cdot} \exp(-t^2)dt + \sqrt{1/\pi} \exp(-l^2 n/2xx'). \qquad (3)$$

Here *p* was the probability of success in a single **Bernoulli** trial, $\mu$, the total number of successes in *n* trials, $q = 1 - p$, *z* was unknown but $|z| < 1$, $x = np + z$, and $x' = nq - z$.

Laplace indicated that his theorem was applicable for estimating the theoretical probability given statistical data, cf. the **Bayes** theorem in § 5.2, but his explanation was not clear, cf. **Todhunter** (1865, pp. 554 – 556). Insufficiently clear is also **Hald's** description (1990, § 24.6).

Already **Daniel Bernoulli** (§ 6.1.1) solved one of Laplace's problem: There are two urns, each containing *n* balls, some white and the rest black; on the whole, there are as many white balls as black ones. Determine the probability *u* that the first urn will have *x* white balls after *r* cyclic interchanges of one ball. The same problem was solved by **Lagrange** (1777/1869, pp. 249 – 251), **Malfatti** (**Todhunter** 1865, pp. 434 – 438) and Laplace (1811; and in the same way in the TAP).

Laplace worked out a partial difference equation, "mutilated it most unsparingly" (**Todhunter** 1865, p. 558), and obtained a partial differential equation with argument *x*

$$u_{r/n} = 2u + 2\mu u_\mu + u_{\mu\mu}, \text{ for } x = (n + \mu\sqrt{n})/2$$

and expressed its solution in terms of functions related to the [**Chebyshev** –] **Hermite** polynomials (Molina 1930, p. 385). Hald (1998, p. 339) showed, however, that Todhunter's criticism was unfounded.

Later **Markov** (1915b) somewhat generalized this problem by considering the cases of $n \to \infty$ for $r/n \to \infty$ and $n \to \infty$ and $r/n = $ const and **Steklov** (1915) proved the existence and uniqueness of the solution of Laplace's differential equation with appropriate initial conditions added whereas **Hald** (2002) described the history of those polynomials. **Hostinský** (1932, p. 50) connected Laplace's equation with the **Brownian** motion and thus with the appearance of a random process (Molina 1936).

Like **Bernoulli**, Laplace discovered that in the limit, and even in the case of several urns, the expected (as he specified on p. 306) numbers of white balls became approximately equal to one another in each of them. He also remarked that this conclusion did not depend on the initial distribution of the balls. Finally, in his *Essai* (1814/1995, p. 42), Laplace noted that nothing changed if new urns, again with arbitrary distributions of balls, were added to the original urns. He declared, apparently too optimistically, that

*These results may be extended to all naturally occurring combinations in which the constant forces animating their elements establish regular patterns of action suitable to disclose, in the very mist of chaos, systems governed by these admirable laws.*

Divine design was absent, cf. **De Moivre's** dedication of his book to **Newton** in § 4.3.

The **Daniel Bernoulli** – Laplace problem coincides with the celebrated **Ehrenfests'** model (1907) which is usually considered as the beginning of the history of stochastic processes so that its formulation is really unnecessary. As I noted in § 6.1.1c, the existence of the limiting state in this problem can be justified by the **Markov** ergodic theorem for Markov chains. See also Sheynin (1972a/1977, pp. 127 – 128).

4) I touch on **Chapter 4** in § 7.2-4. Laplace devoted **Chapter 5** to the detection of constant causes (forces) in nature. Thus, he attempted to estimate the significance of the daily variation of the atmospheric pressure. **K. Pearson** (1978, p. 723) noted that nowadays the **Student** distribution could be applied in such investigations, that some of the assumptions which Laplace made proved wrong, etc. and, in addition, that Laplace had unjustifiably rejected those days during which the variation exceeded 4 *mm*.

Laplace remarked that the *calcul des probabilités* can be applied to medicine and economics. It may be argued that he thought about stochastic analysis of statistical data, see his *Essai* (1814/1995, p. 61). He even inserted a section on the application of the probability calculus to moral sciences.

Concerning geometric probability, Laplace only discussed the **Buffon** problem and stated (p. 365) that geometric probability can be applied to rectify curves "ou carrer leurs surfaces" (nowadays this *carrer* certainly looks strange).To repeat (§ 6.1.6), a needle of length $2r$ falls from above on a set of parallel lines. The distance between adjacent lines is $a \geq 2r$ and the probability $p$ that the needle intersects a line is

$p = 4r/a.$

Without proof Laplace mistakenly stated that, for $a = 1$, $2r = \pi/4$ was the optimal length of the needle for statistically determining $\pi$ although he provided the correct answer, $2r = 1$, in the first edition of the TAP.

Both **Todhunter** (1865, pp. 590 – 591) and Gridgeman (1960) examined this problem by applying stochastic considerations, and both proved that Laplace's previous answer was correct.

5) In **Chapter 6** Laplace solved some problems by means of the **Bayes** approach (see § 5.1) although without referring to him; true, he mentioned Bayes elsewhere (1814/1995, p. 120). Here is one of them.

Denote the unknown probability that a newly born baby is a boy by $x$ and suppose that during some time $p$ boys and $q$ girls were born. Then the probability of that "compound" event will be proportional to

$$y = x^p(1 - x)^q. \tag{4}$$

If $z(x)$ is the prior distribution of $x$, then

$$P(a \leq x \leq b) = \int_a^b yz\, dx \div \int_0^1 yz\, dx,\ 0 < a < b < 1. \tag{5}$$

If, as Laplace nevertheless assumed, $z$ was constant, and if $p$ and $q$ were large, the probability sought will be expressed by an integral of an exponential function of a negative square.

And so, Laplace actually estimated the probability $x$. For the curve (4) the point of its maximum

$$\xi = p/(p + q) \tag{6}$$

seems to be its natural estimator, but $E\xi$, or, more precisely, the expectation of a random variable $\xi$ with distribution

$$x^p(1 - x)^q \div \int_0^1 x^p(1 - x)^q\, dx,$$

does not coincide with (6): the latter is only an asymptotically unbiased estimator of $x$. This expectation is evidently

$$E\xi = \frac{p + 1}{p + q + 2}. \tag{7}$$

The introduction of functions $z(x)$ allowed to assume an equal probability for each value of $x$, but the choice of such functions remained undecided. Of course, nothing more could have been achieved.

Laplace went on to discuss the bivariate case and then solved another problem. Suppose that the inequality $p > q$ persisted during a number of years. Determine the probability that the same will happen for the next hundred years. There is no doubt that Laplace understood that his problem only made sense under invariable social and economic conditions. Here is his answer:

$$P = \int_0^1 x^p(1-x)^q z^{100} dx \div \int_0^1 x^p(1-x)^q dx$$

where $z$ is the sum of the first $n$ terms of the development

$[x + (1-x)]^{2n}$ and $2n = p + q$.

A similar problem in which $q = 0$, $p = m$ and $z = x^n$ led Laplace to the probability of such $z$:

$P = (m + 1) \div (m + n + 1)$.

In the *Essai* Laplace (1814/1995, p. 11) applied this formula, slightly different from his previous formula (7), for solving **Price's** problem about the next sunrise (§ 5.1) but he only mentioned **Buffon** (§ 6.1.4), and, as expected, did not agree with his solution.

I consider this problem in more detail. Recall (§ 5.1) that one of Price's results was doubtful. Fries (1842/1974, p. 7 and 188/140) indicated that $P$ tends to disappear as $n$ and that therefore the phenomenon described by the Laplace formula (he did not mention the sunrise) cannot be a law of nature. He concluded that, when issuing from repeated observations and posterior probability, it was impossible to *guess* the prior probability and he thus hardly recognized the Bernoulli theorem.

Fries decided that the solution of the Price problem just as the Laplacean solution of the D'Alembert – Laplace problem (Note 5 to Chapter 1) is based on philosophical (not mathematical) probabilities and philosophical induction. Finally and without any special study Fries (p. 157/139) declared that the MLSq was *completely subjective*. All this he briefly stated already in his Introduction. It is exactly the Bernoulli theorem that answers Fries. Laplace (1814/1995, p. 10) actually remarked that his formula was necessarily justified by the principle of insufficient reason but with the increase of the number of trials the accepted assumption should have been gradually specified. He himself (§ 7.2-1) recommended a similar procedure.

Zabell (1989) comprehensively studied the Price problem but all but ignored Fries. I note that it was Fries who had effectively introduced philosophical probabilities before Cournot (1843) considered them in detail.

Back to my main subject. Laplace determined the population of France given sampling data, and, for the first time ever, estimated the precision of (his version of) sampling. Suppose that $N$ and $n$ are the known numbers of yearly births in France as a whole and in some of its regions and $m$ is the population of those regions. Laplace naturally assumed that

$M = (m/n)N$.

He then had to estimate its error, the fraction (**Hald** 1998, p. 288)

$$\int_0^1 x^{N+n}(1-x)^{m-n+M-N}dx \div \int_0^1 x^n(1-x)^{m-n}dx.$$

**K. Pearson** (1928a) noted some imperfections in Laplace's reasoning and achieved a reduction of the variance of his result: it should have been multiplied by $[(N-n) \div (N+n)]^{1/2}$. Here are his two main remarks. First, Laplace considered $(m, n)$ and $(M, n)$ as independent samples from the same infinite population whereas they were not independent and the very existence of such a population was doubtful. Second, Laplace chose for the magnitude sought an absolutely inappropriate uniform prior distribution. Pearson also negatively described Laplace's calculation of the incomplete beta-function. However, he (1934, Intro.) also owned that that problem remained very difficult which thus actually exonerated Laplace.

Pearson's first remark had to do with Laplace's supplementary urn problem. Suppose that an urn contains infinitely many white and black balls. After $n$ drawings without replacement $m$ white balls were extracted; a second sample of an unknown volume provided $r$ white balls. Denoting

$k = nr/m + z,$

Laplace derived a limit theorem

$$P(|k - nr/m| < z) = 1 - 2\int \frac{m^3}{\sqrt{\pi S}} \exp(-m^3z^2/S)dz,$$

$S = 2nr(n-m)(m+r).$

The limits of integration, as Laplace formally assumed, were $z$ and .

Later **Markov** (1900b) proved that, for an unknown $m$,

$$P\left[\left|\frac{m}{n} - \frac{r}{k}\right| < \frac{t}{2}\sqrt{(1/k)+(1/n)}\right] > 1 - 1/t^2, t > 0.$$

He (1914a) then specified that all prior probabilities of the appearance of a white ball were equal to one another and proved, in addition, that the same inequality of the **Bienaymé – Chebyshev** type held also for "indefinite" [random] fractions $m/n$ and $r/k$. As though the last of the Mogicanes, Markov consistently refused to use the then new term, *random variable,* see § 14.2-1.

6) In **Chapter 7** Laplace studied the influence of a possible inequality of probabilities assumed equal to each other. For example, when tossing a coin the probability of heads can be $(1 \pm a)/2$ with an unknown *a*. Supposing that both signs were equally probable, Laplace derived the probability of throwing *n* heads in succession

$$P = (1/2) [(1 + a)^n + (1 – a)^n] \div 2^n$$

which for $n > 1$ was greater than $1/2^n$.

Suppose now (the general case) that the probability is not *p*, as assumed, but $(p + z)$, $|z| \leq a$, with density $\varphi(z)$. Then the probability of a "compound" event *y* will be

$$P = \int_{-a}^{a} y(p+z)\varphi(z)\,dz \div \int_{-a}^{a} \varphi(z)\,dz$$

(cf. formula (5) above). In case of an unknown density $\varphi(z)$ it should be replaced by the density of *z* (which he called the probability of *z*), Laplace adds. The appearance of denominators in such formulas seems to be unnecessary.

The special reasoning in this chapter, as also in one of the examples in Chapter 3, can be justified by considering Markov chains. Actually it is tantamount to stating that an infinite shuffling of a deck of cards results in an equal probability of all of their possible arrangements (**Feller** 1950, § 9 of Chapter 15).

Laplace then considers tickets put into an urn. Suppose, he says, that the probabilities of their extraction are not equal to one another. However, the inequalities will be reduced had the tickets been put into the urn not in an assigned order, but according to their random extraction from an auxiliary urn, and still more reduced in case of additional auxiliary urns. Laplace had not proved this statement, he only justified it by a general principle: randomness diminished when subjected to more randomness. This is perhaps too general, but Laplace's example was tantamount to reshuffling a deck of cards.

7) **Chapter 8** was devoted to population statistics, to the mean durations of life and marriages. Laplace did not apply there any new ideas or methods. However, he studied anew the **Daniel Bernoulli** model of smallpox (§ 6.2.3), adopted more general assumptions and arrived at a more general differential equation (**Todhunter** 1865, pp. 601 – 602).

8) In **Chapter 9** Laplace considered calculations made in connection with annuities and introduced the "**Poisson**" generalization of the **Jakob Bernoulli** theorem (as noted by Molina 1930, p. 372). Suppose that two contrary events can occur in each independent (as he clearly indicated) trial *i* with probabilities $q_i$ and $p_i$ respectively, $q_i + p_i = 1$, $i = 1, 2, …, s$, and that these events signify a gain $\lambda$ and a loss µ, respectively. For constant probabilities *q* and *p* the expected gain after all these trials will be $s(q\lambda – p\mu)$, as Laplace for some reason concluded in a complicated way. He then estimated this magnitude for the case of a large *s* by means of his limit theorem (3). Then, generalizing the result obtained to variable probabilities, he introduced the characteristic function of the final gain

$[p_1 + q_1\exp(\alpha_1 i)] [p_2 + q_2\exp(\alpha_2 i)] \ldots [p_s + q_s\exp(\alpha_s i)]$,

applied the inversion formula and obtained the normal distribution, all this similar to the derivation of the law of distribution of a linear function of observational errors (§ 7.2-4, see also his earlier memoirs).

9) In **Chapter 10** Laplace described his thoughts about moral expectation (§ 6.1.1). If the physical capital of a gambler is $x$, his moral capital will be

$y = k\ln x + \ln h$, $h, x > 0$.

Let $x$ take values $a, b, c, \ldots$ with probabilities $p, q, r, \ldots$ Then

$Ey = k[p\ln(x + a) + q\ln(x + b) + \ldots] + \ln h$,

$$E\Delta y < E\Delta x. \qquad (8)$$

In other words, even a just game ($E\Delta x = 0$) is disadvantageous. **Todhunter** (1865, p. 215) proved inequality (8) simpler than Laplace did. However, a more general expression $Ef(x) \leq f(Ex)$ holds for convex functions (**Rao** 1965, § 1e5) so that, if $x > 0$,

$E(-\ln x) \geq -\ln Ex$, $E\ln x \leq \ln Ex < Ex$.

Laplace then proved that the freight in marine shipping should be evenly distributed among several vessels. I provide my own proof (Sheynin 1972a/1977, pp. 111 – 113). Suppose that the capital of the freightowner is $a$, the value of the freight, $A$, the probability of a safe arrival of a vessel, $p$, and $q = 1 - p$. Then

a) If the freight is thus distributed on $n$ vessels, the moral expectation of the freightowner's capital is (here and below $0 \leq k \leq n$ and $k$ is the number of the lost ships)

$$y(n) = \sum C_n^k p^{n-k} q^k \ln\{[A(n-k)/n] + a\}. \qquad (9)$$

b) Independently from $n$ the corresponding expectation is equal to the right side of (9) where, however, the logarithm is replaced by its argument, so that obviously

$a + A \sum [(n-k)/n]p^{n-k} q^k = a + Ap$.

c) For any increasing function $f(x)$ the moral expectation (9) is restricted:

$y(n) = \sum C_n^k p^{n-k} q^k f\{[A(n-k)/n] + a\} < f(A + a)(p + q)^n =$
$\quad f(A + a)$.

d) Let $f(x)$ be continuous and increasing and have a decreasing derivative. Then $y(n)$ increases monotonically but is restricted by the

moral expectation (9). The proof is here rather long and I refer readers to my paper (1972b/1977, pp. 112 – 113).

Many authors after Laplace dwelt on moral expectation (cf. § 6.1.1); **Fourier** (1819) and **Ostrogradsky** (**N. I. Fuss** 1836, pp. 24 – 25), also see Ostrogradsky (1961, pp. 293 – 294), attempted to develop it. The former thought that it should be individually specified for each person whereas Fuss reported that Ostrogradsky

*Did not at all admit Daniel Bernoulli's hypothesis; he expressed the "moral fortune" by an arbitrary function of the physical fortune and he was able to solve the main problems connected with the moral fortune with such breadth and as precisely as could only be desired.*

Note, however, that the logarithmic function also appears in the celebrated **Weber – Fechner** psychophysical law and is applied in the theory of information. Nothing more is known about either of these cases and since Fourier had not published anything more on this subject, and Ostrogradsky himself never published anything, they both possibly had second thoughts.

The connection of marine insurance with moral expectation provided an occasion for Laplace (1814/1995, p. 89) to express himself in favour of insurance of life and even to compare a nation with an association

*Whose members mutually protect their property by proportionally supporting the costs of this protection.*

10) In the **eleventh**, the last, **Chapter,** and, in part, in *Supplement 1* to the TAP, Laplace examined the probability of testimonies. Suppose that an urn contains 1,000 numbered tickets. One of them is extracted, and a witness states that that was ticket number $i$, $1 \leq i \leq 1,000$. He may tell the truth and be deceived or not; or lie, again being deceived or not. Laplace calculated the probability of the fact testified by the witness given the probabilities of all the four alternatives. In accordance with one of his corollaries, the witness's mistake or lie becomes ever more probable the less likely is the fact considered in itself (p. 460).

Laplace next introduced the prior probability of a studied event confirmed by $m$ witnesses and denied by $n$ others. If it is 1/2 and the probability of the truthfulness of each witness is $p$, then the probability of the event was

$$P = \frac{p^{m-n}}{p^{m-n} + (1-p)^{m-n}}.$$

Suppose now that the probabilities of truthfulness are $p_i > 1/2$ and the prior probability of the event is $1/n$. If the event is reported by a chain of $r$ witnesses, then (p. 466)

$$P = (1/n) + [(n-1)/n]\frac{(np_1 - 1)(np_2 - 1)...(np_r - 1)}{(n-1)^r}$$

so that for $n = 2$ and $n \to \infty$

$P = (1/2) + (1/2)(2p_1 - 1)(2p_2 - 1) \ldots (2p_r - 1)$ and $P = p_1 p_2 \ldots p_r$

respectively.

Laplace next examines verdicts brought in by $s$ independent judges (jurors) assuming that each of them decides justly with probability $p > 1/2$. The probability of a unanimous verdict is

$p^s + (1 - p)^s = i/n.$

Here, the right side is known ($n$ is the total number of verdicts of which $i$ were brought in unanimously). For $s = 3$ (p. 470)

$p = 1/2 + [(4i - n)/12n]^{1/2}.$

If $4i < n$, it should have been concluded that Laplace's (very restrictive) assumptions were wrong; he did not however make this remark.

If, in different notation, the probability of a just verdict reached by each judge (juror) was unknown, and $p$ judges condemned, and $q$ of them acquitted the defendant, the probability of a just final verdict was (p. 527)

$$\int_{1/2}^{1} u^p(1-u)^q du \div \int_{0}^{1} u^p(1-v)^q dv$$

(cf. formulas from Laplace's Chapter 6). Laplace stated that the verdicts were independent, but only in passing (1816, p. 523). **Poisson** (1837a, p. 4) indicated that Laplace had considered the defendant innocent unless and until pronounced guilty: his formulas had not included any prior probability of guilt. In Poisson's opinion, this should be assumed to exceed 1/2. I note that his remark had nothing to do with any individual case.

In § 8.9.1 I return to the application of probability in jurisprudence; here, I additionally refer to Zabell (1988a).

### 7.2. Theory of Errors

Laplace's work on the theory of errors can be easily separated into two stages. While treating it in the 18$^{th}$ century, he was applying the comparatively new tool, the density[3], and trying out several rules for the selection of estimators of the true values of the constants sought (cf. § 6.3.3). His equations proved too complicated and he had to restrict his attention to the case of three observations. Later Laplace proved (not rigorously) several versions of the CLT and was able to drop his restriction, but he had to adopt other conditions. Here is **Bienaymé's** precise conclusion (1853/1867, p. 161), also noticed by Idelson (1947, p. 11):

*For almost forty years Laplace had been presenting […] memoirs on probabilities, but […] had not wanted to combine them into a general theory.*

However, Bienaymé continued, the CLT [non-rigorously proved by Laplace] enabled him to compile his TAP.

1) **The Year 1774.** Without substantiation, Laplace assumed that, for any $x_1$ and $x_2$, the sought density $\varphi(x)$ of observational errors satisfied the equation

$$\varphi(x_2)/\varphi(x_1) = \varphi'(x_2)/\varphi'(x_1)$$

and obtained

$$\varphi(x) = (m/2)e^{-m|x|}. \qquad (10)$$

Later, while discussing suchlike decisions, Laplace (1798 – 1825/1878 – 1882, t. 3, p. xi) argued that the adopted hypotheses ought to be "incessantly rectified by new observations" until "veritable causes or at least the laws of the phenomena" be discovered. A similar passage occurred in his *Essai* (1814/1995, p. 116). Cf. Double et al (1835, pp. 176 – 177): the main means for revealing the "vérité" were induction, analogy and hypotheses founded on facts and "incessantly verified and rectified by new observations".

Suppose that the observations are *a*, *b*, and *c*, and $p = b - a$, $q = c - b$. Issuing from the [likelihood function]

$$f(x) = \varphi(x)\varphi(p - x)\varphi(p + q - x) \qquad (11)$$

rather than from the density, Laplace determined the parameter sought, *e* [the median], with respect to curve (11); alternatively, he applied the condition

$$\int |\varepsilon - \xi| f(x) dx = \min, |\xi| < +\infty$$

whence it followed that the integrals of $f(x)$ over $(-\infty; e]$ and $[e; +\infty)$ were equal to each other so that *e*, just the same, was the [median]. Neither function (10) nor (11) contained a location parameter. Eisenhart (1964) noted that Pitman (1939) had made similar use of a function of the type of (11).

For small values of *m* the magnitude $x = e - a \approx (2p + q)/3$ and therefore *e* coincided with the arithmetic mean and function (10) became

$$\varphi(x) = (m/2)(1 - m|x|) \approx m/2 = \text{Const.}$$

Laplace was not satisfied with these corollaries and had thus rejected the [median]. Note that for a [random variable] $\xi$ with density (10) var $\xi = 2/m^2$ so that a small *m* really invites trouble.

He then studied the case of an unknown parameter *m* by applying the principle of inverse probability, that is, by the so-called **Bayes** formula (5.1) with equal prior probabilities, but made a mistake in his calculations (Sheynin 1977a, p. 7). Stigler (1986, pp. 115 – 116) explained its cause but wrongly indicated that, since Laplace had not then read the Bayes memoir, he was unable to borrow the "Bayes formula". Yes, indeed unable, but simply because that formula was lacking in the work of his predecessor. Neither had he mentioned my remark.

2) **The Year 1781.** Laplace again issued from the [likelihood function] of the type of (11) and put forward four possible conditions for determining the real value of the constant sought: the integrals of *f(x)* or *xf(x)* over [– *N*; 0] and [*N*; 0], where *N* was the greatest possible error, should be equal to each other; or, the value of the second integral over [– *N*; *N*] should be minimal; and his final condition was the application of the [maximum likelihood principle]. Recall, however, that the curve (11) did not include a location parameter so that it should have been somehow inserted. Anyway, Laplace decided in favour of his third condition (which coincided with the first one).

So as to select a density, Laplace examined a special problem (§ 7.1-2), and, not really convincingly, obtained a "mean" law of error

$$y = (1/2a) \ln a/|\ |,\ |\ |\quad a. \tag{12}$$

He referred to the principle of insufficient reason and noted that function (12) was even and decreased with $|\ |$, – that is, conformed to the properties of "usual" errors; the restriction $\quad 0$ hardly bothered him.

Next Laplace studied what might be called the multidimensional **Bayes** method. Suppose that observational errors $_i$, $i = 1, 2, …, n$, having "facility" $x_i$ have occurred. Then the probability of the observed system of errors is proportional to

$$P = \frac{x_1 x_2 ... x_n}{\int ... \int x_1 x_2 ... x_n\ dx_1 dx_2 ... dx_n}$$

where the integrals are taken over all possible values of each variable; actually, Laplace considered a more general case in which each $_i$ occurred $k_i$ times. Multiplying the obtained expression by the product of all the differentials, Laplace arrived at the probability element (of the differential) of an *n*-dimensional random vector. It was now possible for him to determine the law of distribution of observational errors provided that the prior information stated above was available.

In connection with density (12) Laplace carried out a special study introducing the **Dirac** delta-function which had already appeared in **Euler's** works, see **Truesdell** (1984, p. 447, Note 4, without an exact reference). One of Laplace's conditions for determining an estimator $x_o$ of the real value of the constant sought given its observations $x_1$, $x_2$, …, $x_n$ was (§ 7.2-1) that the integrals

$$\int (x-x_1)(x-x_2)\ldots(x-x_n)$$

over $[-a; x_o]$ and $[x_o; a]$ should be equal to each other. Laplace indicated, without proof, that in case of an infinite $a$ the arithmetic mean can be obtained from the density law (12). He apparently thought that the function (12) then became constant, cf. his similar derivation in § 7.2-1.

Laplace then went over to a "much more general proposition" for density

$$y = \varphi(x) = \varphi(-x) = q, \text{ if } x = 0, \text{ and } y = 0 \text{ otherwise}; \quad 0.$$

In actual fact, he considered a sequence of functions $\varphi(x)$ such that

$$\varphi(x) = q(\alpha), \quad \alpha = \{\alpha_1; \alpha_2; \ldots \alpha_n; \ldots\} \quad 0.$$

If $\alpha x = t$, then

$$\varphi(t) = q \text{ if } t = 0, ||<+\infty; \quad \varphi(t) = 0 \text{ otherwise (when } t \neq 0, ||=+\infty\text{)},$$

and, obviously,

$$\int_{-\infty}^{\infty} \varphi(t)dt = C (= 1).$$

Laplace had not written these last equalities, but I think that he had actually introduced the Dirac delta-function

$$\varphi(t) = \lim(\alpha/\sqrt{\pi})\exp(-\alpha^2 t^2), \quad .$$

Laplace could have regarded the equalities above as representing a uniform distribution of observational errors having an arbitrarily wide, rather than assigned beforehand domain. His proposition consisted in that the unknown constant $x_o$ was equal to the appropriate arithmetic mean, but it can hardly be proved in the context of generalized functions: Laplace had to consider the integral of

$$[\varphi(x-x_1)][\varphi(x-x_2)]\ldots[\varphi(x-x_n)],$$

which does not exist in their language.

3) **The Years 1810 – 1811.** Laplace (1810a) considered $n$ [independent] discrete random errors (or magnitudes) uniformly distributed on interval $[-h; h]$. After applying a particular case of characteristic functions and the inversion formula, he proved, very carelessly and non-rigorously, that, in modern notation, as $n \to \infty$,

$$\lim P\left(\frac{|\sum \varsigma_i|}{n} \leq s\right) = \frac{\sqrt{3}}{\sigma\sqrt{2\pi}} \int_0^s \exp(-x^2/2\sigma^2)dx, \, i = 1, 2, \ldots, n, \quad (13)$$

where $\sigma^2 = h^2/3$ was the variance of each $\varepsilon_i$. He then generalized his derivation to identically but arbitrarily distributed variables possessing variance. When proving the [CLT] he (1810a/1898, p. 304) made use of an integral of a complex-valued function and remarked that he hoped to interest "géomètres" in that innovation and thus separated himself from (pure) mathematicians, see also similar reservations elsewhere (Laplace 1774/1891, p. 62; 1812/1886, p. 365).

In a supplement to this memoir Laplace (1810b), apparently following **Gauss**, turned his attention to the MLSq, and derived it without making any assumptions about the arithmetic mean (cf. § 9A.2-2), but he had to consider the case of a large number of observations and to suppose that the means of their separate groups were also normally distributed.

Soon enough Laplace (1811) returned to least squares. This time he multiplied the observational equations in one unknown

$$a_i x + s_i = \varepsilon_i, \quad i = 1, 2, \ldots, n,$$

where the right sides were errors rather than residuals, by indefinite multipliers $q_i$ and summed the obtained expressions:

$$[aq]x + [sq] = [\varepsilon q].$$

The estimator sought was

$$x_o = -[sq]/[aq] + [\varepsilon q]/[aq] \equiv -[sq]/[aq] + m.$$

Tacitly assuming that all the multipliers $q_i$ were of the same order, Laplace non-rigorously proved another version of the CLT for $n \to \infty$:

$$P(m = \mu) = \frac{1}{\sigma_m \sqrt{2\pi}} \exp(-\mu^2/2\sigma_m^2), \quad \sigma_m^2 = k \frac{[qq]}{[aq]^2}, \quad k = \int_{-\infty}^{\infty} x^2 \varphi(x)dx$$

where $\varphi(x)$ was an even density of observational errors possessing variance.

Then Laplace determined the multipliers by introducing the condition

$$\int_{-\infty}^{\infty} |z|P(z)dz = \min \qquad (14)$$

which led him to equalities $q = \mu a_i$, and then to the principle of least squares (in the case of one unknown)

$$x = [as] \div [aa].$$

Finally, Laplace generalized his account to the case of two unknowns. He multiplied the observational equations (in two unknowns) by two sets of indefinite multipliers $\{m_i\}$ and $\{n_i\}$[4] and obtained a bivariate normal distribution for independent components and, once more applying the condition of least absolute expectation, arrived at the principle of least squares. Todhunter (1865, pp. 578 – 588), who referred to Ellis (1849), studied the case of three or more unknowns.

And so, the derived principle essentially depended on the existence of the normal distribution. Indeed, the CLT (and therefore a large number of observations) was necessary and the use of conditions of the type of (14) would have otherwise been extremely difficult. No wonder that Laplace's theory had not been enjoying practical success. I adduce a wrong statement formulated on this point by Tsinger (1862, p. 1) who compared the importance of the **Gaussian** and the Laplacean approaches:

*Laplace provided a rigorous* [?] *and impartial investigation* […]; *it can be seen from his analysis that the results of the method of least squares receive a more or less significant probability only on the condition of a large number of observations;* […] *Gauss endeavoured, on the basis of extraneous considerations, to attach to this method an absolute significance* [wrong]. *With a restricted number of observations we have no possibility at all to expect a mutual cancellation of errors and* […] *any combination of observations can* […] *equally lead to an increase of errors as to their diminution.*

With regard to Gauss see § 9A. Here, I note that Tsinger lumped together both justifications of the MLSq due to Gauss and that practice demanded the treatment of a finite (and sometimes a small) number of observations rather than limit theorems. Tsinger's high-handed attitude towards Gauss (and his blind respect for Laplace) was not an isolated occurrence, see §§ 10.8.5 and 13.2-7. This was partly occasioned by Gauss' disrespect of Legendre (§ 9A.1.2) and, consequently, by the (antiscientific) ignorance of Gauss by French mathematicians including Poisson, and partly because they had been unquestionably following Laplace's approach (Tsinger was not the only one so affected).

Even not so long ago Eisenhart (1964, p. 24) noted that the existence of the second formulation of the MLSq

*Seems to be virtually unknown to almost all American users of least squares except students of advanced mathematical statistics.*

4) **Chapter 4 of the TAP.** Laplace non-rigorously proved the CLT for sums and sums of absolute values of independent, identically distributed errors restricted in value as well as for the sums of their squares and for their linear functions. All, or almost all of this material had already been contained in his previous memoirs although in 1811 he only proved the local theorem for linear functions of errors.

In § 23 Laplace formulated his aim: to study the mean result of "observations nombreuses et non faites encore …" This was apparently the first explicit statement concerning general populations; see § 14.2-1 for the appropriate opinion of **Chebyshev** and **Markov** and § 10.8.5 for similar statements in physics.

5) In *Supplement 1* to the TAP Laplace (1816) considered observational equations in (let us say) two unknowns

$$a_i x + b_i y + l_i = v_i, \quad i = 1, 2, \ldots, s.$$

Suppose that $\xi x$ and $\eta y$ are the errors of the least-squares estimators of the unknowns, denote the even density of the observational errors by $\varphi(u/n)$ with $|u| \leq n$, the moments by the letter $k$ with appropriate subscripts, $\alpha = \xi \sqrt{s}$, $\beta = \eta \sqrt{s}$,

$$\sigma^2 = \frac{k}{kk_4 - 2k_2^2}, \quad Q^2 = \sum_{i=1}^{s}(a_i \alpha + b_i \beta)^2 \text{ and } t = \frac{[vv]}{\sqrt{s}} - \frac{2k_2 n^2 \sqrt{s}}{k}.$$

Laplace calculated

$$P(\alpha; \beta) \sim \exp\{-Q^2/(2[vv] - 2t\sqrt{s})\},$$

$$P(t) \sim \exp\{-(\sigma^2/4n^4)[t + (Q^2/s\sqrt{s})]^2\}.$$

It thus occurred that $P(\alpha; \beta; t)$ which he also obtained showed that $t$ was independent of $\alpha$; $\beta$; or, that the sample variance was independent from the estimators of the unknowns, cf. §9A-5; to repeat, the observational errors were assumed to possess an even distribution, – and a normal distribution in the limit. For a proof of Laplace's result see Meadowcroft (1920).

Laplace also considered non-even distributions and recommended, in such a case, to demand that the sum of $v_i$ be zero. Since $[av] = 0$ is the first normal equation written down in another form, this demand is fulfilled for $a_i = $ const (or $b_i = $ const); otherwise, it is an additional normal equation corresponding to a fictitious unknown, the mean systematic error of observations.

Finally, Laplace derived a formula for estimating the precision. Without explanation (which appeared on p. 571 of his *Supplement* 2) he approximated the squared sum of the real errors by the same sum of the residuals and arrived at an estimator of the variance

$$m = \sqrt{\frac{[vv]}{s}}.$$

Without naming anyone **Gauss** (1823b, §§ 37 – 38) remarked that that formula was not good enough, see § 9A.4-6. Interestingly, Laplace (1814/1995, p. 45) stated that

*The weight of the mean result increases like the number of observations divided* [the French word was indeed *divisé*] *by the number of parameters.*

6) In *Supplement 2* to the TAP Laplace (1818) adopted the normal law as the distribution of observational errors themselves and not only as the law for their means. Indeed, the new "repeating" theodolites substantially reduced the error of reading and equated its order with that of the second main error of the measurement of angles in triangulation, the error of sighting, which was not yet sufficient. The error in the sum of the three angles of a triangle (the appropriate discrepancy, or its *closing*) could therefore be also regarded as normally distributed with density

$$\varphi(x) = \sqrt{h/3\pi} \exp(-hx^2/3),$$

where $h = 1/2\sigma^2$ was the measure of precision of an angle.

Tacitly assuming that $h$ was a [random variable], Laplace proved that

$$Eh = 3n/2\Sigma^2, \quad \psi(x) = h^{n/2}\exp[(-h/3)\Sigma^2]$$

were its expectation and density, $\Sigma^2$, the sum of $n$ squares of the triangular discrepancies. He computed the probability of the joint realization of errors obeying the normal law in a triangle and concluded that an equal distribution of the closing of the triangle among its angles was advantageous. The MLSq leads to the same conclusion, and for that matter, irrespective of the normal law. Then, when adjusting a chain of triangles rather than a separate triangle, two additional conditions (never mentioned by Laplace) have to be allowed for, – those corresponding to the existence of two baselines and, possibly, two astronomical azimuths, – and a preliminary distribution of the closings is of course possible but not necessary.

And so, let the observational errors have density

$$\varphi(x) = \sqrt{h/\pi} \exp(-hx^2).$$

Denote the closing of triangle $i$ by $T_i$ and suppose that the errors of the angles $\alpha_i$, $\beta_i$ and $\gamma_i$ already obey the condition

$$\alpha_i + \beta_i + \gamma_i = T_i.$$

Laplace derived the relations

$$P(\alpha_i; T_i) \sim \sqrt{h/3\pi} \exp[-(h/3)T_i^2],$$
$$P(T_1; T_2; \ldots; T_n) \sim (\sqrt{h/3\pi})^{n/2} \exp\{-(h/3)[TT]\},$$

$$P(h) \sim \frac{h^{n/2}e^{-(h/3)[TT]}}{\int_0^\infty h^{n/2}e^{-(h/3)[TT]}dh}, \quad Eh = \int_0^\infty hP(h)dh = \frac{3n+2}{2[TT]} \approx \frac{3n}{2[TT]}.$$

The error involved in the approximation just above can be easily estimated: in his *Supplement 3* (ca. 1819) Laplace took $n_1 = 26$ and $n_2 = 107$. Finally, supposing that $h = Eh$,

$$= 1/\sqrt{2h} = \sqrt{[TT]/3n},$$

not a bad result (improved by the approximation above!).

Laplace next investigates the adjustment of equations in one unknown by the MLSq for normally distributed errors. The interesting point here is that he had not indicated that the distribution of the residuals was also normal; in other words, that that distribution was [stable].

In the same Supplement Laplace discussed the **Boscovich** method of adjusting meridian arc measurements (§ 6.3.2-3). Write the initial equations in a different form,

$$p_i y - a_i + x_i = 0, i = 1, 2, \ldots, n, p_i > 0, a_1/p_1 > a_2/p_2 > \ldots > a_n/p_n.$$

The second unknown is presumed to be eliminated, and $x_i$ are the residual free terms. The Boscovich conditions, or, rather, his second condition, leads to

$$y = a_r/p_r$$

with error $-x_r/p_r$, i.e., to the calculation of this second unknown from one equation only. This latter is determined by inequalities

$$p_1 + p_2 + \ldots + p_{r-1} < p_r + p_{r+1} + \ldots + p_n,$$
$$p_1 + p_2 + \ldots + p_r > p_{r+1} + p_{r+2} + \ldots + p_n.$$

And so, these inequalities determine the sample median of the fractions $a_i/p_i$. Suppose now that the observational errors have an even density $\varphi(x)$ and

$$k = \int_0^\infty x^2 \varphi(x) dx.$$

Then, as Laplace showed, basing his derivation on variances[5] rather than on absolute expectations as before, the Boscovich method was preferable to the MLSq if, and only if,

$$4\varphi^2(0) > 1/(2k).$$

According to **Kolmogorov** (1931), the median is preferable to the arithmetic mean if

$$1/[2\varphi(m)] < 1, \sigma^2 = 2k,$$

and $m$ is the population median.

While translating Laplace's *Mécanique Céleste* into English, **Bowditch** (Laplace 1798 – 1825/1832, vol. 2, § 40, Note) stated:

*The method of least squares, when applied to a system of observations, in which one of the extreme errors is very great, does*

*not generally give so correct a result as the method proposed by Boscovich […]; the reason is, that in the former method, this extreme error* [like any other] *affects the result in proportion to the second power of the error; but in the other method, it is as the first power.*

In other words, the robustness of the Boscovich method is occasioned by its connection with the median.

7) In **Supplement 3 to the TAP** Laplace (ca. 1819) begins by evaluating a chain of 26 triangles (Perpignan – Formentera) which was a part of a much longer independent chain of 107 triangles. For the same normal distribution $\varphi(x)$ he has

$$\eta = \int_{-\infty}^{\infty} |x| \varphi(x) dx, \quad \theta = \int_{-\infty}^{\infty} x^2 \varphi(x) dx, \quad \theta = \pi \eta^2 / 2.$$

The empirical value of $\eta$ for the longer chain was

$(1/107) (|T_1| + |T_2| + \ldots + |T_{107}|) = 1.62$ so that $(1.62^2/2)\pi = 4.13$.

With subscripts 1 and 2 denoting the shorter and the longer chains respectively, Laplace has

$[TT]_1 = 4.13 \cdot 26 = 107.8$; empirical value, $(26/107)[TT]_2 = 108.8$.

This calculation shows that, first, Laplace preferred to evaluate $[TT]_1$ by $[TT]_2$ rather than use its actual value which was hardly correct since the pertinent conditions of observations could well have been different. Second, Laplace has thus qualitatively checked the realization of the normal law.

Next Laplace considers the adjustment of equations

$p_i x = a_i + m_i \varepsilon_i + n_i \varepsilon'_i, i = 1, 2, \ldots, n$

in one unknown, $x$, and independent errors $\varepsilon_i$ and $\varepsilon'_i$ both distributed normally with differing measures of precision; he only mentioned independence later (1827/1904, p. 349). Laplace explained his calculation by referring to his pp. 601 – 603 of the same Supplement which does not help but at least he concluded that the error of $x$ was distributed normally so that after all he knew that the normal law was [stable], cf. § 9A.2-6. However, the variance of the emerged law depended on the application of the MLSq which meant that the result just formulated was not sufficiently general.

Also in 1827 Laplace (1904, p. 343) stated that the MLSq was a particular case of the "most advantageous" method of adjustment (based on the minimal value of the expected absolute error and the presumed normal law, see end of § 7.2-3). Before 1823, he would have been partly in the right, but not afterwards, not since **Gauss'** second justification of the MLSq had appeared.

8) In Note 4 to this Chapter I indicated that Laplace had successfully treated the case of **dependence** between random

variables. Elsewhere, however, he (1827) somehow erred when investigating the atmospheric pressure. Its mean daily variation in Paris during 11 years was 0.763 *mm*, or 0.940 *mm,* if, during the same years, only three months, from February to April, were taken into consideration. When attempting to find out whether the difference between the two values was significant, Laplace had not indicated that they were not independent[6]. He made one more mistake: when solving his equations in two unknowns, the action of the Moon and the time of the maximal pressure, he had not stated that, again, the appropriate free terms were not independent. Without justifying his remark, **K. Pearson** (1914 – 1930, vol. 3A, p. 1) stated that

***Condorcet** often and Laplace occasionally failed because* [the] *idea of correlation was not in their mind*.

Elsewhere, he (1978, p. 658) left a similar remark, again without substantiation; there also, on p. 671, he added that Laplace was "rarely a good collector, or a safe guide in handling [the data]". Pearson exaggerated: on the then possible scientific level, and issuing from observations, Laplace proved that the Solar system will remain stable for a long time and completed the explanation of the motion of its bodies in accordance with the law of universal gravitation.

### 7.3. Philosophical Views

Laplace (1814/1995, p. 2) stated that, for a mind, able to "comprehend" all the natural forces, and to "submit these data to analysis", there would exist no randomness "and the future, like the past, would be open" to it. Nowadays, this opinion cannot be upheld because of the recently discovered phenomenon of chaos (§ 1.2.4); however, other remarks are also in order.

a) Such a mind does not exist and neither is there any comprehensive theory of insignificant phenomena, a fact which Laplace undoubtedly knew. He therefore actually recognized randomness (Dorfman 1974, p. 265).

b) In addition, there exist unstable movements, sensitive to small changes of initial conditions, cf. § 11.2-9.

c) Already previous scholars, for example, **Maupertuis** (1756a, p. 300) and **Boscovich** (1758, §385), kept to the "Laplacean determinism". Both mentioned calculations of past and future ("to infinity on either side", as Boscovich maintained) but, owing to obvious obstacles, see above Item *a*, both disclaimed any such possibility.

In his *Essai* Laplace (1814/1995, p. 37) additionally provided examples of "statistical determinism", – of the stability of the number of dead letters and of the profits made by those who ran lotteries. He explained all this by the action of the LLN (more precisely, by its, then barely known, **Poisson** form, see § 7.1-8). Participation in lotteries only depends on free will, cf. **Quetelet's** similar statement in § 10.5 and **Petty's** opinion (§ 2.1.1)[7].

In his early memoirs, Laplace (e.g., 1776/1891, pp. 144 – 145), like **Newton** (§ 2.2.3), had not recognized randomness and explained it by ignorance of the appropriate causes, or by the complexity of the studied phenomenon. He even declared that the theory of probability, that estimated the degrees of likelihood of phenomena, was indebted for its origin to the weakness of the mind and a similar statement occurred in his *Essai* (1814/1995, p. 3). Thus, probability became for him an applied mathematical discipline servicing natural sciences[8] (cf. § 0.1) and, even for this reason alone, he had not separated mathematical statistics from it, although he (1774/1891, p. 56) noted the appearance of "un nouveau genre de problème sur les hasards", and even (1781/1893, p. 383) of "une nouvelle branche de la théorie des probabilités"[9].

A curious statement deserves to be included (Laplace 1796/1884, p. 504):

*Had the Solar system been formed perfectly orderly, the orbits of the bodies composing it would have been circles whose planes coincided with the plane of the Solar equator. We can perceive however that the countless variations, that should have existed in the temperatures and densities of the diverse parts of these grand masses, gave rise to the eccentricities of their orbits and the deviations of their movement from the plane of that equator.*

The causes mentioned by Laplace could have hardly be called external, and the main relevant explanation of randomness, deviation from the laws of nature, persisted. Leaving aside the planes of the planetary orbits, I deny his (and **Kant's**, and **Kepler's**, see § 1.2.4) opinion concerning eccentricities.

**Newton** theoretically proved that the Keplerian laws of planetary motion resulted from his law of universal gravitation. In my context, it is necessary to stress: he also proved that the eccentricity of the orbit of a given planet is determined by the planet's initial velocity. For some greater values of that velocity the orbit will become parabolic (with its eccentricity equal to unity, rather than less than unity as in the case of ellipses), for other still greater values, hyperbolic (with eccentricities greater than unity). And for a certain value of that velocity an elliptic orbit will become circular. But is it really necessary to refer to a rigorous proof? Indeed, it is difficult to imagine that such changes do not occur gradually, that, consequently, the eccentricity does not vary continuously with the velocity.

All these findings, as Newton proved, persisted for planets (not material points) having uniform density. So what should we think about Laplace's explanation? I believe that the variations of densities (but hardly temperatures) peculiar to a given planet could have somewhat corrupted the eccentricity caused by its initial velocity. I am unable to say whether they could have also caused a corruption of some other kind, but in any case I need not discuss this problem.

So it really seems that Laplace (and Kant) were mistaken (Kepler was obviously ignorant of the law of universal gravitation). I am not sure that Kant had studied Newton attentively enough, but Laplace certainly did, although a bit later, in t. 1 of his *Traité de Méc. Cél.* (1798/1878, Livre 2, chapters 3 and 4).

Witness finally **Fourier's** comment (1829, p. 379) on Laplace's *Exposition*: it "is an ingenious epitome of the principal discoveries". And on the same page, discussing Laplace's "historical works" (to whose province the *Exposition* belonged):

*If he writes the history of great astronomical discoveries, he becomes a model of elegance and precision. No leading fact ever escapes him.* […] *Whatever he omits does not deserve to be cited.*

### 7.4. Conclusions

Laplace collected his earlier memoirs in one contribution which cannot, however, be regarded as a single whole. He never thought about solving similar problems in a similar way (and his *Essai* was not a masterpiece of scientific-popular literature, see Note 1 to this Chapter). Then, many authors complained that Laplace had described his reasoning too concisely. Here, for example, is what **Bowditch** (**Todhunter** 1865, p. 478), the translator of Laplace's *Traité de mécanique céleste* into English, sorrowfully remarked:

*Whenever I meet in La Place with the words 'Thus it plainly appears' I am sure that hours, and perhaps days of hard study will alone enable me to discover how it plainly appears.*

This can also be said about the TAP.

The Laplacean definition of probability (to repeat: first introduced by **De Moivre,** see § 4.3) was of course unsatisfactory, but nothing better had appeared until the advent of the axiomatic theory (or, the **Mises** debatable formula). Here is the testimony of Kamke (1933, p. 14): In 1910, it was said at Göttingen University that probability was a number situated between 0 and 1 about which nothing more was known. Similar statements were due to Mises in 1919, to **Keynes** in 1921, and to **P. Lévy** (who was born in 1886) in his earlier life (**Cramér** 1976, § 2.1) as well as to **Poincaré** (§ 11.2-1) and **Markov** (§ 14.1-5).

But the opinion of **Doob** (1989) was even more interesting. In 1946

*To most mathematicians mathematical probability was to mathematics as black marketing to marketing;* […] *the confusion between probability and the phenomena to which it is applied* […] *still plagues the subject;* [the significance of the **Kolmogorov** monograph] *was not appreciated for years, and some mathematicians sneered that* […] *perhaps probability needed rigor, but surely not <u>rigor mortis</u>;* […] *the role of measure theory in probability* […] *still embarrasses some who like to think that mathematical probability is not a part of analysis.*

All this means that Laplace is here exonerated. However, he had not even heuristically introduced the notion of *random variable* and was therefore unable to study densities or characteristic functions as mathematical objects. His theory of probability quite properly remained an applied mathematical discipline since he made outstanding discoveries in mathematics, astronomy and physics. However, it did not yield to development which necessitated its construction anew. It is opportune to note that **Maxwell** referred to Laplace only twice (§ 10.9.5), and **Boltzmann** did not mention him at all.

Laplace had not regarded himself as a pure mathematician (§ 7.2-3) and I quote the opinion of Fourier (1829, pp. 375 – 376):

*We cannot affirm that it was his destiny to create a science entirely new, like **Galileo** and **Archimedes**; to give to mathematical doctrines principles original and of immense extent, like **Descartes, Newton** and **Leibniz**; or, like **Newton**, to be the first to transport himself into the heavens, and to extend to all the universe the terrestrial dynamics of Galileo: but Laplace was born to perfect everything, to exhaust everything, and to drive back every limit, in order to solve what might have appeared incapable of solution. He would have completed the science of the heavens, if that science could have been completed.*

Laplace introduced partial differential equations and, effectively, stochastic processes into probability, and non-rigorously proved several versions of the CLT by applying characteristic functions and the inversion formula.

On that basis, he constructed his version of the theory of errors, which essentially depended on the existence of a large number of normally distributed observational errors and was therefore unsuccessful. In the not yet existing mathematical statistics Laplace investigated the statistical significance of the results of observation, introduced the method of statistical simulation, studied his version of sampling and extended the applicability of the **Bayesian** approach to statistical problems.

He knew the **Dirichlet** formula (even in a generalized version), introduced the **Dirac** delta-function and integrals of complex-valued functions. He had also indicated (long before the strong law of large numbers became known) that in probability theory *limit* was understood in a special way. Molina (1930, p. 386) quoted his memoir (1786/1894, p. 308) where Laplace had contrasted (although not clearly enough) the "approximations" admitted in the theory of probability with certainty provided in analysis.

It is my belief that a balanced opinion of a classic is needed, and I am now discussing Laplace's mistakes and shortcomings. That his contributions are extremely difficult to understand is generally known, see for example Bowditch's opinion above. Then, Laplace was "extremely careless in his reasoning and in carrying out formal transformations" (**Gnedenko** & Sheynin 1978/2001, p. 224 with examples attached). And here are some of my present comments.

1) See my comments just above the quotation from Doob.
2) Laplace made a mistake when studying the **Buffon** problem (§ 7.1-4).
3) He applied an unsuitable model when calculating the population of France (Ibidem). And he presented his final estimate (1812, pp. 399

and 401) in a hardly understandable way; in any case, **Poisson** (1812) wrongly reported his result.

4) There also, Laplace bravely, and without any reservation, calculated the probability of a certain demographic relation persisting for a hundred years to come.

5) He made a mistake while discussing the treatment of observations (§ 7.2-1).

6) Here is Laplace's opinion (1814/1995, p. 81) about mortality tables: "There is a very simple way of constructing" [them] from the registers of births and deaths". But the main point is to study the plausibility of these registers, to single out possible corruptions and exceptional circumstances etc. Then, the boundaries of the constructed mortality table have to be determined both in time and territory.

7) **Poisson** (§ 8.9.1) noted that in applying probability theory to jurisprudence, Laplace had not allowed for a prior probability of defendants' guilt. Such a probability has nothing in common with the presumption of innocence in each individual case, cf. **Quetelet's** inclinations to crime (§ 10.5).

8) The method of least squares. The proper attitude for Laplace would have been to acknowledge the **Gaussian** demand (§ 9A.2) for studying the treatment of a small number of observations and to restrict therefore the importance of his own results. Instead, he insisted on his own approach and virtually neglected Gauss. Later French scientists including Poisson followed suit and even the most eminent mathematicians (or at least those of them who had not studied attentively the treatment of observations) became confused. When proceeding to prove the CLT, Chebyshev (§ 13.1-4) remarked that it leads to the MLSq!

9) Laplace did not explain, for example in his *Essai* (1814), the dialectic of randomness and necessity.

10) See my comments on the eccentricities of the planetary orbits in § 7.3. The last edition of his book of 1796 during his lifetime appeared in 1813 without any corrections of that subject.

### Notes

**1.** The *Essai* ran through a number of editions and was translated into many languages. It attracted the public to probability, but the complete lack of formulas there hindered its understanding. The appearance of **Quetelet's** superficial contributions written in good style (§ 10.5) had a negative effect on the fate of the *Essai*.

**2.** For a simpler derivation of its equation see **Todhunter** (1865, pp. 545 – 546).

**3.** Laplace applied several pertinent terms. In his TAP, he finally chose *loi de probabilité* or *loi des erreurs.*

**4.** The quantities [ *m*] and [ *n*] which appeared here were not independent. Without indicating this, Laplace correctly solved his problem.

**5.** In his *Supplement 3* Laplace once more applied the variance as the main measure of precision of observations.

**6.** Cf. Note 4. Retaining excessive decimals (see above) was of course traditional, see Note 11 to Chapter 1.

**7. Kant** (1763/1912, p. 111) indicated that the relative number of marriages (which obviously depended on free will) was constant.

**8.** The subjects discussed by Laplace in his *Exposition* (1796/1884) had not demanded stochastic reasoning (see however end of this subsection), but he undoubtedly applied them, for example, in the *Traité* (1798 – 1825), to say nothing about the treatment of observations, and his determinism had not hindered him at all.

Thus, Laplace (1812/1886, p. 361) stated that a certain magnitude, although having been indicated by [numerous] observation[s], was neglected by most astronomers, but that he had proved its high probability and then ascertained its reality (although did not provide his calculations). In general, unavoidable ignorance concerning a single random event becomes a cognizable regularity.

**9. Lagrange**, in a letter to Laplace of 13.1.1775, see t. 14 of his *Oeuvres*, 1892, p. 58, used this latter expression. Inductive stochastic conclusions occurred in the Talmud (§ 1.1.2) and **Arbuthnot's** memoir (§ 2.2.4) and the work of many other authors, especially **Bayes**, which had appeared before Laplace, might be today attributed, at least in part, to mathematical statistics.

## 8. Poisson

Poisson's publications concerning the theory of probability began to appear in 1811 when he published abstracts of two of **Laplace's** early memoirs, and in 1812 he continued this work by an abstract of Laplace's *Théorie analytique*. These abstracts, published in *Nouv. Bull. des Sciences Soc. Philomatique de Paris*, were not really interesting, but who could have done better? Only **Fourier.** See Bru (1981; 2013) for a description of the French mathematical community during Poisson's lifetime and his general and for many years dominant role there. Bru also took care to explain much of Poisson's mathematics.

Like Laplace, Poisson had published a number of memoirs on the theory of probability, then combined them in his monograph (1837a) whose juridical title did not reflect its contents; only its subtitle promised to discuss, as a preliminary, the general principles of the calculus of probability. I describe both this contribution (referring only to its page numbers) and his other works. First, however, I quote Poisson's statement (p. 1) about the place of probability in mathematics and then describe the scope of the *Elements of the Calculus of Probability and Social Arithmetic* as formulated by him (Poisson 1837c, p. 26).

And so, probability became

*Une des principales branches des mathématiques, soit par le nombre et l'utilité de ses applications, soit par le genre d'analyse auquel il a donne naissance* [to which it gave birth].

The *Elements* (1837c) listed: **1**) Topics of probability itself (general principles, the **Bernoulli** theorem, probabilities of future events derived from the probabilities of similar previous events). **2**) Tables of mortality, mean duration of life, smallpox, inoculation and vaccination. Here also, expectation, cf. § 13.3. **3**) Institutions depending on probabilities of events (annuities, insurance, loans). **4**) Mean values of a large number of observations.

The soon forgotten term *social arithmetic* (appearing also below in § 8.9) thus designated population and medical statistics. Now, we would rather say *social statistics*.

### 8.1. Subjective Probability

The aim of the calculus of probability, as Poisson (pp. 35 – 36) maintained, was the determination, in any doubtful "questions", of the ratio of the cases favourable for the occurrence of an event to all possible cases, and its principles should be regarded as "un supplément nécessaire de la logique". He (pp. 30 and 31) remarked that probability changed with experience, and was subjective, but that the chance of an event remained constant. Already **Leibniz** (§ 2.1.1) and then **De Morgan** (1847)[1] and **Boole** (1952) attempted to justify probability by elements of mathematical logic, see also Halperin (1988).

The stressed difference between chance and probability (also recognized by **Cournot**, see § 10.3) is now forgotten, although Poisson attempted to adhere to it. Thus (p. 47), he showed that the subjective probability of extracting a white ball from an urn containing white and black balls in an unknown proportion was equal to 1/2 "as it should have been". This conforms to the principles of the theory of information and he himself was satisfied with the result obtained since it corresponded to "la perfaite perplexite de notre esprit".

Davidov (1854, p. 66) who was well acquainted with foreign literature, noted:

*Vague ideas on probability and an inexact distinction between subjective and objective probabilities are among the main obstacles against a speedy development of practical medicine.*

### 8.2. Two New Notions

Poisson (1829, § 1) defined the distribution function of a discrete random variable as

$$F(x) = P(\xi < x)$$

and (Ibidem) introduced the density as the derivative of $F(x)$. Later he (1837b, pp. 63 and 80) similarly treated the continuous case. **Davidov** (1885) and **Liapunov** (1900) had noted his innovation, but distribution functions only became generally used in the 20$^{th}$ century.

Poisson (pp. 140 – 141) was also the first to introduce the notion of a discrete random variable although he named it by an obviously provisional term, *chose A*[2]. He then (p. 254) considered a random variable with values being multiples of some $\alpha$, assumed that $\alpha \to 0$ and thus went over, in accordance with the tradition of his day, to a continuous variable[3]. As compared with **Simpson**, who studied random observational errors (§ 6.3.1), Poisson's innovation here was a formal heuristic definition of a random variable and its more general (not necessarily connected with the theory of errors) understanding.

### 8.3. The De Moivre – Laplace Limit Theorem

Poisson (1837a, p. 189) provided his own derivation of that theorem by issuing from the probability[4] of the occurrence of contrary events $A$ and $B$ not less than $m$ times (not more than $n$ times) in $\mu = m + n$ **Bernoulli** trials

$$P = p^m \{1 + mq + \frac{m(m+1)}{2!} q^2 + \ldots + \frac{m(m+1)\ldots(m+n-1)}{n!} q^n\} = \quad (1)$$

$$= \int_a^\infty X dx \div \int_0^\infty X dx, \quad X = \frac{x^n}{(1+x)^{\mu+1}}, \quad (2)$$

where $p$ and $q$ were the probabilities of the occurrence of these events in a single trial, $p + q = 1$.

His results were, however, tantamount to formula (7.3), see Sheynin (1978b, pp. 253 – 254). **Montmort** (1708/1713, p. 244), also see

**Todhunter** (1865, p. 9), knew formula (1) and formula (2) occurred in **Laplace's** TAP, Chapter 6.

For small values of $q$ Poisson (p. 205) derived the approximation

$$P \approx e^{-\mu}(1 + \mu + \mu^2/2! + \ldots + \mu^n/n!), \qquad (3)$$

where $mq \approx \mu q = \mu$. He had not provided the expression

$$P(\xi = m) = e^{-\mu}\mu^m/m!.$$

### 8.4. Sampling Without Replacement

Poisson (pp. 231 – 234) examined sampling without replacement from an urn containing $a$ white balls and $b$ black ones ($a + b = c$) and applied the result obtained for appraising a model of France's electoral system. Suppose that the sample contained $m$ white, and $n$ black balls ($m + n = s$). Its probability, as Poisson indicated, was represented by the [hypergeometric] distribution. For large $a$ and $b$ as compared with the sample, Poisson determined an approximate expression for that probability under an additional condition

$$n > m. \qquad (4)$$

If a series of $k$ such samples are made, then

$$s_1 + s_2 + \ldots + s_k = c.$$

After calculating the probability of the condition (4) being fulfilled $j$ times out of $k$, Poisson concluded that, even if $b$ only somewhat exceeded $a$, $j$ will apparently be too large. For $k = 459$, which was the number of electoral districts in France, $c = 200{,}000$, equal to the number of the voters (less than 1% of the population!). Suppose also that each voter is a member of one of the two existing parties; that the voters are randomly distributed over the districts; and that the proportion of party memberships is 90.5:100. Then, as Poisson concluded, remarking, however, that his model was too simplified, the probability of electing a deputy belonging to the less numerous party was very low.

Poisson (1825 – 1826) studied sampling without replacement also in connection with a generally known game. Cards are extracted one by one from six decks shuffled together as a single whole until the sum of the points in the sample obtained was in the interval [31; 40]. The sample is not returned and a second sample of the same kind is made. It is required to determine the probability that the sums of the points are equal. Here is his solution, see Sheynin (1978b, pp. 290 – 292), Gnedenko & Sheynin (1978/2001, pp. 236 – 237) and § 8.5 below.

An urn contains $x_1$ balls numbered one, $x_2$ balls numbered two, …, $x_i$ balls numbered $i$ ($x_1 + x_2 + \ldots + x_i = s$). What is the probability of extracting without replacement $a_1$ balls numbered one, $a_2$ balls numbered two, …, $a_i$ balls numbered $i$ ($a_1 + a_2 + \ldots + a_i = n$) if

$$a_1 + 2a_2 + \ldots + ia_i = x? \qquad (5)$$

Without taking into account this condition Poisson got

$$P = \frac{n!(s-n)!}{s!} \frac{x_1!}{a_1!(x_1-a_1)!} \frac{x_2!}{a_2!(x_2-a_2)!} \cdots \frac{x_i!}{a_i!(x_i-a_i)!} =$$

$$(s+1)\int_0^1 (1-y)^s Y dy,$$

$$Y =$$

$$\frac{x_1!}{a_1!(x_1-a_1)!}\left[\frac{y}{1-y}\right]^{a_1} \frac{x_2!}{a_2!(x_2-a_2)!}\left[\frac{y}{1-y}\right]^{a_2} \cdots \frac{x_i!}{a_i!(x_i-a_i)!}\left[\frac{y}{1-y}\right]^{a_i}.$$

For $i = 2$ the system of probabilities thus calculated defines a hypergeometric distribution, as it is now called. So as to allow for condition (5), Poisson replaced $Y$ by the sum of its values corresponding to such sets of $\{a_1; a_2; \ldots; a_i\}$ for which that restriction is fulfilled.

Considering a generating function

$$(1+\frac{yt}{1-y})^{x_1}(1+\frac{yt^2}{1-y})^{x_2}\ldots(1+\frac{yt^i}{1-y})^{x_i},$$

Poisson noted that the probability sought was equal to the coefficient of $t$ in

$$(s+1)\int_0^1 (1-y+yt)^{x_1}(1-y+yt^2)^{x_2}\ldots(1-y+yt^i)^{x_i} dy.$$

A second sample of extracted balls $\{b_1; b_2; \ldots; b_i\}$ had to be introduced and the probability of the joint occurrence of condition (5) and of restriction

$$b_1 + 2b_2 + \ldots + ib_i = z$$

had to be determined by means of a bivariate generating function, also see § 8.5 below. Only later, when solving the electoral problem (above), Poisson remarked that the result of the second sampling might be (subjectively!) considered independent from the first one if only it, the first one, remained unknown. Poisson's solution was indeed ingenious.

Suppose (Poisson 1837a, pp. 231 – 234) that an urn contains $a$ white balls and $b$ black ones. Two samples without replacement are made, one after the other, and $g$ and $m$ white balls and $h$ and $n$ black ones are extracted respectively, $g + h = r$. The probability of the second sample is

$$P(a; b; m; n) = \sum [P(a - g; b - h; m; n)\, P(a; b; g; h)]$$

where the sum is extended over $g, h = 0, 1, 2, \ldots; g + h = r$ and the letters in parentheses are the appropriate arguments. The right side of the formula does not depend on $r$ which might be therefore assumed to be zero. This remark indeed proves Poisson's statement as well as the finding of another author, Mondesir (1837). This episode and its further history is described in Sheynin (2002b). Here, I only mention **Chuprov**. In a letter of 1921 he (Sheynin 1990/2011, p. 145) stated:

*Not knowing the prior data it is impossible to distinguish a series of numbers obtained when extracting the tickets without replacement from a series obtained according to the usual way of replacing […] the ticket.*

Also see Chuprov (1923, pp. 666 – 667; 1924, p. 490). The first to consider sampling without replacement was **Huygens** (§ 2.2.2).

### 8.5. Limit Theorems for the Poisson Trials

Suppose that contrary events $A$ and $B$ occur in trial $j$ with probabilities $p_j$ and $q_j$ ($p_j + q_j = 1$). Poisson (p. 248) determined the probability that in $s$ trials event $A$ occurred $m$ times, and event $B$, $n$ times ($m + n = s$). He wrote out the generating function of the random variable $m$ (or, the bivariate generating function of $m$ and $n$) as

$$X = (up_1 + vq_1)(up_2 + vq_2) \ldots (up_s + vq_s)$$

so that the probability sought was the coefficient of $u_m v_n$ in the development of $X$. His further calculations (lacking in Chapter 9 of Laplace's TAP) included transformations

$$u = e^{ix},\ v = e^{-ix},\ up_j + vq_j = \cos x + i(p_j - q_j)\sin x = \,_j \exp(ir_j),$$
$$_j = \{\cos^2 x + [(p_j - q_j)\sin x]^2\}^{1/2},\ r_j = \arctan[(p_j - q_j)\tan x].$$

Excluding the case of $p_j$ or $q_j$ decreasing with an increasing $s$, and without estimating the effect of simplifications made, Poisson (pp. 252 – 253) derived the appropriate local and integral limit theorems for large values of $s$. They were, however, complicated and their importance apparently consisted in extending the class of studied random variables.

### 8.6. The Central Limit Theorem

Poisson (p. 254) introduced a [lattice] random variable whose values were multiples of some   on a finite interval and depended on the number of the trial. Applying the appropriate characteristic function and the inversion formula, he determined the probability that the sum of these values $s$ was obeying certain inequalities $a$   $s$   $b$  . He then went over to the sum of continuous variables by assuming that    $0, a, b$   with finite $a$   and $b$   and (p. 268) derived the [CLT] for $s$ under a single (not adequately explained) condition, again without estimating the effect of simplifications made[5]. In accordance with the context, it seems, however, that he supposed that the variances of the terms of $s$ were finite and bounded away from zero.

He (p. 258) also made use of the **Dirichlet** discontinuity factor which he considered known. Dirichlet introduced it in two papers, both published in 1839, see his *Werke*, Bd. 1, 1899, pp. 377 – 410.

Poisson (1824; 1829) earlier proved several versions of the CLT in the same way. He (1824, §§ 4 and 6) introduced then the so-called **Cauchy** distribution and found out that it was [stable]. For a modern exposition see **Hald** (1998, pp. 317 – 327) who took into consideration all the appropriate deliberations of Poisson. Those proofs were methodically bad since the conditions of the theorems were not provided whereas Hald remarked that that defect had then by usual. Laplace (and Poisson) treated probability as a branch of applied mathematics (§ 0.2 and elsewhere), but such omissions are still unforgivable[6]. The proofs themselves were not really studied: Hald somehow thought that they were rigorous.

Poisson (1824, §§ 8 – 10) also considered a linear function

$$E = a_1\xi_1 + a_2\xi_2 + \ldots + a_n\xi_n$$

of discrete and continuous independent random variables $\xi_i$. In the second instance he (1824, p. 288) obtained the appropriate CLT and noted that that theorem did not hold for variables with density

$$\varphi(x) = e^{-2|x|}, \ |x| < +\infty$$

and either $a_i = 1/(i + 1)$ or $1/(2i – 1)$. **Markov** (1899c, p. 42) mentioned these exceptional cases in his debates with **Nekrasov** about the CLT; in the translation of his note, I have inadvertently omitted his exact reference to Poisson.

Poisson also applied the CLT for estimating the significance of discrepancies between empirical indicators obtained from different series of observations. For the **Bernoulli** trials he studied the discrepancies between probabilities of events (p. 224) and between the appropriate frequencies (p. 294), and, for his own pattern (§ 8.7), between the mean values of a random variable (p. 288). **Cournot** (1843, Chapters 7 and 8) borrowed his findings without mentioning him.

### 8.7. The Law of Large Numbers

Here is how he defined this law in his Préambule (p. 7):

*Things of every kind obey a universal law that we may call the law of large numbers. Its essence is that if we observe a very large number of events of the same nature, which depend on constant causes and on causes that vary irregularly, sometimes in one manner sometimes in another, i.e., not progressively in any determined sense, then almost constant proportions will be found among these numbers.*

He went on to state qualitatively that the deviations from his law became ever smaller as the number of observations increased. **Bortkiewicz** (1904, p. 826, Note 13) remarked that the Préambule was largely contained in Poisson's previous work (1835). Poisson (1837a, pp. 8 – 11) illustrated his vague definition by various examples, which, however, did not adequately explain the essence of the law but were interesting indeed. Thus (pp. 9 and 10), the LLN explains the stability of the mean sea level and the existence of a mean interval between molecules. Beginning with 1829, Poisson's contributions had been containing many direct or indirect pronouncements on molecular conditions of substance, local parameters of molecular interactions, etc. sometimes connected with the LLN (Sheynin 1978b, p. 271, note 25) but they remained unnoticed.

Poisson then (pp. 138 – 142) formulated but did not prove three propositions characterizing the LLN. These were based on the standard formula (which Poisson had not written out)

$P(B) = P(A_i) P(B/A_i)$.

In actual fact, he studied the stability of statistical indicators by means of the CLT, see **Hald** (1998, pp. 576 – 582).

Poisson described his law in a very complicated way and its sufficiently detailed exposition is still lacking. No wonder that **Bortkiewicz** (1894 – 1896, Bd. 8, p. 654) declared that

*There hardly exists such a theorem that had met with so many objections as the law of large numbers.*

Here, in addition, is a passage from Bortkiewicz' letter to **Chuprov** of 1897 (Sheynin 1990c/2011, p. 60):

*Or take […] my last three-hour talk with Markov about the law of sm. [small] numbers [§ 15.1.2]. It caused me nothing but irritation. He again demanded that I change the title. With respect to this topic we got into conversation about the law of l. nn. It happens that Markov (like Chebyshev) attributes this term to the case when all the probabilities following one another in n trials are known beforehand. […] In concluding, Markov admitted that perhaps there did exist 'some kind of ambiguity' in Poisson's reasoning, but he believed that it was necessary to take into account the later authors' understanding of the term 'law of l. nn.' …*

It is indeed difficult to examine Poisson's considerations on that point, but at least one of his examples (p. 148ff) is clear. It deals with a throw of many coins of the same denomination and "mode de fabrication". And, although Poisson (p. 147) argued that the probability of (say) heads could be established statistically, it seems that his example had to do with unknown probabilities. Other examples mentioned above (sea level and interval between molecules) certainly included unknown probabilities.

The LLN was not recognized for a long time. In 1855 **Bienaymé** declared that it contained nothing new (§ 10.2) which apparently compelled **Cournot** (1843) to pass it over in silence; his views were certainly known even in 1842 (Heyde & Seneta 1977, pp. 46 – 47). Even much later **Bertrand** (1888a, pp. XXXII and 94) considered it unimportant and lacking in rigour and precision. However, already

**Bessel** (1838a, especially § 9) guardedly called the Poisson law a "principle" of large numbers, **Buniakovsky** (1846, p. 35) mentioned it and **Davidov** (1854; 1857, p. 11) thought it important. It is nevertheless possible (§ 3.2.3) that statisticians had recognized the **Bernoulli**, and the Poisson (and the **Chebyshev**) laws of large numbers only in the qualitative sense.

### 8.8. The Theory of Errors and Artillery

In the theory of errors Poisson offered his proof of the CLT (§ 8.6) and a distribution-free test for the evenness of the density of observational errors (1829, § 10). He (1837b) also applied the theory of probability and the error theory to artillery firing, although mostly in a methodical sense[7]. He recommended the variance as the main estimator of scattering which conformed to **Laplace's** later idea, see § 7.2-6. In one of his problems Poisson (1837b, § 7) determined the distribution of the square of the distance of some point from the origin given the normal distributions of the point's distances from the two coordinate axes. He thus was perhaps the first to treat clearly the densities as purely mathematical objects.

Poisson had not referred to the Gaussian theory of errors, had therefore essentially lessened the importance of his contributions to that discipline. Here, for example, is his opinion (1833, p. 36) about the merits of Legendre:

*It is to him that the sciences of observation owe the rule of calculation which [?] is named the MLSq. Laplace showed its entire probable benefit with respect to the precision of the results*.

Gauss does not exist! And Poisson's (1837a) layman's, as I insist, considerations about measurements and observations are almost useless. Unlike Poisson and other French mathematicians, see end of § 7.2-3, Laplace (1812/1886, p. 353) objectively described the discovery of the MLSq:

*Legendre conceived a simple idea to consider the sum of the squares of observational errors and to render it minimal which immediately leads to as many final equations as there are elements to be corrected. That learned geometer was the first to publish the indicated method, but we ought to acknowledge that Gauss had been invariably applying the same idea for many years before Legendre's publication and communicated it to many astronomers*.

### 8.9. Statistics

In § 6.2 I described the development of statistics in the 18th century and I return to this subject in Chapter 10. Here, I discuss the appropriate pronouncements of Poisson and some other scholars. Recall first of all (§ 8.6) that Poisson investigated the significance of empirical discrepancies. **Quetelet** (1869, t. 1, p. 103), who corresponded with Poisson, testified that the latter had mentioned statisticians, who

*Pretended to substitute their fantasies for the veritable principles of* [their] *science, with derisive severity*.

In a few other cases (and twice in joint papers) Poisson expressed himself more definitely. Thus (Libri-Carrucci et al 1834, p. 535):

*The most sublime problems of the arithmétique sociale can only be resolved with the help of the theory of probability.*

A year later Double et al (1835, p. 174) stated that

*Statistics carried into effect always is, after all, the functioning mechanism of the calculus of probability, necessarily concerning infinite [?] masses, an unrestricted number of facts;* and (p. 176) [with respect to the applicability of mathematics] *the state of the medical sciences is not worse than, not different from the situation with all the physical and natural sciences, jurisprudence, moral and political sciences etc.*

This opinion was, however, questioned. **Poinsot** (Poisson 1836, p. 380) declared that the application of the calculus of probability to "moral things", such as the verdicts of law courts and elections, was a "dangerous illusion", also see § 8.9.1[8]. Double (1837, pp. 362 – 363) sharply objected to the application of statistics in medicine and stated that "each case appear[ed] to me [to him] a new and a separate problem". However, he mistakenly identified statistics with the numerical method (see § 10.8). **Cauchy** (1821/1897, p. V) cautiously pronounced a similar opinion: The only method of natural sciences consisted in subjecting observations to calculus, but the mathematical sciences should not "exceed their bounds". Later he (1845/1896, p. 242), however, expressed himself quite differently: statistics, as he maintained, provided the means, infallible in a sense, for judging doctrines and institutions, and should be applied "with full rigour".

**8.9.1. Jurisprudence.** Poisson (1837a, pp. 1 – 2) thought that the study of the probabilities of verdicts and, in general, of majority decisions, was a most important application of the calculus of probability. He (p. 17) perceived his main goals in that field as an examination of the stability of the rate of conviction, a derivation of the number of jurors and of the majority necessary for conviction so as to lower the probability of miscarriage of justice, a comparison of judicial statistics of different countries and (p. 7) the proof of the applicability of mathematical analysis to "matters that are called moral". On moral sciences see Note 8.

Poisson was mainly interested in studying criminal offences. Unlike **Laplace**, he (pp. 4 and 318) introduced a positive probability ($k$) of the defendant's guilt. One of his formulas (p. 333) determined the probability that the defendant, convicted by ($n – i$) jurors out of $n$, was really guilty:

$$P_i = kt^m/[kt^m + (1 – k)], \; t = u/(1 – u).$$

Here, $m = n – 2i$, and $u$ was the probability of a correct verdict reached by each juror (judge). Poisson noted that the right side did not depend on *n*; however, supposing that *n* was odd (say), *i* could have varied up to its greatest value, $(n – 1)/2$, – and what mattered was the sum of all of the values of $P_i$. Poisson derived a similar formula for a continuous random *u* by introducing its unknown prior density.

One of Poisson's statements (pp. 375 – 376) is debatable: he thought that the rate of conviction should increase with crime.

The application of the theory of probability to jurisprudence continued to be denied. Here are the two most vivid pertinent statements (**Mill** 1843/1886, p. 353; **Poincaré** 1896/1912, p. 20):

1) *Misapplications of the calculus of probability […] made it the real opprobrium of mathematics. It is sufficient to refer to the applications made of it to the credibility of witnesses, and to the correctness of the verdicts of juries.*

2) People "influence each other" and act like the "moutons de Panurge".

Nevertheless, the pertinent work of **Laplace** and Poisson (and of their predecessor, **Condorcet**, § 6.1.5) had undoubtedly attracted the public to the problems of the administration of justice and showed what could be hoped for in the ideal case. I return to Poincaré in § 11.2.

**8.9.2. Medical Statistics.** I mentioned this discipline in § 6.2.3. Now I say that Poisson had certainly contributed to its development. Here is a statement of Gavarret (1840, p. xiii), his former student who later took to medicine:

*Only after long reflection on the lectures and writings of the illustrious geometer, we grasped all the extensiveness of the systematic application of the experimental method in the art of healing.*

In his book, that became very well known, Gavarret explained the normal approximation to the binomial law and the calculation of admissible discrepancies of frequencies in series of Poisson trials and (p. 194) stressed the importance of checking null hypotheses. True, their introduction was a logical consequence of Poisson's deliberations.

In Russia, in the 1850s, **Davidov** (Ondar 1971), who was well acquainted with the work of Poisson and **Cournot** (§ 10.3), popularized the application of the statistical method to medicine, see Ondar (1971). I mention him again in §§ 8.1 and 10.4-8.

Poisson (1837a, p. VI) and Gavarret after him insisted on the need to collect many observations. Here is Poisson:

*Medicine will not become either a science or an art if not based on numerous observations, on the tact and proper experience of the physicians who judge the similarity of cases and take into account exceptional circumstances.*

And so, medicine was then not a science or an art!

Nevertheless, it was difficult to collect many observations. Here is the opinion of a noted physician (Liebermeister ca. 1876, pp. 935 – 940):

*Theoreticians rather often categorically tell us, practical physicians, that all our inferences about the advantages or shortcomings of some methods of treatment, so far as they are based on results which have really taken place, simply remain up in the air if only we do not apply rigorous rules of the theory of probability. […] Until now physicians have applied that theory so seldom not so much because they sometimes did not attach proper significance to it, but mainly since its analytical arsenal was too imperfect and awkward. […] Mathematicians say: If you, physicians, wish to arrive at plausible conclusions, you must invariably work with large numbers; you ought to collect thousands and hundred thousands observations. […] This, however, is impossible for statistics of a general practitioner. And, nevertheless, if this condition is fulfilled, it will often be doubtful whether the theory of probability will be necessary in the same pressing manner. […] Gavarret somewhat arbitrarily presumed, as Poisson also did in several problems, that 0.9953 or 212:213 […] is a sufficient measure of probability. […] Suppose that the successes of two methods of treatment are only as 10:1, would not that be sufficient for preferring the first one?*

He exaggerated, but his criticism was still valid. Then, beginning with 1863 (Chauvenet) and even earlier astronomers and geodesists had begun to offer tests for rejecting outliers quite in vein with Liebermeister's reasoning. Previous practitioners had also made plausible inferences on the strength of scarce data (Bull 1959) whereas **Niklaus Bernoulli** (§ 3.3.2) thought that an absentee ought to be declared dead once his death becomes only twice as probable as the alternative.

Modern statistics including the theory of errors cannot restrict it activity to the case of a large number of observations. Therefore, Liebermeister rather than Gavarret was the pioneer of medical statistics. Freudenthal & Steiner (1966, pp. 181 – 182) mistakenly attached to Gavarret, not to Liebermeister, the transition from unconditional certitude to reasonable degrees of probability. In 1889, Liebermeister's collected medical works had appeared in Leipzig after which he still published a few books. It is worthwhile to mention that Liebermeister in a medical context studied the possibility of distinguishing between equality and inequality of success probabilities in two (small) series of binomial trials. Starting from a **Laplacean** formula based on the existence of uniform prior distribution, and assuming that the two probabilities coincided, he considered the size of the tail probability (of the hypergeometric distribution). His main formula had hardly ever reappeared. See Seneta (1994).

I treat the further history of medical statistics in § 10.8.1, but it is opportune to say now that in the 19[th] century that discipline had still been rather primitive and that the demand advocated by Gavarret (and Poisson) was indeed unrealistic.

## Notes

**1.** When describing his attempt to generalize the normal law, **De Morgan** (1864, p. 421) declared that if the probability of a certain event was 2.5, it meant that the event "must happen twice with an even chance of happening a third time". In 1842, in a letter to **John Herschel** (Sophia De Morgan 1882, p. 147), he stated that "undoubtedly"

$$\sin \varphi = 0, \cos \varphi = 0, \tan \varphi = \mp \sqrt{-1}, \cotan \varphi = \tan \varphi.$$

I did not find Herschel's answer.

**2.** Earlier Poisson (1830, pp. 141 and 146) used the same letter *A* for designating an observed constant, – of "some thing". Consequently, it hardly stood later for "aléatoire".

**3.** Poisson (1837a, p. 274, and earlier (1833, p. 637) corroborated the transition from discrete to continuous by a trick that can be described by **Dirac's** delta-function. When considering density $\varphi(x)$ equal to zero everywhere excepting a finite number of points $c_i$, $i = 1, 2, \ldots, n$, and such that

$$\int_{c_i - v}^{c_i + v} \varphi(x)\, dx = g_i, \quad \sum g_i = 1, \quad v \to 0,$$

Poisson had thus introduced that function of the type of

$$\varphi(x) = \sum g_i \delta(x - c_i).$$

**4.** Poisson was unable to keep to his announced distinction between chance and probability and I am therefore making use of the modern term.

**5.** Poisson referred to his p. 155 and to his memoir (1829, § 8), but, as I see it, the situation remained unclear. Later he (1837a, pp. 312 – 313) repeated the formula of the CLT for the mean value of a random variable without introducing any conditions and even without demanding that its domain was restricted to a finite interval.

**6.** From Laplace to the 1930s the theory of probability had been developing as a branch of applied mathematics (§ 0.1), see also Note 5. And here is Poisson (1837a, § 84):

*There exists a very high probability that these unknown chances little differ from the ratio …*

**7.** From 1812 (and until?) Poisson was "examinateur de l'arme de l'artillerie" (**Arago** 1850, p. 602).

**8.** Laplace (1814/1995, p. 61) urged to

*Apply to the political and moral sciences the method based on observation and the calculus, a method that has served us so successfully in the natural sciences.*

It is difficult to say what exactly is included into moral sciences; see however **Poinsot's** statement below. Beginning at least with **Quetelet**, the study of phenomena depending on free will (although only crimes, suicides, and marriages) was considered to constitute the subject of moral statistics. Then, however, the domain of that branch of statistics essentially broadened and includes now, for example, philanthropy and professional and geographical mobility of the population.

## 9. Gauss, Helmert, Bessel
### 9A. Gauss

This subchapter is mostly devoted to the MLSq[1].

I (1979) have somewhat dwelt on Gauss' investigations in probability proper. He was a tireless collector of statistical data, even of non-essential nature (Biermann 1991), and successfully managed the widows' fund of the Göttingen University. His correspondence and scientific legacy includes a study of the mortality of newly-born and of the members of tontines (of closed societies of mutually insured persons, see § 2.1.3). In the theory of probability, he left the inversion formula for the **Fourier** transform of the density function.

Gauss also solved the first problem in the metric theory of numbers. He considered the expansion of a number *M* (0 < *M* < 1) into a continued fraction with unit numerators and investigated the probability *P*(*n*; *x*) that, beginning with its (*n* + 1)-st convergent, the "tail" of this fraction was less than *x*. If all the permissible values of *M* were equally probable or more or less so, then, as he (*W*-10/1, pp. 371 – 372) explained his problem in a letter of 1812 to **Laplace,** *P*(0; *x*) = *x* and

$$\lim P(n; x) = \frac{\ln(1+x)}{\ln 2}, n \to \infty.$$

Nevertheless, he was not quite satisfied with this solution and asked Laplace to have a look at the problem. He, Gauss, was sure that Laplace will find a *plus complete* solution, – a pre-limiting expression. A phrase from Gauss' *Mathematisches Tagebuch* written in 1800 (p. 552 of the *Werke*) testifies that Gauss had already derived the equality above – and had then been satisfied with his work.

Stäckel (Gauss, *W*-10/1, pp. 554 – 556) and then Kuzmin (1928) proved this equality and the latter also derived an asymptotic expansion for *P*(*n*; *x*).

Here, I also repeat in a few words Gauss' general opinion (*W*-12, pp. 201 – 204) about the applications of the theory of probability as described by **W. E. Weber** in one of his letters of 1841. If only based on numbers, Gauss reasoned, such applications could be greatly mistaken; the nature of the studied subject ought also to be taken into account. However, probability provides clues when nothing except numbers is known, as for example when dealing with annuities; and in jurisprudence and it can determine the desired number of witnesses and jurors [but hardly without allowing for "the nature" of law courts]. See also Gauss' opinion about the rejection of outlying observations in § 9A.5-2. As shown in the sequel, quite a few authors shared Gauss' viewpoint.

### 9A.1. The Method of Least Squares before 1809

It had been indirectly and inaccurately applied from the mid-18$^{th}$ century (§ 6.3.2) and its peculiar version was possibly known even earlier. When some point *P* was graphically intersected from three or more given stations, a triangle, or a polygon of errors appeared on the surveyor's table sheet and it was apparently natural to select the position of *P* by eye in such a manner that the sum of the squares of its distances from the sides of the triangle (of the polygon) was minimal. To a certain extent I can justify my opinion by mentioning an experimental smoothing of a broken line by eye (Tutubalin 1973, p. 27): on the whole, the curves thus drawn were as accurate as if having been determined by the MLSq.

Euler (1778) can be considered as Gauss' predecessor in the heuristic sense (§ 6.3.1), but Gauss possibly did not see that memoir (a commentary on Daniel Bernoulli's memoir of the same year). The regrettably incomplete list of books which Gauss had borrowed from the library of the Göttingen university (Dunnington 1955, pp. 398 – 404) does not include the volume of the contributions of the Petersburg Academy in which those two papers were published. However, in his letters Gauss expressed his surprise at the fact that the principle of least squares was not discovered earlier.

**9A.1.1. Huber.** Many authors, for example Merian (1830, p. 148), stated that somewhat before 1802 the Swiss mathematician and astronomer Huber had discovered the principle of least squares, but that, living far from scientific centres, he had not reported his finding to anyone. However, Dutka (1990), who referred to a forgotten paper (W. Spieß 1939), concluded otherwise. It occurs that Spieß quoted Huber himself who had mentioned "**Legendre's** criterion [Maßstab] of least squares".

**9A.1.2. Legendre** (1805, pp. 72 – 73) introduced the principle of least squares:

*Among all the principles that can be proposed* [for solving redundant systems of linear equations], *I think there is no one more general, more exact, and easier to apply, than that which we have made use of in the preceding researches and which consists in making the sum of the squares of the errors* [of the residuals] *a <u>minimum</u>. In this way there is established a sort of equilibrium among the errors, which prevents the extremes to prevail and is well suited to make us know the state of the system most near to the truth.*

Translation by Hald (1998, p. 119). Legendre also indicated that the absolute values of the extremes [again: of the residuals] should be confined within the shortest possible interval. He had not added that it was the minimax principle (§ 6.3.2) rather than his innovation that ensured his desire.

**9A.1.3. Adrain.** The American mathematician Adrain (1809) justified the principle of least squares and the [normal distribution][2] at about the same time as **Gauss** did and applied it to the solution of several problems, see below (Dutka 1990). He also indicated that the lack of space prevented him to discuss the adjustment of pendulum observations. About ten years later he (1818a) published that study in which he revealed two mistakes in **Laplace's** pertinent calculations (1798 – 1825, t. 2, § 42 of Livre 3). The same year his derivation of the length of the major semi-axis of the Earth's ellipsoid of rotation (1818b) appeared. Incidentally, that length (6378.629 *km*) was sufficiently close to a determination of 1940 by **F. N. Krasovsky** (6378.245 *km*)**.**

Adrain's main paper was first mentioned much later (C. Abbe 1871) but his second article had become known to **Olbers** (W/Erg-4, No. 2, p. 711) who informed **Gauss** about it. An American author, wrote Olbers to Gauss, had mentioned his previous paper and "ascribed" the MLSq to himself. Gauss hardly made any comment; the priority strife with **Legendre** was apparently enough for him.

Here are Adrain's derivations of the normal distribution.

a) Lines *a* and *b* are measured in the field with errors *x* and *y* respectively and

$$x/a = y/b \qquad (1)$$

and the total error is fixed:

$$x + y = c. \qquad (2)$$

Introducing the density of the observational errors $\varphi$ and tacitly assuming their independence, Adrain applied the principle of [maximum likelihood]

$$\varphi(x; a)\varphi(y; b) = \max$$

so that, after allowing for conditions (1) and (2),

$$[\varphi'(x; a)/\varphi(x; a)]dx + [\varphi'(y; b)/\varphi(y; b)]dy = 0, \quad \varphi'(x; a)/\varphi(x; a) = mxa, \text{ etc.}$$

b) Suppose that for linear measurements

$$x^2 + y^2 = r^2,$$

then

$$W = \varphi(x)\varphi(y) - \lambda(x^2 + y^2) = \max, \quad \varphi'(x)\varphi(y) - 2\lambda x = 0,$$
$$\varphi(x)\varphi'(y) - 2\lambda y = 0,$$
$$\varphi'(x)/x\varphi(x) = \varphi'(y)/y\varphi(y) = c, \text{ etc.}$$

Adrain then wrote out the joint distribution of both these errors and indicated that the appropriate contour lines of equal probabilities were ellipses (ellipses of errors, as they were later called in the theory of errors).

Conditions (1) and (2) hardly conform to reality; thus, the former describes the action of systematic errors. Also arbitrary is the condition applied in the second justification. Nevertheless, **John Herschel** (1850), **Maxwell** (1860), **Thomson** & **Tait** (1867, p. 314) and **Krylov** (1950, Chapt. 8) repeated that demonstration without any references (Sheynin 1965). Later on **Kac** (1939) and **Linnik** (1952) weakened the condition of independence.

Adrain was now able to prove quite simply that the arithmetic mean of direct measurements was optimal; this, of course, conformed to the principle of least squares in case of several unknowns. Finally, Adrain showed how to adjust a traverse (a polygon with measured sides and bearings) by the principle of least squares and, what is also remarkable, he calculated corrections to directly measured magnitudes rather than to their functions which were not independent from each other, see Sheynin (1965), a paper with which in general I am now dissatisfied.

**9A.1.4. Gauss.** He (1809a; 1809b, § 186) applied the principle of least squares from 1794 or 1795. In the second instance, he called it "our principle": *Unser Princip, dessen wir uns seit dem Jahre 1795 bedient haben*, and in both cases he mentioned Legendre. Much later Gauss (1823a) once more mentioned Legendre, although not in his main memoir of 1823.

**Legendre** (letter to Gauss 31.5.1809, see Gauss, *W*-9, p. 380) was badly offended. He rightfully stated that priority is only established by publication. Gauss did not answer him, and Legendre (1820, pp. 79 – 80) charged him with appropriating the MLSq.

Many authors commented on this episode, and I quote May (1972, p. 309) and Biermann (1966, p. 18):

*Gauss cared a great deal for priority.* […] *But to him this meant being first to discover, not first to publish; and he was satisfied to establish his dates by private records, correspondence, cryptic remarks in publications.* […] *Whether he intended it so or not, in this way he maintained the advantage of secrecy without losing his priority in the eyes of later generations.*

*What is forbidden for usual authors, ought to be allowed for Gausses and in any case we must respect his* [Gauss'] *initial considerations.*

It seems that Legendre could have stated in 1820 that he was the inventor of the MLSq and that, in spite of Gauss' claim, everyone will agree with him. Regrettably, this did not happen. Laplace (end of § 8.8) objectively described the discovery of the MLSq but did not add that Legendre had not really substantiated it. Then, he offered his own version of the theory of errors only suitable for a large number of observations and for other conditions necessary for the CLT, and other

French mathematicians including Poisson followed him and Legendre and did not even mention Gauss. This seriously diminished the importance of their work.

In his letter to Gauss, Legendre also noticed that Euler rather than Laplace, as Gauss (1809, § 177) stated, calculated the integral of the exponential function of the negative square. Later Gauss (*Monatliche Correspondenz*, Bd. 21, p. 280) explained that he noticed his mistake when the printing of his book was almost completed and that it was Laplace who provided the final formulation of the result. The Editors of Gauss (1887), who noticed that explanation did not supply the year of that Bd. 21.

On 30 Jan. 1912 Gauss (W-10/1, p. 373) only answered Laplace: he had applied the MLSq long before 1805 but had no desire to publish a fragment. Gauss' preliminary reports had appeared in the *Göttingische gelehrte Anzeigen* and Legendre hardly saw the report (1809a). Later it was reprinted (W-6, pp. 59 – 60), and here is an excerpt:

*The author has been applying the main principles which are here considered for fourteen years now and long ago had communicated to his astronomical friends. They lead to the same method which Legendre* […] *published a few years ago*.

I (Sheynin 1999a; 1999c) described the possible cases in which Gauss could have applied the MLSq before 1805 and named many of his colleagues and friends to whom he had communicated his discovery. Unexpectedly, it occurred that **von Zach**, who allegedly refused to testify to Gauss' priority, had not until 1805 known the formulation of the principle of least squares, and, furthermore, that he (1813, p. 98n) indirectly agreed with the latter's statements by repeating them without any qualification remarks:

*The celebrated Dr Gauss was in possession of that method since 1795 and he advantageously applied it when determining the elements of the elliptical orbits of the four new* [minor] *planets as it can be seen in his excellent work* [*Theoria motus*].

Regrettably, "it" is not seen there.

This passage is even more important than Zach's editorial acceptance of Gauss' priority (noticed by Dutka 1996, p. 357). In 1809, Zach's periodical, *Monatliche Correspondenz*, carried an anonymous review of Gauss' *Theoria motus*, and there, on p. 191, Gauss' pertinent claim was repeated.

Gauss' early use of the MLSq is not generally accepted, see for example Marsden (1995, p. 185) who nevertheless had not mentioned the opposite opinion of Brendel (1924) and Galle (1924, p. 9) or of Gauss' contemporaries[3]. In any case, Gerardy (1977), drawing on archival sources, discovered that Gauss, in 1802 – 1807, had participated in land surveying (in part, for his own satisfaction) and concluded, on p. 19 (note 16) that Gauss started using the method not later than in 1803. Regrettably, Gerardy concentrated on describing

Gauss' simple calculations and his statement mentioned just above was not quite definite.

There are many other instances including that mentioned by von Zach (above) in which Gauss could have well applied his invention at least for preliminary, trial calculations, or short cuts. For him, the MLSq was not a cut and dry procedure, see § 9A.5-3. Then, weighing the observations as well as the presence of possible mistakes in the data could have hampered justification.

As to the communication of his discovery, I proved that among those whom Gauss had informed before 1805 were **Bessel** and **Wolfgang Bolyai** (the father of the cofounder of the non-Euclidean geometry, **Janos** or **Johann Bolyai**), and **Olbers** about whom it was known long ago. In 1812, Olbers promised Gauss to state publicly that he came to know about the MLSq from Gauss a few years before Legendre's publication, but he only did so in 1816. In 1812 – 1815, Olbers just did not publish anything suitable for inserting such a remark, see the appropriate volume of the *Catalogue of Scientific Literature* of the Royal Society.

### 9A.2. *Theoria Motus* (1809b)

In accordance with the publisher's demand, this book appeared in Latin. Its German original is lost and Gauss' correspondence (letter from **Olbers** of 27.6.1809, see W/Erg-4, No. 1, p. 436) proves that, while translating, he essentially changed its text. The treatment of observations occupies only a small part of the book.

1) The **Boscovich** method (see § 6.3.2-3). Suppose that $n$ equations (1.2) in $m$ unknowns ($n > m$) are adjusted by that method. Then, as Gauss (§ 186) remarked, equation (6.13) meant that exactly $m$ residual free terms will be zero. Somewhat below, in the same § 186, Gauss qualified his statement by taking into account the other **Boscovich** equation (6.12) but mistakenly attributed it to **Laplace**. In § 174 he stated that the formulated corollary was undesirable although in §§ 188 – 189 he apparently agreed that the Boscovich method (§ 6.3.2 might ensure a first approximation. His remark, that can be easily proved, means that he knew an important theorem in linear programming.

2) The [normal distribution] (§§ 175 – 177). Gauss (§ 177) assumed "as an axiom" that the arithmetic mean of many observations was the most probable value of the measured constant "if not absolutely precisely, then very close to it". He (§ 175) derived the density of observational errors (although not introducing any new term) believing that it was [unimodal] and "in most cases" even; this, then, was his understanding of the properties of random errors. Finally, in order to justify the principle of [maximal likelihood], Gauss (§ 176) proved the "fundamental principle" of inverse probability, see my § 7.1-1, for the case of equal probabilities of the various hypotheses. However, the principle of the arithmetic mean (above) already implied this restriction (**Whittaker** & **Robinson** 1924/1949, p. 219).

And so, if the observations are denoted by $x_i$, $i = 1, 2, …, n$, then, according to the principle of maximal likelihood,

$$[\ (x_1 - a)/\ (x_1 - a)] + [\ (x_2 - a)/\ (x_2 - a\ )] + ... +$$

$$[\sum(x_n - a)/\varphi(x_n - a)] = 0,$$

where *a* is the estimator sought, coinciding, as stipulated, with the arithmetic mean $x_o$. If

$$x_i = x_1 - nN, \ i = 2, 3, \ldots, n,$$

then

$$x_1 + (x_2 + x_3 + \ldots + x_n) = x_1 + (n-1)x_1 - n(n-1)N,$$
$$N = (x_1 - x_o)/(n-1),$$
$$\varphi(x_1 - x_o)/\varphi(x_1 - x_o) = (1-n)\varphi(-N)/\varphi(-N) = -(1-n)\varphi(N)/\varphi(N),$$
$$[\varphi N(n-1)]/\{(1-n)\varphi[N(n-1)]\} = -\varphi(N)/\varphi(N),$$
$$\varphi(x)/x\varphi(x) = \text{Const},$$

$$\varphi(x) = (h/\sqrt{\pi})\exp(-h^2x^2), \ h > 0. \tag{3}$$

Gauss (§ 178) called *h* the "measure of precision" (*gradus praecisionis*). It might be supposed that, from the very beginning, he was not satisfied with his derivation. His wording of the principle of the arithmetic mean and of the properties of the density of observational errors contained qualification remarks whereas the obtained principle of least squares (see below § 9A.2-3) occurred to be an axiom. Again, it is difficult to believe that Gauss was pleased with the appearance of a universal law of error. Later he (1821/1887, pp. 193 and 194; 1823a/1887, p. 196) remarked that his derivation had depended on a hypothetically assumed distribution. And here is **Bertrand's** opinion (1888a, p. XXXIV): Gauss had not claimed to establish the "vérité", he attempted to search for it. Bertrand (pp. 180 – 181) also remarked that the mean of the values of some function did not coincide with the mean value of its arguments, which, in his opinion, testified against the principle of arithmetic mean. Gauss, however, considered direct measurements. Note also that he (his letter to **Encke** of 1831; *W*-8, pp. 145 – 146) "not without interest" acquainted himself with the attempt of his correspondent to justify the arithmetic mean by deterministic analytical axioms. Many authors made similar efforts and Zoch (1935) concluded that, although they were unsuccessful, the postulate of the arithmetic mean can nevertheless be established without stochastic considerations. His finding was unrelated to the theory of errors, but the pertinent investigations apparently served as the point of departure for the theory of invariant statistical hypotheses and estimators (Lehmann 1959, Chapter 6). Encke (1832, p. 74), as it occurred, was not satisfied with either substantiations of the MLSq.

Gauss (1845/1873, p. 143) left a lesser known statement about the arithmetic mean. He remarked that the random variations corrupting observations mostly compensate one another so that the mean becomes ever more reliable as the number of observations increases. This is "generally absolutely right", and often led to "splendid results" in natural sciences. However, Gauss continued, an important condition, often overlooked and difficult to check, was that the disordered variations ought to be entirely independent from each other, cf. § 9A.4-4.

3) The principle of least squares (§ 179) followed immediately. Gauss, however, added that, similar to the principle of the arithmetic mean, it should be considered an axiom [considered as a corollary of an axiom?]. A special point here is that, instead of the real errors the principle of least squares was formulated with regard to residual free terms. **Helmert** (1872, p. 75) indicated this fact but paid scant attention to it and had not mentioned Gauss. Apparently he had not realized that the normal law was [stable], cf. §§ 7.2.7.

4) The precision of the arithmetic mean. Gauss, naturally, restricted his attention to the case of the [normal distribution]. Later he (§ 9A.4) certainly abandoned this restriction.

5) The precision of a random sum (marginal note to § 183, included not in the German translation, but in Gauss' *W*-7). Suppose that

$$x = a + b + c + ...,$$

then

$$h_x = 1 \div [(1/h_a^2) + (1/h_b^2) + (1/h_c^2) + ... ]^{1/2}.$$

Gauss did not explain his note; it might be supposed that the terms above were normally distributed since he only introduced *h* for that law. However, he may well have derived this formula in the general case. Incidentally, it shows that at least in the case of the normal distribution the random error of the sum increases as the square root of the number of its terms.

6) The precision of the [estimators of the] unknowns (§ 182; 1811, § 13). Suppose that these estimators are determined by solving a system of normal equations in accordance with the Gauss method of successive eliminations. Then, assuming that the precision of a direct measurement is unity, the precision of the estimator of the last unknown is equal to the root of its coefficient in the last reduced equation. Also see my § 9A.4-5.

### 9A.3. *Determining the Precision of Observations* (1816)

1) The precision of the measure of precision *h* in formula (3). Suppose that the errors of *m* [independent] observations are , , , … Then the most probable value of that magnitude is determined by the condition

$$h^m \exp[- h^2 (\ ^2 +\ ^2 +\ ^2 + ...)] = \max$$

and is therefore equal to

$$h_o = \{m/[2(\ ^2 +\ ^2 +\ ^2 + ...]\}^{1/2} = 1/\ \sqrt{2}.$$

In the last expression, which is my own,   is the mean square error of an observation. Gauss also indicated that

$$P(h_o - \delta \le h \le h_o + \delta) = \Phi(\delta\sqrt{m}/h_o),\ \Phi(t) = (2/\sqrt{\pi})\int_0^t \exp(-z^2)dz,$$

so that, for $P = 1/2$, $\delta = \rho h_o/\sqrt{m}$, $\rho \approx 0.477$. In addition, for distribution (3),

$$P(|\xi| \le h) = 1/2,\ \text{and}\ r = \rho/h$$

is the probable error formally introduced by **Bessel** (1816, pp. 141 – 142).

Let

$$S_n = |\xi_1|^n + |\xi_2|^n + |\xi_3|^n + ...,\ K_n = \int_{-\infty}^{\infty} x^n \varphi(x)dx,$$

then, for large values of $m$,

$$P(-\beta \le S_n - mK_n \le \beta) = \Phi\{\beta/[2m(K_{2n} - K_n^2)]^{1/2}\}, \qquad (4)$$

where $mK_n$ is the most probable [the mean] value of $S_n$. In actual fact, Gauss treated absolute moments and the formula for $K_n$ should be corrected. Formula (4) was proved by **Helmert** (1876b) and then by Lipschitz (1890), but **Cramér** (1946, § 28.2) noted that it was a particular case of the CLT.

Finally, Gauss derived a formula for the absolute moments of the normal law

$$mK_n = S_{n0} = m\ \Gamma[(n-1)/2]/h^n\sqrt{\pi},\ \Gamma(x) = \Pi(x+1),$$

so that $h$ (and therefore $r$) could have been estimated by $\bar{S}_n$, the mean value of $S_n$. Comparing the probable intervals of $r$ for different $n$, Gauss concluded that $n = 2$ secured its best estimator.

In one of his letters of 1825 Gauss (*W-8*, p. 143) objected to the probable error as "depending on a hypothesis" [on the law of distribution]. Still, again in his correspondence (Sheynin 1994a, p. 261), and even in a paper (1828b), he applied it quite a few times. Natural scientists, for example **Mendeleev** (§ 10.9.3) and **Newcomb** (§ 10.9.4), followed suit and Bomford (1971, pp. 610 – 611) "reluctantly" changed from probable to mean square error in that edition of his book. However, L. O. Struve (1887, last, unnumbered, page) proposed to abandon the probable error.

2) Denote $1/h\sqrt{2} = \sigma$ and let $n = 2$. Then

$$[m(K_4 - K_2^2)]^{1/2} = \sigma^2\sqrt{2m}$$

and, in accordance with formula (4), the sum of squares $S_2$ is distributed normally $N[m\sigma^2;\ \sigma^2\sqrt{2m}]$. This is the asymptotic chi-squared distribution, cf. **Cramér** (1946, § 20.2).

### 9A.4. *The Theory of Combinations* (1823 – 1828)

I consider the main part of this memoir in which Gauss provided his definitive justification of the MLSq by the principle of maximum weight [of minimal variance], and I add a few words about its supplement (1828a).

1) Random errors and the density of observational errors. Gauss (§§ 1 – 3) distinguished between random and systematic errors but had not provided their formal definition. He (§ 4) then repeated (see my § 9A.2-2) the definition of density and listed its properties. The mean value of the errors (§ 5) was equal to zero; otherwise, as Gauss additionally remarked, it determined the action of constant errors.

2) The measure of precision. Gauss (§ 6) introduced a measure of precision [the variance]

$$m^2 = \int_{-\infty}^{\infty} x^2 \varphi(x)dx$$

calling it the mean error to be feared, – *des mittleren zu befürchtenden Fehler, errorum medium metuendum* (1821/1887, p. 194; 1823b, § 7).

In his letters to Encke of 23 Aug. 1831 (W-8, pp. 145 – 146), to **Bessel** of 28.2.1839 (Ibidem, pp. 146 – 147), to Schumacher of 25 Nov. 1844 (reported by Helmert in his Introduction to Gauss 1887) and in § 7 of the *Theor. motus* Gauss stressed that an integral measure of precision was preferable to a local measure. He (1823b, § 6) also indicated that the quadratic function was the simplest [from integral measures], and in 1821 he (1887, p. 192) dwelt on his choice in more detail: it was also connected with

*Some other, extremely important advantages which no other function possesses. However, any other even degree could have been selected as well …*

Could have been chosen in spite of the advantages of the variance? **Bienaymé** (1853/1867, pp. 167 – 169) proved that a formula of the type of (5), see below, was not valid for any other even exponent; a clear exposition of this proof is due to Idelson (1947, pp. 269 – 271). Therefore, Bienaymé continued, the choice of the variance was unavoidable. I doubt, however, that, as he believed (p. 169), Gauss was here mistaken. The sample variance (see § 9A.4-6) is distribution-free.

3) An inequality of the **Bienaymé – Chebyshev** type. Gauss (§ 9) examined the probability

$$\mu = P(|\xi| \leq \lambda m) = \int_{-\lambda m}^{\lambda m} \varphi(x)dx$$

for a [unimodal] density of observational errors $\varphi$ having variance $m^2$ and proved (§ 10) that

$$\lambda \geq \sqrt{3} \text{ for } \mu \leq 2/3 \text{ and } \lambda \geq \frac{2}{3\sqrt{1-\mu}} \text{ for } \mu \geq 1.$$

**Cramér** (1946, § 15.7 and Example 4 to Chapters 15 – 20) more easily proved this "remarkable" theorem, as Gauss called it, whereas Seal (1967/1970, p. 210) indicated, that Gauss had wished to abandon the universality of the normal distribution since it occurred that, anyway, $P(|\ | \ 2m) \ 0.89$. But should we forget his own, although indirect arguments and doubts?

4) Independence. In § 18 Gauss offered his definition, although not quite formal, of independent functions of observations: they should not contain common observations. In § 19 he specified that those functions were linear; otherwise his statement would have contradicted the Student – Fisher theorem on the independence of the sample variance and the arithmetic mean.

Therefore, if some observation was common for two functions of observational results, the errors of these functions will not be independent from one another and the mean value of their product will not therefore vanish[4]. In one of his examples, Gauss calculated the variance of a linear form of independent random variables.

Gauss (1809b, § 175; 1823b, § 15) mentioned independence even earlier but without explanation, and, later he (1826/1887, p. 200; 1828, § 3) described the mutual dependence of magnitudes known from observation by the existence of functional connections between them. This meant, for example, that the adjusted angles of a triangle, since their sum was equal to 180° plus the spheroidal excess, were dependent on one another. See also end of § 9A.2-2.

In mathematical statistics the definition of independence is different. An orthogonal transformation of independent and normally distributed magnitudes leads to their as though "adjusted" values, – to their linear forms of a certain type, which are nevertheless independent (the **Fisher** lemma; **Cramér** (1946, § 29.2)). Here is **K. Pearson's** appropriate statement (1920/1970, p. 187) which I do not however understand: for Gauss

*The observed variables are independent, for us* [they] *are associated or correlated. For him the non-observed variables are correlated owing to their known geometrical relations with observed variables; for us,* [they] *may be supposed to be uncorrelated causes, and to be connected by unknown functional relations with the correlated variables.*

According to Krengel (2011), the modern notion of independence of events is due to Bohlmann whom Kolmogorov in 1933 had not mentioned. Kolmogorov introduced independence of events and random variables.

5) The principle of maximum weight for [unbiased] estimators. Gauss described this subject ponderously. For that matter, **Helmert** (1872) and **Idelson** (1947) are in general much better understood. Suppose that, without loss of generality, the initial equations are

$a_i x + b_i y = G_i = g_i + \ _i, i = 1, 2, …, n$

where $\varepsilon_i$ is the error of the free term $g_i$. The estimators of the unknowns might be represented by linear forms, for example by $x = [\alpha G]$ with unknown coefficients $\alpha_i$ so that

$$m_x^2 = [\alpha\alpha] m^2 \qquad (5)$$

where $m^2$ is the variance of an observation.

It is easy to prove that $[a\alpha] = 1$, $[b\alpha] = 0$ and the condition of maximal weight will be

$$W = [\alpha\alpha] - 2Q_{11}[a\alpha] - 2Q_{12}[b\alpha] = \max$$

where $Q_{11}$ and $Q_{12}$ are the **Lagrange** multipliers. Similar considerations, and, in particular, an estimation of precision resembling formula (5), are also possible for the other unknowns. It occurs that the estimators of the unknowns are determined from the normal equations and their weights are calculated by means of the Lagrange multipliers of the type of $Q_{ii}$ which, like the other multipliers $Q_{ij}$, are determined from the same normal equations with partly unit and partly zero free terms. Thus, in formula (5) $[\alpha\alpha] = Q_{11}$.

According to the above, it follows that such formulas can be made use of even before observation; the general layout of the geodetic network and the crude values of its angles obtained during reconnaissance make it possible to calculate the $Q_{ij}$. And (what Gauss had not known) these multipliers are connected with covariations; thus, $Q_{12} = E(xy)$.

6) The estimator of the sample [variance]. Gauss (§§ 37 – 38) proved that, for $n$ observations and $k$ unknowns, the unbiased sample variance and its estimator were, respectively,

$$m^2 = E[vv]/(n - k), \quad m_o^2 = [vv]/(n - k) \qquad (6a, b)$$

where $v_i$ were the residual free terms of the initial equations. Instead of the mean value, the sum of squares $[vv]$ itself has to be applied. Coupled with the principle of maximal weight (of least variance), formulas (6) provide effective estimators, as they are now called. Gauss (1823a/1887, p. 199) remarked that the acceptance of his formula (6b) instead of the previous expression (§ 7.2-5), whose denominator was equal to $n$, was demanded by the "dignity of science".

Gauss stressed that the estimator (6) was unbiased; however, the practically applied estimator is not $m^2$, but the biased $m$. Furthermore, unbiased estimators do not exist in every case and some bias is allowed (Sprott 1978, p. 194). Finally, I note that **Czuber** (1891, p. 460) testified that Helmert had thought that $\text{var} m_o^2/m_o^2$ was more important than $\text{var} m_o^2$ by itself and **Eddington** (1933, p. 280) expressed the same opinion. Czuber also proved that, for the normal distribution, that relative error was minimal for the estimator (6b).

7) The precision of the estimator of the sample variance. Gauss (§ 40) directly calculated the boundaries of the $\text{var } m_o^2$ by means of the fourth moment of the errors and indicated that for the normal distribution

$$\mathrm{var} m_o^2 = 2m^4/(n-k). \tag{6c}$$

He somehow erred in calculating the abovementioned boundaries, see (15); in addition, his formulas should have included the unknown magnitude E $\varepsilon_i^2$ ($\varepsilon_i$ were the observational errors) rather than $m^2$. Formula (6c) shows that $m_o^2$ is a consistent estimator of the sample variance; this property persists in the general case, see formulas (15).

Many years later Bertrand (1888a) criticized the Gauss formula (6c). Tacitly assuming the normal distribution, he provided an example in which his own estimate of $\sigma^2$ was less than that provided by Gauss. He forgot, however, that formula (6c) provided an unbiased estimate whereas his own estimate was biased. Then, he calculated $\sigma^2$ forgetting the Gauss formula for the case of normal distribution. It was this episode that led Czuber to the discussion described in § 9A.4-6.

8) *Other topics*. Gauss also determined the variance of a linear function of the estimators of the unknowns (which are not independent) and provided expedient procedures for further calculations after additional data become known or after the weights of some observations have to be changed.

9) *Another manner of adjusting observations*. In the supplement (1828a) to his memoir Gauss described the adjustment of observations by the MLSq according to the pattern of conditional observations. In geodetic practice, it is often expedient to issue from the directly measured magnitudes and *conditional equations* rather than from observational equations (1.2). Sometimes both kinds of equations are made use of at the same time, but I leave this case aside and consider now a (later) typical chain of, say, 10 triangles of triangulation. Each angle is measured as are the lengths of two extreme sides (baselines) whose directions (azimuths) are determined by astronomical observations. The observational errors are such that both the baselines and the azimuths might be considered exact; only the angles are adjusted. Each measured angle $q_i$ provides an equation

$$x_i - q_i = v_i \tag{7}$$

where the first term is the true value of the angle and the right side is the sought correction. Now, the condition of closing the first triangle (I disregard its excess) is

$$x_1 + x_2 + x_3 - 180° = 0. \tag{8}$$

Extremely simple is also the condition that demands that the azimuth of the first baseline plus the algebraic sum of the appropriate angles be equal to the azimuth of the second baseline. The sine theorem is however needed for the transition from the first baseline to the second one, but a first approximation is achieved by introducing the measured angles so that the required trigonometric equation is linearized. It follows that all the conditions can be written as

$$[av] + w_1 = 0, [bv] + w_2 = 0, \textit{etc.} \qquad (9)$$

Formed by means of equations (7), they should be exactly fulfilled and the number of the terms in the square brackets is either three, as in equations of the type of (8), or more, depending on the number of the triangles in the chain. The adjustment proper consists in determining the conditional minimum of $[vv]$ with the usual application of the **Lagrange** multipliers and the corrections $v_i$ are determined through these multipliers. Strangely enough, only **Helmert** (1872, p. 197) was the first to provide such an explanation.

10) **A new exposition of the memoir** (Sheynin 2012a; 2014). Gauss could have derived formulas (6) in the very beginning of the memoir, just after he introduced the variance (§ 9A.4-2) since the required conditions (linearity of the initial equations, (physical) independence of their free terms and unbiasedness of the estimators of the unknowns) were not connected with the further exposition. And the MLSq directly followed from those formulas (6). Hundreds of textbooks had been describing the MLSq as justified by Gauss in 1809: such an approach was incomparably easier. Now we see that it is quite possible to follow the memoir of 1823. Its very existence had been barely known, see Eisenhart (1964, p. 24) in § 7.1-3.

The most eminent scientists (Boltzmann 1896/1909, p. 570; Chebyshev, see § 13.2-7) had been barely acquainted with the work of Gauss.

Many authors beginning with Gauss himself had derived the formula (6) which was not difficult. The main point, however, is that the proof does not depend on the condition of least squares. On the contrary, this condition can now be introduced at once since it means minimum variance. The formulas derived by Gauss for constructing and solving the normal equations and calculation of the weights of $\hat{x}, \hat{y},...$ and of their linear functions will still be useful.

Gauss had actually provided two justifications of least squares (of which I only left the second one), but why did not he even hint at this fact? I can only quote Kronecker (1901, p. 42) and Stewart (Gauss 1823b – 1828/1995, p. 235):

*The method of exposition in the "Disquisitiones [Arithmeticae", 1801] as in his works in general is Euclidean. He formulates and proves theorems and diligently gets rid of all the traces of his train of thoughts which led him to his results. This dogmatic form was certainly the reason for his works remaining for so long incomprehensible.*

*Gauss can be as enigmatic to us as he was to his contemporaries.*

Gauss himself actually said so. His eminent biographer, Sartorius von Waltershausen (1856/1965, p. 82) testified: He had *used to say* that, after constructing a good building, the *scaffolding* should not be seen. And he had often remarked that his method of description *strongly hindered* readers *less experienced* in mathematics.

Finally, I note Gauss' words (letter to W. Olbers 30.7.1806): *Meine Wahlspruch* [motto] *ist aut Caesar, aut nihil*.

The second substantiation of the MLSq can be accomplished by applying the notions of multidimensional geometry (**Kolmogorov** 1946; **Hald** 1998, pp. 473 – 474). Nevertheless, the new exposition of the memoir of 1823 is essential, and it appeared more than 200 years after its publication!

Kolmogorov (p. 64) also believed that the formula for $m^2$ (6a) should, after all, be considered as its definition. Much earlier Tsinger (1862, § 33) stated that it already "concealed" the MLSq which, however, was only a hint at the real possibility of understanding Gauss. Harter (1977, p. 28) stated almost the same.

### 9A.5. Additional Considerations

Having substantiated the MLSq, Gauss nevertheless deviated from rigid rules; one pertinent example is in § 6.3.2. Here, I have more to say.

1) The number of observations. In his time, methods of geodetic observations were not yet perfected. Gauss himself was successfully developing them and he understood that a formal estimation of precision can only describe the real situation after all the conditions (§ 9A.4-9) were allowed for, i.e., only after all the field work was done. It is no wonder, then, that Gauss continued to observe each angle at each station until being satisfied that further work was useless, see the very end of Chapter 6.

2) Rejection of outliers. This delicate operation does not yield to formal investigation since observations are corrupted by systematic errors, and, in general, since it is difficult to distinguish between a blunder and a "legitimate" large error. Statistical tests, that had appeared in the mid-19[th] century, have not been widely used in the theory of errors. Here is the opinion of the authors (Barnett & Lewis 1978, p. 360) of a book on this subject:

*When all is said and done, the major problem in outlier study remains the one that faced the very earliest workers […] – what is an outlier and how should we deal with it?*

I still refer to Laplace (1918, p. 534) and Gauss (letter to **Olbers** of 1827, W-8, pp. 152 – 153):

*For a successful application of the calculus of probability to geodetic observations one should honestly report about all of his own observations which he admitted and not reject any of them solely because they are somewhat* (!) *remote from the rest*.

*When the number of observations was not very large, and a sound knowledge of the subject was lacking, rejection was always doubtful, and in any case nothing should be concealed so that others will be able to decide otherwise.*

3) Calculations. Without even an arithmometer, Gauss was able to carry out difficult calculations; once he solved a system of 55 normal

equations (letter to **Olbers** of 1826; W-9, p. 320). For other examples see Sheynin (1979, p. 53). His preparatory work (station adjustment; compilation of the initial equations, see § 9A.4-9, and of the normals themselves) had to be very considerable as well.

Sometimes Gauss applied iterative calculations (letter to **Gerling** of 1823; W-9, pp. 278 – 281), also see Forsythe (1951) and Sheynin (1963). The first to put on record this fact, in 1843, was Gerling himself. Then, Gauss (1809b, § 185) left an interesting qualitative remark stating that "it is often sufficient" to calculate approximately the coefficients of the normal equations. The American astronomer Bond (1857) had applied Gauss' advice and **Newcomb** (1897a, p. 31) followed suit.

As a calculator of the highest calibre (Maennchen 1918/1930, p. 3),

*Gauss was often led to his discoveries by means of mentally agonizing precise calculations* […]; *we find* [in his works] *substantial tables whose compilation would in itself have occupied the whole working life of some calculators of the usual stamp.*

I ought to add that Gauss made some mistakes in his computations possibly because, first, he had not invariably checked them, see for example Gerardy (1977) or his own methodological note (1823c) where the signs of *dx* and *dy* were wrong. Second, Gauss calculated "unusually fast" (Maennchen 1918/1930, p. 65ff) which caused mistakes and additional difficulties in proving that he had applied the MLSq before 1805.

Maennchen did not study Gauss' geodetic calculations possibly because in his time the solution of systems of linear equations had not yet attracted the attention of mathematicians.

For my part, I note that, when compiling a certain table of mortality, Gauss (W-8, pp. 155 – 156) somehow calculated the values of exponential functions $b^n$ and $c^n$ for $n = 3$ and 7(5)97 with $\lg b = 0.039097$ and $\lg c = -0.0042225$.

Here, now, is **Subbotin's** conclusion (1956, p. 297) about the determination of the orbits of celestial bodies but applicable to my subject as well: **Lagrange** and **Laplace**

*Restricted their attention to the purely mathematical aspect* [of the problem] *whereas Gauss had thoroughly worked out his solution from the point of view of computations taking into account all the conditions of the work of astronomers and* [even] *their habits.*

4) Estimation of precision (Sheynin 1994a, pp. 265 – 266). In his letters to **Bessel** (in 1821) and **Gerling** (in 1844 and 1847) Gauss stated that the estimation of precision based on a small number of observations was unreliable. In 1844 he combined observations made at several stations and treated them as a single whole, cf. **Laplace's** attitude (§ 9A.2-7). And in 1847 Gauss maintained that, lacking sufficient data, it was better to draw on the general knowledge of the situation.

### 9A.6. More about the Method of Least Squares

1) In spite of Gauss' opinion, his first justification of the MLSq became generally accepted (§ 9A4-10 and Sheynin 1995c, § 3.4), in particular because the observational errors were (and are) approximately normal whereas his mature contribution (1823b) was extremely uninviting; and the work of **Quetelet** (§ 10.5) and **Maxwell** (§ 10.8.5) did much to spread the idea of normality. Examples of deviation from the normal law were however accumulating both in astronomy and in other branches of natural sciences as well as in statistics (Sheynin 1995c, § 3.5; again Quetelet and **Newcomb**, see § 10.8.4). See also § 9C in which I describe the surprising attitude of Bessel who had actually concealed such deviations.

And, independently from that fact, several authors came out against the first substantiation. **Markov** (1899a), who referred to Gauss himself (to his letter to **Bessel**, see my § 9A.4-2), is well known in this respect but his first predecessor was **Ivory** (§ 10.9-1).

The second justification was sometimes denied as well. Thus, **Bienaymé** (1852, p. 37) declared that, unlike Laplace, Gauss had provided considerations rather than proofs (?); see also **Poincaré's** opinion in § 11.2-7.

2) When justifying the MLSq in 1823 in an essentially different way, Gauss called the obtained estimators *most plausible* (*maxime plausibiles*, or, in his preliminary report (1821), *sicherste*, rather than as before, *maxime probabile*, *wahrscheinlichste*. For the case of the normal distribution, these are jointly effective among unbiased regular estimators[5].

3) Mathematicians had not paid due attention to Gauss' work on the MLSq (§§ 9A.2-3 and 13.2-7), and neither did statisticians, see the Epigraph to this book which apparently complements the following passage (Eisenhart 1978, p. 382):

*When **Karl Pearson** and **G. Udny Yule** began to develop the mathematical theory of correlation in the 1890s, they found that much of the mathematical machinery that Gauss devised […] was immediately applicable. […] Gauss' contributions to the method of least squares embody mathematics essential to statistical theory and its applications in almost every field of science today.*

I really think that K. P. and Yule only discovered Gauss at a late stage of their work.

### 9B. Helmert

It was Helmert who mainly completed the development of the classical **Gaussian** theory of errors; furthermore, some of his findings were interesting for mathematical statistics. With good reason Schumann (1917, p. 97) called him "Master of both the *niedere* [surveying and applied] geodesy and higher [triangulation etc., gravimetry, figure of the Earth] geodesy". Until the 1930s, Helmert's treatise (1872) remained the best source for studying the error theory and the adjustment of triangulation.

Indeed, its third, posthumous edition of 1924 carried a few lines signed by a person (H. Hohenner) who explained that, upon having been asked by the publishers, he had stated that the treatise still remained the best of its kind. His opinion, he added, convinced the publishers.

Helmert (1886, pp. 1 and 86) was the first to consider appropriate geodetic lines rather than chains of triangulation, and this innovation, developed by **Krasovsky**, became the essence of the method of adjustment of the Soviet primary triangulation (see Note 19 to Chapter 6 and Sakatow 1950, § 91). Another of his lesser known contributions (Helmert 1868) was a study of various configurations of geodetic systems. Quite in accordance with the not yet existing linear programming, he investigated how to achieve necessary precision with least possible effort, or, to achieve highest possible precision with a given amount of work. Some equations originating in the adjustment of geodetic networks are not linear, not even algebraic; true, they can be linearized (§ 9A.4-9), and perhaps some elements of linear programming could have emerged then, in 1868, but this had not happened. Nevertheless, Helmert noted that it was expedient to leave some angles of a particular geodetic system unmeasured, but this remark was only academic: all angles have always been measured at least for securing a check upon the work as a whole.

I describe now Helmert's stochastic findings.

1) The chi-square distribution. I (1966) noted that **E. Abbe** (1863), see also **M. G. Kendall** (1971), derived it as the distribution of the sum of the squares of normally distributed errors. He wished to offer a test for revealing systematic errors, and he required, in particular, the distribution of the abovementioned function of the errors since it was indeed corrupted by those errors. Exactly his test rather than the distribution obtained was repeatedly described in the geodetic literature whereas **Linnik** (1958/1961, pp. 109 – 113) introduced a modified version of the **Abbe** test.

Helmert (1876b) provided his own derivation of the $\chi^2$ distribution which he first published without justification (1875a). Neither then nor much later (see Item 2) did he mention Abbe. Actually, he continued after **Gauss** (1816), see § 9A.3, by considering observational errors $\varepsilon_1, \varepsilon_2, \ldots, \varepsilon_n$ and the sum of their powers $\varepsilon_i^n$ for the uniform and the [normal] distributions and for an arbitrary distribution as $n \to \infty$. In the last instance, he proved the Gauss formula (4) and then specified it for the abovementioned distributions. He derived the $\chi^2$ distribution by induction beginning with $n = 1$ and 2; **Hald** (1952, pp. 258 – 261) provided a modernized derivation.

2) Much later Helmert (1905) offered a few tests for revealing systematic influences in a series of errors which he wrote down as

$$v_1 \varepsilon_1 + v_2 \varepsilon_2 + \ldots + v_n \varepsilon_n$$

with $v_i = 1$ or $-1$ and $\varepsilon_i > 0$. He issued from the formula

$$P(|\varepsilon - \bar\varepsilon| \le m) \approx 0.68 \qquad (10)$$

where $m$ was the mean square error of $\varepsilon$ (and thus restricted his attention to the normal law): if the inequality in the left side of (10) did not hold, then, as he thought, systematic influences were present.

When deriving his tests, Helmert considered $v_i$, $|v_i|$, runs of signs of the $v_i$ and functions of the errors $\varepsilon_i$ themselves and in this last-mentioned case he provided a somewhat modified version of the Abbe test.

3) The Peters formula (1856) for the mean absolute error. For $n$ normally distributed errors it was

$$\vartheta = |v_i|/\sqrt{n(n-1)}, \quad 1 \le i \le n \tag{11}$$

with $v_i$ being the deviations of the observations from their arithmetic mean. Helmert (1875a) derived this formula anew because Peters had tacitly and mistakenly assumed that these deviations were mutually independent. Passing over to the errors $\varepsilon_i$, Helmert calculated the appropriate integral applying for that purpose the **Dirichlet** discontinuity factor. However, since the normal distribution is stable, it is possible to say now at once (H. A. David 1957) that formula (11) is correct because

$$E|v_i| = \sqrt{n(n-1)}/h$$

where $h$ is the appropriate parameter [measure of precision] of the initial normal distribution and, as it should be, $\vartheta = 1/h\sqrt{\pi}$.

Helmert also attempted to generalize the Peters formula by considering indirect measurements with $k$ unknowns ($k > 1$). He was unable to derive the appropriate formula but proved that a simple replacement of $(n-1)$ in formula (10) by $(n-k)$ resulted in underestimating the absolute error.

4) Helmert (1876b) calculated the variance of the estimator (11). His main difficulty here was the derivation of $E|v_i v_j|$, $i < j$, but he was able to overcome it and obtained

$$\{\pi/2 + \arcsin[1/(n-1)] - n + \sqrt{n(n-2)}\}/\pi n h^2.$$

**Fisher** (1920, p. 761) derived this formula independently.

5) In the same paper Helmert investigated the precision of the **Gauss** formula (6b). For direct measurements it can be replaced by the expression for the mean square error

$$m = \sqrt{\frac{[vv]}{n-1}}.$$

Helmert derived it for the normal distribution by the principle of maximum likelihood, but had not remarked that the estimator obtained (which, however, directly followed from (6a) and was always applied in practice in geodesy) was, unlike the Gauss formula, biased.

Denote the observational errors by $\varepsilon_i$ and their mean by $\bar\varepsilon$, then

$$v_i = \varepsilon_i - \bar\varepsilon$$

and the probability that these errors had occurred, as Helmert indicated in the context of his proof, was equal to

$$P = n(h/\sqrt{\pi})^n \exp[-h^2([vv] + n\bar{v}^2)] \, dv_1 \, dv_2 \ldots dv_{n-1} \, d\bar{v}. \qquad (12)$$

This formula shows that, for the normal distribution, $[vv]$, – and, therefore, the variance as well,– and the arithmetic mean are independent. Helmert had thus proved the important **Student – Fisher** theorem although without paying any attention to it.

A special feature in Helmert's reasoning was that, allowing for (6c), he wrote down the **Gauss** formula 6b) for the case of direct measurements (and, to repeat, for the normal distribution) as[6]

$$m_o^2 = \frac{[vv]}{n-1} [1 \pm \frac{\sqrt{2}}{\sqrt{n-1}}]; \qquad (13)$$

that is, he considered the variance together with its mean square error. Gauss (1816, §§ 6 and 8) sometimes, but not always, acted the same way (although without applying that term, *mean …*, which I found in the mid-19$^{th}$ century in the works of Chebyshev and Russian artillerists).

In addition, Helmert noted that for small values of $n$ the var$m_o^2$ did not estimate the precision of formula (6b) good enough and derived the following formula

$$E[m - \frac{[vv]}{\sqrt{n-1}}]^2 = (1/h^2)\{1 - 2\frac{\Gamma(n/2)}{\Gamma[(n-1)/2]\sqrt{n-1}}\}. \qquad (14)$$

He issued from the probability of the values of $v_i$, $i = 1, 2, \ldots, (n-1)$,

$$P = \sqrt{n}(h/\sqrt{\pi})^{n-1} \exp(-h^2[vv]) \, dv_1 \, dv_2 \ldots dv_{n-1}$$

that follows from formula (12), noted that the probability $P(\sigma \leq [vv] \leq \sigma + d\sigma)$ was equal to the appropriate integral, and introduced new variables

$$
\begin{aligned}
t_1 &= \sqrt{2}(v_1 + 1/2 v_2 + 1/2 v_3 + 1/2 v_4 + \ldots + 1/2 v_{n-1}), \\
t_2 &= \sqrt{3/2}\,(v_2 + 1/3 v_3 + 1/3 v_4 + \ldots + 1/3 v_{n-1}), \\
t_3 &= \sqrt{4/3}\,(v_3 + 1/4 v_4 + \ldots + 1/4 v_{n-1}), \ldots, \\
t_{n-1} &= \sqrt{n/(n-1)}\, v_{n-1}.
\end{aligned}
$$

Note that $[vv] = [tt]$ where, however, the first sum consisted of $n$ terms and the second one, of $(n-1)$ terms, and the **Jacobian** of the transformation was $\sqrt{n}$. The derivation of formula (14) now followed immediately since Helmert knew the $\chi^2$ distribution. Taken together, the transformations from $\{\xi\}$ to $\{v\}$ and from $\{v\}$ to $\{t\}$ are called after him.

**Kruskal** (1946) transformed formula (12) by introducing a bivariate "Helmert distribution" with variables

$$s = \sqrt{[vv]/n}, \; u = x - \mu,$$

where $x$ was the arithmetic mean of $n$ normally distributed observations $N(\mu; \;)$, and replaced $h$ by . He mentioned several authors who had derived that new distribution by different methods, determined it himself by induction and indicated that the **Student** distribution followed from it, see **Hald** (1998, p. 424).

Finally, Helmert corrected the boundaries of the estimator (6b). As indicated by **Gauss** they were

$$2(\nu_4 - 2s^4)/(n-k); \; [1/(n-k)](\nu_4 - s^4) + (k/n)(3s^4 - \nu_4)$$

where $\nu_4$ was the fourth moment of the errors and $s^2 = Em^2$. Helmert had discovered that the lower boundary was wrong and **Kolmogorov** et al (1947) independently repeated his finding. Here is the final result; Maltzev (1947) proved that the lower bound was attainable. For non-negative and then non-positive $(v_4 - 3s^4)$, the product $(n-k)\text{var } m_o^2$ as it occurred, was contained within, respectively, boundaries

$$[(\nu_4 - s^4) - (k/n)(\nu_4 - 3s^4); \; (\nu_4 - s^4)], \qquad (15a)$$
$$[(\nu_4 - s^4); \; (\nu_4 - s^4) + (k/n)(3s^4 - \nu_4)]. \qquad (15b)$$

### 9C. Bessel

Friedrich Wilhelm Bessel (1784 – 1846) was an outstanding astronomer and an eminent mathematician. He and Gauss were the originators of a new direction in practical astronomy and geodesy which demanded a thorough examination of the instruments and investigation of the plausibility of observational methods. Each instrument was *accused* of every possible defect and was only exonerated after proving itself irreproachable. See an appropriate quotation from Newcomb in Schmeidler (1984, pp. 32 – 34).

I mentioned Bessel in §§ 6.3 and 8.7. His achievements in astronomy and geodesy are well known; I name the determination of astronomical constants; the first determination of a star's parallax (and thus the definitive establishment of the heliocentric system); the discovery of the personal equation; the development of a method of adjusting triangulation (see however below); design and use of metallic bars for measuring baselines, and the derivation of the parameters of the Earth's ellipsoid of rotation which enjoyed international recognition (Strasser 1957, p. 39).

The personal equation is the systematic difference of the moments of the passage of a star through the cross-hairs of an astronomical instrument as recorded by two observers. When studying this phenomenon, it is necessary to compare the moments fixed by two astronomers at different times and, consequently, to take into account the correction of the clock (if they use only one). Bessel (1823) had indeed acted appropriately, but in one case he failed to do so, and his pertinent observations proved useless. He made no such comment;

Bessel (1838a, §§ 1 and 2) determined the densities of two functions of a continuously and uniformly distributed [random variable], and, unlike **Laplace**, he clearly formulated this problem. Nevertheless, he erred in his computations of the pertinent variances and probable errors. He also determined the density of the total observational error made up of many heterogeneous components but a rigorous solution of such problems became only possible mush later (§ 13.1-4)[7].

I have discovered 33 mistakes in arithmetic and elementary algebra (except those noticed by the Editor) in his *Abhandlungen* (1876). They did not influence his conclusions but they throw doubt on his more serious calculations. Here is just one of them (1876, Bd. 2, p. 376):

 4:  5 = 1/1.409; actually, however, 1/1.118.

One more example, this time concerning Bessel's reasoning (1818; 1838a). He presented three series of Bradley's observations, 300, 300 and 470 in number, and stated that their errors almost precisely obeyed normal distributions. Actually, he was wrong and it is difficult to believe that he was mistaken (especially see below). Moreover, he thus missed the opportunity to discover an example of long series not quite normally distributed errors of precise observations. Later, scientists gradually discovered such series, in the first place see Newcomb (1886).

Bessel's contribution (1838a) included a proof of a version of the CLT (rigorously proved only by Liapunov and Markov). Bessel stated that, given more observations, the deviation from normality will disappear. Did not he notice that he thus undermined the essence of that theorem? Did not he formulate his wrong conclusion to save that proof?

It became customary to measure each angle of a chain of triangulation an equal number of times and, which was more important, to secure their mutual independence so as to facilitate the treatment of the observations, – to separate the station adjustment from the adjustment of the chain as a whole. Bessel, however, did not keep to the abovementioned condition (and had to adjust all the observations at once). There are indications that the actual rejection of his method annoyed him[8].

I have since discovered other examples of Bessel's misleading statements in his popular writings. True, at least one of them pertains to the time of his fatal illness, but I venture to suppose that a very ill person should all the more try to avoid mistakes.

1. Bessel (1843). This is his report of the same year read out to the physical section of the Königsberg physical-economic society in which he had been very active. Schumacher published the texts of these reports (1848b), and Bessel (1848a), about which I say a few words below, is included in that collection.

And so, Bessel (1843) described the life and work of William Herschel. Among other things, he properly discussed Herschel's hunt for double stars and his attempt at counting the stars in the Milky Way, but he did not explain that there are two types of double stars nor did he say that the Milky Way is only one of the countless galaxies.

Herschel came to understand that his telescope did not penetrate to the boundaries of the Milky Way (F. G. W. Struve 1847, p. 34; Hoskin 1959) whereas Bessel (p. 474, left column) stated quite the opposite. Another mistake concerned the discovery of the planet Uranus. Contrary to Bessel's statement (p. 469, left column), Herschel discovered a moving body and decided that it was a comet. Finally, Bessel (p. 470, right column) mentioned Caroline, the sister of William, and remarked that she was still alive and assisted her brother. Actually, Caroline died several decades later than he.

2. Bessel (1845). This is a newspaper article which had nothing to do with astronomy. Bessel stated that under such parameters as territory, climate etc. (political system not mentioned) only mental development of the population determined its acceptable maximal number. However, a territory becomes more or less populated when people turn from hunting to farming (Bessel's own example), but are farmers more mentally developed than hunters?

Then Bessel turned his attention to the United States and provided his own data about the population of Native Americans taken out of thin air and damnably wrong.

3. Bessel (1848a). The date of the report is unknown. Bessel mentioned Delambre's *Astronomy* which was not quite definite, but sufficient for stating that the report was read in 1821 or later.

The significance of Jakob Bernoulli's law of large numbers was not discussed, Lambert's preference of maximum-likelihood estimators over the arithmetic mean (p. 401) was mostly imagined and Laplace's *Essai philosophique* of 1816 was not even mentioned. Population statistics studied, for example, by De Moivre, Nicolaus and Daniel Bernoulli, was completely left out. It is difficult to conclude that Bessel's quite elementary exposition had satisfied his listeners.

In his correspondence, Gauss several times indicated Bessel's shortcomings.

1. Gauss – Olbers, 2 Aug. 1817. Bessel had overestimated the precision of some of his measurements. On 2 Nov. 1817 Olbers *confidentially* informed Bessel about Gauss' opinion.

2. Gauss – Schumacher, between 14 July and 8 Sept. 1826. He stated the same about Bessel's investigation of the precision of the graduation of a limb.

3. Gauss – Schumacher 27 Dec. 1846. He negatively described some of Bessel's posthumous manuscripts. In one case he was *shocked* by Bessel's *carelessness*[8].

I am at a loss: how was it possible to pass these statements over? And, again, how was it possible for Bessel to be at once a great scholar and a happy-go-lucky scribbler? Cf. Goethe (*Faust*, pt. 1, Sc. 2): *Two souls are living in his breast*.

### Notes

**1.** This term should only be applied to the method as substantiated by Gauss in 1823; until then, strictly speaking, the *principle* of least squares ought to be thought of.

**2.** Adrain included his work in a periodical published by himself for the year 1808; however, its pertinent issue appeared only in 1809 (Hogan 1977). Adrain's library included a copy of **Legendre's** memoir (Coolidge 1926) in which, however, the normal distribution was lacking; furthermore, it is unknown when had Adrain obtained that memoir. The term *normal distribution* appeared in 1873 (**Kruskal** 1978) and was definitively introduced by **K. Pearson** (1894).

**3.** Their opinion should not be forgotten. Here is another example. **Encke** (1851, p. 2) believed that Gauss had applied the MLSq when determining the orbit of Ceres, the first-discovered minor planet (Gauss did not comment). In Note 19 to Chapter 6 I mentioned an inadmissible free and easy manner adopted by a certain author (Stigler 1986) with respect to **Euler**. His attitude towards Gauss was not better. Here are his statements: **Legendre** "immediately realized the method's potential" (p. 57), but "there is no indication that [Gauss] saw its great potential before he learned of Legendre's work" (p. 146); then (p. 143), only **Laplace** saved Gauss's argument [his first justification of the MLSq] from joining "an accumulating pile of essentially ad hoc constructions"; and, finally (p. 145), Gauss "solicited reluctant testimony from friends that he had told them of the method before 1805". I (Sheynin 1999a; 1999c) had refuted these astonishing declarations (see also the very end of § 9A.1 about Olbers) which Stigler (1999), the first ever slanderer of the memory of the great man, repeated slightly less impudently. Regrettably, no one supported me; on the contrary, Stigler's first book met with universal approval although he, in addition, left aside the ancient history as well as such scholars as **Kepler**, **Lambert** and **Helmert**. **Hald** (1998, p. xvi), whose outstanding contribution deserves highest respect, called Stigler's book "epochal". I am unable to understand suchlike opinions. The attitude of the scientific comunity towards Stigler proves that it is seriously ill.

**4.** It is not amiss to add that the primary triangulation of the Soviet Union consisted of chains independent one from another in the Gauss' sense. This, together with other conditions, enabled the geodesists to estimate realistically the precision of the whole great net (Sakatow 1950/1957, pp. 438 – 440). And in general, geodesists, not necessarily mentioning Gauss, were keeping to his opinion. I also note that **Kapteyn** (1912), who had not cited Gauss and was unsatisfied with the then originating correlation theory, proposed to estimate quantitatively the dependence between series or functions of observations by issuing from the same notion of independence, see Sheynin (1984a, § 9.2.1). His article went unnoticed.

**5.** Concerning this rarely mentioned concept see **Cramér** (1946, § 32.6).

**6.** In the theory of errors, the application of the mean square error with a double sign became standard (§ 9B-4).

**7.** In 1839 Gauss informed **Bessel** (*W*-8, pp. 146 – 147) that he had read the latter's memoir with interest although the essence of the problem had been known to him for many years.

**8.** In 1825, Gauss had a quarrel with **Bessel** but no details are known (Sheynin 2001d, p. 168). Even in 1817 **Olbers** (Erman 1852, Bd. 2, p. 69) regretted that the relations between Bessel and Gauss were bad. In 1812, in a letter to Olbers, Bessel (Ibidem, Bd. 1, p. 345) had called Gauss "nevertheless" the inventor of the MLSq, but in 1844, in a letter to **Humboldt** (Sheynin 2001d, p. 168), he stressed **Legendre's** priority.

## 10. The Second Half of the 19th Century

Here, I consider the work of several scholars (§§ 10.1 – 10.6), statistics (§ 10.7), and its application to various branches of natural sciences (§ 10.8). The findings of some natural scientists are discussed in § 10.9 since it proved difficult to describe them elsewhere.

### 10.1. Cauchy

Cauchy published not less than 10 memoirs devoted to the treatment of observations and the theory of probability. Eight of them (including those of 1853 mentioned below) were reprinted in t. 12, sér. 1, of his *Oeuvres complètes* (1900). In particular, he studied the solution of systems of equations by the principle of minimax (§ 6.3.2) and proved the theorem in linear programming known to **Gauss** (§ 9A.2-1). He had also applied the method of averages (§ 6.3.2) and **Linnik** (1958/1961, § 14.5), who cited his student L. S. Bartenieva, found out that the pertinent estimators were unbiased and calculated their effectiveness for the cases of one and two unknown(s). I briefly describe some of Cauchy's findings.

Cauchy (1853b) derived the even density of observational errors demanding that the probability for the error of one of the unknowns, included in equations of the type of (1.2), to remain within a given interval, was maximal. Or, rather, he derived the appropriate characteristic function

$$\varphi(\ ) = \exp(-c\ ^{\mu+1}), \quad c,\ > 0, \mu \text{ real} \qquad (1)$$

and noted that the cases $\mu = 1$ and $0$ led to the [normal law] and to the "Cauchy distribution", see § 8.6. The function (1) is characteristic only when $1 < \mu \leq 1$ and the appropriate distributions are [stable].

In two memoirs Cauchy (1853c; 1853d) proved the [CLT] for the linear function

$$A = [m\ ] \qquad (2)$$

of [independent] errors $\varepsilon_i$ having an even density on a finite interval. In both cases he introduced characteristic functions of the errors and of the function (2), obtained for the latter

$$\varphi(\ ) = \exp(-s\ ^2)$$

where $2s$ was close to $\sigma^2$, the variance of (2), and, finally, arrived at

$$P(|\ | \leq \ ) \approx \frac{\sqrt{2}}{\sqrt{\ }} \int_0^{\ } \exp(-x^2/2\sigma^2) dx.$$

It is important that he had also estimated the errors due to assumptions made and **Freudenthal** (1971, p. 142) even declared that his proof was rigorous by modern standards; see, however Heyde & Seneta (1977, pp. 95 – 96).

Cauchy devoted much thought to interpolation of functions, and, in this connection, to the MLSq, but, like **Poisson**, he never cited **Gauss**.

In one case he (1853a/1900, pp. 78 – 79) indicated that the MLSq provided most probable results only in accordance with the **Laplacean** approach [that is, only for the normal distribution] and apparently considered this fact as an essential shortcoming of the method.

## 10.2. Bienaymé

Heyde & Seneta (1977) described his main findings; I follow their account and abbreviate their work as HS. Bru et al (1997) published two of Bienaymé's manuscripts and other relevant archival sources.

1) A limit theorem (Bienaymé 1838; HS, pp. 98 – 103). Bienaymé had "essentially" proved the theorem rigorously substantiated by **Mises** (1919; 1964a). The abovementioned adverb appeared in Mises (1964b, p. 352). Suppose that $n$ trials are made with some event $A_i$ from $m$ mutually exclusive events ($i = 1, 2, …, m$) occurring in each trial with probability $p_i$ and that $x_i$ is the number of times that $A_i$ happened, $\sum x_i = n$. Treating the probabilities $p_i$ as random variables, Bienaymé studied the distribution of their linear function in the limiting case $x_i, n \to \infty$, $x_i/n = C_i$. As a preliminary, he had to derive the posterior distribution of the $p_i$ given $x_i$ tacitly assuming that the first ($m - 1$) of these probabilities were random variables with a uniform prior distribution. Actually, Bienaymé proved that the assumption about the prior distributions becomes insignificant as the number of the multinomial trials increases.

Note that **Nekrasov** (1890) had forestalled **Czuber**, whom **Mises** named as his predecessor. Assuming some natural restrictions, he proved a similar proposition concerning the **Bernoulli** trials.

2) The **Liapunov** inequalities (Bienaymé 1840b; HS, pp. 111 – 112). Without proof, Bienaymé[1] indicated that the absolute initial moments of a discrete [random variable] obeyed inequalities which could be written as

$$(E|\xi|^m)^{1/m} \geq (E|\xi|^n)^{1/n}, \quad 0 \leq m \leq n.$$

Much later Liapunov (1901a, § 1) proved that

$$(E|\xi|^m)^{s-n} < (E|\xi|^n)^{s-m} < E(|\xi|^s)^{m-n}, \quad s > m > n \geq 0.$$

He applied these inequalities when proving the [CLT].

3) The law of large numbers. Bienaymé (1839) noted that the fluctuation of the mean statistical indicators was often greater than in accordance with the **Bernoulli** law, and suggested a possible reason: some causes acting on the studied events, as he thought, remained constant within a given series of trials but essentially changed from one series to the next one. **Cournot**, **Lexis** and other "Continental" statisticians took up this idea without citing Bienaymé (Chapter 15) but it was also known in the theory of errors since systematic errors can behave in a similar way. Bienaymé, in addition, somehow interpreted the Bernoulli theorem as an attempt to study suchlike patterns of the action of causes. He (1855/1876) repeated this statement and, on p. 202, he mistakenly reduced the **Poisson** LLN to the case of variable probabilities whose mean value simply replaced the constant probability of the **Bernoulli** trials, also see HS, § 3.3.

4) The Bienaymé – **Chebyshev** inequality (Bienaymé 1853; HS, pp. 121 – 124; **Gnedenko** & Sheynin 1978/2001, pp. 258 – 262). This is the name of the celebrated inequality

$$P(|\xi - E\xi| < \varepsilon) > 1 - \text{var}\,\xi/\varepsilon^2, \quad \varepsilon > 0. \qquad (3)$$

Differing opinions were pronounced with regard to its name and to the related method of moments. **Markov** touched on this issue four times. In 1912, in the Introduction to the German edition of his *Treatise* (1900a/1908), he mentioned "the remarkable **Bienaymé – Chebyshev** method". At about the same time he (1912b, p. 218) argued that

*Nekrasov's statement* [that Bienaymé's idea was exhausted in Chebyshev's works] *is refuted by indicating a number of my papers which contain the extension of Bienaymé's method* [to the study of dependent random variables].

Then, Markov (1914b/1981, p. 162) added that the "starting point" of Chebyshev's second proof of **Poisson's** LLN "had been […] indicated by […] Bienaymé" and that in 1874 **Chebyshev** himself called this proof "a consequence of the new method that Bienaymé gave". Nevertheless, Markov considered it "more correct" to call the method of moments after both Bienaymé and Chebyshev, and "sometimes" only after the latter, since "it only acquires significance through Chebyshev's work" [especially through his work on the CLT]. Finally, Markov (*Treatise*, 1924, p. 92) stated that Bienaymé had indicated the main idea of the proof of the inequality (3), although restricted by some conditions, whereas Chebyshev was the first to formulate it clearly and to justify it.

Bienaymé (1853/1867, pp. 171 – 172) considered a random sum, apparently (conforming to the text of his memoir as a whole) consisting of identically distributed terms, rather than an arbitrary magnitude $\xi$, as in formula (3). This is what Markov possibly thought of when he mentioned some conditions. HS, pp. 122 – 123, regarded his proof, unlike Chebyshev's substantiation [§ 13.1-3], "short, simple, and […] frequently used in modern courses …" Yes, **Hald** (1998, p. 510) repeated it in a few lines and then got rid of the sum by assuming that it contained only one term. **Gnedenko** (1954/1973, p. 198) offered roughly the same proof but without citing Bienaymé.

Bienaymé hardly thought that his inequality was important (Gnedenko & Sheynin 1978/2001, p. 262; Seneta 1998, p. 296). His main goal was to prove that only the variance was an acceptable estimator of precision in the theory of errors (see § 9A.4-2) and, accordingly, he compared it with the fourth moment of the sums of random [and independent] errors. Consequently, and the more so since he never used integrals directly, I believe that Chebyshev (1874), see also Gnedenko & Sheynin 1978/2001, p. 262) overestimated the part of his predecessor in the creation of the method of moments. Here are his words:

*The celebrated scientist presented a method that deserves special attention. It consists in determining the limiting value of the integral [...] given the values of the integrals...*

The integrand in the first integral mentioned by Chebyshev was $f(x)$ and the limits of integration were $[0; a]$; in the other integrals, $xf(x)$, $x^2 f(x)$, ... and the limits of integration, $[0; A]$, $f(x) > 0$ and $A > a$.

5) Runs up and down (Bienaymé 1874; 1875; HS, pp. 124 – 128). Suppose that $n$ observations of a continuous random variable are given. Without proof Bienaymé indicated that the number of intervals between the points of extremum (almost equal to the number of these points) is distributed approximately normally with parameters

mean ...$(2n – 1)/3$, variance ... $(16n – 29)/90$.  (4)

He maintained that he knew this already 15 or 20 years previously. HS states that these findings were discovered anew; nevertheless, the authors derive formulas (4), the first of them by following **Bertrand**, also see Moore (1978, p. 659). Bienaymé checked the agreement between several series of observations and his findings. Some of the data did not conform to his theory and he concluded that that happened owing to unrevealed systematic errors. I return to his test in § 10.3.

6) The method of least squares (Bienaymé 1852; HS, pp. 66 – 71). Bienaymé correctly remarked that least variance for each estimator separately was not as important as the minimal simultaneous confidence interval for all the estimators. Keeping to the **Laplacean** approach to the MLSq (see his remark in my § 9A.6-1), he restricted his attention to the case of a large number of observations. Bienaymé also assumed that the distribution of the observational errors was known and made use of its first moments and even introduced the first four cumulants and the multivariate **Gram – Charlier** series (Bru 1991, p. 13; **Hald** 2002, pp. 8 – 9). He solved his problem by applying the principle of maximum likelihood, introducing the characteristic function of the [vector of the] errors and making use of the inversion formula. True, he restricted his choice of the [confidence] region; on the other hand, he derived here the $\chi^2$ distribution. Bienaymé's findings were interesting indeed, but they had no direct bearing on the theory of errors. Furthermore, his statement (pp. 68 – 69) that both the absolute expectation and variance were unreliable estimators of precision was certainly caused by his adoption of the method of maximum likelihood and is nowadays forgotten.

7) A branching process (Bienaymé 1845; HS, pp. 117 – 120). Bienaymé had formulated the properties of criticality of a branching process while examining the same problem of the extinction of noble families that became attributed to **Galton**. **D. G. Kendall** (1975) reconstructed Bienaymé's proof and reprinted his note and Bru (1991) quoted a passage from **Cournot's** contribution of 1847 who had solved a stochastic problem in an elementary algebraic way and had indicated that it was tantamount to determining the probability of the duration of the male posterity of a family,– to a problem in which "Bienaymé is engaged". Bru thought it highly probable that Cournot had borrowed his study from Bienaymé.

8) An approach to the notion of sufficient estimator (Bienaymé 1840a; HS, pp. 108 – 110). When investigating the stability of statistical frequencies (see also Item 3), Bienaymé expressed ideas that underlie the notion of sufficient estimators. For $m$ and $n$ **Bernoulli** trials with probability of success $p$ the number of successes has probability

$$P(\mu_i = k) = C_s^k p^k (1-p)^{s-k}$$

with $s = m$ and $s = n$ respectively and $i = 1$ and 2 denoting these series of trials. It is easy to ascertain that the probability $P(\mu_1 = r, \mu_2 = a - r | \mu_1 + \mu_2 = a)$ does not depend on $p$. Bienaymé thought that this property, that takes place when the totality of the trials is separated into series, could prove the constancy of the laws of nature. However, statisticians (**Fourier**, whom he mentioned; **Quetelet** 1846, p. 199 ) pragmatically considered such a separation as the best method for revealing variable causes. HS additionally noted that Bienaymé should have understood that all the information about the unknown probability $p$ [if it were constant] was provided by the totality of the trials, that Bienaymé had wrongly calculated the variance of the hypergeometric distribution, and that he made use of a particular case of the CLT for checking a null hypothesis.

### 10.3. Cournot

Cournot intended his main contribution (1843) for a broader circle of readers. However, lacking a good style and almost completely declining the use of formulas, he hardly achieved his goal. Recall also (the end of § 8.7) that Cournot passed over in silence the LLN. I describe his work as a whole; when referring to his main book, I mention only the appropriate sections.

1) The aim of the theory of probability. According to Cournot (1875, p. 181), it was "The creation of methods for assigning quantitative values to probabilities". He thus moved away from **Laplace** (§ 7.3) who had seen the theory as a means for revealing the laws of nature. Cf. **Chebyshev's** opinion in § 13.2-1.

2) The probability of an event (§ 18): this is the ratio of the extent (*étendue*) of the favourable chances to the complete extent of all the chances[2]. The modern definition replaced "extent" by a clear mathematical term, "measure". I stress that Cournot's definition included geometric probability, which until him had been lacking any formula, and thus combined it with the classical case. Cournot (§ 113) also introduced probabilities unyielding to measurement and (§§ 233 and 240-8) called them *philosophical*. They might be related to expert estimates whose treatment is now included in the province of mathematical statistics. Cf. Fries in § 7.1-5.

3) The term *médiane*. This is due to Cournot (§ 34).

4) The notion of randomness. In a book devoted to exonerating games of chance La Placette (1714) explained (not clearly enough) randomness as an intersection of independent chains of determinate events, thus repeating the statements of many ancient scholars (see § 1.1.1). Cournot (§ 40) expressed the same idea, and, in § 43, indirectly connected randomness with unstable equilibrium by remarking that a right circular cone, when stood on its vertex, fells in a "random" direction. This was a step towards **Poincaré's** viewpoint (§ 11.2-9). Cournot (1851, § 33, Note 38; 1861, § 61, pp. 65 – 66) also recalled **Lambert's** attempt to study randomness (see my § 6.1.3)[3], and (1875, pp. 177 – 179) applied **Bienaymé's** test (§ 10.2-3) for investigating whether the digits of the number were random. He replaced its first 36 digits by signs plus and minus (for example, 3; 1; 4; 1 became – ; +; –) and counted 21 changes in the new sequence. Comparing the fractions 21/36 = 0.583 and [(2 36 + 1)/3 36] = 0.667, see the first of the formulas (3), Cournot decided that the accordance was good enough (?), but reasonably abstained from a final conclusion.

5) A mixture of distributions. Given the densities of separate groups of $n_i$, $i = 1, 2, …, m,$ observations, Cournot (§ 81) proposed the weighted mean density as their distribution. He had not specified the differences between the densities, but in § 132 he indicated that they might describe observations of different precision and in § 135 he added, in a Note, that observational errors approximately followed the [normal law]. I describe the attempts to modify the normal law made by astronomers in § 10.8.4.

6) Dependence between the decisions of judges and/or jurors. Cournot (1838; 1843, §§ 193 – 196 and 206 – 225) gave thought to this issue. Suppose (§ 193) that in case of two judges the probabilities of a correct verdict are $s_1$ and $s_2$. Then the probability that they agree is

$$p = s_1 s_2 + (1 - s_1)(1 - s_2) \qquad (5)$$

so that, if $s_1 = s_2 = s > 1/2$, $s = \dfrac{1}{2} + \dfrac{1}{2\sqrt{2p-1}}$.

Statistical data provided the value of $p$; it should have obviously exceeded 1/2. If the data made it also possible to ascertain $s_1$ and $s_2$, and equations of the type of (5) will not be satisfied, it will be necessary to conclude that the verdicts were not independent. Cournot's study was hardly successful in the practical sense, but it at least testified to an attempt to investigate dependence between some events.

7) A critical attitude towards statistics; a description of its aims and applications. Cournot (§ 103) declared that statistics had *blossomed exuberantly* and that [the society] should be on guard against its *premature and wrong* applications which might discredit it for some time and delay the time when it will underpin all the theories concerning the "organization sociale". Statistics, he (§ 105) continued, should have its theory, rules, and principles, it ought to be applied to natural sciences, to physical, social and political phenomena; its main goal was (§ 106) to ascertain "the knowledge of the essence of things", to study the causes governing the phenomena of the physical world and social life (§ 120). The theory of probability was applicable to statistics (§ 113) and the "principe de **Bernoulli**" was its only pertinent sound foundation (§ 115).

These statements were not at all unquestionable (§§ 6.2.1 and 9B). Cournot, however, went further: the theory of probability might be successfully applied in astronomy, and the "statistique des astres, if such an association of words be permitted, will become a model for every other statistics" (§ 145). He himself, however, statistically studied the parameters of planetary and cometary orbits, but not the starry heaven, and his statement was inaccurate: first (§ 10.8), by that time statistics had begun to be applied in a number of branches of natural sciences; second (Sheynin 1984a, § 6), stellar statistics had by then already originated.

8) Explanation of known notions and issues. Cournot methodically explained the notion of density (§§ 64 – 65) and the method of calculating the density of a function of a (of two) [random variable(s)] (§§ 73 – 74). He also described how is or should statistics be applied in natural sciences and demography, discussed the treatment of data when the probabilities of the studied events were variable, etc.

Taken as a whole, Cournot made a serious contribution to theoretical statistics. **Chuprov** (1905/1960, p. 60), bearing in mind mathematics as well as philosophy and economics, called him a man of genius. Later he (1909/1959, p. 30) stated that Cournot was

*One of the most profound thinkers of the 19<sup>th</sup> century, whom his contemporaries failed to appreciate, and who rates ever higher in the eyes of posterity*.

Lastly, Chuprov (1925a/1926, p. 227) characterized the French scientist as "the real founder of the modern philosophy of statistics". All this seems to be somewhat exaggerated and in any case I do not agree with Chuprov's opinion about Cournot's "real substantiation" and "canonical" proof of the LLN (1905/1960, p. 60; 1909/1959, pp. 166 – 168). Up to 1910, when he began corresponding with **Markov**, Chuprov was rather far from mathematical statistics. He had not remarked that Cournot did not even formulate that law, and that his "Lemma", as Chuprov called it [rare events do not happen, see Cournot (1843, § 43) interpreted by Chuprov as "do not happen often"], was not new at all, see my §§ 2.1.2, 2.2.2 and 3.2.2 concerning moral certainty and § 6.1.2 with regard to **D'Alembert** who formulated the same proposition.

**Cournot** obviously never had anything in common with precise measurements (observations) and his considerations about them are hardly useful. He had not, but should have known that in 1817 Humboldt introduced isotherms (§ 10.8.3). He had not mentioned Daniel Bernoulli's study of 1766 of smallpox epidemics and his description of the tontines was wrong.

### 10.4. Buniakovsky

Several European mathematicians had attempted to explicate the theory of probability simpler than **Laplace** did. **Lacroix** (1816), the mathematical level of whose book was not high, **Cournot** (§ 10.3), whose main contribution was translated into German in 1849, and **De Morgan** (1845) might be named here; concerning the last-metioned author see however Note 1 to Chapter 8. In Russia, **Buniakovsky** (1846) achieved the same aim; his treatise was the first comprehensive Russian contribution so that **P. B. Struve** (1918, p. 1318) called him "a Russian student of the French mathematical school". I discuss the main issues considered by him both in his main treatise and elsewhere, see also Sheynin (1991b), this being a general description of Buniakovsky's work. A practically complete list of his contributions is in *Materialy* (1917).

1) The theory of probability. In accordance with its state in those times, Buniakovsky (1846, p. I) attributed it to applied mathematics. He (Ibidem) also maintained that

*The analysis of probabilities considers and quantitatively estimates even such phenomena […] which, due to our ignorance, are not subject to any suppositions.*

This mistaken statement remained, however, useless: Buniakovsky never attempted to apply it; furthermore, he (p. 364; 1866a, p. 24) went back on his opinion.

2) Moral expectation (see § 6.1.1). Independently from **Laplace**, Buniakovsky (1846, pp. 103 – 122) proved **Daniel Bernoulli's** conclusion that an equal distribution of a cargo on two ships increased the moral expectation of the freightowner's capital as compared with transportation on a single ship. Later he (1880) considered the case of unequal probabilities of the loss of each ship. Then, he (1866a, p. 154) mentioned moral expectation when stating that the statistical studies of the productive population and the children should be separated, and he concluded with a general remark:

*Anyone, who does not examine the meaning of the numbers, with which he performs particular calculations, is not a mathematician.*

For a long time statisticians shunned mathematics (end of § 10.7) since they wrongly understood the essence of that science. And even in the first half of the 20th century Soviet statisticians shunned it since they correctly understood the ensuing danger to Marxism (Note 7 to Chapter 15).

Buniakovsky passed over in silence **Ostrogradsky's** attempt to generalize the concept of moral expectation (§ 7.1-9).

3) Geometric probabilities (§ 6.1.4). Buniakovsky (1846, pp. 137 – 143) generalized the **Buffon** problem by considering the fall of the needle on a system of congruent equilateral triangles. His geometric reasoning was, however, complicated and his final answer, as **Markov** (*Treatise*, 1900/1924, p. 270) maintained, was wrong. Markov himself had been solving an even more generalized problem concerning a system of congruent scalene triangles, but his own graph was no less involved and, as it seems, no one has checked his solution.

Buniakovsky also investigated similar problems earlier (1837) and remarked then that their solution, together with [statistical simulation], might help to determine the values of special transcendental functions. In the same connection, **Laplace** only mentioned the number .

4) "Numerical" probabilities. Buniakovsky (1836; 1846, pp. 132 – 137) solved an elementary problem on the probability that a quadratic equation with coefficients "randomly" taking different integral values had real roots. Much more interesting are similar problems of later origin, for example on the reducibility of fractions (**Chebyshev**, see § 13.2-8) or those concerning the set of real numbers.

5) A random walk. Buniakovsky (1846, pp. 143 – 147) calculated the probability that a castle, standing on square $A$ of a chessboard, reached square $B$ (possibly coinciding with $A$) in exactly $x$ moves if its movement was "uniformly" random. Before that, random walks (here, it was a generalized random walk) had occurred only indirectly, when studying a series of games of chance.

Buniakovsky's problem was, however, elementary. The castle can only be in two states,– it can reach $B$ either in one move, or in two moves; the case $A \quad B$ belongs to the latter, but might be isolated for the sake of expediency. Buniakovsky divided the squares in three groups: square A itself; 14 squares lying within reach of the very first move from A; and the rest 49 squares and he formed and solved a system of three pertinent difference equations for the number of cases leading to success. It turned out that the mean probability (of its three possible values) was equal to 1/64 and did not depend on $x$. He had not interpreted his result, but properly indicated that it was also possible to solve the problem in an elementary way, by direct calculation. Note that the first $n$ moves ($n \quad 1$), if unsuccessful, do not change anything, and this circumstance apparently explains the situation.

6) Statistical control of quality. Buniakovsky (1846, Addendum; 1850) proposed to estimate probable military losses (actually, mean losses, cf. beginning of § 11.2) in battle by sample data in each arm of the engaged forces, – to use stratified sampling, as it is now called. His study was hardly useful, the more so since he applied the **Bayesian** approach assuming an equal prior probability of all possible losses, but he (1846, pp. 468 – 469) also indicated that his findings might facilitate the acceptance "of a very large number of articles and supplies" only a fraction of which was actually examined.

Statistical control of quality was then still unknown although even **Huygens** (§ 2.2.2) had solved a pertinent urn problem. **Ostrogradsky** (1848), possibly following Buniakovsky[4], picked up the same issue. He stated, on p. 322, that the known solutions [of this problem] "sont peu exactes et peu conformes aux principes de l'analyse des hasards". He did not elaborate, and his own main formula (p. 342) was extremely involved (and hardly checked by anyone since).

Two other authors had preceded Buniakovsky who considered an equivalent problem: **Simpson** (1740, Problem 6) and Öttinger (1843, p. 231). Here is the former:

*There is a given Number of each of several Sorts of Things* […] *as (a) of the first Sort, (b) of the second, etc. put promiscuously together; out of which a given Number (m) is to be taken, as it happens; To find the probability that there shall come out precisely a given Number of each Sort…*

7) The history of the theory of probability. Buniakovsky was one of the first (again, after **Laplace**) to consider this subject and a few of his factual mistakes might well be overlooked. His predecessors were Montucla (1802) and Lubbok & Drinkwater (1830). In his popular writings Buniakovsky showed interest in history of mathematics in general and in this field he possibly influenced to some extent both **Markov**, see Sheynin (1989 , § 3), and the eminent historian of mathematics, **V. V. Bobynin** (1849 – 1919).

8) Population statistics. Buniakovsky (1846, pp. 173 – 213) described various methods of compiling mortality tables, studied the statistical effect of a weakening or disappearance of some cause of death (cf. § 6.2.3), calculated the mean and the probable durations of marriages and associations and, following Laplace, solved several other problems.

After 1846, Buniakovsky actively continued these investigations. He compiled mortality tables for Russia's Orthodox believers and tables of their distribution by age (1866a; 1866b; 1874) and estimated the number of Russian conscripts ten years in advance (1875b). No one ever verified his forecast and the comments upon his tables considerably varied. **Bortkiewicz** (1889; 1898b) sharply criticized them, whereas **Davidov** (Ondar 1971), who, in 1886, published his own study of mortality in Russia, noted a serious methodical mistake in their compilation but expressed an opposite opinion. Finally, Novoselsky (1916, pp. 54 – 55) mainly repeated Davidov's criticism, but indicated that Buniakovsy's data were inaccurate and incomplete (as Buniakovsky himself had repeatedly stressed) and called his tables "a great step forward":

*A new period in the study of mortality in Russia started with* […] *the demographic investigations* […] *made by Buniakovsky* […] *and especially* […] *with his <u>Essay</u> (1866a). His contributions on population statistics represent an outstanding and remarkable phenomenon not only in our extremely poor demographic literature but in the rich realm of foreign writings as well, and particularly of his time.* […] *Due to the clearness, depth and nicety of his analysis, Buniakovsky's works fully retain their importance for the present day …*

*His tables constitute <u>a great step forward</u> but do not represent sufficiently well the picture of Russian mortality. This is explained by the inaccuracy of the main initial materials*, lack […] *of many necessary data and a defect in the very method of compiling these tables.*

On Buniakovsky's method of compiling mortality tables see articles on his method and on Mortality tables in *Demograficheskii* (1985).

In 1848 Buniakovsky published a long newspaper article devoted to a most important subject, to the dread of cholera. However, he likely had not paid due attention to this work. Much later **Enko** (1889) provided the first mathematical model of an epidemic (of measles). It is now highly appreciated (Dietz 1988; Gani 2001) and it might be regretted that Buniakovsky did not become interested in such issues.

9) Among other studies, I mention Buniakovsky's solution of a problem in the theory of random arrangements (1871) that, however, hardly found application, and an interesting urn problem (1875a) connected with partition of numbers. An urn contains *n* balls numbered from 1 through *n*. All at once, *m* balls ($m < n$) are extracted; determine the probability that the sum of the numbers drawn was equal to *s*. This problem is due to **Laurent** (1873) who referred to Euler (1748, Chapter 16).

The problem demanded from Buniakovsky the calculation of the coefficient of $t^m x^s$ in the development of

$$(1 + tx)(1 + tx^2) \ldots (1 + tx^n).$$

He solved this problem by means of an involved partial difference equation and only for small values of *m*; he then provided a formula for the transition from *m* to $(m + 1)$. **Laplace** (§ 7.1-2) solved the same problem in a different way.

10) Approximate stochastic summing and an automatic abacus. Buniakovsky (1867) solved a few problems on the approximate summing of the values of functions (e. g., square and cubic roots for consecutive natural numbers), or air pressure during six months. He had not stated the aim of this latter summing, but it is evident: for deriving mean values.

Buniakovsky reported about his abacus at a sitting of the mathematical and physical class of the Petersburg Academy, then published an appropriate memoir (1867). There, he described the application of such abacuses for treating observations of meteorological elements. Prudnikov (1954, p. 81) suggested that Chebyshev had arrived at the idea about the structure of his arithmometer under the influence of that memoir.

11) From other subjects, apart from participation in compiling explanatory dictionaries of Russian language, I name the solution of a (barely useful) problem on the theory of random arrangements, linguistics (see below), the LLN, treatment of observations, and testimonies, elections and verdicts. The last-mentioned item deserves notice because of later events. I (Sheynin 1989 , p. 340) noted that in 1912 Markov applied to the Most Holy Synod for excommunication from the Russian orthodox Church. Explaining his request, he mentioned, in part, his disagreement with Buniakovsky who had stated that some events ought to be simply believed. Markov, however, had not known the most clear statement made by the latter (1866b, p. 4) to the effect that "truths cognized by revelation" ought to be separated from everything else. See also § 14.1.

In a popular article Buniakovsky (1847) dealt in particular with linguistics. He mentioned his previous unpublished work on this topic and explained the aims of statistically studying linguistics. Regrettably, nothing more is known about him here. On early statistical work in this field see Knauer (1955).

For several decades Buniakovsky's treatise (1846) continued to influence strongly the teaching of probability theory in Russia; in spite of the work of previous Russian authors, he it was who originated the real study of the theory beyond Western Europe.

Several authors expressed their opinion about Buniakovsky (1846). Vasiliev (1921, p. 36) decided that

*This thorough and clearly written book is one of the best from the European mathematical literature on the theory of probability. It much assisted the dissemination of the interest in this science among Russian mathematicians and increased the significance of teaching probability here as compared with the universities of other countries.*

Indeed, Buniakovsky was vice-president of the Petersburg Academy! See also in § 10.5 a statement about the instruction in probability in Belgium.

Markov (1914b/1981, p. 162) considered Buniakovsky's treatise (1846) *a beautiful work* and Steklov (1924, p. 177), president of the Russian Academy of Sciences, believed that *for his time* it was *complete and outstanding*.

I am duty bound, however, to remark that Buniakovsky did not pay attention to the work of **Chebyshev**; after 1846, he actually left probability for statistics. Then, Buniakovsky devoted more than 60 pages to the treatment of observations, but had not thrown proper light on the achievements of Gauss and his exposition was old-fashioned. This is all the more regrettable since even Shiyanov (1836) did more justice to Gauss (but still had not described the second justification of the MLSq).

### 10.5. Quetelet

At the beginning of his scientific career Quetelet visited Paris and met leading French scientists. Some authors indicated that he was much indebted to **Laplace** but I think that the main inspiration to him was **Fourier** (1821 – 1829).

Quetelet tirelessly treated statistical data and attempted to standardize statistics on an international scale. He was co-author of the first statistical reference book (Quetelet & Heuschling 1865) on the population of Europe (including Russia) and the USA that contained a critical study of the initial data; in 1853, he (1974, pp. 56 – 57) served as chairman of the *Conférence maritime pour l'adoption d'un système uniforme d'observation météorologiques à la mer* and the same year he organized the first *International Statistical Congress*.

**K. Pearson** (1914 – 1930, 1924, vol. 2, p. 420) praised Quetelet for "organizing official statistics in Belgium and […] unifying international statistics". The latter (1846, p. 364) complained that *different states are apparently pleased to prevent any rapprochement* [of statistical materials].

Mouat (1885, p. 15) testified that about 1831 – 1833 Quetelet had suggested "the formation of a Statistical Society in London". It was indeed established in 1834 and now called *Royal Stat. Soc*.

Quetelet popularized the theory of probability. His own writings on this subject intended for the general reader are not really interesting, but Mansion (1904, p. 3) claimed that there were "peu de pays" where instruction in that discipline had been on a par with Belgian practice which fact he attributed to the lasting influence of Quetelet.

Quetelet's writings (1869; 1871) contain many dozens of pages devoted to various measurements of the human body, of pulse and respiration, to comparisons of weight and stature with age, etc. and he extended the applicability of the [normal law] to this field. Following **Humboldt's** advice he (1870), introduced the term *anthropometry* and thus curtailed the boundaries of anthropology. He was possibly influenced by **Babbage** (1857), an avid collector of biological data. In turn, Quetelet impressed **Galton** (1869, p. 26) who called him "the greatest authority on vital and social statistics". While discussing Galton (1869), **K. Pearson** (1914 – 1930, vol. 2, 1924, p. 89) declared:

*We have here Galton's first direct appeal to statistical method and the text itself shows that* [the English translation of Quetelet (1846)] *was Galton's first introduction to the* […] *normal curve.*

In those days the preliminary analysis of statistical materials was extremely important first and foremost because of large systematic corruptions, forgeries and incompleteness of data. Quetelet came to understand that statistical documents were only probable and that, in general, *tout l'utilité* of statistical calculations consisted in estimating their trustworthiness (Quetelet & Heuschling 1865, p. LXV).

Quetelet (1846) left recommendations concerning the compilation of questionnaires and the preliminary checking of the data; maintained (p. 278) that too many subdivisions of the data was a *charlatanisme scientifique*, and, what was then understandable, opposed sampling (p. 293). This contribution contained many reasonable statements. In 1850, apparently bearing in mind its English translation, **Darwin** (1887, vol. 1, p. 341) noted:

*How true is a remark* […] *by Quetelet,* […] *that no one knows in disease what is the simple result of nothing being done, as a standard with which to compare homeopathy, and all other such things.*

Almost a null hypothesis! Quetelet (1846, p. 259) also declared that

*The plants and the animals have remained as they were when they left the hands of the Creator.*

Lamarck was the first who attempted to construct a theory of evolution, and Quetelet's statement possibly testifies that his thoughts had been more or less discussed. But Quetelet never mentioned Lamarck or Wallace or Darwin. His attitude partly explains why the statistical study of the evolution of species had begun comparatively late (and only in the Biometric school). I note that **Knapp** (1872b), while discussing Darwin's ideas, had not mentioned randomness and said nothing about a statistical study of biological problems.

Quetelet collected and systematized meteorological observations[5] and described the tendency of the weather to persist by elements of the theory of runs, see § 10.8.3. Keeping to the tradition of political arithmetic (§ 2.1.4), he discussed the level of postal charges (1869, t. 1, pp. 173 and 422) and rail fares (1846, p. 353) and recommended to study statistically the changes brought about by the construction of telegraph lines and railways (1869, t. 1, p. 419). His special investigation (1836, t. 2, p. 313) was a quantitative description of the changes in the probabilities of conviction of the defendants depending on their personality (sex, age, education, first noted in 1832) and **Yule** (1900/1971, pp. 30 – 32) favourably, on the whole, commented upon his work as the first attempt to measure dependence (association).

Quetelet is best remembered for the introduction of the Average man (1832a, p. 4; 1832b, p. 1; 1848a, p. 38), inclinations to crime (1832b, p. 17; 1836, t. 2, p. 171 and elsewhere) and marriage (1848a, p. 77; 1848b, p. 38), – actually, the appropriate probabilities, – and statements about the constancy of crime (1829, pp. 28 and 35 and many other sources). In spite of his shortcomings (below), the two last-mentioned items characterized Quetelet as the originator of moral statistics with Süssmilch as his predecessor.

The Average man, as he (1848a, p. 38) thought, had mean physical and moral and intellectual features, was the alleged type of the nation and even of the entire mankind. And he (1848a, pp. 91 – 92) properly related the mean inclinations to the Average man. Reasonable objections were levelled against this concept to which I add that Quetelet had not specified the notion of *average* as applied here. Sometimes he had in mind the arithmetic mean, in other cases (1848a, p. 45), however, it was the median, and he (1846, p. 216) only mentioned the **Poisson** LLN in connection with the mean human stature. **Cournot** (1843, p. 143) stated that the Average man was physiologically impossible (the averages of the various parts of the human body [and of weight and height] were inconsistent one with another), and Bertrand (1888a, p. XLIII) ridiculed Quetelet:

*In the body of the average man, the Belgian author placed an average soul*. [The average man] *has no passions or vices* [wrong]*, he is neither insane nor wise, neither ignorant nor learned.* […] [He is] *mediocre in every sense. After having eaten for thirty-eight years an average ration of a healthy soldier, he has to die not of old age, but of an average disease that statistics discovers in him.*

Nevertheless, the Average man is useful even now at least as an average producer and consumer; **Fréchet** (1949) replaced him by a closely related "typical" man.

Quetelet (1848a, p. 82; 1869, t. 2, p. 327) indicated that the real inclination to crime of a given person might well differ considerably from the apparent mean tendency. It seems, however, that he had not sufficiently stressed these points; a noted statistician (Rümelin 1867, p. 25) forcibly denied any criminal tendency in himself.

**Chuprov** (1909/1959, p. 23) neatly summed up the arguments of Quetelet's followers and opponents:

*Their* [his worshippers' *zealous beyond reasoning*] *naïve admiration for "statistical laws"; their idolizing of stable statistical figures; and their absurd teaching that regarded everyone as possessing the same "mean inclinations" to crime, suicide and marriage, undoubtedly provoked protests. Regrettably, however, the protests were hardly made in a scientific manner.*

Quetelet (1836, t. 1, p. 10) declared that the [relative] number of crimes was constant, and that

*Each social regime presupposes […] a certain number and a certain order of crimes, these being merely the necessary consequences of its organization.*

However, he had not justified his statement by statistical data. The alleged constancy did not take place (Rehnisch 1876): Quetelet had not studied criminal statistics attentively enough. Then, constancy of crime could only happen under constant social conditions, but this consideration had only indirectly followed from his statements.

A special point concerns laws of distribution. Quetelet (1848a, p. 80; 1869, t. 2, pp. 304 and 347) noticed that the curves of the inclinations to crime and to marriage plotted against ages were exceedingly asymmetric. He (1846, pp. 168 and 412 – 424) also knew that asymmetric densities occurred in meteorology and he (1848a, p. viii) introduced a mysterious "loi des causes accidentelles" whose curve could be asymmetric (1853, p. 57)! And still he (1853) explained this fact by special causes and anomalies.

In short, he had revealed here and elsewhere, see above and my paper (1986a), his general attitude which **Knapp** (1872a, p. 124) explained by his "spirit, rich in ideas, but unmethodical and therefore unphilosophical". His was a polite appraisal indeed; in actual fact, Quetelet often acted in a happy-go-lucky manner.

Nevertheless, Quetelet had been the central figure of statistics in the mid-19$^{th}$ century. **Freudenthal** (1966, p. 7) correctly concluded that before Quetelet there were statistical bureaux and statisticians, but no statistics [as a discipline].

### 10.6. Galton

Being influenced by his cousin, **Darwin**, Galton began to study the heredity of talent and published an important treatise (1869) on that subject; incidentally, he introduced the term *eugenics*. In a letter of 1861 Darwin (1903, p. 181) favourably mentioned that contribution. He (1876, p. 15) also asked Galton to examine his investigation of the advantages of cross-fertilization as compared with spontaneous pollination. Galton (Ibidem) compared the two processes with regard to their characteristics and, in particular, to the ordered heights of the seedlings. In the latter instance, it occurred that the signs of almost all the differences between the corresponding heights coincided. See a similar study by Seidel in § 10.8.1.

Galton (1863) devised an expedient system of symbols for weather charts and immediately discovered the existence of previously unknown anticyclones. From the point of view of statistics, he had thus reasonably studied his initial data. Galton (**K. Pearson** 1914 – 1930, vol. 2, Chapter 12) also invented *composite photographs* of kindred persons (of people of a certain nationality or occupation, or criminals), all of them taken on the same film with an appropriately shorter exposure. In any case, his innovation was much more justified than **Quetelet's** Average man.

Galton(1892) became the main inventor of fingerprinting. Because of its reliability, it did not demand statistical analysis and superseded the previous system of identification developed by Alph. Bertillon (1893). This latter procedure was partially based on anthropometry and made use of from the 1890s to the beginning of the 20$^{th}$ century. Another of Galton's invention (1877) was the so-called *quincunx*, a device for visually demonstrating the appearance of the normal distribution as the limiting case of the binomial law (Stigler 1986, pp. 275 –281). A special feature of that device was that it showed that the normal law was stable. Galton's main statistical merit consisted, however, in the introduction of the notions of regression and correlation. The development of correlation theory became one of the aims of the Biometric school (§ 14.2), and Galton's close relations with **Pearson** were an important cause of its successes.

Recall (§§ 1.1.1 and 1.1.3) that reasoning in the spirit of qualitative correlation was not foreign to ancient scholars which was in conformity with the qualitative nature of the science of those days. And what about modernity? In the 1870s, several scientists (C. Meldrum, in 1872 and 1875; J. N. Lockyer, in 1873; H. F. Blanford, in 1880, see Sheynin (1984a, p. 160)) took notice of the dependence between solar activity and elements of terrestrial magnetism and on meteorological phenomena but not a word did they say about developing a pertinent quantitative theory. And even though Seidel, in 1865 – 1866 (§ 10.8.1), quantitatively studied the dependence between two, and then three factors, he did not hint at generalizing his findings. Galton was meritorious indeed! For the sake of comprehensiveness I repeat (Note 4 to Chapter 9) that in 1912 **Kapteyn** provided an "astronomical" version of the correlation coefficient.

### 10.7. Statistics

Here, I discuss the situation in statistics in the 19$^{th}$ century. Related material is in §§ 6.2 and 10.5.

The *Staatswissenschaft* held its ground for many decades. In France, **Delambre** (1819, p. LXVII) argued that statistics was hardly ever engaged in discussions or conjectures and did not aim at perfecting theories, and that political arithmetic ought to be "distinguished" from it. Under statistics he understood geodetic, meteorological and medical data, mineralogical descriptions and even art expositions. I believe however that the two last-mentioned items were rather soon excluded from such lists[6].

The newly established London Statistical Society declared that statistics "does not discuss causes, nor reason upon probable effects" (Anonymous 1839, p. 1). True, they denied that "the statist [!] rejects all deductions, or that statistics consists merely of columns of figures" and stated that "all conclusions shall be drawn from well-attested data and shall admit of mathematical demonstration". This announcement was thus ambiguous; the Society attempted to adhere to its former statement, but in vain. Anyway, Woolhouse (1873, p. 39) testified that

*These absurd restrictions have been necessarily disregarded in [...] numerous papers*.

Indeed, that statistics should explain the present state of a nation by considering its previous states was declared a century before that (Gatterer 1775, p. 15). And the very title of Dufau (1840) called statistics "The theory of studying the laws according to which the social events are developing"[7].

During the 19$^{th}$ century the importance of statistics had been considerably increasing. **Graunt** (1662/1939, p. 79) was not sure whether his work would be "necessary to many, or fit for others, than the Sovereign, and his chief Ministers ..." and the classical investigations of the sex ratio at birth (§§ 2.2.4, 3.3.4, 4.4, 6.1.1) had not found direct applications. However, by the mid-19$^{th}$ century it became important to foresee how various transformations will influence society and **Quetelet** (§ 10.5) repeatedly stressed this point. Then, at the end of the 19$^{th}$ century censuses of population, answering an ever widening range of questions, began to be carried out in various countries. It is nevertheless instructive to compare the situation at that time with what is happening nowadays[8].

1) Public opinion was not yet studied, nor was the quality of mass production checked by statistical methods, cf. § 10.4-6.

2) Sampling had been considered doubtful[9]. **Cournot** (1843) passed it over in silence and **Laplace's** determination of the population of France based on sampling (§ 7.1-5) was largely forgotten even in spite of the inexorable increase in statistical materials. Indeed, already the beginning of the century witnessed "legions" of new data (Lueder 1812, p. 9) and the tendency to amass sometimes useless or unreliable data revealed itself in various branches of natural sciences (§ 10.8). Note however that the need for observations (in natural science?) increases with knowledge (Descartes 1637/2012, p. 63).

**Quetelet** (§ 10.5) opposed sampling. Even much later **Bortkiewicz** (1904, p. 825) and **Czuber** (1921, p. 13) called sampling "conjectural calculation" and **Chuprov** (1912) had to defend that procedure vigorously.

3) The development of the correlation theory began at the end of the 19$^{th}$ century (§§ 10.6, 15.2), but even much later Kaufman (1922, p. 152) declared that

*The so-called method of correlation adds nothing essential to the results of elementary analysis*.

See, however, § 14.1-4.

4) Apart from the error theory the variance began to be applied in statistics only after **Lexis** (§ 15.1), but even later **Bortkiewicz** (1894 – 1896, Bd. 10, pp. 353 – 354) stated that the study of precision was an accessory goal, a luxury, and that the statistical flair was much more important, cf. the opinion of **Gauss** in § 9A.5-1. This point of view had perhaps been caused by the presence of large systematic corruptions in the initial materials.

5) Not just a flair, but a preliminary data analysis (which, however, does not call off the definitive estimation of the plausibility of the final results and which received general recognition only a few decades ago) is necessary, and should be the beginning of the statistician's work. Splendid examples of such analysis had occurred much earlier and to these I attribute the introduction of contour lines (**Halley**, in 1701**,** see § 2.1.4, drew lines of equal magnetic declinations over North Atlantic, also see § 10.8.3, the discoveries made by Humboldt in § 10.8.3 and by Galton in § 10.6).

6) Econometrics originated only in the 1930s.

I list now the difficulties, real and imaginary, of applying the theory of probability to statistics.

7) The absence of "equally possible" cases whose existence is necessary for understanding the classical notion of probability. Statisticians repeatedly mentioned this cause, also see § 3.2.3. True, **Cournot** (§ 10.3-7 and -8) explained that equipossibility was not necessary (and, in the first place, mentioned the "**Bernoulli** principle"), but his advice was hardly heard. **Lexis** (report of 1874/1903, pp. 241 – 242; 1886, pp. 436 – 437; 1913, p. 2091) also cited equipossibility. In the second case, in a paper devoted to the application of probability theory to statistics, he even added that the introduction of equipossibility led to the subjectivity of the theory of probability. Elsewhere, however, Lexis reasoned differently; he had no integral viewpoint. Thus, statistics is mainly based on the theory of probability (1877, p. 5); if the statistical probability tends to its theoretical counterpart, equally possible cases might be assumed (Ibidem, p. 17); and the "pattern" of the theory of probability is the highest scientific form in which statistics might be expressed (1874/1903, p. 241).

8) Disturbance of the constancy of the probability of the studied event and/or of the independence of trials. I repeat (§ 3.2-3) that before **Lexis** statisticians had only recognized the **Bernoulli** trials; and even much later Kaufman (1922/1928, pp. 103 – 104) declared that the theory of probability was applicable only to these trials, and, for that matter, only in the presence of equally possible cases. He mentioned several allegedly likeminded authors including **Markov** and **Yule**, but did not supply the exact references and I am inclined to believe that these authors wished to investigate whether or not the given statistical trials were Bernoullian. As to the equally possible cases, see Item 7 above. After all, Kaufman's opinion about correlation theory (Item 3) can be understood, but this time he proved himself hopelessly lagging behind life. See however § 14.1-4.

9) The abstract nature of the (not yet axiomatized) theory of probability. As I noted in § 6.3.3, the history of mathematics testifies that the more abstract it became, the wider had been the range of its applicability. Nevertheless, statisticians had not expected any help from the theory of probability. Block (1878/1886, p. 134) thought that it was too abstract and should not be applied "too often", and **Knapp** (1872, p. 115) called it difficult and hardly useful beyond the sphere of games of chance and insurance.

There also, on pp. 116 – 117, he added a sound statement:

*Placing coloured balls in Laplacean urns is not enough for shaking scientific statistics out of them.*

For his part, **Chuprov** (1922/1960, pp. 415 – 416) formulated a related opinion:

*Only statisticians armed with mathematics can defeat mathematicians playing with statistics.*

Statisticians had not been *armed* and did not trust mathematicians. Even in 1911 G. von Mayr, an eminent statistician of the old school stated that mathematical formulas were unnecessary for statistics and then privately confessed that he cannot endure mathematics (Bortkevich & Chuprov 2005, Letter 109). And mathematicians, if not playing with statistics, had not always understood statistics (§ 8.9-2).

It is not amiss to mention here the pioneer attempt to create mathematical statistics (Wittstein 1867). He compared the situation in statistics with the *childhood* of astronomy and stressed that statistics (and especially population statistics) needed a Tycho and a Kepler to proceed from reliable observations to regularities. Specifically, he noted that statisticians did not understand the essence of probability theory and never estimated the precision of the results obtained. The term *mathematical statistics* is apparently due to him.

### 10.8. Statistics and Natural Sciences

In the 19th century, the statistical method gave rise to a number of disciplines and I discuss the relevant situation in several branches of natural sciences. First, I cite the opponents of that method; as described in § 10.5, their stance can be explained by ignorance, by attaching mean indicators to individuals. Thus, Comte (1830 – 1842, t. 3/1893, No. 40, p. 329): the use of statistics in medicine is a "profonde dégénération directe de l'art médicale". Also see similar statements in § 8.9.

To begin in earnest, however, I note the existence of the so-called numerical method usually attributed to the French physician **Louis** (1825) who introduced it by calculating the frequencies of the symptoms of various diseases so as to facilitate diagnosing. No stochastic considerations or estimations of the reliability of the conclusions were involved. Louis (pp. xvii – xviii) even thought that, given observations, any physician ought to make the same conclusion. His method remained in vogue for a few decades.

He and his adherents attempted to replace qualitative descriptions by directly obtained statistical data (cf. **Petty's** statement in § 2.1.4). Bouillaud (1836), who inserted numerous passages from **Laplace's** *Essai philosophique* (1814) in his book, favourably described the numerical method and called it a supplement to other methods (pp. 190 – 191), and (p. 187) added only a few words about the "calcul approximatif ou des probabilités". That method, as he stated, was almost always the only means for generalizing the results obtained; and the advantages of this "kind of calculus" are such that its discussion was not necessary. Also see Proctor's statement in § 10.8.4. Bouillaud (pp. 186 – 187) thought that the calculus of probability was *approximate*, but almost always the only means for generalizing the obtained results, whereas medicine, if only based on probability, will remain a game of chance of sorts since physicians ought to consider the personality of their patients. And he also stated that medical statistics was yet in its cradle, but is still applied with some success and will considerably develop, cf. § 10.8.1.

And here is **D'Alembert** (1759/1821, p. 163):

*Systematic medicine seems to me* […] *a real scourge of mankind. Observations, numerous and detailed, sensibly corresponding one to another, – this is* […] *to what reasoning in medicine ought to reduce itself.*

On D'Alembert see also Note 7 to Chapter 6. Another author who advocated compilation of medical observations is mentioned at the beginning of § 6.2.3.

Unlike Bouillaud, Gavarret did not sidestep this issue (§ 8.9.2), and he (1840, p. x) reasonably remarked that the numerical method was not in itself scientific and was not based on "general philosophy". There also existed a wrong opinion (D'Amador 1837, p. 12) attributing the numerical method to probability theory. Guy (1852, p. 803), who was an eminent physician and a statistician of note, applied the term "numerical or statistical method", but only to prevent statistics from being confused with Staatswissenschaft (§ 6.2.1). There also (pp. 801 – 802) he claimed that medicine had a special relation with statistics:

*There is no science which has not sooner or later discovered the absolute necessity of resorting to figures as measures and standards of comparison; nor is there any reason why physiology and medicine should claim exemption.* […] *On the contrary, they* […] *may hope to derive the greatest benefit from the use of numbers.* […] *Without statistics a science bears to true science the same sort of relation which tradition bears to history.*

However, the numerical method (not necessarily in medicine) can be traced back to the 18th century (see below and § 6.2.3) and my description (§§ 10.8.1 – 10.8.4) shows that it continued in existence for many decades. Furthermore, empiricism had been a feature of the

Biometric school (§ 15.2). Greenwood (1936, p. 139) indirectly praised it, perhaps excessively:

*Some heart-breaking therapeutic disappointments in the history of tuberculosis and cancer would have been avoided if the method of Louis had been not merely praised but generally used during the last fifty years.*

In statistics proper, **Fourier's** fundamental *Recherches* (1821 – 1829) concerning Paris and the Département de la Seine might be here mentioned. This contribution almost exclusively consisted of statistical tables with data on demography, industry, commerce, agriculture and meteorology. True, empiricism was not sufficient even for compiling tables. Then, the abundance of materials led to the wrong idea that a mass of heterogeneous data was better than a small amount of reliable observations (§ 10.8.1).

In actual fact, the numerical method originated with **Anchersen** when statisticians have begun to describe states in a tabular form (and thus facilitated the use of numbers), see § 6.2.1. Recall (§ 2.1.4), moreover, that **Leibniz** recommended compilation of *Staatstafeln*.

**10.8.1. Medicine.** See also § 8.9.2. In 1835, Double et al (§ 8.9) indicated that statistics might be applied in medicine. Surgery occurred to be the first branch of medicine to justify their opinion. Already in 1839 there appeared a (not really convincing) statistical study of the amputation of limbs. Soon afterwards physicians learned that the new procedure, anaesthesia, could cause complications, and began to compare statistically the results of amputation made with and without using it. The first such investigation (**J. Y. Simpson** 1847 – 1848/1871, p. 102) was, however, unfortunate: its author had attempted to obtain reliable results by issuing from materials pertaining to several English hospitals during 1794 – 1839:

*The data I have adduced […] have been objected to on the ground that they are collected from too many different hospitals, and too many sources. But, […] I believe all our highest statistical authorities will hold that this very circumstance renders them more, instead of less, trustworthy.*

I ought to add, however, that Simpson (Ibidem, p. 93) stated that only a statistical investigation could estimate the ensuing danger.

He (1869 – 1870/1871, title of contribution) also coined the term *Hospitalism* which is still in vogue. He compared mortality from amputations made in various hospitals and reasonably concluded, on the strength of its monotonous behaviour, that mortality increases with the number of beds; actually (p. 399), because of worsening of ventilation and decrease of air space per patient. Virchow (1868 – 1869/1879, Bd. 2, p. 21) noted that the pernicious influence of bad air had been known about six hundred years previously.

Justification of conclusions based on monotonic behaviour of some indications was not restricted to medicine; **Quetelet's** table of probabilities of conviction of defendants depending on their personality showed a monotonous increase of those probabilities (cf. § 10.5).

At about the same time **Pirogov** introduced anaesthesia in military surgery and began to compare the merits of the conservative treatment of the wounded versus amputation. Much later he (1864, p. 690) called his time "transitional":

*Statistics shook the sacred principles of the old school, whose views had prevailed during the first decades of this century, – and we ought to recognize it,– but it had not established its own principles.*

Pirogov (1849, p. 6) reasonably believed that the application of statistics in surgery was in "complete agreement" with the latter because surgical diseases depended incomparably less on individual influences. However, he repeatedly indicated that medical statistics was unreliable. Thus (1864/1865 – 1866/1920, p. 20):

*Even a slightest oversight, inaccuracy or arbitrariness makes* [the data] *far less reliable than the figures founded only on a general impression with which one is left after a mere but sensible observation of cases.*

Later he (1879/1882, p. 40) singled out an important pertinent cause:

*Extremely different circumstances separate the entire mass of information in too insignificant and very dissimilar groups which does not allow any correct conclusion about the worth of a certain amputation.* [In 1849] *I* […] *had not yet known all the false ways to which the number sometimes leads us.*

In essence, he advocated attentive allowance for all circumstances and minimal statistical technique which was in accordance with his time and especially so with the originating military surgery (of which he was the founder).

Pirogov was convinced in the existence of regularities in mass phenomena. Thus (1850 – 1855/1961, p. 382), each epidemic disease as well as each "considerable" operation had a constant mortality rate, whereas war was a "traumatic epidemic" (1879/1882, p. 295). This latter statement apparently meant that under war conditions the sickness rate and mortality from wounds obeyed statistical laws. Then (1854, p. 2), the skill of the physicians [but not of witch doctors] hardly influenced the total result of the treatment of many patients. Here is his highly relevant opinion (1871, pp. 48 – 49):

*On what does the success of treatment or the decrease of mortality in the army depend? Surely not on therapy and surgery by themselves. Without an efficient administration* [of medicine] *little can the masses expect from therapy and surgery even in time of peace, much less during such a catastrophe as war.*

Note finally Pirogov's possibly correct statement (1864, pp. 5 – 6): without the not yet existing doctrine of individuality, real progress in medical statistics is impossible.

Pirogov participated in the Crimean war, in which **Florence Nightingale**, on the other side, showed her worth both as a medical nurse and a statistician, cf. **Pearson's** relevant statement in § 2.2.3. She would have wholeheartedly approved of Pirogov's conclusion (above) concerning the success of treatment.

Such new disciplines as epidemiology and public hygiene, the forerunner of ecology, appeared within medicine in the 19$^{th}$ century. I discussed the inoculation of smallpox in § 6.2.3 and mentioned **Enko's** essential finding at the end of § 10.4. In 1866, **Farr** (Brownlee 1915) preceded Enko; his study of cattle plague only methodically belonged to epidemiology, and, interestingly enough, Brownlee published his note in a medical journal. Farr indicated that he had also investigated the visitations of cholera and diphtheria of 1849 and 1857 – 1859 respectively.

Here is his reasoning. Denote the number of attacks of the plague during a period of four weeks by *s*. He noted that the third differences of ln*s* were constant, so that

$$s = C\exp\{\ t[t + m)^2 + n]\}, C > 0,\ < 0.$$

It was Brownlee who supplied this formula because Farr was unable to insert it in his newspaper letter. Farr's calculated values of *s* did not agree with actual figures, but at least he correctly predicted a rapid decline of the epidemic.

It seems that epidemiology was properly born when cholera epidemics had been ravaging Europe. The English physician **Snow** (1855) compared mortality from cholera for two groups of the population of London, – for those whose drinking water was either purified or not. He ascertained that purification decreased mortality by eight times, and he thus discovered how did cholera epidemics spread, and proved the essential applicability of the first stage of the statistical method (§ 0.4). **Pettenkofer** (1886 – 1887) published a monstrous collection of statistical materials pertaining to cholera, but he was unable to process them. He (1865, p. 329) stressed that no cholera epidemic was possible at a certain moment without a local "disposition" to it and he attached special importance to the level of subsoil water. His view does not contradict modern ideas about the necessary threshold values. However, Pettenkofer did not believe in contemporary bacteriological studies and opposed Snow. For an estimate of his views see Winslow (1943/1967, p. 335).

**Seidel** (1865) arranged the years 1856 – 1864 in the decreasing order of the first, and then of the second series of numbers describing two phenomena. The conformity between the two series was, in his opinion, striking so that Seidel thus, like Galton later on, applied rank correlation.

He (1865 – 1866) investigated the dependence of the monthly cases of typhoid fever on the level of subsoil water, and then on both that level and the rainfall. It occurred that the signs of the deviations of these figures from their mean yearly values coincided twice more often than not and Seidel quantitatively (although indirectly and with loss of information) estimated the significance of the studied connections. His work remained, however, completely forgotten and Weiling (1975) was likely the first to recall it.

Already **Leibniz** (§ 2.1.4) recommended to collect and apply information concerning a wide range of issues, which, as I add now, pertained to public hygiene. **Condorcet** (1795/1988, pp. 316 and 320) described the aims of "mathématique sociale" [social statistics] and mentioned the study of the influence of temperature, climate, properties of soil, food and general habits on the ratio of men and women, birth-rate, mortality and number of marriages. Much later, M. Lévy (1844) considered the influence of atmosphere, water and climate as well as of the suitable type of clothes and appropriate food on man.

From its origin in the mid-19$^{th}$ century, public hygiene began statistically studying a large number of problems, especially those caused by the Industrial Revolution in England and, in particular, by the great infant mortality. Thus, in Liverpool only 2/3 of the children of gentry and professional persons lived to the age of five years (Chadwick 1842/1965, p. 228).

**Pettenkofer** (1873) estimated the financial loss of the population of Munich ensuing from such diseases as typhoid fever and his booklet can be attributed to this discipline. In Russia his student **Erismann** (1887) published a contribution on sanitary statistics.

At the turn of the 19$^{th}$ century Jenner introduced smallpox vaccination instead of inoculation. It had not, however, preclude the need for solving statistical problems about the estimation of various possible versions of vaccination and inoculation was not ruled out at once. See also § 6.2.3.

**10.8.2. Biology.** The attempts to connect the appearance of leaves, flowers and fruits on plants of a given species with the sums of mean daily temperatures began in the 18$^{th}$ century (**Réaumur** 1738) and **Quetelet** (1846, p. 242) proposed to replace those sums by the sums of squares, but he was still unable to compare both procedures quantitatively. Also in the 19$^{th}$ century, vast statistical materials describing the life of plants were published (**Aug. P. DeCandolle** 1832), and **Babbage** (1857) compiled a statistical questionnaire for the class of mammalia. I mentioned his similar work in § 10.5. In Russia, **Baer** (1860 – 1875) with associates conducted a large-scale statistical investigation of fishing.

**Humboldt** created the geography of plants (Humboldt & Bonpland 1815; Humboldt 1816) which was based on collection and estimation of statistical data. Alph. DeCandolle (1855, t. 1, p. vi) however also mentioned Linné, his own father Aug. DeCandolle and Brown; Darwin (1903, vol. 2, p. 26, letter of 1881) only mentioned Humboldt. **Darwin** had to study various statistical problems, for example on cross-fertilization of plants (§ 10.6), the life of earthworms (§ 12-2) and on the inheritance of a rare deformity in humans (1868/1885, vol. 1, p. 449). In the last-mentioned case **Stokes** provided the solution (apparently by applying the **Poisson** distribution) at his request. Statistical tables and summaries with qualitative commentaries occur in a number of Darwin's writings and he also collected statistical data.

Being the main author of the hypothesis of the origin of species, he made use of such terms as variation and natural selection without defining any of them. And, when reasoning about randomness, he understood it in differing ways. In the problem concerning the deformity Darwin decided that it was not random (not merely possible, cf. Kolmogorov's statement in § 1.1). In two other cases in which he discussed the hypothesis of evolution he understood randomness as ignorance of causes (1859/1958, p. 128), cf. **Laplace** (§ 7.3), and, in 1881, as lack of purpose (1903, p. 395), cf. **Aristotle** (§ 1.1.1). It is also remarkable that Darwin (1859/1958, p. 77) actually described randomness as the effect of complicated causes(cf. **Poincaré** § 11.2-9):

*Throw up a handful of feathers, and all fall to the ground according to definite laws; but how simple is the problem where each shall fall compared with problems in the evolution of species.*

The stochastic essence of the evolution hypothesis was evident both for its partisans and the opponents; **Boltzmann**, however, was an exception (§ 10.8.5-4).

I reconstruct now **Darwin's** model of evolution. Introduce an $n$-dimensional (possibly with $n = $  ) system of coordinates, the body parameters of individuals belonging to a given species (males and females should, however, be treated separately), and the appropriate **Euclidean** space with the usual definition of distances between its points. At moment $t_m$ each individual is some point of that space and the same takes place at moment $t_{m+1}$ for the individuals of the next generation. Because of the "vertical" variation, these, however, will occupy somewhat different positions. Introduce in addition point (or subspace) $V$, corresponding to the optimal conditions for the existence of the species, then its evolution will be represented by a discrete stochastic process of the approximation of the individuals to $V$ (which, however, moves in accordance with the changes in the external world) and the set of individuals of a given generation constitutes the appropriate realization of the process. Probabilities describing the process (as well as estimates of the influence of habits, instincts, etc.) are required for the sake of definiteness, but they are of course lacking.

The main mathematical argument against Darwin's hypothesis was that a purposeful evolution under "uniform" randomness was impossible; see end of § 6.1.3 with regard to the difficulties of generalizing the notion of randomness. Only **Mendel's** contributions (1866; 1866 – 1873, publ. 1905), forgotten until the beginning of the 20[th] century, allowed to answer such criticisms. True, great many objections and problems still remain, but at the very least **Darwin** had transformed biology as a science. In addition, his work was responsible for the appearance of the Biometric school (§ 15.2).

An essentially new stage in the development of the Darwinian ideas had occurred at the end of the 18$^{th}$ century (De Vries 1905) and somewhat later (Johannsen 1922/1929). De Vries stressed the importance of *sports* (although did not explain their relation to mutations, a term which he himself introduced somewhat earlier) which considerably strengthened the theory of evolution.

Andersson (1929), who briefly discussed the work of Johannsen, quoted him, regrettably without indicating the source:

*The science of evolution has turned into an Augeas stables which really ought to be mucked.*

And here is Johannsen himself.

p. 355. *Galton's statistics of heredity were quite erroneous – a combination between collective measurements of unsorted rough material and biological analysis of the real units of certain populations.*

p. 356. *These statistical researches in heredity are naturally of importance from the sociological point of view and of practical interest in insurance calculations and so on. But they do not reach the biological problems of heredity.*

p. 357. *The* [genotypic] *differences are <u>discontinuous</u> […] rather contrary to Darwinism.*

p. 359. *We cannot do without statistics!*

From the mathematical point of view, Mendel did nothing except for an elementary application of the binomial distribution, but his memoir marked the origin of a new direction in biology, of genetics, and provided an example of a fundamental finding achieved by elementary means. Mendel had based his conclusions on experiments, and these became the object of many discussions with regard to his initial data and to his subjective and objective honesty. Such scholars as **Fisher** (1936) and **van der Waerden** (1968) participated in the debates, and finally all doubts have possibly blown over the more so since Mendel's life and his meteorological observations and investigations unquestionably testify in his favour. It is thought that Mendel was born in a mixed Czech-German family; actually, however, he was German, and in 1945 – 1946 the descendants of his relatives were driven out of the then Czechoslovakia[10].

**10.8.3. Meteorology.** The material pertaining to the 18$^{th}$ century is in § 6.2.4. **Humboldt** (1818, p. 190) maintained that

*To discover the laws of nature* [in meteorology] *we ought to determine the mean state of the atmosphere and the constant type*[s] *of its variations before examining the causes of the local perturbations*[11].

He (1845 – 1862, Bd. 1, pp. 18 and 72; Bd. 3, p. 288) conditioned the investigation of natural phenomena by examination of mean states. In the latter case he mentioned "the sole decisive method [in natural sciences], that of the mean numbers" which (1845, Bd. 1, p. 82)

"show us the constancy in the changes". He himself (1817, p. 466) introduced isotherms and climatic belts (known to ancient scholars who had only possessed qualitative knowledge of temperature). He thus separated climatology from meteorology and he also noted the existence of local corruptions of the temperature.

Much later he (1845 – 1862, Bd. 4, p. 59) added that he had borrowed the idea of contour lines from **Halley** (§ 2.1.4) [and had therefore also applied a splendid particular instance of exploratory data analysis].

Humboldt (1817, p. 532) also recommended the application of contour lines for winter and summer. It is somewhat strange that, when offering a definition of climate, he (1831, p. 404) had not directly linked it with mean states but later scholars have formulated this tie ever more explicitly (Körber 1959, p. 296). **Chuprov** (1922b/1960, p. 151), for example, had identified climate with a system of certain mean values.

Humboldt (1843, t. 1, p. 83) also formulated a statement which could have shown that the time for passing from mean values to distributions had not yet come:

*Exact sciences only developed to the extent that they considered physical phenomena in their ensemble and gradually ceased to attach too much importance either to culminating points isolated among a line of facts, or to the extreme temperatures during some days of the year.*

He (Ibidem, p. 405) provided a splendid example: large nuggets cannot be a sure indication of the mean content of gold in their vicinity.

**Köppen** (1874, p. 3) dwelt on the same subject. He believed that the "introduction of the arithmetical mean in meteorology was the most important step", but that it was not sufficient all by itself.

**Dove** (1837, p. 122) came out against "the domination" of mean values; largely following Humboldt (see above), he (1839, p. 285) formulated the aims of meteorology as the "determination of mean values [of temperature], derivation of the laws of [its] periodic changes and indication of rules for [determining its] irregular changes", and he attached no less importance to the spatial scatter of the temperature. Later Dove (1850, p. 198) introduced monthly isotherms.

**Buys Ballot** (1850, p. 629) was even more straightforward: the study of deviations from mean values (mean states) constituted the second stage in the development of meteorology. He (1847, p. 108) noted that a similar process was going on in astronomy (planetary orbits had been assumed elliptical, then perturbations began to be studied). The same, he added, was the situation concerning all sciences that did not admit experimentation. He could have cited geodesy and the ever more precise determination of the figure of the Earth.

In the 18$^{th}$ and apparently in the early 19$^{th}$ century the mean monthly and even the mean yearly temperatures had been considered as the means of their extreme values (Cotte 1788, p. 9) so that the occurred introduction of the arithmetic mean of all of their values was an essential step forward.

Meteorological observations multiplied, and they had been published without being of use to the general readership of scientific periodicals. **Biot** (1855, pp. 1179 – 1180), for example, had opposed that practice and **Mendeleev** (1876/1946, p. 267) remarked that the prevailing "collecting" school of meteorologists needed nothing but "numbers and numbers". Later he (1885/1952, p. 527) optimistically decided that a new meteorology was being born and that "little by little" it had begun, basing its work on statistical data, to "master, synthesize, forecast".

**Lamont** (1867, p. 247) maintained that the irregular temporal changes of the atmosphere were not random "in the sense of the calculus of probability" and (p. 245) recommended his own method of studying, instead, simultaneous observations made at different localities. **Quetelet** (1849, t. 1, Chapter 4, p. 53) remarked that the differences of such observations conformed to accidental errors, but he did not elaborate. Much earlier Lamont (ca. 1839, p. 263) indicated that the air pressure and temperature are very changeable and their mean values are barely reliable. Again, he did not elaborate and stated without proof that a year of simultaneous observations is tantamount to 30 years of usual observations.

**Quetelet** (1846, p. 275) resolutely contended that meteorology was alien to statistics: unlike the "physicist", the statistician wishes to know, first of all, everything that can influence man and contribute to his welfare. In addition to meteorology, he cited other "alien" sciences, such as physical geography, mineralogy, botany. His statement is correct only insofar as statistical meteorology, stellar statistics etc. belong to the appropriate sciences.

The study of densities of the distributions of meteorological elements began in the mid-19$^{th}$ century; Quetelet, for example, knew that these densities were asymmetric (§ 10.5). At the end of the century Meyer (1891, p. 32), when mentioning that fact, stated that the theory of errors was not applicable to meteorology. However, mathematical statistics does not leave aside the treatment of asymmetric series of observations, and already **K. Pearson** (1898) made use of Meyer's material for illustrating his theory of asymmetric curves.

Lamarck was one of the first scholars to note the dependence of the weather on its previous state, see for example t. 5, pp. 5 and 8 and t. 11, p. 143 of his *Annuaries* (1800 – 1811). In t. 11, p. 122 he essentially repeated his first pronouncement and formulated it as an aphorism:

*The entire state of things in the atmosphere* […] *results not only from an ensemble of causes which tend to operate, but also from the influence of the previous state*.

Quetelet (1852; 1853, p. 68; 1849 – 1857/1857, pt. 5, pp. 29 and 83) repeatedly mentioned lasting periods of fair or foul weather. He

(1852; 1857) analysed that phenomenon applying elementary stochastic considerations and concluded that the chances of the weather persisting (or changing) were not independent. Köppen's analysis (1872) was more mathematically oriented, see Sheynin (1984b, p. 80).

Quetelet was also praiseworthy for compiling and systematizing meteorological observations. Köppen (1875, p. 256) noted that observations made since the early 1840s at the "entire network of stations" in Belgium "proved to be the most lasting [in Europe] and extremely valuable". And **Faraday** praised Quetelet's observations of atmospheric electricity (Note 5).

**10.8.4. Astronomy.** Here, I only discuss the events of the 19$^{th}$ century. Already **Daniel Bernoulli** (§ 6.1.1) and **Laplace** (§ 7.1-2) stochastically studied regularities in the Solar system. In actual fact, they considered the planets as elements of a single population, and this approach was vividly revealed in the later investigations of the asteroids. **Newcomb** (1861a and elsewhere) repeatedly compared the theoretical (calculated in accordance with the uniform distribution) and the real parameters of their orbits; true, he was yet unable to appraise quantitatively his results. **Poincaré** (§ 11.2-5) stochastically although worthlessly estimated the total number of the small planets.

Of special interest are Newcomb's considerations (1862) on the distribution of the asteroids, likely based on his later published and even more interesting statement (1881). His former contribution makes difficult reading mostly because of its loose style. As I understand him, Newcomb intuitively arrived at the following proposition: a large number of independent points $A_1 = (B_1 + b_1 t)$, $A_2 = (B_2 + b_2 t)$, … where $t$ denoted time, and the other magnitudes were constant, will become almost uniformly distributed over a circumference. In 1881 Newcomb remarked that the first pages of logarithmic tables wore out "much faster" than the last ones and set out to derive the probability that the first significant digits of empirically obtained numbers will be $n_1$, $n_2$, … Without any proof he indicated that, if numbers $s_1$, $s_2$, …, $s_n$ were selected "at random", the positive fractional parts of the differences $(s_1 - s_2)$, $(s_2 - s_3)$, … will tend, as $n$ , to a uniform distribution over a circumference, and that the empirical magnitudes, to which these differences conform, will have equally probable mantissas of their logarithms. Newcomb's reasoning heuristically resembles the **Weyl** celebrated theorem that states that the terms of the sequence $\{nx\}$, where $x$ is irrational, $n = 1$, 2, …, and the braces mean "drop the integral part", are uniformly distributed on a unit interval. In the sense of the information theory, Newcomb's statement means that each empirical number tends to provide one and the same information. Several authors, independently one from another, proved that Newcomb was right. One of them (Raimi 1976, p. 536) called his statement an "inspired guess" and reasonably noted that it was not, however, universally valid.

By the mid-century, after processing observations made over about a century, a rough periodicity of the number of sunspots was established (cf. § 1.2.3). Newcomb (1901), who studied their observations from 1610 onward, arrived at $T = 11.13$ years. The present-day figure is $T \approx 11$ years but a strict periodicity is denied. In any case, it might be thought that the numbers of sunspots constitute a time series, an object for stochastic studies. I note that Newcomb considered the maxima and the minima of that phenomenon as well as half the sums of the numbers of the sunspots "corresponding to the year of minimum and the following maximum, or vice versa" (p. 4). He determined the four appropriate values of $T$ and their mean without commenting on the possible dependence between them.

The variation of the terrestrial latitudes is known to be caused by the movement of the pole about some point along a curve resembling a circumference with period 1.2 years. Newcomb (1892) checked the then proposed hypothesis that the movement was periodic with $T = 1.17$ years and he assumed that the pole moved uniformly along a circumference. Some of his calculations are doubtful (and not sufficiently detailed, a feature peculiar to many of his works), but he correctly concluded that the hypothesis was [apparently] valid.

In 1767 **Michell** (§ 6.1.6) attempted to determine the probability that two stars were close to each other. By applying the **Poisson** distribution, **Newcomb** (1859 – 1861, 1860, pp. 437 – 439) calculated the probability that some surface with a diameter of 1° contained $s$ stars out of $N$ scattered "at random" over the celestial sphere and much later **Fisher** (**Hald** 1998, pp. 73 – 74) turned his attention to that problem. **Boole** (1851/1952, p. 256) reasoned on the distinction between a uniform and any other random distribution:

*A 'random distribution' meaning thereby a distribution according to some law or manner, of the consequences of which we should be totally ignorant; so that it would appear to us as likely that a star should occupy one spot of the sky as another. Let us term any other principle of distribution an indicative one.*

His terminology is now unsatisfactory, but his statement shows that Michell's problem had indeed led to deliberations of a general kind. See also Newcomb (1904a) who thought about the subjective difference between those distributions.

Newcomb (1861b) also solved a related problem in which he determined the probability of the distance between the poles of two great circles randomly situated on a sphere. Issuing from other initial considerations, **Laplace** (1812, p. 261) and **Cournot** (1843, § 148) earlier provided solutions differing both from each other and from Newcomb's answer (Sheynin 1984a, pp. 166 – 167).

F. G. W. Struve (1827, pp. xxxvii – xxxix) determined the probability that two or three stars are situated near to each other. Bertrand (1888, pp. 170 – 171) noted however, that likely/unlikely relative situations of two stars can be considered in different ways (not only by distances between them) and concluded (pp. 4 – 7) that the Michell problem was unsolvable (that there are different ways to describe randomness, see § 11.1-1.

About 1784 **William Herschel** started counting the number of stars situated in different regions of the sky. He thought that his telescope was able to penetrate right up to the boundaries of the Milky Way and hoped to determine its configuration but later understood his mistake (§ 9C).

In one of its sections he (1784/1912, p. 158) counted the stars in six fields selected "promiscuously" and assumed the mean number of them as an estimate for the entire section. Much later Herschel (1817) proposed a model of a uniform spatial distribution of the stars. He fixed the boundaries for the distances of the stars of each magnitude but allowed the stars to be randomly distributed within these boundaries.

Herschel thus provided an example of randomness appearing alongside necessity; cf. **Poincaré's** statement in § 1.2.4. When estimating the precision of his model for the stars of the first seven magnitudes, Herschel calculated the sum of the deviations of his model from reality. For the first four magnitudes that sum was small although the separate deviations were large. Recall (§ 6.3.2) that, when adjusting observations, **Boscovich** applied a similar test with respect to absolute deviations. Herschel himself (1805) made use of it when determining the direction of the Sun's movement (cf. Note 17 to Chapter 6).

Herschel (1817/1912, p. 579) indicated that

*Any star promiscuously chosen* […] *out of* [14,000 stars of the first seven magnitudes] *is not likely to differ much from a certain mean size of them all*.

He certainly did not know that, with regard to size, the stars are incredibly different; its mean value is a worthless quantity, and, in general, stochastic statements, made in the absence of data, are hardly useful. A formal check in accordance with the **Bienaymé – Chebyshev** inequality would have revealed Herschel's mistake. But in any case it occurred that the stars, even earlier than the asteroids, had been considered as elements of a single population (in the last-mentioned instance, wrongly).

Stellar statistics really originated in the mid-19[th] century with the study of the proper motions of hundreds of stars (until 1842, when astronomers started to use the **Doppler's** invention, only in the directions perpendicular to the appropriate lines of sight).

Argelander (1837, p. 581) considered 560 stars with perceptible proper motions and determined the Sun's motion more reliably than it was achieved previously. Otto Struve (1842; 1844), then F. G. W. Struve (1852, pp. clxxxii – clxxxv) made the next steps. For the first of his studies the Royal Astronomical Society awarded O. Struve its gold medal (Airy 1842).

When studying the Sun's motion, astronomers beginning with Herschel thought that the peculiar motions of the stars (their motions relative to the Sun) were [random variables] and Kapteyn (1906a, p. 400) called the random distribution of the direction of the peculiar motions *a fundamental hypothesis*.

The calculated mean proper motions for stars of a given magnitude proved, however, almost meaningless since the magnitudes depended on distances. Beginning with **W. Herschel**, astronomers thought that the proper motions were random, but they understood randomness in different ways. **Newcomb** (1902a) assumed that their projections on an arbitrary axis were normally distributed. He derived, although without providing his calculations, the density laws of their projections on an arbitrary plane and their own distribution. Both these laws were connected with the $\chi^2$ distribution.

The general statistical study of the starry heaven became more important than a precise determination of the parameters of some star (Hill & Elkin 1884, p. 191):

*The great Cosmical questions to be answered are not so much what is the precise parallax of this or that particular star, but – What are the average parallaxes of those of the first, second, third and fourth magnitude respectively, compared with those of lesser magnitude?* [And] *What connection does there subsist between the parallax of a star and the amount and direction of its proper motion or can it be proved that there is no such connection or relation*?

Then, **Kapteyn** (1906b; 1909) described a stochastic picture of the stellar universe by the laws of distribution of the (random!) parameters, parallaxes and peculiar motions, of the stars. He (1906a) also initiated the study of the starry heaven by [stratified] sampling; here is a passage from a letter that he received in 1904 on this subject from one of his colleagues, Edward Pickering, and inserted on his p. 67:

*As in making a contour map, we might take the height of points at the corners of squares a hundred meters on a side, but we should also take the top of each hill, the bottom of each lake*, […], *and other distinctive points.*

In statistics, sampling became recognized at about the same time, although not without serious resistance (You Poh Seng 1951) and its most active partisan was Kiaer, also see § 10.7-2.

**Newcomb** (1902b, pp. 302 and 303) offered a correct estimate of Kapteyn's work:

*In recent times what we may regard as a new branch of astronomical science is being developed.* […] *This is what we now call the science of stellar statistics. The statistics of the stars may be said to have commenced with Herschel's gauges of the heavens* […]. *The outcome of Kapteyn's conclusions is that we are able to describe the universe as a single object* …

The compilation of vast numerical materials (catalogues, yearbooks) was also of a statistical nature. Moreover, sometimes this direction of work had been contrasted to theoretical constructions. Thus, **Proctor** (1872) plotted 324 thousand stars on his charts attempting to leave aside any theories on the structure of the stellar system, but the development of astronomy proved him wrong.

Calculation and adjustment of observations, their reasonable comparison has always been important for astronomy. Here, I again ought to mention, in the first place, **Newcomb**. Benjamin (1910) and many other commentators stated that he had to process more than 62 thousand observations of the Sun and the planets and that his work included a complete revision of the constants of astronomy. He necessarily discussed and compared observations obtained at the main observatories of the world but he hardly had any aids except for logarithmic tables. In addition he published some pertinent theoretical studies. He was of course unable to avoid the perennial problem of the deviating observations. At first he regarded them with suspicion, then (1895, p. 186), however, became more tolerant. If a series of observations did not obey the normal law, Newcomb (1896, p. 43) preferred to assign a smaller weight to the "remote" observations, or, in case of asymmetric series, to choose the median instead of the arithmetic mean. He had not mentioned **Cournot** (§ 10.3-3), and, in two memoirs published at the same time, he (1897a; 1897b) called the median by two (!) other, nowadays forgotten, terms.

**Mendeleev** (§ 10.9.3) objected to combining different summaries of observations; Newcomb, however, had to do it repeatedly, and in such cases he (1872) hardly managed without subjective considerations and assigned weights to individual astronomical catalogues depending on their systematic errors. Interestingly enough, he then repeated such adjustments with weights, depending on random errors.

After determining that the normal law cannot describe some astronomical observations necessarily made under changing conditions, Newcomb (1886) proposed for them (and, mistakenly, for all astronomical observations altogether) a generalized law, a mixture of normal laws with differing measures of precision occurring with certain probabilities. The measure of precision thus became a discrete random variable, and the parameters of the proposed density had to be selected subjectively. Newcomb noted that his density led to the choice of a generalized arithmetic mean with weights decreasing towards the "tails" of the variational series. I (§ 6.3.1) have remarked, however, that that mean was hardly better than the ordinary arithmetic mean.

He had also introduced some simplifications, and Hulme & Symms (1939, p. 644) noted that they led to the choice of the location parameter by the principle of maximum likelihood. Newcomb hardly knew that his mixture of normal laws was not normal (**Eddington** 1933, p. 277). In turn, two authors generalized Newcomb's law (Lehmann – Filhès 1887; K. F. Ogorodnikov, three English papers in 1928 – 1929), see Sheynin (1995c, pp. 179 – 182), but their work was of little practical importance.

Like Mendeleev (§ 10.9.3), Newcomb (1897b, p. 165) thought that the discrepancy between two empirical magnitudes was essential if it exceeded the sum of the two appropriate probable errors, and it seems that this rigid test had been widely accepted in natural sciences. Here is **Markov's** relevant pronouncement from a rare source (Sheynin 1990b; pp. 453 – 454): he

> *Like*[d] *very much Bredikhin's rule according to which 'in order to admit the reality of a computed quantity, it should at least twice numerically exceed its probable error'. I do* [he does] *not know, however, who established this rule or whether all experienced calculators recognized it.*

In other words, the difference between zero and a "real" non-zero magnitude must twice exceed its probable error, a statement that conformed to **Mendeleev's** and **Newcomb's** opinion. But still, Newcomb several times indicated that some magnitude *a* determined by him had mean square error *b* even when the latter much exceeded the former including the case (1901, p. 9) of $a = 0.05$ and $b = 0.92$! Dorsey & Eisenhart (1969, p. 50) recommended another rule in which the probable errors characterized separate observations.

Repeatedly applying the MLSq, Newcomb sometimes deviated from strict rules; see one such example in § 9A.5-3. In another case he (1895, p. 52) thought that small coefficients in a system of normal equations might be neglected, but he had not provided any quantitative test. Newcomb realized that, when forming normal equations, the propagation of round-off errors could result in their interdependence, and he reasonably concluded that in such cases the calculations should be made with twice as many significant digits. This is what he (1867) did when studying the calculations of the Kazan astronomer **Kowalski**, who had noted that, out of the four normal equations which he formed, only two were independent. It is now known that ill-conditioned observational equations should rather be processed without forming normal equations, – for example, by successive approximations.

Newcomb's calculation (1874, p. 167) presents a special case. Having 89 observational equations in five unknowns, he formed and solved the normal equations. Then, however, he calculated the residual free terms of the initial equations and somehow solved these equations anew (providing only the results of both solutions). He apparently wished to exclude systematic influences as much as possible, but how?

Newcomb (1895, p. 82; 1897b, p. 161) mistakenly stated, although mentioning earlier the definitive **Gaussian** justification of the MLSq, that the method was inseparable from the normal law. I note also his unfortunate reasoning (Newcomb & Holden 1874, p. 270) similar to the one made by Clausius (§ 10.8.5): for systematic error *s* and random errors $r_1$ and $r_2$, as he proved, and only for the normal law by considering the appropriate double integral, that

$$E[(s + r_1)(s + r_2)] = s^2.$$

It might be concluded that Newcomb necessarily remained more or less within the boundaries of the classical theory of errors and simple stochastic patterns. At the same time, the extant correspondence between him and **K. Pearson** during 1903 – 1907 (Sheynin 2002a, § 7.1) testifies that he wished to master the then originating mathematical statistics. Here is a passage from his archival letter of 1903 to Pearson (I have supplied the reference):

*You are the one living author whose production I nearly always read when I have time and can get at them, and with whom I hold imaginary interviews while I am reading.*

I mention finally Newcomb's statistical contribution (1904b) in which he examined the classical problem of the sex ratio at birth (see §§ 2.2.4, 3.3.4, 4.4 and 6.1.1). He assumed that there existed three kinds of families numbered, say, *m*, *n*, and *n*, for whom the probabilities of the birth of a boy were *p*, *p* +  and *p* –  respectively and he studied, in the first place, the births of twins. The sex of the embryo, as he thought, became established only after the action of a number of successive causes made it ever more probable in either sense.

**10.8.5. Physics.** 1)The kinetic theory of gases originated in mid-19$^{th}$ century as the result of the penetration of the statistical method into physics. Truesdell (1975) discussed its early history; thus (p. 28), it was Waterson who, in 1843, introduced the mean free path of a molecule, cf. § 8.7, but his innovation was not published. **Clausius** likely published the first memoir (1849) which belonged to physics (but did not deal with the molecular hypothesis) and contained ideas and methods of the theory of probability.

After **Poisson's** death that theory sank into oblivion (§§ 0.1 and 7.4). No wonder that Clausius (1889 – 1891, p. 71) made a point to prove the equality $E( /E ) = 1$ for the velocity  of a molecule. See a similar case in § 10.8.4.

**Maxwell** twice mentioned Laplace (Sheynin 1985, pp. 364 and 366n), although without providing any definite references, whereas Boltzmann, who cited many scholars and philosophers in his popular writings, never recalled him. **Khinchin** (1943/1949, p. 2) maintained that Maxwell and Boltzmann applied

*Fairly vague and somewhat timid probabilistic arguments* that *do not pretend here to be the fundamental basis, and play approximately the same role as purely mechanical considerations. […] Far reaching hypotheses are made concerning the structure and the laws of interaction between the particles […]. The notions of the theory of probability do not appear in a precise form and are not free from a certain amount of confusion which often discredits the mathematical arguments by making them either devoid of any content or even definitely incorrect. The limit theorems […] do not find any applications […]. The mathematical level of al these investigations is quite low, and the most important mathematical problems which are encountered in this new domain of application do not yet appear in a precise form.*

His statement seems too harsh, written from the standpoint of statistical mechanics of the mid-20$^{th}$ century. Then, I believe that it was partly occasioned by Boltzmann's verbose style of writing. Third, physicists certainly applied the LLN indirectly. Fourth, Khinchin said nothing about positive results achieved in physics (formulation of the ergodic hypothesis, use of infinite general populations, Maxwell's indirect reasoning about randomness). My first remark is indeed essential; here is an extract from Maxwell's letter of 1873 (Knott 1911, p. 114):

*By the study of Boltzmann I have been unable to understand him. He could not understand me on account of my shortness, and his length was and is an equal stumbling block to me.*

And Boltzmann (1868/1909, p. 49) indeed owned that it was difficult to understand Maxwell's "Deduktion" (1867) "because of its extreme brevity". I emphasize that statistical mechanics could not have appeared unless and until the kinetic theory with its mathematical shortcomings had been established (Truesdell 1975; Brush 1976).

2) Clausius. He (1857/1867, pp. 238 and 248) asserted that molecules moved with essentially differing velocities. Even **Boscovich** (1758, § 481) stated something similar but perhaps presumed that the differences between these velocities were not large: The "points" [atoms] of "a particle" [of light, as in § 477, or of any body, as in § 478] move "together with practically the same velocity", and the entire particle will "move as a whole with the single motion that is induced by the sum [the mean] of the inequalities pertaining to all its points".

Clausius used a single mean velocity such as to make the entire kinetic energy of a gas equal to its actual value. Later he (1862/1867, p. 320) maintained that the velocities of molecules randomly differed one from another.

And he (1858/1867, p. 268) studied the length of the free path of a molecule. Denote the probability of a unit free path by $a$, then

$$W = a^x = (e^{-x})^\alpha, \alpha > 0$$

will be the probability of its being equal to $x$; here, $\alpha$ is derived from the molecular constants of the substance. Similar considerations are in other works of Clausius (1862/1867, § 29; 1889 – 1891, pp. 70 – 71 and 119).

He (1889 – 1891, pp. 70 – 71) also calculated the mean free path of a molecule. Actually, without writing it out, he considered free paths of random length $\xi$ and calculated the expected free path as an integral over all of its possible values from 0 to .

Suppose now that

$$\xi = \xi_1 + \xi_2 + \ldots + \xi_m$$

where $m$ is an arbitrary natural number. Then, according to Clausius' assumptions, $\xi_k$, $k = 1, 2, \ldots, m$, will not depend on
$(\xi_1 + \xi_2 + \ldots + \xi_{k-1})$ and the characteristic function for $\xi_k$ will be equal to the product of these functions for the previous $\xi$'s. In this instance, all these functions are identical, and *F(s)*, the integral distribution function of ξ, is therefore infinitely divisible.

Distribution functions had first appeared (or were easy to be derived from the correspondence of) Huygens (§ 2.2.2), from De Moivre (§ 4.2) and Davidov (1885) and directly in Poisson (§ 8.2).

Clausius' achievements were interesting, but he did not attempt to construct the kinetic theory of gases on a stochastic basis. Nevertheless, his role in this direction, at least from the viewpoint of probability theory, had not yet ben studied. See Schneider (1974) who reviewed his concrete physical research by stochastic considerations. His great merit is seen in Maxwell's statement (1875/1890, p. 427):

*Clausius opened up a new field […] by showing how to deal mathematically with moving systems of innumerable molecules.*

3) **Maxwell** (1860) established his celebrated distribution of the velocities of monatomic molecules

$$\varphi(x) = \frac{1}{\sqrt{\pi}} \exp(-x^2/\alpha^2).$$

He tacitly assumed that the components of the velocity were independent; later this restriction was weakened (§ 9A.1.3). He then maintained that the average number of particles with velocities within the interval [*v*; *v* + *dv*] was proportional to

$$f(x) = \frac{4}{3\sqrt{\pi}} v^2 \exp(-v^2/\alpha^2) dv.$$

This can be justified by noting that the probability of such velocities can also be represented as

$$\int_0^{2\pi} d\varphi \int_0^{\pi} \sin\vartheta \, d\vartheta \int_v^{v+dv} t^2 \exp(-t^2/\alpha^2) dt.$$

It is presumed here that the components of the velocity in each of the three dimensions have the same distribution. Maxwell left interesting statements about the statistical method in general, and here is one of them (1873b/1890, p. 374):

*We meet with a new kind of regularity, the regularity of averages, which we can depend upon quite sufficiently for all practical purposes, but which can make no claim to that character of absolute precision which belongs to the laws of abstract dynamics.*

The drafts of the source just mentioned (Maxwell 1990 – 2002, 1995, pp. 922 – 933) include a previously unpublished and very interesting statement (p. 930): abandoning the "strict dynamical method" and adopting instead the statistical method "is a step the philosophical importance of which cannot be overestimated".

And here is his definition (not quite formal) of the statistical method which heuristically resembles the formulation provided by **Kolmogorov** & **Prokhorov** (§ 0.2): it consisted in "estimating the average condition of a group of atoms" (1871/1890, p. 253), in studying "the probable [not the average!] number of bodies in each group" under investigation (1877, p. 242).

Maxwell gave indirect thought to randomness. Here is his first pronouncement (Maxwell 1859/1890, vol. 1, pp. 295 – 296) which was contained in his manuscript of 1856 (1990 – 2002, 1990, p. 445), and it certainly describes his opinion about that phenomenon:

*There is a very general and very important problem in Dynamics* […]. *It is this – Having found a particular solution of the equations of motion of any material system, to determine whether a slight disturbance of the motion indicated by the solution would cause a small periodic variation, or a derangement of the motion* […].

Maxwell (1873a, p. 13) later noted that in some cases "a small initial variation may produce a very great change […]". Elsewhere he (report read 1873, see Campbell & Garnett 1882/1969, p. 440) explained that in such instances the condition of the system was unstable and prediction of future events becomes impossible. He (Ibidem, p. 442) provided an example of instability of a ray within a biaxial crystal and prophetically stated (p. 444) that in future physicists will study "singularities and instabilities".

In a manuscript of the same year, 1873 (p. 360), Maxwell remarked that

*The form and dimensions of the orbits of the planets* […] *are not determined by any law of nature, but depend upon a particular collocation of matter. The same is the case with respect to the size of the earth*.

This was an example illustrating **Poincaré's** statement concerning randomness and necessity (§ 1.2.4), but I ought to add that it was not sufficiently specific; the eccentricities of planetary orbits depend on the velocities of the planets, cf. end of § 7.3.

And here is Maxwell's position (1875/1890, p. 436) concerning randomness in the atomic world:

*The peculiarity of the motion of heat is that it is perfectly irregular*; […] *the direction and magnitude of the velocity of a molecule at a given time cannot be expressed as depending on the present position of the molecule and the time*.

At the very end of his life Maxwell (1879/1890, pp. 715 and 721) introduced a definition for the probability of a certain state of a system of material particles:

*I have found it convenient, instead of considering one system of* […] *particles, to consider a large number of systems similar to each other*

*[…]. In the statistical investigation of the motion, we confine our attention to the number of these systems which at a given time are in a phase such that the variables which define it lie within given limits.*

*Boltzmann (1868, § 3) defines the probability of the system being in a phase […] as the ratio of the aggregate time during which it is in that phase to the whole time of the motion […].*

4) If the classical definition of probability is included here, we can say that Boltzmann used three formulations. Maxwell (§ 10.8.5-2) mentioned one of them, and another reference can be added: Boltzmann (1868/1909, Bd. 1, p. 50). Still another was that applied by Maxwell (see same subsection) although sometimes Boltzmann (1878/1909, p. 252) did not indicate which one he was employing. He (1872/1909, p. 317) apparently thought that these posterior probabilities were equivalent.

In other words, with respect to separate molecules Boltzmann introduced the time average probability, – and maintained that it was equivalent to the "usual" phase average probability, also see § 12-2. When studying polyatomic gases, Boltzmann (1871) defined the probability of its state as a product such as $fd\tau$ where $f$ was some function, varying in time, of the coordinates and velocities of the separate molecules and $d\tau$, the product of the differentials of those parameters. For stochastic processes, such functions determine the distribution of a system of random variables at the appropriate moment. Zermelo (1900, p. 318) and then Langevin (1913/1914, p. 3) independently stressed the demand to provide a "définition correcte et claire de la probabilité" (Langevin).

Like Maxwell, Boltzmann (1887/1909, p. 264; 1899, Bd. 2, p. 144) used the concepts of fictitious physical systems and infinite general population.

**Boltzmann** (1896/1909, p. 570) stated that the [normal law] followed from equal probabilities of positive and negative elementary errors of the same absolute value. His was of course an unworthy formulation of the CLT.

I ought to add that Boltzmann respected the theory of probability. Thus (1872/1909, p. 317)

*An incompletely proved theorem whose correctness is questionable should not be confused with completely proved propositions of the theory of probability. Like the results of any other calculus, the latter show necessary inferences made from some premises.*

And again (1895/1909, p. 540): the theory of probability "is as exact as any other mathematical theory" if, however, "the concept of equal probabilities, which cannot be determined from the other fundamental notions, is assumed".

From 1871 onward **Boltzmann** had been connecting the proof of the second law of thermodynamics with stochastic considerations. Thus, he (1872/1909, pp. 316 – 317) declared that the problems of the mechanical theory of heat are also problems of probability theory.

Then, however, he (1886/1905, p. 28) indicated that the 19$^{th}$ century will be the age of "mechanical perception of nature, the age of **Darwin**", and (1904a/1905, p. 368) that the theory of evolution was understandable in mechanical terms, that (1904b, p. 136) it will perhaps become possible to describe electricity and heat mechanically.

The possible reason for his viewpoint was that he did not recognize objective randomness. Another reason valid for any scholar was of course the wish to keep to "abstract dynamics", see Maxwell's statement on the "new kind of regularity" (§ 10.8.5-3) and the opinion of Hertz (1894, Vorwort): "Physicists are unanimous in that the aim of physics is to reduce the phenomena of nature to the simple laws of mechanics". And here is a lucid description of this point as far as Boltzmann was considered (Rubanovsky 1934, p. 6): in his works

*Randomness […] struggles with mechanics. Mechanical philosophy is still able […] to overcome randomness and wins a Pyrrhic victory over it but recedes undergoing a complete ideological retreat.*

### 10.9. Natural scientists

**10.9.1. Ivory.** In a letter to **Olbers** of 1827, **Gauss** (*W*/Erg-4, No. 2, pp. 475 – 476) called Ivory an "acute" mathematician, but indicated that the "spirit" of the MLSq was alien to him. In 1825 – 1830 Ivory published 11 papers in one and the same periodical [the last of these was Ivory (1830)] devoted to the derivation of the flattening of the Earth's ellipsoid of rotation by means of pendulum observations. It is not amiss to add that his main contributions pertained to the theory of the figure of the Earth (attraction of a material point by an ellipse).

In accordance with the **Clairaut** theorem, the Earth's flattening (see Note 10 to Chapter 6) is determined by two observations of [the acceleration of] gravity at different latitudes; however, unavoidable errors and local irregularities in the figure of the Earth necessitate the use of redundant observations. See the same remark in § 6.3. To strengthen **Gauss'** remark, I state that Ivory was simply ignorant of the MLSq and without justification called it not good enough. He denied it in words but applied the MLSq, perhaps not even realizing it at once. Thus, starting from equations of the type (1.1) with $a_i = 1$ he (1826b, pp. 244 – 245) stated that the condition $v_i = 0$, unlike the requirement of the MLSq $[av] = [bv] = … = 0$, see equations (1.5), was expedient. He failed to notice that in his case the *expedient* condition coincided with the demand that $[av] = 0$.

Then, having at his disposal 5 – 7 observations, only one of which was made at a southern station, he (1826a, p. 9) combined it with each of the others (so as to have pairs with a large latitudinal difference between stations) and calculated the flattening from the thus obtained pairs. The weight of the equatorial observation became absurdly great and its error corrupted all the pairs in the same way. An utterly unworthy manner of treatment, as Gauss stated but it seems that the answer to this problem depends on the unknown magnitude of the systematic errors.

Only later did Ivory remark that the local anomalies of gravity can essentially influence the end result, – and rejected a large part of the available observations, – up to 31%, see Ivory (1826b, p. 242), – and even began doubting whether it was possible to derive a single flattening. Local anomalies are indeed extremely troublesome (also see § 10.8.3 where I indirectly mentioned local perturbations of temperature) but Ivory attempted to get rid of them too radically. Finally, when estimating the precision of his results, he had not applied the variance.

Before adjusting pendulum observations it is possible to replace stations having almost the same latitude by one fictitious mean station (which Ivory had not done). If, however, the longitudes of such stations are also almost the same, all of them could have been corrupted by the same local gravimetric anomaly so that the weight of the mean station should not be increased as compared with either of them.

I ought to add two remarks. First, his final result (1828, p. 242) was sufficiently close to the flattening of the **Krasovsky** ellipsoid of 1940 (Sakatov 1950, p. 364): $e = 0.00333 – 0.00338$ and $0.00335$, respectively. In addition, Ivory (1825, p. 7), without, however, mentioning **Gauss**, maintained that the MLSq should be substantiated by the principle of maximum weight. Second, Ivory actually wished to solve two problems at once: to find out whether the observations were consistent with an ellipsoidal Earth, and to adjust them. It is the minimax method (§ 6.3.2) rather than the MLSq that is best for solving the first problem.

**10.9.2 Fechner.** He (1860) was the founder of psychophysics and therefore became one of the first to introduce the statistical method into physics, although not in the crucial direction. He (1860, Bd. 1, p. 8, see also 1877, p. 213) defined it as an "exact doctrine on the functional correspondence or interdependence of body and soul".

And here is a modern down to earth definition (*New Enc. Brit.*, 15th ed., vol. 9, 1997, p. 766): psychophysics is a

*Study of quantitative relations between psychological events and physical events, or, more specifically, between sensations and the stimuli that produce them.*

There also Fechner's book just mentioned is called a classic.

Fechner (1855 and 1864) missed the opportunity to comment on the developing kinetic theory of gases. Moreover, he (1874b, pp. 7 and 9; 1897, p. 15) repeatedly treated physics on a par with (practical) astronomy by stating that both these branches of natural sciences had to do with symmetric distributions and true values of magnitudes sought. His mathematical tools and the very approach were primitive and almost everything he achieved had to be repeated at a much higher level. Ebbinghaus (1908, p. 11) possibly having in mind Fechner's non-scientific writings called him "a philosopher full of fantasies" and, at the same time, "a most strict physicist" who had "put […] together psychophysics as a new branch of knowledge". As chance would have it, Fechner published those writings under the pen name *Dr Mises*.

Fechner himself (1877, p. 215) figuratively estimated his own work (and made known the existence of opposition to it):

*The Tower of Babel was not completed because the workers were unable to explain to each other how should they build it. My psychophysical structure will probably survive because the workers cannot see how they might demolish it.*

Being the co-author of the logarithmic **Weber – Fechner** law connecting stimuli with sensations, Fechner extended the range of its application after making a great number of experiments (1860; 1887). He studied the methods of experimentation, and made reasonable pertinent statements. The modern method of paired comparisons (H. A. David 1963) owes much to him.

In the theory of errors Fechner had been mentioning **Gauss**, but he also attempted, sometimes unsuccessfully, to introduce his own innovations, or to repeat unknown to him previous findings but he nevertheless somewhat furthered that theory. Thus, issuing from elementary but apparently non-rigorous considerations, he (1874a, p. 74) provided a formula for estimating the precision of observations which coincided with the Peters formula (9B.11) but was applicable to any distributions. Then, proceeding from the Gaussian formulas (§ 9A.3), he compared two competing expressions connecting the magnitudes of the stars with their brightness, solved redundant systems of equations by the method of pairwise combinations (§ 6.3.2), and remarked, without substantiation (and hardly correctly), that that method asymptotically tended to the MLSq (1887, p. 217).

Fechner's main innovation was, however, the collective, – actually, the set of observed values of a random variable. He (1897) proposed to study collectives by applying several mean values, their mutual arrangement, and their deviations (both absolute and normed) from the observations. He mostly paid attention to asymmetric collectives and even attempted to discover a universal asymmetric distribution for errors in natural sciences (cf. § 10.5). Fechner especially examined the double normal law (two different normal laws for the smaller and the larger values of observations in the variational series respectively, turning into one another at the point of maximal probability, i.e., at the mode), and the double lognormal law. Fechner also attempted, although not very successfully, to separate the real and the apparent (caused by an insufficient number of observations) asymmetry.

Finally, Fechner (1897, pp. 365 – 366) studied the dependence of the successive daily air temperatures by comparing their course with the arrangement of winning (numbered) tickets of a reputed German lottery. When examining the results of the lottery, he achieved an interesting result pertaining to the runs up and down (cf. § 10.2-5). Furthermore, Fechner even introduced a measure of dependence which varied from 0 to 1, but only described "positive" dependences. His contribution appeared posthumously, after the **Galton** correlation theory had emerged.

**Mises** (1928/1972, pp. 26 and 99) highly appraised Fechner's efforts and owned (p. 99) that Fechner's "constructions prompted, at least me [Mises], to adopt a new viewpoint". Two more passages by **K. Pearson** (1905, p. 189) and **Freud** (1925/1963, p. 86) are in order:

*All the leading statisticians from Poisson to Quetelet, Galton, Edgeworth and Fechner […] have realized that asymmetry must be in some way described before we can advance in our theory of variation* [in biology].

*I was always open to the ideas of Fechner and have followed that thinker upon many important points.*

**10.9.3. Mendeleev.** From 1893 to 1907 Mendeleev was Director of Russia's Main Board of Measures and Weights which was established on his initiative on the basis of its modest predecessor (where, from 1892, Mendeleev was scientific curator). Actually, he was an extremely versatile scientist, studied statistics of population and industry and thought that the range of statistics is unrestricted.

Mendeleev processed observations both as a chemist and a metrologist. He (1872b/1951, p. 101) distrusted data "obtained under differing conditions, by different methods and observers" as compared with those "achieved by precise methods and experienced persons". And (1887/1934, p. 82) "disadvantageous" data ought to be rejected "by a clear critical appraisal", otherwise "a realistic result" is impossible to get.

No wonder that Mendeleev (1872a/1939, p. 144) preferred "to make a few but precise and repeated measurements" and objected to amassing observations; true, this attitude was partly due to his wish to avoid calculations, cf. **Boyle's** statement in § 1.2.2. Mendeleev effectively repeated this reasoning (1875a/1939, p. 256) and, for example in the first case (an attempt at refining the **Boyle – Mariotte** law), he added that the measurements ought to be made at "significantly different pressures". He was apparently concerned with systematic errors which could have been otherwise almost the same.

Mendeleev (1875b/1950, p. 209) thought that an observational series should be "harmonious", that is, that its median should coincide with its arithmetic mean, or (his second definition) that the mean of its middlemost third should coincide with the mean of the means of its extreme thirds. In the first case, he mistakenly added that the coincidence meant that the appropriate distribution was normal. He had not said how to treat observations which did not obey his wish; cf. end of § 10.8.3.

The deviation of the arithmetic mean from the median, normed in a certain way, is nowadays recognized as a measure of asymmetry of the appropriate distribution (**Yule** & **Kendall** 1958, p. 161). For Mendeleev, the probable error was the main estimator of precision and he (1860/1947, p. 46) assumed that the admissible deviation between two means was the sum of their probable errors (cf. § 10.8.4). Suppose that these errors are equal to each other; then, for the normal distribution, their sum is 1.35 where is the standard deviation (or the mean square error) of each mean. On the other hand, the standard deviation of the difference between the means is 2 and it thus occurs that the studied difference is essential when it is equal to its standard deviation. Mendeleev's (or **Bredikhin's**, or **Newcomb's**) rule seems to be too rigid. A different rule was recommended more recently (Dorsey & Eisenhart 1969, p. 50) in which the probable errors of individual measurements were involved instead.

Mendeleev had not mentioned the second **Gaussian** justification of the MLSq and made a few mistakes in his theoretical considerations. One of them was an excessive belief in the arithmetic mean (1856/1937, p. 181; 1877/1949, p. 156; 1895/1950, p. 159): he thought that that estimator ought to be chosen if nothing was known about the precision of the individual results; cf. § 7.2-6.

### Notes

**1.** He published many very short notes and insufficiently described his findings, sometimes, like in this case, without proof.

**2. Cournot** only explained his understanding of the distinction between chance and probability in § 48 and not clearly enough, see his §§ 12 and 240/3. True, in his Préface he published **Poisson's** letter of 1836 where its author had indicated that with regard to that point they were unanimous, cf. § 8.1.

**3.** No one apparently recalled it before **Cournot**; and only **Chuprov** (1909/1959, p. 188) mentioned it afterwards.

**4.** He read his work in 1846, but the existing materials testify that already then he could have known **Buniakovsky's** book, and exactly him did **Ostrogradsky** apparently criticize (see next sentence). He barely busied himself with probability, but he (1858) made the calculations necessary for the work of a society of mutual insurance. In § 7.1-9 I mentioned his attempt to generalize the notion of moral expectation. On Ostrogradsky see **Gnedenko** (1951).

**5.** Over the years, **Faraday** (1991 – 2008) several times expressed his high opinion about Quetelet's measurements of atmospheric electricity, and I especially note two of his letters (to Quetelet, No. 1367 of 1841, 1996, p. 42, and to Richard Taylor, the then Editor of the *Lond., Edinb. and Dublin Phil. Mag.*, No. 2263 of 1850, 1999, p. 270):

1) *You are indeed a worthy example in activity & power to all workers in science and if I cannot imitate your example I can at least appreciate & value it.*

2) Faraday approvingly remarked on the absence of "imagination or hypothesis" in Quetelet's work and added: *Such was the true method by which advances in science in this very difficult part could be really made.*

See also the opinion of **Köppen** about the meteorological observations in Belgium at the very end of § 10.8.3.

**6.** Achenwall (Schlözer 1804, § 3) declared that statistics of a given state is a collection of its remarkable features and it is unclear how Delambre's differing understanding of statistics related to political arithmetic. One point is however clear: the refusal to study causes and effects was quite foreign to Petty and Graunt.

**7.** The Paris statistical society was established a bit earlier than 1804 (Schlözer 1804, § 7). It attempted to collect and treat data on many diverse matters, in the first place concerning France.

**8.** The *Handbook of social indicators* published in 1989 by the UN, see W. F. M. De Vries (2001), listed several hundred indicators separated into 13 groups. They help to trace the range of problems of modern statistics. Several papers on the newest goals of statistics in the "information society" are collected in the *International Statistical Review*, vol. 71, No. 1, 2003. It is opportune to add that Leibniz (Sheynin 1977, pp. 222 – 227) had interesting thoughts about the aim of statistics.

**9.** It had been applied in England from the $12^{th}$ century onward for checking the quality of batches of new coins (Stigler 1977). Ptukha (1961) described its usage in Russia from the $17^{th}$ century. In the beginning of the $18^{th}$ century Marshall Vauban estimated France's agricultural production by sampling, but apparently unsuccessfully (Moreau de Jonnès 1847, pp. 53 – 54).

Lagrange (1796), see Pearson (1978, pp. 628 – 635), published an essay on political arithmetic in which he had to use sampling.

**10.** Private communication (2003) by Prof. Walter Mann, a grandson of **Mendel's** nephew, Alois Schindler, and a typographic text of the latter's manuscript (of his report of 1902). That manuscript was, however, published (Krizenecky 1965, pp. 77 – 100).

**11.** Still earlier the problem of allowing for local anomalies presented itself when pendulum observations began to be used for determining the flattening of the Earth's ellipsoid of rotation, see § 10.9.1.

## 11. Bertrand and Poincaré

Passing now to **Bertrand**, I disturb the chronology of description, but not its logic: he was not interested in the work of **Chebyshev**. This is also true with regard to **Poincaré** who never mentioned Chebyshev's followers in probability (**Markov** and **Liapunov**) either.

### 11.1. Bertrand

In 1855 Bertrand had translated **Gauss'** works on the MLSq into French[1], but his own work on probability began in essence in 1887 – 1888 when he published 25 notes in one and the same periodical as well as his main treatise (1888a), written in great haste and carelessly, but in a very good literary style. I take up its main issues and state right now that it lacks a systematic description of its subject.

1) "Uniform" randomness. By several examples Bertrand proved that the expression "at random", or even "uniformly" random, was not definite enough. Thus, he maintained that the **Michell** problem (§ 6.1.6) should have been generalized: remarkable was not only a small distance between stars, but some other features of their mutual arrangement as well. One of his examples (p. 4) became classical. Determine the probability, Bertrand asked, that a randomly drawn chord of a given circle was longer than the side of an equilateral triangle inscribed in the circle. He listed three possible answers:

a) One endpoint of the chord is fixed; $p = 1/3$.
b) The chord's direction is fixed; $p = 1/2$.
c) The location of the centre of the chord in any point of the circle is equally probable; $p = 1/4$.

A curious statement about this problem is due to **Darboux** (1902/1912, p. 50):

*In accord with considerations which seem equally plausible, he* [Bertrand] *derived two different values for the probability sought, 1/2 and 1/3. He investigated this question and found its solution, but left its discovery to the readers.*

In failing to mention the third solution he possibly followed **Poincaré** (§ 12-4). I return to this problem in Chapter 12.

2) Statistical probability and the **Bayesian** approach. Heads appeared $m = 500{,}391$ times in $n = 10^6$ tosses of a coin (p. 276). The statistical probability of that event is $p = 0.500391$; it is unreliable, not a single of its digits merits confidence. After making this astonishing declaration, Bertrand compared the probabilities of two hypotheses, namely, that the probability was either $p_1 = 0.500391$, or $p_2 = 0.499609$. However, instead of calculating

$$[p_1^m p_2^n] \div [p_2^m p_1^n],$$

he applied the **De Moivre – Laplace** theorem and only indicated that the first probability was 3.4 times higher than the second one. So what should have the reader thought? And his choice of the two probabilities was unfortunate.

As I understand him, Bertrand (p. 161) "condemned" the **Bayes** "principle" only because the probability of the repetition of the occurrence of an event after it had happened once was too high (cf. the problem about the sunrise in § 5.1). This conclusion was too hasty, and the reader was again left in suspense: what might be proposed instead? Note that Bertrand (p. 151) mistakenly thought that the De Moivre – Laplace theorem precisely described the inverse problem, the estimation of the theoretical probability given the statistical data.

3) Statistics of population. Bertrand indicated that there existed a dependence between trials (or their series) and that the probabilities of the studied events could change. He referred only to **Dormoy** and had not provided any concrete examples, but he (p. 312) noted that, when studying the sex ratio at birth, both **Laplace** and **Poisson** had assumed without justification that the probability of a male birth was constant in time and space. Yes, but their mistake was only methodological since they could not have failed to understand this circumstance (as I mentioned in § 7.1-5).

**Bortkiewicz** (1930, p. 53) concluded that Dormoy was much less important than **Lexis**; see however end of § 15.1.1.

4) Mathematical treatment of observations. Bertrand paid much attention to this issue, but his reasoning was amateurish and sometimes wrong. Even if, when translating **Gauss** (see above), he had grasped the essence of the MLSq, he barely remembered that subject after more than 30 years. Thus, he (pp. 281 – 282) attempted to prove that the sample variance (9.6b) might be replaced by another estimator of precision having a smaller variance. He failed to notice, however, that, unlike the Gauss statistic, his new estimator was biased. Furthermore, when providing an example, Bertrand calculated the variance of (9.6b) for the case of the normal distribution instead of applying the Gauss formula (9.6c).

At the same time Bertrand formulated some sensible remarks. He (p. 248) expressed a favourable opinion about the second Gauss justification of the MLSq but indicated (p. 267) that, for small errors, the even distribution

$$(x) = a + bx^2$$

could be approximately represented by an exponential function of a negative square, – that the first substantiation of the method was also approximately valid. Keynes (1921/2014, p. 240) stated the same Finally, Bertrand provided an argument against the postulate of the arithmetic mean, see § 9A.2-2.

5) Several interesting problems in addition to that described in § 11.1-1 dwell on a random composition of balls in an urn; on sampling without replacement; on the ballot problem; and on the gambler's ruin.

a) White and black balls are placed in the urn with equal probabilities and there are $N$ balls in all. A sample made with replacement contained $m$ white balls and $n$ black ones. Determine the most probable composition of the urn (pp. 152 – 153). Bertrand calculated the maximal value of the product of the probabilities of the sample and of the hypotheses on the composition of the urn.

b) An urn has *sp* white balls and *sq* black ones, $p + q = 1$. Determine the probability that after *n* drawings without replacement the sample will contain ($np – k$) white balls (p. 94). Bertrand solved this problem applying the [hypergeometric distribution] and obtained, for large values of *s* and *n*, an elegant formula

$$P = \frac{1}{\sqrt{2\pi pqn}} \sqrt{s/(s-n)} \exp[-k^2 s/2pqn(s-n)].$$

He (1887b) published this formula earlier without justification and noted that the variable probability of extracting the balls of either colour was "en quelque sorte un régulateur".

c) Candidates *A* and *B* scored *m* and *n* votes respectively, $m > n$ and all the possible chronologically differing voting records were equally probable. Determine the probability *P* that, during the balloting, *A* was always ahead of *B* (p. 18). Following André (1887), who provided a simple demonstration, Bertrand proved that

$$P = (m - n)/(m + n), \qquad (1)$$

see also **Feller** (1950/1968, § 1 of Chapter 3). Actually, Bertrand (1887a) was the first to derive formula (1) by a partial difference equation. This *ballot problem* has many applications (Feller, Ibidem). Takácz (1982) traced its history back to **De Moivre** (§ 4.1-5). He indicated that it was extended to include the case of $m \geq \mu n$ for positive integral values of μ and that he himself, in 1960, had further generalized that extended version.

d) I select one out of the few problems on the gambler's ruin discussed by Bertrand (pp. 122 – 123). Gambler *A* has *m* counters and plays with an infinitely rich partner. His probability of winning any given game is *p*. Determine the probability that he will be ruined in exactly *n* games ($n > m$). Bertrand was able to solve this problem by applying formula (1). Calculate the probability that *A* loses $(n + m)/2$ games and wins $(n – m)/2$ games; then, multiply it by the probability that during that time *A* will never have more than *m* counters, that is, by *m/n*.

In a brief chapter he largely denied almost everything done in the "moral applications" of probability by **Condorcet** (but did not refer to **Laplace** or **Poisson**).

In two of his notes Bertrand (1888b; 1888c) came close to proving that for a sample from a normal population the mean and the variance were independent. Heyde & Seneta (1977, p. 67n) indicated this fact with respect to Bertrand's second note; see §§ 7.2-5 and 9B-5 for the previous findings of Laplace and **Helmert**.

Taken as a whole, Bertrand's treatise is impregnated with its non-constructive negative (and often unjustified) attitude towards the theory of probability and treatment of observations. And at least once he (pp. 325 – 326) wrongly alleged that **Cournot** had supposed that judges decided their cases independently one from another, see § 10.3-

6. I ought to add, however, that Bertrand exerted a strong (perhaps too strong) influence upon **Poincaré**, and, its spirit and inattention to Laplace and **Bienaymé** notwithstanding, on the revival of the interest of French scientists in probability (Bru & Jongmans 2001).

## 11.2 Poincaré

In the theory of probability, Poincaré is known for his treatise (1896); I refer to its extended edition of 1912. I note first of all that he had passed over in silence not only Russian mathematicians, but even **Laplace** and **Poisson**, and that his exposition was imperfect. Commenting on the first edition of his treatise, **Bortkiewicz** (Bortkevich & Chuprov 2005, Letter 19 of 1897) noted:

*The excessively respectful attitude towards […] Bertrand is surprising. No traces of a special acquaintance with the literature on probability are seen. The course is written in such a way as though Laplace and Poisson, especially the latter, never lived.*

Following **Bertrand**, Poincaré (p. 62) called the expectation of a random variable its probable value; denoted the measure of precision of the normal law either by $h$ or by $h$; made use of loose expressions such as "$z$ lies between $z$ and $z + dz$" (p. 252).

Several times Poincaré applied the formula

$$\lim \frac{\int \varphi(x) \Phi^n(x) dx}{\int \psi(x) \Phi^n(x) dx} = \frac{\varphi(x_0)}{\psi(x_0)}, n \to \infty \qquad (2)$$

where $(x)$ was a restricted positive function, $x_0$, the only point of its maximum, and the limits of integration could have been infinite (although only as the result of a formal application of the **Bayesian** approach). Poincaré (p. 178) only traced the proof of (2) and, for being true, some restrictions should perhaps be added. To place Poincarè's trick in the proper perspective, see Erdélyi (1956, pp. 56 – 57). I discuss now some separate issues mostly from Poincaré's treatise.

1) The theory of probability. Poincaré (p. 24) reasonably stated that a satisfactory definition of prior probability was impossible, cf. § 7.4. Strangely enough, he (1902/1923, p. 217) declared that "all the sciences" were nothing but an "unconscious application" of the calculus of probability, that the theory of errors and the kinetic theory of gases were based on the LLN and that the calculus of probability "will evidently ruin them" (*les entrainerait évidemment dans sa ruine*). Therefore, as he concluded, the calculus was only of practical importance[2]. Another strange pronouncement is in his treatise (1896, p. 34). As I understand him, he maintained that a mathematician is unable to understand why forecasts concerning mortality figures come true.

In a letter of ca. 1899 partly read out at the hearing of the notorious **Dreyfus** case (*Le procès* 1900, t. 3, p. 325; Sheynin 1991a, pp. 166 – 167) Poincaré followed **Mill** (§ 8.9.1), even generalized him to include all "moral sciences" and declared that the appropriate findings made by **Condorcet** and **Laplace** were senseless. And he objected to a stochastic study of handwriting for identifying the author of a certain document.

The interest in application of probability to jurisprudence is now revived. Heyde & Seneta (1977, p. 34) had cited several pertinent sources published up to 1975; to these I am adding Zabell (1988a), Gastwirth (2000) and Dawid (2005) who emphasized the utmost importance of interpreting background information concerning stochastic reasoning, cf. § 3.1.2.

2) Poincaré (1892a) had published a treatise on thermodynamics which **Tait** (1892) criticized for his failure to indicate the statistical nature of this discipline. A discussion followed in which Poincaré (1892b) stated that the statistical basis of thermodynamics did not satisfy him since he wished to remain "entirely beyond all the molecular hypotheses however ingenious they might be"; in particular, he therefore passed the kinetic theory of gases over in silence. Soon he (1894/1954, p. 246) made known his doubts: he was not sure that that theory could account for all the known facts. In a later popular booklet Poincaré (1905/1970, pp. 210 and 251) softened his attitude: physical laws will acquire an "entirely new aspect" and differential equations will become statistical laws; laws, however, will be shown to be imperfect and provisional.

3) Geometric probability. On its previous history see § 6.1.6; its further development is described in Chapter 12. Here, I only indicate that Poincaré explained the paradoxical nature of the **Bertrand** problem (§ 11.1-1).

4) The binomial distribution. Suppose that $m$ **Bernoulli** trials with probability of success $p$ are made and the number of successes is $\mu$. Poincaré (pp. 79 – 84), in a roundabout and difficult way, derived (in modern notation) $E(\mu – mp)^2$ and $E|\mu – mp|$. In the first case he could have calculated $E\mu^2$; in the second instance he obtained

$$E|\mu – mp| \approx 2mpq\, C_m^{mp} p^{mp} q^{mq},\ q = 1 – p.$$

5) The **Bayesian** approach: estimating the total number ($N$) of the asteroids. Poincaré (pp. 163 – 168) assumed that only $M$ of them were known and that, during a certain year, $n$ minor planets were observed, $m$ of which were known before. Introducing a constant probability $p = n/N$ of observing an asteroid during a year and applying the Bayesian approach, he obtained

$$EN \approx n/p.$$

He was not satisfied with this pseudo-answer and assumed now that $p$ was unknown. Again applying the Bayesian approach and supposing that $p$ took with equal probability all values within the interval [0; 1], he derived instead

$$EN = (M/m)n.$$

He could have written this formula at once; in addition, it was possible to recall the **Laplace** problem of estimating the population of France by sample data (§ 7.1-5). It is nevertheless interesting that Poincaré considered the unknown number of the minor planets as a random variable.

6) Without mentioning **Gauss** (1816, § 5), Poincaré (pp. 192 – 194) derived the moments of the [normal] distribution

$$\varphi(y) = \sqrt{h/\pi} \exp(-hy^2) \qquad (3)$$

obtaining

$$Ey^{2p} = \frac{(2p)!}{h^p p! 2^{2p}} \qquad (4)$$

and proved, by issuing from formula (2), that the density function whose moments coincided with the respective moments of the [normal] law was [normal]. This proposition was, however, due to **Chebyshev** (1887a), see also **Bernstein** (1945/1964, p. 420).

Then Poincaré (pp. 195 – 201) applied his investigation to the theory of errors. He first approximately calculated $E\bar{y}^{2p}$ for the mean $\bar{y}$ of a large number $n$ of observations having $Ey_i = 0$ and $Ey_i^2 =$ Const, equated these moments to the moments (4) and thus expressed $h$ through $Ey_i^2$. This was a mistake: $\bar{y}$, being a mean, had a measure of precision $nh$ rather than $h$. Poincaré (p. 195) also stated that **Gauss** had calculated $E\bar{y}^2$; actually, Gauss (1823b, §15) considered the mean value of $\sum y_i^2/n$.

The main point here and on pp. 201 – 206, where Poincaré considered the mean values of $(y_1 + y_2 + \ldots + y_n)^{2p}$ with identical and then non-identical distributions of the terms and $Ey_i = 0$, was a non-rigorous proof of the CLT: for errors of *sensiblement* the same order and constituting *une faible part* of the total error, the resulting error follows *sensiblement* the Gauss law (p. 206)[3].

Also for proving the normality of the sum of errors Poincaré (pp. 206 – 208, only in 1912) introduced characteristic functions which did not conform to their modern definition. Nevertheless, he was able to apply the **Fourier** formulas for passing from them to densities and back. These functions were

$$f(\alpha) = \sum p_x e^{\alpha x}, f(\alpha) = \int \varphi(x) e^{\alpha x} dx \qquad (5)$$

and he noted that

$$f(\alpha) = 1 + \alpha Ex/1! + \alpha^2 Ex^2/2! + \ldots \qquad (6)$$

**Markov** (1898/1951, p. 269) referred, but had not commented on Poincaré (1896, pp. 169 – 186 = 1912, pp. 189 – 206). I repeat that there, on p. 173/194, Poincaré had applied his formula (2).

7) Homogeneous [**Markov** chains]. Poincaré provided interesting examples that might be interpreted in the language of these chains.

a) He (p. 150) assumed that all the asteroids moved along one and the same circular orbit, the ecliptic, and explained why they were uniformly scattered across it. Denote the longitude of a certain minor planet by $l = at + b$ where $a$ and $b$ are random and $t$ is the time, and, by $\varphi(a; b)$, the continuous joint density function of $a$ and $b$. Issuing from the expectation

$$E e^{iml} = \iint \varphi(a; b) e^{im(at + b)} da db$$

(which is the appropriate characteristic function in the modern sense), Poincaré not very clearly proved his proposition that resembled the celebrated **Weyl** theorem (beginning of § 10.8.4). The place of a planet in space is only known with a certain error, and the number of all possible arrangements of the asteroids on the ecliptic might therefore be assumed finite whereas the probabilities of the changes of these arrangements during time period $[t; t + 1]$ do not depend on $t$. The uniform distribution of the asteroids might therefore be justified by the ergodic property of homogeneous **Markov** chains having a finite number of possible states.

b) The game of roulette. A circle is alternately divided into a large number of congruent red and black sectors. A needle is whirled with force along the circumference of the circle, and, after a great number of rotations, stops in one of the sectors. Experience proves that the probabilities of *red* and *black* coincide and Poincaré (p. 148) attempted to justify that fact. Suppose that the needle stops after travelling a distance $s$ ($2\pi < s < A$). Denote the corresponding density by $\varphi(x)$, a function continuous on $[2\pi; A]$ and having a bounded derivative on the same interval. Then, as Poincaré demonstrated, the difference between the probabilities of *red* and *black* tended to zero as the length of each red (and black) arc became infinitesimal (or, which is the same, as $s$ became infinitely large). He based his proof on the method of arbitrary functions (**Khinchin** 1961, No. 2, pp. 88 – 89/2004, pp. 421 – 422; von Plato 1983) and himself sketched its essence. Poincaré also indicated that the rotation of the needle was unstable: a slight change in the initial thrust led to an essential change in the travelled distance (and, possibly, to a change from *red* to *black* or vice versa).

c) Shuffling a deck of cards (p. 301). In an extremely involved manner, by applying hypercomplex numbers, Poincaré proved that after many shuffling all the possible arrangements of the cards tended to become equally probable. See end of § 7.1-6.

8) Mathematical treatment of observations. In a posthumously published *Résumé* of his work, Poincaré (1921/1983, p. 343) indicated that the theory of errors "naturally" was his main aim in the theory of probability and that statement reflected the situation in those times. In his treatise he (pp. 169 – 173) derived the normal distribution of observational errors mainly following **Gauss**; then, like **Bertrand**, changed the derivation by assuming that not the most probable value of the estimator of the [location parameter] coincided with the arithmetic mean, but its mean value. He (pp. 186 – 187) also noted that, for small absolute errors $x_1, x_2, \ldots, x_n$, the equality of some $f(z)$ to the mean value of $f(x_i)$, led to the equality of $z$, the estimate of the true value of the constant sought, and the arithmetic mean of $x_i$. It seemed to him that he thus corroborated the **Gauss** postulate[4]. Finally, Poincaré (p. 188) indicated that the [variance] of the arithmetic mean tended to zero with the increase in the number of observations and referred to Gauss (who nevertheless had not stated anything at all about the case of $n \to \infty$, cf. § 9A.4-7). "Nothing", however, followed since other linear means had the same property, as **Markov** (1899a/1951, p. 250) stated mentioning a wrong remark made by Maievsky. Poincaré himself (1912, pp. 196 – 201 and 217) twice proved the [consistency] of the arithmetic mean. In the second case he issued from a characteristic function of the type of (5) and (6) and passed on to the characteristic function of the arithmetic mean. He noted that, if that function could not be represented as (6), the consistency of the arithmetic mean was questionable, and he illustrated that fact by the **Cauchy** distribution. Perhaps because of all this reasoning on the mean Poincaré (p. 188) declared that **Gauss'** rejection of his first substantiation of the MLSq was "assez étrange" and corroborated this conclusion by remarking that the choice of the [parameter of location] should not be made independently from the distribution (and thus directly contradicting Gauss' mature approach which presumed unknown distributions). This statement corresponds to modern thoughts, but the distributions are hardly known.

Poincaré (pp. 217 – 218) also stated that very small errors made it impossible to obtain absolute precision even as $n \to \infty$. If so, these errors originate from the non-evenness of the law of distribution (**Bayes**; see Stigler (1986, pp. 94 – 95) and **Cournot** (1843, § 137)), the variability of that law (again Cournot) and, I would add, some interdependence of the observations.

9) Randomness. Poincaré discussed randomness both in his treatise and in his popular scientific booklets, but his various interpretations of chance were not definitely enough compared one with another.

a) Instability of equilibrium or movement. Some of the statements made by **Aristotle** (§ 1.1.1) and **Galen** (§ 1.1.3) meant that small causes might lead to considerable consequences, and **Maxwell** (§ 10.8.5-3), who described this phenomenon, apparently thought about randomness but did not mention it.

Poincaré (p. 4) was the first to state directly that randomness was instability of equilibrium or movement and he (pp. 4 – 5) provided a few examples: the instability of a cone stood on its vertex formerly mentioned by Cournot (§ 10.3-4); the roulette; the scattering of the asteroids; unstable states of the atmosphere. His third example, just like **Newton's** reasoning on the irregularities in the Solar system (§ 2.2.3), was nevertheless connected with great intervals of time. Poincaré also argued that **Laplace** (§ 7.3), whom he did not name, was wrong: forecasts of the future were impossible because of the instability of motion. I have not found any connections between the just described explanation of randomness and Poincaré's study of stability in mathematics or astronomy.

b) Complicated causes. Already **Leibniz** (§ 3.1.2) heuristically explained randomness by the complexity of causes. Laplace (end of § 7.4) wrongly explained the existence of the eccentricities of planetary orbits by the action of great many complicated causes. Maxwell (§ 10.8.5-3) assumed that the distribution of the velocities of molecules set in after a great number of collisions among a great number of particles, but he did not mention randomness. And once more Poincaré was the first to do so. He (pp. 7 – 8) maintained that the molecular motion was random because of the combined action of instability and complexity of causes, but he then mentioned the shuffling of cards, the mixing of liquids and powders and (p. 15) "even" of molecules in the kinetic theory of gases.

c) Small causes leading to small consequences. Poincaré (p. 10) provided only one example, that, furthermore, did not belong to natural sciences: small causes led to small errors of measurement; he also indicated that these errors were considered random because their causes were too complicated.

d) Intersection of chains of determinate events. I mentioned this explanation in §§ 1.1.1 and 10.3-4. Poincaré (p. 10) allowed it, but his first two explanations were his main ones; and he apparently forgot here about the third one.

e) Randomness and necessity. Poincaré (see § 1.2.4) formulated a highly proper idea on the combined action of randomness and necessity. Regrettably, he had not mentioned the appearance of necessity in mass random phenomena.

For Poincaré, the theory of probability remained an accessory subject, and his almost total failure to refer to his predecessors except **Bertrand** testifies, as mentioned by Bortkievicz, that he was not duly acquainted with their work. Furthermore: in 1912 he was already able to, but did not apply the arsenal of **Markov** chains. At the same time, however, he became the author of a treatise that for about 20 years had remained the main writing on probability in Europe. Le Cam's declaration (1986, p. 81) that neither **Bertrand**, nor Poincaré "appeared to know" the theory was unjust: he should have added that, at the time, Markov was apparently the only one who did master probability. On Bertrand see end of § 11.1.

### Notes

**1.** The title-page of the French translation carried a phrase "Translated and published avec l'autorisation de l'auteur", but Bertrand himself (*C. r. Acad. Sci. Paris*, t. 40, 1855, p. 1190) indicated that **Gauss**, who had died the same year, was only able to send him "quelques observations de détail".

**2.** Poincaré always applied the term "calcul" rather than "théorie" of probability. It is hardly amiss to note that in 1882 – 1891 **Markov** had published five mimeographed editions of his lectures called *Theory of probability*, but that he called his treatise (1900 and later editions) *Calculus of probability* (in both cases, in Russian). Another point: at least in 1892 Poincaré was not prepared to believe in the statistical nature of the second law of thermodynamics; in addition to § 11.2-2 above, see Sheynin (1991a, p. 141).

**3.** For Poincaré the theory of probability remained a branch of applied mathematics (§ 0.1 and Note 6 to Chapter 8).

**4.** In the same context Poincaré (p. 171) argued that everyone believed that the normal law was universal: experimentalists thought that that was a mathematical fact and mathematicians believed that it was experimental. Poincaré referred to the oral statement of Lippmann, an author of a treatise on thermodynamics.

## 12. Geometric Probability

On the development of the notion of geometric probability in the 18$^{th}$ century and earlier see § 6.1.6, and on its definition by **Cournot** see § 10.3-4; I described the **Bertrand** problem on the length of a random chord in § 11.1-1. Here, I discuss the further history of the same notion.

1) Cournot (1843, § 74) applied geometric probability for deriving the distribution of a function of several random arguments. Here is one of his examples. The arguments of the function $u = |x - y|$ are uniformly distributed on segment [0; 1]. After calculating the areas of the appropriate figures, he concluded that

$$P(u \leq a) = (1 - a^2), 0 \leq a \leq 1.$$

The determination of the probability of the contrary event would have led Cournot to the once popular encounter problem (**Laurent** 1873, pp. 67 – 69): two persons are to meet at a definite spot during a specified time interval, their arrivals are independent and occur "at random". The first one to arrive waits for a certain time and then leaves. Determine the probability of the encounter.

2) Most eminent natural scientists of the 19$^{th}$ century tacitly applied geometric probability. **Boltzmann** (§ 10.8.5-1) defined the time average probability that the velocity of a molecule was contained in an interval [$c$; $c + dc$] as the ratio of the time during which that event took place to the total time of observation. I do not dwell on an earlier definition of probability in physics or the further considerations concerning the ergodic hypothesis. **Maxwell** (§ 10.8.5-2) applied geometric probability while deriving his celebrated law.

When studying the life of earthworms, **Darwin** (1881/1945, pp. 52 – 55) strewed paper triangles over some ground. They were dragged away by the worms but he recovered most of them and found out that the worms had not seized "indifferently by chance any part" of the triangles. He thought about several possibilities of [uniform] "chance", and, in particular, he decided, in actual fact, that the number of times a worm would have seized "by chance" any side of a triangle was proportional to its length[1].

3) Seneta et al (2001) described the pertinent investigations of **Sylvester**, **Crofton** and **Barbier** which led to the appearance of integral geometry. I only mention Sylvester's remarkable problem: To determine the probability that four points taken "at random" within a finite convex domain will form a convex quadrilateral. See **Czuber** (1908/1968, pp. 99 – 102) for a few particular cases of that problem.

4) **Czuber** (1884, p. 11): Two points, M and N, are randomly situated on segment AB = $a$. Determine the probability of MN > NA. We have

$$P(x \leq MA \leq x + dx) = dx/a, P(MA > NA | x \leq MA \leq x + dx) = x/a.$$

$$P(MN > NA) = \int_0^a x dx/a^2 = 1/2.$$

Indeed, suppose that the points A, M, N, B are situated in that order from *left* to *right*. (The other possible case can be considered likewise.) Then, if MN > NA, move M towards N until MN < NA, which is always possible, QED.

5) **Poincaré** (1896, p. 97; 1912, p. 118) noted that the probability that a point (*x*; *y*) was situated within some figure was equal to the appropriate integral

$$\iint (x; y)\, dx\, dy$$

where should be somehow specified. He then went over to the **Bertrand** problem (§ 11.1-1) but mentioned only two of its solutions and provided his own reasoning tacitly assuming that 1. The chord can be fixed with respect to the centre of the circle *O* and the polar axis passing through, and beginning in *O*, by two parameters, and , – the polar angles of *A*, an endpoint of the chord, and of *P*, its centre; or, by two other parameters, and , – the polar coordinates of *P*. Now, the integrals over the given circle

$$\iint d\, d \qquad \iint d\, d$$

and this, as Poincaré stated, explained the paradoxical nature of the problem.

He also studied the probability that rotated figures satisfy certain conditions but he did not state that this investigation was connected with the Bertrand paradox, cf. § 12-7 below.

6) Czuber (1908/1968, pp. 107 – 108) discovered three more natural solutions of the Bertrand problem.

   a) One endpoint of the chord is fixed, and the chord passes through any point of the circle; $p = 1/3 + \ 3/2 \ \ 0.609$.

   b) Both endpoints of the chord are chosen randomly; this case coincided with Bertrand's first version.

   c) Two points of the chord situated inside the circle are chosen randomly; $p = 1/3 + 3\ 3/4 \ \ 0.746$.

7) It turned out (De Montessus 1903) that the Bertrand problem had an uncountable set of answers. Suppose that *Ox* is the *x*-axis and mark points *D* and *C* on its positive half, – its intersections with concentric circumferences with common centre in point *O* and radii $OD = 1/2$ and $OC = 1$. Arbitrary points $M_2(x)$ and $M_3(x)$ are situated on the same semiaxis, between the two circles and beyond the larger of them respectively. Tangents $A_2B_2$ and $A_3B_3$ to the smaller circumference pass through $M_2$ and $M_3$ respectively, and $M_3T$ is the tangent to the larger circumference with point of contact *T*. Finally, $M_1(x)$ is an arbitrary point on the same semiaxis inside the smaller circle.

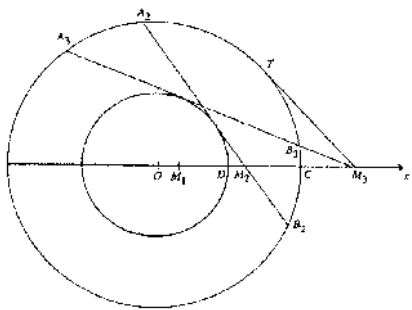

Fig. 1. De Montessus (1903). A point moves along the axis from *D* to infinity, and, correspondingly, the probability sought in the Bertrand problem is seen to have an uncountable set of values. $OD = 1/2$, $OC = 1$.

For points $M_2$ and $M_3$ the probability sought is, respectively,

$P_2$ = angle $A_2M_2O/$ = [2arcsin(1/2x)]/ ,
$P_3$ = angle $A_3M_3O$/angle $TM_3O$ = [arc sin(1/2x)]/[arcsin(1/x)],

with $1/2$ x 1 and x 1 respectively.

When moving from point *O* in the positive direction (say), the probability $P_2$ decreases from 1 at point *D* to 1/3, and, from point *C* to infinity, probability $P_3$ increases from 1/3 to 1/2. It is rather difficult to prove that $P_3$ increases monotonically (and De Montessus had not done it), but already for $x = 1.01$ and 1.1 it is 0.36 and 0.41 respectively and it reaches value $(1/2 – 1/1,600)$ at $x = 10$.

Note that the coincidence of points $M_2$ or $M_3$ with *D* leads to **Bertrand's f**irst solution and the movement $M_3$ provides his second case. His third solution concerned a point rather than a straight line and was thus different.

De Montessus calculated the general mean probability of the studied event. However, it was hardly proper to include in the calculation, as he did, points such as $M_1$ for which the stipulated condition was certainly satisfied. More important, while calculating the mean probability for the continuous case, De Montessus first determined a finite sum, and, when adding together the appropriate fractions, he added separately their numerators and their denominators. Nevertheless, his mean probability ($P = 1/2$) was correct and could have been established by noting that the studied interval beyond the circle was infinite.

8) **Borel** (1909/1950, p. 132) solved the encounter problem (§ 12-1) and on pp. 148 – 149 he indicated, without referring to anyone that most of the natural methods of solving the Bertrand problem led to $P = 1/2$.

9) **Schmidt** (1926) issued from **Poincaré's** considerations and indicated that the probability sought should persist under translation and rotation of the coordinate system (invariance under reflection is now also included). Accordingly, he proved that this condition is only fulfilled for the ( ; ) coordinate system, see Item 5, and when transforming that system into another one (with the appropriate **Jacobian** being of course allowed for)[2].

And so, commentators have finally agreed that the probability sought was 1/2 which was tantamount to ignorance and could have been stated from the very beginning.

For a modern viewpoint on geometric probability see **M. G. Kendall** & **Moran** (1963); in particular, following authors of the 19th century (e.g., **Crofton** 1869, p. 188), they noted that it might essentially simplify the calculation of integrals. Then, Ambartzumian (1999) indicated that geometric probability and integral geometry are connected with stochastic geometry.

### Notes

**1. Darwin** thus considered several possibilities of a "random" dragging of triangles and in that sense his study forestalled the **Bertrand** problem on the length of a random chord. Darwin attempted to ascertain whether or not the worms acted somewhat intelligently, and he concluded that they had not seized the triangles indifferently.

**2. Prokhorov** (1999b) believed that, from the geometrical point of view, the most natural assumption in the **Bertrand** problem was that   and   were independent and uniformly distributed, 0   2 , 0   1.

## 13. Chebyshev
### 13.1. His Contributions

1) His Master's dissertation (1845). It was intended as a manual for students of the Demidov Lyceum in Yaroslavl and Chebyshev gave there an account of the theory of probability barely applying mathematical analysis; for example, he replaced integration by summing. Already then, however, he consistently estimated the errors of pre-limiting relations. The dissertation apparently had an addendum published somewhat later, see § 13.1-2.

Incidentally, I doubt that this manual was a proper textbook. A general survey of the theory of probability and its applications would have been more useful the more so since Chebyshev's reasoning was necessarily burdensome.

2) The **Poisson** LLN (Chebyshev 1846); see **Prokhorov** (1986) for a detailed exposition. Chebyshev solved the following problem. In $n$ [independent] trials the probability of success was $p_1, p_2, \ldots, p_n$. Determine the probability that the total number of successes was $\mu$. By clever reasoning he obtained the formula

$$P(\mu \quad m) \quad \frac{1}{2\sqrt{n}} \frac{\sqrt{m(n-m)}}{m-ns} \left(\frac{ns}{m}\right)^m \left(\frac{n(1-s)}{n-m}\right)^{n-m+1}$$

where $m > ns + 1$ and $s$ was the mean probability of success.

This result was interesting in itself and, in addition, it enabled Chebyshev to prove the **Poisson** theorem, cf. § 8.7. He did not fail to indicate the necessary number of trials for achieving a stipulated probability of the approximation of the frequency of success $\mu/n$ to $s$. As stated in the title of the memoir, Chebyshev had indeed not applied any involved mathematical tools, but his transformations were burdensome. His proof was rigorous (although he had not indicated that the trials were independent) and he (1846) had the right to reproach Poisson whose method of derivation did not provide the limits of the error of his approximate analysis. See however the very end of § 13.3. Later Chebyshev (1879 – 1880/1936, pp. 162 – 163, translation p. 152) explicated one of his intermediate transformations more clearly, also see **Bernstein** (1945/1964, p. 412).

3) The **Bienaymé – Chebyshev** inequality (cf. § 10.2-4). Chebyshev (1867) considered discrete [random variables] with a finite number of possible values; without loss of generality I simplify his derivation by assuming that each of the $n$ variables takes an equal number of values. Chebyshev showed that

$$P\{|\ (_i - E\ _i)| < \sqrt{\sum [E\ _i^2 - (E\ _i)^2]}\ \} > 1 - 1/\ ^2, \quad > 0. \quad (1)$$

Unlike Heyde & Seneta (§ 10.2-4) I believe that Chebyshev derived this inequality in about the same way as **Bienaymé** did, only in much more detail. True, he restricted his attention to discrete variables whereas Bienaymé, without elaborating, apparently had in mind the continuous case; his memoir was devoted to the mathematical treatment of observations. Modern authors, whom I mentioned in § 10.2-4, repeat the derivation for the latter instance; actually, already Sleshinsky (1893) had done it.

Chebyshev immediately derived a corollary, which, in somewhat different notation, was

$$\lim P(\ [|\sum u_i - E\sum u_i|/n] < \varepsilon) = 1, n \to \infty$$

and he (1879 – 1880/1936, pp. 166 – 167, translation pp. 155 – 156) specified this formula for the case in which the [random variables] coincided one with another. Chebyshev thus obtained a most important and very simple corollary: the arithmetic mean was a [consistent] estimator of the expectation of a random variable. Both corollaries assume that the expectations and variances of the appropriate variables are uniformly restricted and Chebyshev had indeed indicated this restriction (in another language). In the last-mentioned source, and even earlier in another context, he (1867, p. 183) introduced indicator variables (taking values 0 and 1 with respective probabilities) but not the term itself.

4) [The central limit theorem]. The title of the appropriate memoir (1887b) mentions two theorems the first of which was the proposition on the arithmetic mean (see Item 3) and Chebyshev only repeated its formula. He then went on to the CLT noting that it "leads" to the MLSq, – leads, as I comment, in accordance with the **Laplacean** approach.

Chebyshev first of all referred to his inequalities for an integral of a non-negative function whose moments up to some order coincided with the same moments of the appropriate, in a definite sense, normal distribution. He (1874) had published these inequalities without proof and **Markov** (1884) and then **Stieltjes** substantiated them. Chebyshev himself also justified them afterwards but without mentioning his predecessors. A detailed history of these inequalities is due to **Krein** (1951). Later Stieltjes (1885) expressed regret that he had overlooked that Markov paper.

Chebyshev considered [random variables] $u_1, u_2, \ldots, u_n$ having densities $\varphi_i(x)$ and moments

$$Eu_i = 0, |Eu_i^2| < C, |Eu_i^3| < C, \ldots$$

These conditions are not sufficient. The random variables ought to be independent, and Chebyshev certainly thought so, but he had not indicated the restriction

$$\lim \frac{\sum Eu_i^2}{n} \neq 0, i = 1, 2, \ldots, n, n \to \infty. \qquad (2)$$

On the other hand, it was not necessary to demand that the moments were uniformly bounded and Chebyshev possibly did not express such a restriction. Here is **Liapunov's** indirect testimony (1901b, p. 57): it occurs that Chebyshev sometimes used the singular form instead of the plural. Liapunov provided a few examples and, in particular, quoted Chebyshev's expression "the absolute value of the mathematical expectations" from his formula of the CLT.

Chebyshev noted that the density $f(x)$ of the fraction

$$x = \sum u_i / \sqrt{n} \qquad (3)$$

can be determined by means of the multiple integral

$$f(x) = \int \int \ldots \int \varphi_1(u_1) \varphi_2(u_2) \ldots \varphi_n(u_n) du_1 du_2 \ldots du_n \qquad (4)$$

extended over the values of the variables at which the fraction above is situated within the interval $[x; x + dx]$. He multiplied both parts of (4) by $e^{sx}$ where $s$ was some constant and integrated them over $(-\infty; +\infty)$ so that the right side became separated into a product of $n$ integrals with the same limits of integration. Chebyshev then developed both parts in powers of $s$ (the right side, after taking its logarithm) and equated the coefficients of the same powers of that magnitude to each other. Thus the integrals

$$\int f(x)dx, \int xf(x)dx, \int x^2 f(x)dx, \ldots$$

or, in other words, the moments of magnitude (3), were determined up to some order $(2m - 1)$. It occurred that, as $n \to \infty$, again with the same limits of integration,

$$\int e^{sx} f(x)dx = \exp(s^2/2q^2) \qquad (5)$$

where $1/q^2$ was the arithmetic mean of the second moments of $u_i$ and it is here that the condition (2) was needed. Applying his previously mentioned estimates of the integral of a non-negative function, Chebyshev now completed his proof:

$$\lim P(\alpha \leq \frac{\sum u_i}{\sqrt{2 \sum E u_i^2}} \leq \beta) = (1/\sqrt{\pi}) \int_\alpha^\beta \exp(-x^2)dx, \, n \to \infty. \qquad (6)$$

For finite values of $n$ the same probability, as Chebyshev indicated without a rigorous demonstration, was determined by a development in polynomials now called after him and **Hermite**.

**Bernstein** (1945/1964, pp. 423 – 425) indicated that the abovementioned expansion in powers of $s$ diverged at $|s| \neq 0$ and that **Markov** (1898/1951, p. 268), when proving the Chebyshev theorem anew, without explaining the situation had therefore introduced an additional restriction,– not (2), but

$$\lim E u_n^2 \neq 0, \, n \to \infty.$$

In addition, Markov wrote out the expression that Chebyshev had actually applied in his investigation:

$$\lim \left( \frac{\sum u_i}{\sqrt{2 \sum \mathrm{E} u_i^2}} \right)^m = \frac{1}{\sqrt{\pi}} \int_{-\infty}^{\infty} t^m \exp(-t^2) dt, \quad n \quad . \qquad (7)$$

Issuing from the Chebyshev inequalities, Markov also proved that these expressions meant that the appropriate density tended to the normal density whereas Chebyshev had apparently thought it evident. For a detailed discussion see **Kolmogorov's** commentary in Chebyshev (1944 – 1951, vol. 3, 1948, pp. 404 – 409).

Sleshinsky (1892) turned his attention to the Chebyshev demonstration of the CLT even before Markov did. He issued from **Cauchy's** findings (§ 10.1) and, as he stated on his p. 204, aimed at simplifying (not specifying) his predecessor. In spite of **Freudenthal's** opinion and following Heyde & Seneta (1977, pp. 95 – 96), I think that the Cauchy investigation was nevertheless imperfect, see end of § 13.1. And, once more repeating the last-mentioned commentators, I note that Sleshinsky apparently proved the CLT rigorously although only for a linear function of observational errors having an even density. **Liapunov** (1900, p. 360) remarked that these conditions were too restrictive. Like Chebyshev (§ 0.3, Note 1), Sleshinsky maintained that his findings justified the MLSq [once more: only in the **Laplacean** sense].

### 13.2. His Lectures

Chebyshev delivered lectures on the theory of probability at Petersburg University from 1860 to 1882. In 1936, **A. N. Krylov** published those read in 1879/1880 as recorded by **Liapunov** and I refer to his publication by mentioning only the page numbers of this source and of its translation of 2004.

In his Foreword Krylov declared that Liapunov had reproduced the lectures "exactly as they were read, including all the fine points of his accompanying remarks". Prudnikov (1964, p. 183), however, thought differently:

*It was hardly possible to write down Chebyshev's lectures minutely and it is natural that their extant record is fragmentary.*

This seems to be at least partly true. Krylov also indicated that he had rewritten the Liapunov manuscript

*In accordance with the new system of spelling at the same time checking all the derivations …*

I translated this book correcting perhaps a hundred (I repeat: a hundred) mathematical misprints but I do not claim that I revealed all of them. Ermolaeva (1987) briefly described a more detailed record of Chebyshev's lectures read during September 1876 – March 1878, discovered by herself but still unpublished. She had not indicated whether the newly found text essentially differed from the published version and she is behaving like a dog in the manger.

The lectures were devoted to definite integrals, the theory of finite differences and the theory of probability. I discuss only their last section but I begin by several general comments. Chebyshev attempted to apply the simplest methods; for example, he used summing, and, if necessary, went on to integration only at the last moment; he introduced characteristic functions only in the discrete case; as I mentioned above; he did not specify that he considered independent events or variables; he was not interested in the philosophical aspect of probability[1]; and, among the applications of the theory of probability, he almost exclusively discussed the mathematical treatment of observations.

1) The main notions. Chebyshev (p. 148, translation p. 141) declared that the aim of the theory of probability was

*To determine the chances of the occurrence of a certain event*, and he continued: *the word 'event' means anything whose probability is being determined*, and probability *serves to denote some magnitude that is to be measured*.

The use of "chance" and "probability" in the same sentence was perhaps an elegant variation; in essence, however, Chebyshev made a small heuristic step towards an axiomatic theory[2]. I also adduce a modern formula (**Prokhorov** & **Sevastianov** 1999, p. 77): the theory of probability studies mathematical models of random events and

*Given the probabilities of some random events, makes it possible to determine the probabilities of other random events somehow connected with the first ones.*

Chebyshev (p. 160/150) introduced an unusual and hardly useful generalized definition of expectation, – of the expectation of the occurrence of one out of several incompatible events. The sum of the products of the type $p_i a_i$, as he stated, described these events by their probabilities and the magnitudes "measuring" them. Note that he mainly discussed discrete variables.

Tacitly following **Laplace** (§ 7.4), Chebyshev (p. 165/155) indicated that the concept of limit in probability theory differed from that in analysis. Still, I am unable to agree with such equalities (or misprints?) as (pp. 167, 183, 204/156, 171, 190)

$$\lim m/n = p. \tag{8}$$

2) The limit theorem for **Poisson** trials (p. 167 and 201ff/156, 187ff). Determine the probability $P_{n,m}$ that in $n$ trials an event having probabilities $p_i$, $i = 1, 2, …, n$, respectively, occurred $m$ times. Applying a little known formula from the first section of his lectures (p. 59/63)

$$A_m = (1/2\pi) \int_{-\pi}^{\pi} f[e^{i\varphi}]e^{-m i\varphi}d\varphi, \qquad (9)$$

for the coefficients of the series

$$f(x) = A_0 + A_1 x + A_2 x^2 + \ldots + A_m x^m + \ldots,$$

Chebyshev obtained ($q_i = 1 - p_i$)

$$P_{n,m} = (1/2\pi) \int_{-f}^{f} [p_1 e^{i\varphi} + q_1][p_2 e^{i\varphi} + q_2]\ldots[p_n e^{i\varphi} + q_n]e^{-m i\varphi}d\varphi.$$

After some transformations when considering only small values of $\varphi$ it occurred that

$$P_{n,m} = (1/\pi) \int_0^{\varphi} \exp(-nQ\varphi^2/2)\cos[(np - m)\varphi]d\varphi$$

where $p$ was the mean probability of success and $Q = \sum p_i q_i / n$. Assuming for large values of $n$ an infinite upper limit in the obtained integral, Chebyshev finally got

$$P[|m/n - p| < t\sqrt{2Q/n}] = (2/\sqrt{\pi}) \int_0^t \exp(-z^2)dz$$

(without the sign of limit!) and noted that formula (8), or, as he concluded, the **Poisson** LLN, followed from it. He naturally did not here admonish his predecessor.

3) The **Bernoulli** pattern (pp. 168 – 175/157 – 163). Chebyshev wrote out the generating function of the binomial distribution (in usual modern notation)

$$\sum P_{n,m} t^m = (pt + q)^n, m = 0, 1, 2, \ldots, n \qquad (10)$$

and calculated the appropriate expectation and variance by the modern method (by differentiating this function once and twice etc.) although without indicating its generality. He then repeated this derivation otherwise. Assuming that in equation (10) $t = e^\tau$, Chebyshev multiplied both its parts by $e^{-\tau pn}$, developed the exponential functions in powers of $\tau$ and equated the coefficients of $\tau$ and then of $\tau^2$. In concluding, he (pp. 179 – 183/167 – 171) derived the local and the integral **De Moivre** – **Laplace** limit theorems and (pp. 183 – 186/171 – 175) paid attention to the calculation of the integral of the exponential function of a negative square. I note his unusual manner which, in this case, becomes evident when he stated that the abovementioned integral with the limits of integration being $[u; +\infty)$ was equal to the value of the integrand at the lower limit multiplied by some proper fraction, – rather than by a real number situated in the interval (0; 1).

4) A limit theorem for the multinomial distribution (pp. 205 – 207 and 214 – 218/190 – 193, 198 – 203). Chebyshev considered $n$ trials in each of which occurred one and only one event out of $A_1, A_2, \ldots, A_k$ with event $A_i$ meaning that some function took the value $i$. All the events were equally probable so that each had probability $1/k$. Suppose that event $A_i$ happened $m_i$ times, and

$$m_1 + m_2 + \ldots + m_k = n, \ P(m_1 + 2m_2 + \ldots + km_k = s) = P_s,$$

then

$$P_s t^s = t^n(t^k - 1)^n / k^n (t - 1)^n \qquad (11)$$

after which Chebyshev determined $P_s$.

When considering the limiting case he expressed the right side of (11) as

$$f(t) = A_o + A_1 t + A_2 t^2 + \ldots + A_s t^s + \ldots$$

and made use of the expression (9) so as to obtain

$$k^n P_s = (1/2\pi) \int_{-\pi}^{\pi} e^{i\varphi(n-s)} \{[e^{k i\varphi} - 1] \div [e^{i\varphi} - 1]\}^n d\varphi ,$$

where the $n$-th power of the fraction was equal to

$$e^{n(k-1)i\varphi/2} [\sin(k\varphi/2) \div \sin(\varphi/2)]^n$$

so that

$$P_s = \frac{1}{\pi} \int_0^\pi \cos\{([n(k-1) - 2s]\frac{\varphi}{2}\}\{[\sin k\varphi/2] \div [k\sin \varphi/2]\}^n d\varphi ,$$

and, for large values of $k$, again without the sign of limit,

$$P(|s - kn/2| < ku\sqrt{n/6}) = (2/\sqrt{\pi}) \int_0^u \exp(-t^2)dt.$$

5) [The central limit theorem] (pp. 219 – 223/203 – 206). At the time, Chebyshev had not yet known its rigorous proof. I only note his pronouncement (p. 224/206): the formula that he obtained was not derived

*In a rigorous way […]. We have made various assumptions but did not determine the boundary of the ensuing error. In its present state, mathematical analysis cannot derive this boundary in any satisfactory fashion.*

6) Statistical inferences. Chebyshev solved two problems which, however, were considered before him. In the first of these he (pp. 187 – 192/175 – 180) derived the **Bayes** limit theorem (§ 5.2). He ascertained that, given the length $2w$ of the segment $[b, c]$ and large numbers $p$ and $q$ (notation as in formula (5.2)), probability $P$ takes its maximal value at

$$b = p/(p + q) - w, \ c = p/(p + q) + w$$

which is natural. After that Chebyshev (p. 142) had indeed derived formula (5.3), but in his *Lectures* he had not mentioned Bayes at all.

In the second problem he (pp. 193 – 201/181 – 187) studied the probability of a subsequent result in **Bernoulli** trials. An event occurred $m$ times in $n$ trials; determine the probability that it will happen $r$ times in $k$ new trials. Only guiding himself mostly by the **Stirling** theorem, Chebyshev non-rigorously derived an integral limit theorem similar to that obtained by **Laplace** (§ 7.1-5). His formula (again without the sign of limit) was

$$P\left(|\frac{r}{k} - \frac{m}{n}| < t\sqrt{2\frac{m}{n}(1-\frac{m}{n})(\frac{1}{n} + \frac{1}{k})}\right) = (2/\ \ )\int_0^t \exp(-z^2)dz. \quad (12)$$

Later **Markov** (1914a) indicated the same formula, although correctly written down as a limit theorem. He hardly remembered its occurrence in Chebyshev's lectures.

**7)** Mathematical treatment of observations (pp. 224 – 252/207 – 231). Chebyshev (p. 227/209) proved that the arithmetic mean was a [consistent] estimator of the unknown constant. Unlike **Poincaré** (§ 11.2-7), he (pp. 228 – 231/209 – 212) justified its optimality by noting that, among linear estimators, the mean ensured the shortest probable intervals for the ensuing error. The variance of the arithmetic mean was also minimal (Ibidem); although Chebyshev had not paid special attention to that estimator of precision. In principle, he based his reasoning on the definitive **Gaussian** substantiation of the MLSq (§ 9A.4).

But at the same time Chebyshev (pp. 231 – 236/212 – 216) derived the normal distribution as the universal law of error in about the same way as Gauss did in 1809 (§ 9A.2). "The Gauss method", Chebyshev (p. 250/229) maintained, bearing in mind exactly that attempt later abandoned by Gauss, was based on the doubtful "law of hypotheses", – on the "**Bayes** theorem" with equal prior probabilities. Chebyshev several times censured that "law" when discussing the Bayesian approach in his lectures. In this case, it is opportune to recall **Whittaker** & **Robinson's** remark in § 9A.2-2. I also note that Chebyshev (p. 249/228) wrongly thought that the Gauss formula (9.6) for the sample variance had only appeared "recently" and that it assumed a large number of observations. Chebyshev indicated that he considered only random errors having zero expectations, but he did not mention that the Gauss formula provided an unbiased estimation.

It might be concluded that the treatment of observations hardly interested him.

8) *Cancellation of a fraction* (pp. 152 – 154/144 – 146). Determine the probability $P$ that a "random" fraction $A/B$ cannot be cancelled. Denote the probability that a prime number $m$ cannot be cancelled out of $A/B$ by $p_m$. Then

$$P = p_2 \, p_3 \, p_5 \, \ldots \, p_m.$$

Since the probability that $A$ or $B$ is divisible by $m$ is $1/m$ (this was an essential assumption, see comment below!),

$$p_m = 1 - 1/m^2,$$
$$P = (1 - 1/2^2)(1 - 1/3^2)(1 - 1/5^2) \ldots (1 - 1/m^2) \ldots, \qquad (13)$$
$$1/P = 1 + 1/2^2 + 1/3^2 + 1/4^2 + \ldots = \pi^2/6, \qquad (14)$$

$$P = 6/\pi^2.$$

Chebyshev did not explain the transition from product to series, but it was known to **Euler** (1748, Chapter 15, §§ 275 – 277). Chebyshev determined the sum (14) by two different methods. One of them consisted in equating the coefficients of $x^2$ in two different expansions of $\ln(\sin x/x)$ of which at least the second one was again known to Euler (Ibidem, Chapter 9, § 158):

$$\ln(1 - x^2/6 + x^4/120 - \ldots) = \ln[(1 - x^2/\pi^2)(1 - x^2/4\pi^2)(1 - x^2/9\pi^2) \ldots].$$

Chebyshev also remarked that if a fraction cannot be reduced by 2, 3, or 5, then $1/19 < 1 - P < 1/20$, which testifies once again that he paid due attention to practical considerations[3]. **Markov**[4] remarked that **Kronecker** (1894, Lecture 24) solved the same problem but indicated **Dirichlet's** priority. Kronecker had not supplied an exact reference and I was unable to check his statement; he added that Dirichlet had determined the probability sought "if it existed at all".

**Bernstein** (1928/1964, p. 219) refuted Chebyshev's solution by noting that his assumption led to contradiction. He also adduced further considerations, and, in particular, indicated, on p. 220, that the theory of numbers dealt with regular number sequences whose limiting or asymptotic frequencies of numbers of some class, unlike probabilities, "we will never determine experimentally", might be studied. See Postnikov (1974) on the same problem and on the stochastic theory of numbers.

### 13.3. Some General Considerations

And so, Chebyshev argued that the propositions of the theory of probability ought to be rigorously demonstrated and its limit theorems should be supplemented by estimation of the errors of pre-limiting relations (**Kolmogorov** 1947, p. 56). He himself essentially developed the LLN and, somewhat imperfectly, proved for the first time the CLT; on the study of these two issues depended the "destiny" of the theory of probability (**Bernstein** 1945/1964, p. 411). His students,

**Markov** and **Liapunov** in the first place, also contributed to the theory (§§ 14.1, 14.2, 14.4).

Kolmogorov continued:

*Chebyshev was the first to appreciate clearly and use the full power of the concepts of random variable and* [its] *expectation.*

I take issue with *the full power*. Indeed, Chebyshev had not made use of **Poisson's** heuristic definition of random variable (§ 8.2), had not applied this term[5] and did not study densities or generating functions as mathematical objects. Then, the entire development of the theory of probability might be described as an ever fuller use of the power of the abovementioned concepts; thus, it had since began to study dependent random variables, their systems and chains.

Here also is Bernstein's conclusion (1945/1964, p. 432):

*The genius of Chebyshev and his associates, who, in this field* [theory of probability]*, have left mathematicians of Western Europe far behind, have surmounted the crisis of the theory of probability that had brought its development to a stop a hundred years ago.*

"Crisis" may be understood as a dangerous and unstable state; in this case, as the theory's extremely unfavourable state as compared with the main branches of mathematics then rapidly developing in Europe. Two circumstances ought to be mentioned here. First,

*In spite of his splendid analytical talent, Chebyshev was a pathological conservative.*

This is the opinion of **Novikov** (2002, p. 330) who corroborated it by referring to **V. F. Kagan** (1869 – 1953), an eminent geometrician. The latter,

*When being a young Privat-Docent*, had listened to Chebyshev's scornful statement on the *trendy disciplines like the **Riemann** geometry and complex-variable analysis*.

Even **Liapunov** (1895/1946, pp. 19 – 20), who (**Bernstein** 1945/1964, p. 427)

*Understood and was able to appreciate the achievements of the West European mathematicians, made in the second half of the* [19th] *century, better than the other representatives of the* [Chebyshev] *Petersburg school,*

called Riemann's ideas "extremely abstract"; his investigations, "pseudo-geometric" and sometimes, again, too abstract and having nothing in common with **Lobachevsky's** "deep geometric studies". Liapunov obliquely recalled **Klein**, but disregarded him. Klein had in 1871 presented a unified picture of the non-Euclidean geometry in which the findings of Lobachevsky and Riemann appeared as particular cases.

On the other hand, Tikhomandritsky (1898, p. IV) testified that in 1887 he had shown Chebyshev his "course" and that the latter "stated that […] it is necessary to transform the entire theory of probability". It is difficult to say what exactly did he mean. His words must have been known; I found later references to them (Maciejewski 1911, p. 87; **Gnedenko** & **Gikhman** 1956, p. 487). On the Petersburg school of the theory of probability see also **Bernstein** (1940).

But why neither Poisson nor Poincaré had not attempted to prove rigorously their statements? With respect to the CLT Poisson apparently was unable to act differently, but I have indicated the main reason in Notes 5 and 6 to Chapter 8 and in § 7.4. Neither densities nor characteristic functions had hardly been treated as mathematical objects.

And the *destiny* of probability theory not less depended on the proper handling of those two objects which had only occurred in the 1920s (P. Lévy). True, Poisson (1837a, p. 1) indicated that in the 18$^{th}$ century did the calculus (not the theory) of probability

*Wholly extend and became one of the main branches of mathematics owing both to the number and utility of its applications and the type of the analysis which it beget*.

Still, this does not contradict the statements above. And here is a suitable statement of Markov from his report of 1921 (Sheynin 2006a, p. 152):

*It was usual to consider the theory of probability as an applied science in which mathematical rigour was superfluous.*

A special comment on Chebyshev is warranted. For many yedars had had been associated with problems in artillery (Prudnikov 1950) and certainly with, and essentially influenced Maievsky, the most eminent specialist in ballistics, and *the founder of modern external ballistics,* see his book of 1872 (Mandryka 1954, p. 162).

Chebyshev published a note (1870) in the memoirs of the Belgian academy which was favourably discussed by mostly French authors. It apparently had not contained any essentially new material as compared with his previous publications but nevertheless it should have been, but was not mentioned in his collected works (1899 – 1907; 1944 – 1951). I have almost completely reprinted it in one of my collections (**S, G,** 78) as Mayevski (1872) also did long before me. However, I have also described the comments mentioned above. In somewhat more detail see the same in my Russian collection **S, G,** 82.

As a sideline I add that both Chebyshev and Mayevski had been applying the term *mean square error* so that it seems that it had thus at least found its way into Russian geodesy.

And here are a few lines from Maievsky (1872, p. XII):

*Chebyshev provides formulas for interpolating by the method of least squares. We have applied them for determining the projection of the path* [of the shell] *on the vertical plane of the firing by issuing from its results*.

## Notes

**1.** Prudnikov (1964, p. 91) quoted a paper of V. A. Latyshev, an educationalist and a student of Chebyshev, published in 1893:

*One of the most distinguished* [Russian] *mathematicians* […] *had the habit of expressly telling his students that he did not advise* [them] *to engage in the philosophical aspect of mathematics since this was not very helpful for acquiring the knowledge of mathematics, and even rather harmful.*

Prudnikov added that Latyshev had certainly meant Chebyshev. Recall (§ 5.1) that Chebyshev formulated the problem on the next sunrise in everyday language.

**2.** Chebyshev (1845/1951, p. 29) provided a similar definition of the aims of the theory of probability much earlier. It is hardly amiss to remark that for **Laplace** the theory served for discovering the laws of nature. **Boole** (1851/1952, p. 251) expressed ideas similar to those formulated by Chebyshev:

*The object of the theory of probabilities may be thus stated: Given the separate probabilities of any propositions to find the probability of another proposition.*

He (1854, p. 288) was also the first to argue that the theory should be axiomatized, see **Supplement**. On Boole's probability see Hailperin (1976) who does not, however, dwell on axiomatization.

**3.** Note however Chebyshev`s unqualified statement (1879 – 1880/1936, p. 214, translation p. 198): Different lotteries are equally fair if the expected gains are the same and equal to the [equal] stakes. This contradicts the reasonable opinion of both **D'Alembert** and **Buffon** (§§ 6.1.2 and 6.1.4) that a low probability of a single [favourable] event be disregarded.

**4.** I refer to the German translation of his *Treatise* (1912, p. 148) from the Russian edition of 1908 and to p. 241 of its last edition of 1924.

**5.** The term "random quantity" appeared at the end of the 19$^{th}$ century (Vasiliev 1885, pp. 127 – 131; **Nekrasov** 1888, p. 77) whereas the English expression "random magnitude" was possibly introduced later (Whitworth 1901, p. 207). I had not, however, seen the previous editions of that book. The present English term is of course *random variable*.

## 14. Markov, Liapunov, Nekrasov

I consider here the work of three outstanding scholars; with regard to **Markov** and **Nekrasov**, however, qualification remarks will follow.

### 14.1. Markov: General Scientific Issues

I consider his main results in § 14.2; here, I briefly treat some additional matters, but see § 15.1.3 for his study of statistical series. Finally, in § 14.3 I describe Markov's personality.

1) History of the theory of probability. Markov undoubtedly paid attention to it. He investigated the **Bernoulli** LLN (§ 3.2.3); in 1913 he initiated a jubilee meeting of the Petersburg Academy of Sciences celebrating the bicentenary of that law, as well as the publication of a Russian translation of pt. 4 of the *Ars Conjectandi* (see § 3). Markov left several statements about the history of the **Bienaymé** – **Chebyshev** inequality (§ 10.2-4), somehow argued for the second **Gauss** justification of the MLSq (as mentioned in § 9A.6-1), introduced an apt term, "**De Moivre** – **Laplace** limit theorem" (*Treatise* 1924, p. 53) and stressed De Moivre's part in establishing the "**Stirling** formula". This last edition of his *Treatise* includes many interesting historical remarks. As far as his number-theoretic papers collected in his *Selected works* (1951) are concerned, I can at least add that they contain many references to his predecessors.

2) Insurance of life. In his *Treatise* (1900a), Markov described the pertinent theory but did not add any new findings. However, he actively collaborated with pension funds scrupulously considering all practical details of their work (Sheynin 1997d), and in 1906 he published two newspaper articles destructively criticizing a proposed scheme for insuring children (reprinted in the same later translated article).

3) Calculations. "Markov liked calculating and was good at it" (**Linnik** et al 1951, p. 615). In the theory of probability, most important is his table of the normal distribution (1888) giving it to 11 digits for the argument $x = 0 \ (0.001) \ 3 \ (0.01) \ 4.8$ and adduced the differences of all the necessary orders (for example, with the first three differences for $x \ 2.649$). According to a reputed reference book (Fletcher et al 1962), two tables of the normal distribution, one of them Markov's, and the other, published ten years later, remained beyond compare up to the 1940s. In an indirect way, Markov (1899b, p. 30) made known his attitude toward calculation:

*Many mathematicians apparently believe that going beyond the field of abstract reasoning into the sphere of effective calculations is humiliating.*

4) Correlation theory. In § 10.6-3 I indicated that statisticians had doubted it. The same was true with regard to Markov. **Slutsky** (1912a) had collected and generalized the relevant findings of the Biometric school, and even a few decades later **Kolmogorov** (1948) called his book important and interesting. Markov, however, did not duly estimate it. He mentioned it in three letters to **Chuprov**, all written in 1912 (Ondar 1977/1981, pp. 53 – 58), and he (p. 53) stated that it interested, but did not "attract" him, and (p. 58) did not "like it very much".

Also in 1912, Slutsky exchanged a few letters with Markov and, in particular (Sheynin 1990c/2011, p. 64) stated:

*I believe that the shortcomings of* **Pearson's** *exposition are temporary and of the same kind as the known shortcomings of mathematics in the 17$^{th}$ and 18$^{th}$ centuries. A rigorous basis for the work of the geniuses was built only post factum, and the same will happen with Pearson. I took upon myself to describe what was done. Sometime A. A. Chuprov will set forth the subject of correlation from the philosophical and logical point of view, and describe it as a method of research. An opportunity will present itself to a ripe mathematical mind to develop the mathematical basis of the theory.*

In a few years Markov (1916a/1951, p. 533) critically mentioned the correlation theory:

*Its positive side is not significant enough and consists in a simple usage of the method of least squares to discover linear dependences. However, not being satisfied with approximately determining various coefficients, the theory also indicates their probable errors and enters here the region of fantasy, hypnosis and faith in such mathematical formulas that, in actual fact, have no sound scientific justification.*

Now, discovering dependences, even if only linear, is indeed important; and the estimation of plausibility of the results obtained is an essential part of any investigation. True, at the time such estimation had not been done properly. Considering a paper of a contemporary Russian author, Markov (Ibidem pp. 534 – 535) pointed out an obvious senselessness: the calculated correlation coefficient was 0.09 with probable error 0.14. In addition, these figures greatly changed when Markov left aside some of the observations made use of. However (**Linnik**; see his comment on that paper (Markov 1951, p. 670), without knowledge of the distribution of the population, the sample correlation coefficient cannot properly estimate the general coefficient.

Nevertheless, I quote Bernstein (1928/1964, p. 231):

*Excluding biological applications, most of its* [of the correlation theory] *practical usage is based on misunderstanding*.

He had not regrettably elaborated.

5) Principles of the theory of probability. In essence, Markov left that issue aside. Thus, in the German edition of his *Treatise* (1912, translated from the Russian edition of 1908, p. iii) he declared that he did not discuss it in detail. At about the same time, he (1911c/1981, pp. 149 – 150) pessimistically estimated suchlike efforts:

*I shall not defend these basic theorems connected to the basic notions of […] equal probability, of independence of events, and so on, since I know that one can argue endlessly on the basic principles even of a precise science such as geometry.*

Markov (*Treatise*, 1908, p. and 1924, c. 2) also stated, somewhat indefinitely, that

*Various concepts are defined not by words, each of which can in turn demand definition, but rather by* [our] *attitude towards them ascertained little by little.*

Apparently, some (not *various*) concepts must be admitted without definition; there are no modern definition either of *point* or *straight line*.

It ought to be added, however, that, except for the axiomatic approach, only **Mises** was able to abandon the classical definition of probability (see also § 7.4). Note, however, that Markov, like a student of **Chebyshev**, underrated the then originating axiomatic direction of probability as well as the theory of functions of a complex variable (**A. A. Youshkevich** 1974, p. 125).

On p. 10 of his *Treatise* (1924) Markov formulated the following *axiom*: If there are several equally possible events, some of them favourable, the others not, with regard to event *A*, then, after *A* occurs, the unfavourable events "fall through" whereas the others remain equally possible. I do not see how can it be otherwise. Then, on pp. 13 – 19 Markov proved the addition and the multiplication theorems (in a rather complicated way and mentioning his axiom) and on p. 24 concluded that these theorems along with his axiom serve

*As an unshakeable base for the calculus of probability as a chapter of pure mathematics.*

So here we are! His axiom (many pages apart, not displayed in the text and difficult to find!), never mentioned by any later author, allegedly transformed the theory of probability[1] …

6) Mathematical statistics. In 1910 Markov (Ondar 1977/1981, p. 5) had denied **Pearson**, but by the end of his life he somewhat softened his attitude. Here is a passage from **Chuprov's** letter, written apparently in 1924, to Isserlis (Sheynin 1990c/2011, p. 76):

*Markov regarded Pearson, I may say, with contempt. Markov's temper was no better than Pearson's, he could not stand even slightest contradictions either* [§ 14.3]. *You can imagine how he took my persistent indications to the considerable scientific importance of Pearson's works. My efforts thus directed were not to no avail as proved by* [Markov 1924]. *After all, something* [Pearsonian] *was included in the field of Markov's scientific interests.*

Chuprov (1925b/1981, p. 155) also published a review of the mentioned edition of Markov's *Treatise*. Here, I only cite his reasonable criticism of Markov's treatment of correlation theory:

*The choice of questions on which attention is concentrated is fortuitous, their treatment within the bounds of the chapter on the method of least squares is incomplete, the connection made between the theory of correlation and the theory of probability is inadequate…*

Now, what statistical innovations had Markov included in this last edition of his *Treatise*? A study of statistical series and of the Pearsonian correlation theory (§ 14.2-1). He considered linear correlation and applied the MLSq for determining the parameters of the lines of regression and discussed the case of [random variables] with densities of their distributions being quadratic forms. Markov also included a general reference to **Slutsky** (1912a) and certainly did not repeat his earlier harsh words about imagination, hypnotism, etc.

Below (§ 14.2-1) I shall add that Markov paid no attention either to the chi-squared test or to the Pearsonian curves.

7) Teaching probability theory in school. In 1914 **Nekrasov** made an attempt to introduce probability into the school curriculum. Markov, who could not stand him at all, either as a man, or as a mathematician, was not invited to the pertinent discussion by correspondence, but he voiced his opinion in an ad hoc paper (1915a). He sharply protested against the concrete school programme proposed by Nekrasov, but, as I understand him, did not object to the very principle. In 1914, he published a relevant newspaper article (reprinted in Sheynin 1993a, p. 200), and in 1916 he was member of the Commission established by the Academy of Sciences to study Nekrasov's proposal. Its review was extremely negative (*Report* 1916) both with respect to Nekrasov's programme and to his understanding of the main issues of mathematical analysis, cf. Nekrasov's statement about the concept of limit in § 14.5.

8) Methodological issues. Many authors praised the methodological value of Markov's contributions. **Bernstein** (1945/1964, p. 425) stated that Markov's *Treatise* and memoirs were "specimens of preciseness and lucidity of exposition". **Linnik** et al (1951, p. 615) maintained that Markov's language was distinct and clear, and that he thoroughly trimmed the details. A striking example proving the opposite is Markov's failure to discuss the adjustment of direct conditional observations (§ 9A.4-9) in his *Treatise*. And I do not trust **Chuprov** (1925b/1981, p. 154) who thought that the exposition in Markov's *Treatise* was "transparently clear".

Except Markov himself (§ 14.2-1) the only author with whom I agree is Idelson (1947, p. 101). He remarked that the chapter on the MLSq in the *Treatise* (1924) was ponderously written. Indeed, Markov's general rule was to rewrite his formulas rather than to number, and then to refer to them. Thus, on pp. 328 – 330 of the *Treatise* a long equality appeared five times in succession! Then, he disregarded demonstrative pronouns. On p. 328, for example, he wrote: "The choice of coefficients [a displayed line of them followed] is at our disposal. We shall subject the coefficients [the same line was repeated] to two conditions …"

Then, Markov refused to apply the term *random magnitude* (as it has been called in Russia), see § 14.2-1, and the expressions *normal law* and *coefficient of correlation* were likewise absent in his works. As the Russian saying goes: "The whole company broke step, the lieutenant alone is in step". Markov's literary style was pedestrian and sometimes hardly understandable (1906/1951, p. 341) and, from one edition to another, the structure of his *Treatise* became ever more complicated. A few more remarks are in §14.2-1.

## 14.2. Markov: His Main Investigations

1) Mathematical treatment of observations. **Linnik** et al (1951, p. 637) believed that, when substantiating the MLSq, Markov had "in essence" introduced concepts "equivalent" to the modern concepts of unbiased and effective statistics for estimating parameters of the laws of distribution. Markov, however, only indirectly estimated parameters (he never used such an expression), and, which is more important, it is just as possible to attribute those concepts to **Gauss**. Nor do I agree with Idelson (1947, p. 14) who mentioned the Gauss method, developed by Markov "up to the highest logical and mathematical perfection". In § 9A.6-1 I mentioned Markov's resolute stand for the second Gauss substantiation of the MLSq; this (and his remark that the [consistency] of the arithmetic mean was inadequate, see § 11.2-8) is all that he really accomplished here. **Neyman** (1934, p. 595) erroneously attributed that justification to Markov and **F. N. David** & Neyman (1938) repeated that statement, but then Neyman (1938/1952, p. 228) recognized his mistake. Nevertheless, the mysterious *Gauss – Markov theorem* is still alive and kicking. It was Scheffé (1959, p. 14) who introduced that term (H. A. David, see H. A. David & Edwards 2001, p. 218) although Plackett (1949) had noted that mistake.

And I add now that Markov (1899a/1951, p. 246) had at the same time denied any optimality of the MLSq, so why did he feel necessary to substantiate it at all?

In his *Treatise* (1900) Markov in essence combined the treatment of observations with the study of correlation (§ 14.1-4), statistical series and interpolation; this, perhaps, reflected his attempt to include the MLSq into the then originating mathematical statistics, but his innovation was methodically doubtful.

The discussion of statistical series was rather involved for an educational aid and did not mention **Chuprov's** relevant papers (1916; 1918 – 1919) the first of which Markov himself had communicated to the *Izvestia* of the Petersburg Academy of Sciences. Note that Chuprov (1925b/1981, pp. 154 and 155) politely remarked that Markov had left out the works of other authors not belonging to the "stream" of his own contributions. I would say that this criticism was too mild.

In connection with statistical series Markov (*Treatise* 1924, pp. 349 – 353) considered Weldon's experiment with 26,306 throws of 12 dice (**K. Pearson** 1900) and decided, after applying the CLT and the **Bayes** theorem with transition to the normal law, that the probability of a 5 or a 6 was higher than 1/3. Unlike Pearson, he had not used the chi-squared test and he could have left an impression that (although suitable for a small number of trials as well) it was not needed at all.

As to interpolation, the only point of contact with the MLSq was the calculation of the empirical coefficients according to the principle of maximal weight.

Markov passed over in silence the Pearsonian curves perhaps owing to their insufficient substantiation. However, again in 1924, he reprinted the Introduction to the edition of 1913 where he had stated that the use of approximate methods in applied mathematics was unavoidable even when an estimation of their error was impossible, and (1915a, p. 32) maintained that "empirical" formulas did not demand theoretical proof. I think that Markov followed here (as he did with his chains which he left without any applications to natural sciences) his own rigid principle hardly worthy of exact imitation (Ondar 1977/1981, Letter 44 to Chuprov of 1910): "I shall not go a step out of that region where my competence is beyond any doubt".

The explication of the MLSq proper was involved, and Markov himself knew it. In a letter of 1910 to Chuprov he (Ibidem, p. 21) wrote: "I have often heard that my presentation is not sufficiently clear". In 1893, his former student, Koialovitch (Sheynin 2006b, p. 85), writing to Markov, formulated some puzzling questions:

*As far as I understand you, you consider each separate observation as a value of a possible result. Thus, a series of results [...] is possible for each measurement, and one of them is realized. I am prepared to understand all this concerning one observation. However, if there are, for example, two observations, then I cannot understand the difference between the series of all the possible results of the first observation [...] and the similar series for the second measurement [...] The problem will certainly be solved at once if you say that the probabilities of the same error in these two series are different, but you will hardly want to introduce the notion of probability of error in your exposition.*

See also Ermolaeva (2009) and my Russian contribution on Koialovich in **S, G,** 85. The situation had not improved with time. Contrary to what he himself (*Treatise* 1924, pp. 323 and 373) stated while following **Chebyshev** (1879 – 1880/1936, p. 227, translation p. 208), Markov (1924, pp. 327 and 374) maintained that only one possible observation corresponded to each actually made. He never clearly explained that observational errors were [random variables] and that a series of observations was a [random] sample and had a density function. I note that Laplace (§ 7.2-4) mentioned *numerous and not yet made observations*…
Wherever possible, Markov (Ondar 1977/1981, Letter 53 to Chuprov of 1912) excluded "the completely undefined expressions *random* and *at random*". Instead, he added, he introduced an appropriate explanation in each particular case. However, at least sometimes he simply wrote *indefinite* which was much worse; incidentally, the translators of Ondar (1977) modernized Markov's letters by rendering *indefinite* as *random*. Cf. Note 5 to Chapter 13.

The chapter on the MLSq in Markov's *Treatise* was hardly inviting either for mathematicians or geodesists. Both would have been disappointed by the lack of discussion of **Pearson's** work whereas the latter, in addition, had not needed interpolation or investigation of statistical series but would have wished to see much more about correlation. And the absence of the **Gauss** brackets (Note 19 in my Chapter 1) as well as the appearance of the long-ago dated term *practical geometry* instead of *geodesy* (p. 462) would have annoyed them.

I also mention that Markov destructively criticized a paper (**Galitzin** 1902) devoted to the study of the solidity of glass tubes. His review was extant as a manuscript and I (Sheynin 1990b) published it. Markov had not applied any new method, but he thoroughly treated Galitzin's data and allowed for every possible circumstance. It was in connection with the discussion of Galitzin's paper that Markov stated his opinion about the "**Bredikhin** rule" (§ 10.8.4).

2) The LLN. Markov (1906/1951, p. 341) noted that the condition

$$\lim E \frac{1}{n^2} [(\xi_1 + \xi_2 + \ldots + \xi_n) - (E\xi_1 + E\xi_2 + \ldots + E\xi_n)]^2 = 0, n \to \infty \qquad (1)$$

was sufficient for the sequence $\xi_1, \xi_2, \ldots, \xi_n, \ldots$ of [random variables] to obey the LLN; or, in accordance with his formula, to comply with the condition

$$\lim P\{|(\xi_1 + \xi_2 + \ldots + \xi_n) - (E\xi_1 + E\xi_2 + \ldots + E\xi_n)| < \varepsilon\} = 1, n \to \infty.$$

Then Markov (Ibidem, pp. 342 – 344; *Treatise*, 1913, pp. 116 – 129) derived a few relevant sufficient conditions for sequences of independent, and, especially, dependent random variables and (1906/1951, p. 351; *Treatise* 1913, p. 119; *Treatise*, 1924, p. 174) provided examples of sequences not obeying the law, and, in addition (*Treatise*, 1913, p. 129), proved that apart from condition (1), independent variables obeyed the LLN if, for every *i*, there existed the moments

$$E\xi_i = a_i, E|\xi_i - a_i|^{1+\delta} < C, 0 < \delta < 1.$$

In connection with his investigations of the LLN Markov (*Treatise*, 1900; p. 86 in the edition of 1924) had proved that, for a positive random variable $\xi$,

$$P(\xi \leq t^2 E\xi) > 1 - 1/t^2$$

and **Bortkiewicz** (1917, p. 36) and **Romanovsky** (1925a; 1925b) called this inequality after Markov.

3) [The CLT]. As I mentioned at the end of § 13.1, Markov specified the conditions of theorem (13.6) proved by **Chebyshev**. He (1898/1951, p. 268) considered independent [random variables] $u_i$ with zero expectations[2] and, following Chebyshev, supposed that, for finite or[3] infinite values of $k$,

$$\lim |E u_n^k| < +\infty, \quad n \to \infty. \tag{2}$$

In addition, Markov, however, demanded that

$$\lim E u_n^2 \to 0, \quad n \to \infty. \tag{3}$$

He several times returned to the CLT.

a) For equalities (13.7) to hold, he (1899a/1951, p. 234) assumed condition (3) and, for the transition to the new theorem he (p. 240) additionally introduced two restrictions: as $n \to \infty$,

$$\lim E[(u_1 + u_2 + \ldots + u_n)^2] = \infty,$$
$$\lim [E(u_1 + u_2 + \ldots + u_n)^2 / n] \to \infty. \tag{4; 5}$$

b) Later Markov (1907, p. 708) again proved formula (13.7). Referring to his papers (1898; 1899a), he now introduced conditions (2) (for finite values of $k$) and (5) but did not restrict the values of $u_i$. On his next page Markov abandoned condition (5) "if only"

$$\lim E u_n^2 = \infty, \quad n \to \infty \tag{6}$$

and the values of $u_i$ remained finite. Restrictions (4) and (6) certainly coincided.

Finally, Markov (1908a) essentially extended the applicability of the method of moments by replacing his conditions by **Liapunov's** single restriction (1901a/1954, p. 159)

$$\lim \frac{\sum E |u_i|^{2+\delta}}{\left(\sum \text{var } u_i^2\right)^{1+\delta/2}} = 0, \quad \delta > 0, \quad n \to \infty. \tag{7}$$

In 1913 Markov included a modified version of his last-mentioned study in his *Treatise* (Markov 1924; reprint: 1951, pp. 319 – 338).

In connection with condition (3) Markov (1899c, p. 42) mentioned the example provided by **Poisson** (1824, § 10). The latter proved that the limiting distribution of the linear form

$$L = \xi_1 + 1/3 \, \xi_2 + 1/5 \, \xi_3 + \ldots$$

of [random variables] $\xi_i$ with density $e^{-2|\xi|}$ was

$$\lim P(|L| \le c) = 1 - (4/\pi)\arctan e^{-2c}, \quad n \to \infty.$$

In this example

$$\lim \text{var}[\xi_n/(2n-1)] = 0, \quad n \to \infty.$$

Markov himself (1899a/1951, pp. 242 – 246) also provided an example in which the condition (3) had not held and the CLT did not take place and he mentioned Poisson without adducing the exact reference.

The appearance of condition (5) remains, however, unclear. **Nekrasov** (1900 – 1902, 1902, pp. 292 and 293) introduced it for independent variables instead of restriction (3). **Liapunov** (1901a/1954, p. 175) maintained that it was not sufficient (but was he then acquainted with the third part of Nekrasov's contribution?) and mentioned Markov's examples. Seneta (1984, p. 39) indicated, however, that Markov's published papers had not contained such examples and that condition (3) was necessary and sufficient for the CLT in the case of uniformly restricted variables.

4) Markov chains. This term is due to **Bernstein** (1926, §16); Markov himself (1906/1951, p. 354) called them simply *chains*. He issued from a paper by **Bruns** of the same year, but the prehistory of Markov chains is much richer. Here are the main relevant issues.

a) The **Daniel Bernoulli** – **Laplace** urn problem, the predecessor of the **Ehrenfests'** model (§ 7.1-3).

b) The study of the **Brownian** motion (Brush 1968).

c) The problem of the extinction of families (§ 10.2-7).

d) The problem of random walks (Dutka 1985).

e) Some of **Poincaré's** findings (§ 11.2-6).

f) The work of **Bachelier**, for example, Bachelier (1900) on financial speculations, also see Courtault et al (2000) and Taqqi (2001).

Markov (1906/1951, pp. 345 and 354) considered simple homogeneous chains of random events and discrete [random variables] and proved that the LLN was applicable both to the number of successes and to the sequences of these variables. Later he (1910/1951, p. 476) extended the first of these findings to simple nonhomogeneous chains.

Markov proved the CLT for his chains by means of condition (13.7). He considered simple homogeneous chains of events (1906) and of [random variables] (1908b); simple nonhomogeneous (1910) and complex homogeneous (1911a; 1911b) chains of [random variables]; simple homogeneous chains of indirectly observed events (1912a). While studying the chains, Markov established important ergodic theorems but had not paid them any special attention; in this connection, I mentioned one of his solved problems in § 7.1-3.

It is difficult to imagine that Markov had not grasped the essential importance of the chains for various applications, but he did not say anything about that, and his only relevant and mostly methodical example (*Treatise*, 1913) was a study of the alternation of consonants and vowels in the Russian language, see Petruszewycz (1983). In 1910 Markov himself, in his letters to **Chuprov**, remarked more than once that he was restricting his field of work by what was well known to him, see § 14.2-1. I note also that at the time physics was not yet duly studied in Russia (**Kolmogorov** 1947, p. 59).

5) In concluding, I ought to add that Markov widely applied the method of moments, and only he who repeats some of his

investigations (for example, of his study of the limiting behaviour of the terms obtained by decomposing algebraic fractions), will be able to estimate the obstacles which he overcame. **Bernstein** (1945/1964, p. 427) contrasted Markov and **Liapunov**. The latter had applied the classical transcendental analysis as developed by that time whereas the method of moments, Bernstein maintained,

*Did not facilitate the problem [of proving the CLT] but rather transferred all its difficulties elsewhere*.

Soon Bernstein (1947, p. 44) remarked that Chebyshev had expressed a negative opinion about the method of moments for characteristic functions which in those times had not yet answered the requests of mathematical rigour.

It might be imagined, however, that Markov wished to ascertain how powerful was the method of moments; he himself (*Treatise*, 1913, p. 322) indicated that Liapunov had "shaken" the importance of the method of moments and that he, Markov, therefore decided to prove the CLT anew (see above).

Liapunov (1900/1954, p. 125) called the Markov proof of the CLT *too complicated and unwieldy* because of its connection with *a special theory*. However, Krein (1951, pp. 8 – 9) noted that the method of moments had not yet *become useless*. It is still being applied (A. D. Soloviev 1997/2008, p. 355):

*Considerably improved since Markov's time*, [it] *is still used* […] *for solving such problems in which the moments of random variables are derived much easier than their distributions*.

### 14.3. Markov: His Personal Traits

The sequel will substantiate my decision to discuss Markov's personality at the end of this section. For his biography see **Markov Jr** (1951), a noted mathematician in his own right, and Grodzensky (1987). These authors describe many cases of his principled stand on burning social and political issues whereas Grodzensky published many of his newspaper letters[4]. In particular, Markov struggled against anti-Semitism (a subject which Markov Jr was unable to describe in 1951!) and denounced the Russian Orthodox Church, see also Sheynin (1989a, pp. 340 – 341; 2007c). The Press used to call him *Militant academician* (Nekrasov 1916a, p. 9) and *Andrew the Furious* (Neyman 1978). Other sources include Gnedenko & Sheynin (1978), and Seneta (2001). See Markov's correspondence in **S, G,** 16. I have mentioned one of his letters in §14.1-7.

In 1901 Tolstoy was excommunicated from the Church. Then, during his last days, the Most Holy Synod discussed whether he should be "admitted to the bosom of the Church" and decided against it (Anonymous 1910; archival sources of the Synod). This goes to show that in 1912 Tolstoy's excommunication was likely well remembered.

Yes, in 1912 Markov submitted a request to the Synod for excommunication. He quoted his *Treatise* to the effect that

> *We should regard stories about incredible events allegedly having occurred in bygone times with extreme doubt*

and added that he did not

> *Sympathise with religions which, like Orthodoxy, are supported by, and in turn lend their support to fire and sword.*

His request was not granted; the Synod resolved that Markov "had seceded from God's church" (Emeliakh 1954, pp. 400 – 401 and 408). This paper provides, in particular, the text of an internal letter of the Synod which stated that excommunication would have been too honourable for Markov. I think that it was the notorious Beilis case (a blood libel) which prompted Markov's request.

Markov's last and certainly useless protest took place in 1921 (Grodzensky 1987, p. 137) when 15 professors of the Petrograd University declared that applicants ought to be chosen according to their knowledge rather than to class or political considerations; Markov was the first to sign their statement.

Markov's attitude towards other scholars had not rarely been wrong. The excessive sharpness of his statements is generally known and here is a passage from a letter of **Zhukovsky**, the then President of the Moscow Mathematical Society, of 23.11.1912 to **Markov** (**S, G,** 16, section 5.1, Letter 47):

> *I cannot fail to reproach you for the expressions concerning the honourable **Sergei Alekseevich Chaplygin** in your letter. They can hardly be called proper.*

Chaplygin (1869 – 1942) was cofounder of aero hydrodynamics (and an active member of the Society). Markov's letter reflected the polemic between him and **Nekrasov** in the Society's periodical, *Matematichesky Sbornik.*

The second and last example (**K. A. Andreev's** letter of 1915 to **Nekrasov**; Sheynin 1994f, p. 132): Markov

> *Remains to this day an old and hardened sinner in provoking debate. I had understood this long ago, and I believe that the only way to save myself from the trouble of swallowing the provocateur's bait is a refusal to respond to any of his attacks …*

Andreev had published a posthumous manuscript of **V. G. Imshenetsky** and Markov severely criticized its incompleteness. Nevertheless, Markov himself, soon before his death, agreed to publish his last, and also incomplete manuscript (Besikovitch 1924, p. XIV).

Markov had been sharply and even provocatively behaving towards Nekrasov whom he *overpowered by many rude postcards* (Nekrasov's complaint of 1915 to the permanent secretary of the Academy of Sciences (Nekrasov 1916a, pp. 56 – 62; **S, G,** 16).

In his own scientific work, Markov had been too rigid, see §§ 14.1-4, 14.1-5 (the statement of A. A. Youshkevich), 14.2-1 and 14.2-4, which negatively influenced his work. During his last years, in spite of extremely difficult conditions of life in Russia and his worsened health, he completed (?) the last edition of his *Treatise* (published posthumously) but hardly sufficiently described the new findings of the Biometric school. Such scholars as **Yule** and **Student** (**Gosset**) were not mentioned, **Fisher** was also ignored.

### 14.4 Liapunov

The theory of probability remained an episode in his scientific work. He (1900; 1901a) proved the [CLT] assuming a single condition (7). I briefly repeat (**Bernstein** 1945/1964, pp. 427ff) that a characteristic function determines the sought law of distribution independently from the existence of the relevant moments and that the expansion in powers of *s* which **Chebyshev** (§ 13.1-4) made use of did not anymore lead to difficulties after replacing that argument by *is*. Liapunov proved that under his condition the characteristic function of a centred and normed sum of random variables tended to the characteristic function of a normed normal law. I also mention **Lindeberg** (1922b, p. 211) whose proof of the CLT was simpler and became better known[5]. He referred to his previous paper (1922a) and continued:

*I see now that already Liapunov had explicated general findings which not only surpass the results achieved by* **Mises** [a reference to his article of 1919 followed] *but which make it possible to derive most of what I have established.* […] *The study of Liapunov's work prompted me to check anew the method that I have applied.*

A special point is connected here with the CLT for large deviations. **Chebyshev** thought that the limits of integration,   and   in formula (13.6) describing that theorem, were "any". **Nekrasov** (1911, p. 449) arbitrarily interpreted that expression as "variable". I discuss Nekrasov in § 14.5; here, I say that he could have well indicated that, on the contrary, he had generalized the Chebyshev theorem. In his previous polemic paper Liapunov (1901b, p. 61) declared that he had assumed that these limits were given beforehand and that otherwise the probability, written down in the left side of formula (13.6), could have no limit at all, – but nevertheless be asymptotically expressed by the normal law of distribution[6].

### 14.5. Nekrasov

Nekrasov's life and work are clearly separated into two stages. From 1885 and until about 1900 he had time to publish remarkable memoirs both in Russia and Germany and to become Professor and Rector of Moscow University; I mentioned him in § 10.2. In 1898 he sketched the proof of the [CLT] for sums of [lattice random variables]. Then, however, his personality changed. His writings (only on probability and statistics) became unimaginably verbose, sometimes obscure and confusing, and inseparably linked with ethical, political and religious considerations. Here is a comparatively mild example (1906, p. 9): mathematics accumulated

*Psychological discipline as well as political and social arithmetic or the mathematical law of the political and social development of forces which depend on mental and physiological principles.*

Furthermore, Nekrasov's work began to abound with elementary mathematical mistakes and senseless statements. Thus (Nekrasov 1901, p. 237): it is possible to assume roughly that $x^n$, $n > 0$, is the limit of sin $x$ as | |   0, and

*The conclusions made by* [**Chebyshev**, **Markov** and **Liapunov**] *never differ from such an understanding of limit*.

I provide a second and last out of many possible illustrations from Nekrasov's letter of 20.12.1913 to Markov (Archive, Russian Acad. Sci., Fond 173, inventory 1, 55, No. 5; **S, G,** 16, section 5.1, Letter 17):

*I distinguish the viewpoints of **Gauss** and **Laplace** [on the MLSq] by the moment with regard to the experiment. The first one is posterior and the second one is prior. It is more opportune to judge à posteriori because more data are available, but this approach is delaying, it lags behind, drags after the event.*

The attendant reasons for such a change were Nekrasov's religious upbringing (before entering Moscow University he graduated from a Russian Orthodox seminary), his work from 1898 onward as a high official at the Ministry of People's Education[7], and his reactionary views. At least once Nekrasov (A. V. Andreev 1999, p. 103) mentioned the *Integral Knowledge* of the religious philosopher **V. S. Soloviev** (1853 – 1900) and it is opportune to quote Soloviev's pronouncement (Radlov 1900, p. 787) with which, in actual fact, Nekrasov became absorbed: "veritable knowledge is a synthesis of theology, rational philosophy and positive science". Andreev indeed maintains that Nekrasov became split between mathematics and such philosophy. **Bortkiewicz** (1903) notes that Nekrasov "especially often mentioned Soloviev in vain", – and sometimes justifiably, as I am inclined to believe.
And here is Liapunov's suitable remark (1901b, p. 63):

*All of Nekrasov's objections are based on misunderstandings* […]. *Some of them are simply unjustified statements* […], *others either have no relations to the papers which he criticizes or are extremely indefinite*.

Concerning Nekrasov's social and political views I turn to his letter of 1916 to **P. A. Florensky** (Sheynin 1993a, p. 196): "the German – Jewish culture and literature" pushes "us" to the crossroads. World War I was then going on, but that fact only to some extent exonerates Nekrasov. Florernsky was an eminent theologian, a philosopher of mathematics and an out-and-out anti-Semite. I have seen in the internet a story about his statement to the effect that, had he been a Jew, he would himself killed the Christian boy. He spoke about the notorious blood libel of Beilis (who was later acquitted).

I shall now dwell on some concrete issues.

1) Teaching the theory of probability. In § 14.1-7 I mentioned Nekrasov's proposal for teaching probability in school and the rejection of the curriculum drawn up by him. I add now that already in 1898 Nekrasov made a similar proposal concerning the Law Faculty of Moscow University, also rejected or at least forgotten; however (Sheynin 1995a, p. 166), during 1902 – 1904 the theory of probability was not taught there even at the Physical and Mathematical Faculty, and hardly taught during 1912 – 1917.

2) The MLSq. Nekrasov (1912 – 1914) mistakenly attributed to **Legendre** an interpolation-like application of the method and (1914) acknowledged his failure to notice, in 1912, the relevant work of Yarochenko (1893a; 1893b), but still alleged (wrongly) to have considered the issue in a more general manner. Yarochenko justified the arithmetic mean and the MLSq in general by a reference to Chebyshev's memoir (1867), – that is, by the **Bienaymé – Chebyshev** inequality (§ 9A.4-3). Note that the first such statement appeared simultaneously with the Chebyshev memoir (Usov 1867). Recall also Nekrasov's strange pronouncement about **Laplace** and **Gauss** quoted above.

3) [The CLT]. It was Nekrasov who had considered the CLT for large deviations, – for the case that began to be studied only 50 years later. Suppose that independent [lattice] random variables (linear functions of integral variables) $\xi_i$, $i = 1, 2, …, n$, have finite mean values $a_i$ and variances $\sigma_i^2$ and

$$m = \xi_1 + \xi_2 + … + \xi_n.$$

Denote

$$|x(m)| = |m - \sum a_i|/(\sum \sigma_i^2)^{1/2}.$$

Nekrasov restricted his attention to the case in which $|x| < n^p$, $0 < p < 1/6$ and stated that, for all values of $m_1$ and $m_2$ which conformed to that condition,

$$P(m_1 < \xi_1 + \xi_2 + … + \xi_n < m_2) \sim [1/\sqrt{2\pi}] \int \exp(-t^2/2)dt.$$

The limits of integration were $x(m_1)$ and $x(m_2)$ respectively.

In all, Nekrasov (1898) formulated six theorems and proved them later (1900 – 1902). Neither Markov, nor Liapunov had sufficiently studied them; indeed, it was hardly possible to understand him and **A. D. Soloviev** (1997/2008, p. 359) reasonably stated:

*I am firmly convinced that no contemporaneous mathematician, or later historian of mathematics had* [has] *ever studied it* [the memoir (1900 – 1902)] *in any detail*.

He himself was only able to suggest that Nekrasov had indeed proved his theorems and he reminded his readers that **Markov** had indicated some mistakes made by Nekrasov. Furthermore, Soloviev (pp. 356 – 357) remarked that Nekrasov had wrongly understood the notion of lattice variables (not like I described above). In his general conclusion Soloviev (p. 362) stated that Nekrasov had imposed on the studied variables an excessively strict condition (the analyticity of the generating functions in some ring, which was much stronger than presuming the existence of all of the moments) and that it was generally impossible to check his other restrictions. Both Soloviev, and the first of the modern commentators, Seneta (1984, §6), agree in that Nekrasov's findings had not influenced the development of the theory of probability[8]. This regrettable outcome was certainly caused both by Nekrasov's inability to express himself intelligibly and by the unwieldiness of his purely analytical rather than stochastic approach (Soloviev, p. 363).

Chuprov (Bortkevich & Chuprov 2005, letter to Bortkiewicz of 22 Nov. 1896) actually proved that Nekrasov did not study mathematical statistics at all or even sufficiently mastered the theory of probability:

*Noting the word <u>dispersia</u> [variance] in my [student] composition, he timorously asked me:<u> So do you really apply the theory of probability to the dispersion of light</u>?*

### Notes

**1.** Still, Markov had a forerunner: Donkin (1851) expressed a quite similar statement apparently based on the principle of insufficient reason which Keynes (1921/1973, p. 44) later renamed *the principle of indifference*. Boole (1854b/2003, p. 163) thought that *perhaps* [Donkin's principle] *might* […] *be regarded as axiomatic*. Neither Donkin, nor Keynes mentioned that that principle or axiom transforms the theory of probability.

**2.** Until he began to study his chains, Markov had always introduced these two conditions. In one case (1899a/1951, p. 240) he apparently had not repeated them from his p. 234.

**3.** A misprint occurred in the Russian translation of the French original.

**4.** Grodzensky regrettably had not adduced an index of the letters he discovered and did not indicate which of them had indeed been published at once.

**5.** Thus, **Gnedenko** (1954/1973, pp. 254 – 259) proves the theorem under the **Lindeberg** condition and then explains that the **Liapunov** restriction leads to the former.

**6.** Liapunov's correspondence with **K. A. Andreev** in 1901 (Sheynin 1989b) testifies that he had initially wished to publish his note in the *Matematichesky Sbornik*, that the leadership of the Moscow Mathematical Society (**Bugaev**, **Nekrasov** (!)) opposed his desire, and that he essentially expanded his first draft on Andreev's advice.

**7.** Here is **K. A. Andreev's** opinion (letter of 1901 to **Liapunov**; Gordevsky 1955, p. 40): Nekrasov

*Reasons perhaps deeply, but not clearly, and he expresses his thoughts still more obscurely. I am only surprised that he is so self-confident. In his situation, with the administrative burden weighing heavily upon him, it is even impossible, as I imagine, to have enough time for calmly considering deep scientific problems, so that it would have been better not to study them at all.*

Here is an example of Nekrasov's deep thought (1916b, p. 23): he mentioned almost all the problems of the yet not existing theory of catastrophes and introduced that very term, *catastrophe*.

**8.** It might be added, however, that **Markov** (1912b, p. 215) sometimes considered the refutation of **Nekrasov's** mistaken statements as one of the aims of his work. A similar explanation is contained in one of his letters of 1910 to **Chuprov** (Ondar 1977/1981, p. 5).

## 15. The Birth of Mathematical Statistics
## 15.1. The Stability of Statistical Series

By the end of the 19th, and in the beginning of the 20th century, statistical investigations on the Continent were mostly restricted to the study of population. In England, on the contrary, the main field of application for statistical studies at the time had been biology. It is possible to state more definitely that the so-called Continental direction of statistics originated as the result of the work of **Lexis** whose predecessors had been **Poisson**, **Bienaymé**, **Cournot** and **Quetelet**. Poisson and Cournot (§ 8.6) examined the significance of statistical discrepancies "in general" and assuming a large number of observations, – without providing concrete examples. Cournot (§ 10.3-6) also attempted to reveal dependence between the decisions reached by judges (or jurors). Bienaymé (§ 10.2-3) was interested in the change in statistical indicators from one series of trials to the next one and Quetelet (§ 10.5) investigated the connections between causes and effects in society, attempted to standardize statistical data worldwide and created moral statistics.

All this had been occurring against the background of statements that the theory of probability was only applicable to statistics if, for a given totality of observations, "equally possible cases" were in existence, and the appropriate probability remained constant (§ 10.7).

**15.1.1. Lexis.** He (1879) proposed a distribution-free test for the equality of probabilities in different series of observations; or, in other words, a test for the stability of statistical series. Suppose that there are $m$ series of $n_i$ observations, $i = 1, 2, \ldots, m$, and that the probability of success was constant throughout and equal to $p$. If the number of successes in series $i$ was $a_i$, the variance of these magnitudes can be calculated by two independent formulas (Lexis 1879, § 6)

$$\sigma_1^2 = pqn, \quad \sigma_2^2 = [vv]/(m-1) \qquad (1; 2)$$

where $n$ was the mean of $n_i$, $v_i$, the deviations of $a_i$ from their mean, and $q = 1 - p$. Formula (2) was due to **Gauss**, see (9.6b); he also knew formula (1), see a posthumously published note: W-8, 1900, p. 133. The frequencies of success could also be calculated twice. Note however that Lexis applied the probable error rather than the variance.

Lexis (§ 11) called the ratio

$$Q = \sigma_2/\sigma_1 \qquad (3)$$

the *coefficient of dispersion* perhaps choosing the letter $Q$ in honour of Quetelet. In accordance with his terminology, the case $Q = 1$ corresponded to normal dispersion (with some random deviations from unity nevertheless considered admissible); he called the dispersion supernormal, and the stability of the totality of observations subnormal if $Q > 1$ (and indicated that the probability $p$ was not then constant); and, finally, Lexis explained the case $Q < 1$ by dependence between the observations, called the appropriate variance subnormal, and the stability, supernormal. He did not, however, pay attention to this last-mentioned case.

But how exactly, in Lexis' opinion, could the probability vary? No universal answer was of course possible. He (1876, pp. 220 – 221 and 238) thought that the variations followed a normal law, but then he (1877, § 23) admitted less restrictive conditions (evenness of the appropriate [density function]) and noted that more specific restrictions were impossible. Moreover, I am not sure that Lexis had firmly broken off with previous traditions, see § 10.7-7.

Then, Lexis (1879) discussed this issue once more, and even mentioned "irregular waves" (§ 22), but it is very difficult to follow him. Time and time again he interrupted himself by providing statistical examples and never gave precise formulations. Recall (my § 10.8.4) that **Newcomb** introduced a mixture of normal laws instead of one single distribution and that at the very least his exposition was quite definite. At the same time Lexis made a common mistake by believing that the relation between the mean square error and the probable error remained constant irrespective of the appropriate distribution.

He had not calculated either the expectation, or the variance of his coefficient (which was indeed difficult), neither did he say that that was necessary. Recall (§ 9A.4) that **Gauss**, after introducing the sample variance, indicated that it was [unbiased] and determined its variance. Lexis' main achievement was perhaps his attempt to check statistically some stochastic model; it is apparently in this sense that **Chuprov's** remark on the need to unite him and **Pearson** (§ 15.2) should be understood.

I note finally that Lexis (1879, § 1) qualitatively separated statistical series into several types and made a forgotten attempt to define stationarity and trend.

A French actuary **Dormoy** (1874; 1878) preceded Lexis, but at the time even French statisticians (who barely participated in the development of the Continental direction of statistics) had not noticed his theory. It was Lexis who first discovered Dormoy (Chuprov 1909/1959, p. 236) and Chuprov (1926, p. 198, Russian translation 1960, p. 228, English translation 2004, p. 78) argued that the Lexian theory of dispersion ought to be called after Dormoy and Lexis. He thus opposed a later opinion of **Bortkiewicz** (1930, p. 53) who ranked Dormoy far below Lexis. Be that as it may, later statisticians had only paid attention to Lexis.

**15.1.2. Bortkiewicz.** I mentioned him in § 8.7 in connection with the LLN and, in § 10.7-4, I dwelt on his statement about the estimation of precision of statistical inferences. Of Polish descent, Vladislav Iosifovich Bortkevich was a lawyer by education. He was born and studied in Petersburg, but at the end of the 19$^{th}$ century he continued his education in Germany (he was **Lexis'** student) and in 1901 secured a professorship in Berlin and remained there all his life as Ladislaus von Bortkiewicz. In 1912 the Russian statistician P. D. Azarevich (Fortunatov 1914, p. 237) mentioned him thus: "Each time I see him, I feel sorry that he was lost to Russia. There's a genuine man of science". In a letter of 1905 to **Chuprov** (Sheynin 1990c/2011,

p. 56) Bortkiewicz indicated that in Germany he felt himself "perfectly well", whereas a cataclysm was possible in Russia. Bortkiewicz had indeed published most of his contributions in German (which he knew hardly worse than Russian), but he did not lose his ties with Russia. He (1903) sharply criticized **Nekrasov** for the latter's statements that the theory of probability can soften "the cruel relations" between capital and labour (p. 215) and (p. 219) exonerate the principles of firm rule and autocracy as well as for Nekrasov's "sickening oily tone" (p. 215) and "reactionary longings" (p. 216)[1]. Then, **Slutsky** (1922) referred to a letter received from Bortkiewicz, see also Sheynin (2007f), and, finally, at least during his last years he was connected with the then existing in Berlin Russian Scientific Institute and Russian Scientific Society (Sheynin 2001f, p. 228; Bortkevich & Chuprov, 2005, pp. 9 – 12).

Bortkiewicz achieved interesting findings and his example is extremely instructive since he was not initially acquainted with mathematics. In 1896, in a letter to **Chuprov** (Sheynin 1990c/2011, p. 58), he declared that the differentiation of an integral with respect to its (lower) limit was impossible. Many authors deservedly praised him for his scientific work. Thus (Woytinsky 1961, pp. 452 – 453), he was called "the statistical Pope" whereas H. Schumacher (1931, p. 573) explained Bortkiewicz' attitude towards science by a quotation from the Bible (Exodus 20:3): "You shall have no other gods before me".

Here is how Chuprov's student and the last representative of the Continental direction, **Anderson**[2] (1932, p. 243/1963, Bd. 2, p. 531) described his achievement in studying statistical series:

*Our (younger) generation of statisticians is hardly able to imagine that mire in which the statistical theory had got into after the collapse of the Queteletian system, or the way out of it which only Lexis and Bortkiewicz have managed to discover.*

Yes, but neither the former (§ 15.1.1) nor the latter (see below) had successfully overcome the occurring difficulties. Later Anderson included Chuprov.

Bortkiewicz' work is insufficiently known mostly because of his pedestrian style and excessive attention to details but also since German statisticians and economists of the time (Bortkiewicz was also a celebrated economist) had been avoiding mathematics. He did not pay attention to improving his style. Chuprov (Bortkevich & Chuprov 2005, Letter 35 of 1898) criticized his friend for a methodically unfortunate paper on Pareto but Bortkiewicz refused to mend his ways. And Winkler (1931, p. 1030) quoted a letter from Bortkiewicz (date not given) who was glad to find in him one of the five expected readers of his work of 1923 – 1924 on indexes! Here is what Anderson (1932, p. 245/1963, Bd. 2, p. 533) had to say:

*Bortkiewicz did not write for a wide circle of readers* […] *and was not at all a good exponent of his own ideas. In addition, he made very high demands on the readers' schooling and intellect. With* *stubbornness partly caused by his reclusive life, […] he refused to follow the advice of […] Chuprov …*

His Italian papers of 1908 – 1909 written as a defence of his law of small numbers (see below) are virtually unknown. The initial German text of one of them is kept at Uppsala and copied in **S, G,** 25.

Bortkiewicz had determined $EQ$ and $EQ^2$. Chuprov several times mentioned this fact (Sheynin 1990c/2011, pp. 87, 93, 139) and in 1916 **Markov** (Ondar 1977/1981, p. 93) stated that Bortkiewicz' "research […] while not fully accurate, is significant" and (Markov 1911c/1981, p. 153) that "some" of his relevant studies "deserve greater attention".

It is most interesting that Bortkiewicz introduced his law of small numbers (1898a) for studying the stability of statistical series[3]. He argued that a series consisting of independent observations with differing probabilities of the occurrence of a rare event might be considered as a sample from a single totality. This fact, or, more precisely, the decrease of the pertinent coefficient of dispersion to unity with the decrease of the number of observations he had indeed called the law of small numbers.

From the very beginning his publication aroused debates (Sheynin 1990c/2011, pp. 59 – 62). **Chuprov** (Bortkevich & Chuprov 2005, Letter 2 of 1896) advised Bortkiewicz to refer to **Poisson**, but in 1909 – 1911, in his letters to Chuprov, Bortkiewicz stressed the distinction between his law and the Poisson formula. The low value of probability, as he argued, was not his main assumption; the rarity of the event might have been occasioned by a small number of observations. Incidentally, this explanation raises doubts about the applicability here of the Poisson distribution, and, for that matter, Bortkiewicz had never comprehensively explained his law. Here is what Chuprov wrote to **Markov** in 1916 (Sheynin 1990c/2011, pp. 91 – 92):

*It is difficult to say to what extent the law of small numbers enjoys the recognition of statisticians since it is not known what, strictly speaking, should be called the law of small numbers. Bortkiewicz did not answer my questions formulated in the note on p. 398 of the second edition of the <u>Essays</u>* [Chuprov 1909/1959, p. 285] *either in publications or in written form; I did not question him orally at all since he regards criticisms of the law of sm. numb. very painfully.*

Mathematicians now simply dismiss the law of small numbers as another term for the **Poisson** limit theorem (**Kolmogorov** 1954, without explaining the situation). However, the first to deny that law was L. Whittaker (1914) who showed her standpoint in the very title of her paper. **Bortkiewicz** sharply objected to it in his polemic paper (1915a).

**Markov** repeatedly discussed that law in his letters of 1916 to **Chuprov** (Ondar 1977/1981); he indicated that Bortkiewicz had wrongly combined his data and (p. 108) "chose material that was pleasing to him"[4] and that (pp. 81 and 108) for small numbers the coefficient of dispersion cannot be large. He also publicly repeated his last-mentioned statement (1916b, p. 55). In 1916, in an answer to Markov, Chuprov (Sheynin 1990c/2011, p. 91) apparently disagreed that Bortkiewicz had wrongly combined his materials and reported that Yastremsky (1913) had also proved Markov's main statement. Then, Quine & Seneta (1987) described Bortkiewicz' law and indicated more definitely that for small independent and integral random variables a large value of $Q$ was unlikely.

I ought to add that after 60 years of its neglect Bortkiewicz was the author who picked up **Poisson's** law (also Newcomb, see § 10.8.4, and Chebyshev) and that for a long time his contribution (1898a) had remained the talk of the town. Thus, **Romanovsky** (1924, book 17, p. 15) called Bortkiewicz' innovation "the main statistical law".

Bortkiewicz invariably defended his law but it is difficult to study it thoroughly; commentators whom I mentioned above did not examine it comprehensively (and neither did other authors), and I refer readers to my paper (2008). I have minutely studied his contribution and decided that he did nothing except (as Kolmogorov stated) recalling Poisson.

No one ever noted that the Bortkiewicz coefficient $Q_1$ differed from Q (so that some comments had been wrong) and that it was the ratio of two dependent random variables, call tem   and   : $Q_1 = E /E$ , and, consequently, no one noted that Bortkiewicz wrongly assumed that E /E  = E( / ). Chuprov was an exception, but he said nothing in public whereas Bortkiewicz arbitrarily answered him that that wrong equality holds approximately.

Another point here is that the decrease of the coefficient of dispersion with the decrease of the number of observations tells us nothing about the underlying probability (probabilities) of the studied event. Chuprov (1909/1959, p. 277) noted this fact but only mentioned Lexis, not Bortkiewicz, and did not repeat his remark in any later Russian or German paper.

**15.1.3. Markov and Chuprov.** In his letters of 1910 to Chuprov, Markov (Ondar 1977) proved that **Lexis'** considerations were wrong. Thus, it occurred that the dispersion could also be normal when the observations were dependent. In addition, he constructed an example of independent observations which, when being combined into series in different ways, were characterized either by super- or subnormal dispersions. However, later Chuprov, in a letter of 1923 to his former student Chetverikov (Sheynin 1990c/2011, p. 139), remarked that stability was only determined for concrete series.

Also in 1910, Chuprov, in a letter to Markov, provided examples of dependences leading to super- and sub-normality of dispersion; in 1914 he even decided that the coefficient of dispersion should be "shelved" to which **Bortkiewicz** strongly objected (Sheynin 1990c/2011, p. 140). Then, in 1916 both Markov and Chuprov proved that $EQ^2 = 1$ (see details Ibidem, pp. 140 – 141). Finally, Chuprov (1918 – 1919; see Ibidem, pp. 113 – 114) definitively refuted the applicability of the coefficient of dispersion.

I cannot understand why, even in 1921 Chuprov (2009, p. 88) wrote to Gulkevich (a noted Russian diplomatist who refused to return to the

Soviet Union and became the assistant of F. Nansen at the League of Nations):

*One of the most important doctrines of the theory of statistics, which I had until now completely acknowledged and professed, the Lexian theory of the stability of statistical numbers, is to a large extent, as it occurs, based on a mathematical misunderstanding. This knocks out one of the foundations of the theory whose central part is now hanging on air. I do not want to resign myself to this fact without providing a replacement, but I cannot manage. My attempts encounter that same objection and I have almost concluded that the obstacle is essentially insurmountable. […] That's how it is turning out for the time being.*

All this was almost forgotten. Thus, Särndal (1971, pp. 376 – 377), who briefly described the work of **Lexis** and noted that it prompted **Charlier** "to look […] into questions of non-normality of data", did not mention either it or any later developments.

Here is what Anderson (1926/1963, Bd. 1, p. 31) stated, and what significantly complements his own opinion cited in § 15.1.2: Chuprov

*Essentially changed the traditional Lexis' doctrine […] from which actually pretty little is now left.*

In the same contribution Chuprov (1918 – 1919, p. 205) proved, in a most elementary way, a general formula for the variance:

$$(1/n)E\left(\sum_{i=1}^{n}(x_i - \sum_{i=1}^{n}E\xi_i)^2\right) =$$
$$(1/n^2)\sum_{i=1}^{n} E(\xi_i - E\xi_i)^2 + (1/n^2)\sum_{i=1}^{n}\sum_{j\neq i} [E(x_i x_j) - E\xi_i E\xi_j].$$

Included here were $n$ random variables $\xi_i$ anyhow dependent on each other and the results of a single observation $x_i$ of each of them. I note that Chuprov partly issued from his manuscript (1916 or 1917) There, he determined $EQ^2$ anew and provided qualitative considerations concerning the distribution of the coefficient of dispersion. Chuprov sent his manuscript to Markov, and it is mentioned or implied in their correspondence of 1917 (Sheynin 1990c/2011, pp. 94 – 96) and, possibly, Ondar 1977/1981, p. 116ff).

**Romanovsky** (1923) published a very favourable review of Chuprov's work but did not indicate that the latter's notation was too involved and impeded understanding. A notorious case is presented in his paper (Chuprov 1923). There, for example on p. 472, he applied two-storeyed superscripts and two-storeyed subscripts in the same formula. Hardly has any other author (certainly including Bortkiewicz) allowed himself to take such liberties, to expect his readers to understand suchlike gibberish.

A few additional words are here in order. While studying the stability of statistical series, Chuprov achieved really interesting results, see Seneta (1987). On the other hand, since he considered problems of the most general nature, he inevitably derived awkward formulas, and the same Romanovsky (1930, p. 216) noted that Chuprov's formulas, although "being of considerable theoretical interest", were "almost useless" due to complicated calculations involved. And, on the next page: the estimation of the empirical coefficient of correlation for samples from arbitrary populations was possible almost exclusively by means of Chuprov's formulas, which however were "extremely unwieldy […], incomplete and hardly investigated".

### 15.2. The Biometric School

That name itself recalls the periodical *Biometrika* whose first issue appeared in 1902 with a subtitle *Journal for the Statistical Study of Biological Problems*. Its first editors were Weldon (a widely educated biologist who died in 1906), **Pearson** and Davenport[5] "in consultation" with **Galton**. The editorial there contained the following passage:

*The problem of evolution is a problem in statistics […] [***Darwin*** established] *the theory of descent without mathematical conceptions*[6] *[…] [but] every idea of Darwin – variation, natural selection […] – seems at once to fit itself to mathematical definition and to demand statistical analysis. […] The biologist, the mathematician and the statistician have hitherto had widely differentiated fields of work. […] The day will come […] when we shall find mathematicians who are competent biologists, and biologists who are competent mathematicians …*

Much later Pearson (1923, p. 23) once more mentioned Darwin:

*We looked upon Darwin as our deliverer, the man who had given new meaning to our life and the world we inhabited.*

Here is a passage from a note which **Pearson** had compiled (and apparently sent around) in 1920 and which his son, **E. S. Pearson** (1936 – 1937, vol. 29, p. 164), quoted: The aim of the Biometric school was

*To make statistics a branch of applied mathematics, […] to extend, discard or justify the meagre processes of the older school of political and social statisticians, and, in general, to convert statistics in this country from being the playing field of dilettanti and controversialists into a serious branch of Science. […] Inadequate and even erroneous processes in medicine, in anthropology* [anthropometry]*, in craniometry, in psychology, in criminology, in biology, in sociology, had to be criticized […] with the aim of providing those sciences with a new and stronger technique.*

Note that almost all the disciplines mentioned above were included in Pearson's main field of interests and that he had not found a single kind word for Continental statisticians.

Here is a peculiar passage (Pearson 1907, p. 613):

*I have learned from experience with biologists, craniologists, meteorologists, and medical men (who now occasionally visit the biometricians by night!) that the first introduction of modern statistical methods into an old science by the layman is met with characteristic scorn; but I have lived to see many of them tacitly adopting the very processes they began by condemning.*

The immediate cause for establishing *Biometrika* seems to have been scientific friction and personal disagreement between Pearson and Weldon on the one hand, and biologists, especially Bateson, on the other hand, who exactly at that time had discovered the unnoticed Mendel. It was very difficult to correlate Mendelism and biometry: the former studied discrete magnitudes while the latter investigated continuous quantitative variations.

*Between the Biometric and Mendelian schools there was nothing fundamentally incompatible; in pursuit of the same objective they followed different lines of approach which were essentially complementary rather than antagonistic*. E. S. Pearson (1936, p. 227).

On pp. 169 – 170 E. S. P. noted that the study of social problems by the Mendelian approach had been superficial.

The rapid success of the new school was certainly caused by the hard work of its creators, but also by the efforts of their predecessor, **Edgeworth**. **Chuprov** (1909/1959, p. 27 – 28) provided his correct characteristic. A talented statistician (and economist), he was excessively original and had an odd style; he was therefore unable to influence strongly his contemporaries. However, in his native country he at least paved the way for the perception of mathematical-statistical ideas and methods. His works have appeared recently in three volumes (1996). See also **Schumpeter** (1954/1955, p. 831) and **M. G. Kendall** (1968/1970, pp. 262 – 263). I have quoted them as well as Chuprov elsewhere (Sheynin 2006/2009, NNo. 537 – 539).

Pearson, perhaps at once, became the main editor of *Biometrika,* and among his authors were Chuprov and **Romanovsky**[7].

A list of more than 600 of Pearson's publications is in Morant et al (1939) and Merrington et al (1983), and his son, Egon Pearson (1936 – 1937) described his life and work. Many of his earlier papers are reprinted (K. Pearson 1948) and his manuscripts are held in University College London.

Hardly known are his thoughts about physics which he studied until 1893. Thus (1891, p. 313) *negative matter* exists *in the universe*, and physical *variation effects* were *perhaps due to the general construction of our space*, see Clifford (1885/1886, p. 202). He did not however mention Riemannian spaces whereas it is nowadays thought that the curvature of space–time is caused by forces operating in it. Lastly (Pearson 1887, p. 114), *all atoms in the universe* of *whatever kind appear to have begun pulsating at the same instant*.

Unexpectedly, at least for me, Tee (2003, 15 May) noted that Pearson's *Grammar of Science* (1892) denied the existence of atoms and only repented in 1911, in the third edition of that book. Tee (from University of Auckland) is an extremely careful reviewer; I had not checked his statement which nevertheless seems at least partly doubtful (*all atoms* […] *have begun* etc., see a bit above).

Pearson (1857 – 1936) was an applied mathematician and a philosopher, but in the first place a co-founder of biometry. The beginning of his scientific work can be connected with his *Grammar of science* which earned him the brand of a

*Conscientious and honest enemy of materialism* and *one of the most consistent and lucid* **Machians**.

That was **Lenin's** conclusion (1909/1961, pp. 190 and 274). Note that the latter term is tantamount to Mach's philosophy, i.e., to a variety of subjective idealism. It is, however, difficult to imagine that Pearson evaded reality. But at the same time Mach's followers define the aim of science as description rather than study of phenomena and Pearson separated experience (statistical data) from theory (from the appropriate stochastic patterns), although he did not at all keep to the tradition of the *Staatswissenschaft* (§ 6.2.1).

I note that, in turn **Pearson** (1978, p. 243) had mentioned **Lenin**: Petersburg

*Has now for some inscrutable reason been given the name of the man who has practically ruined it*.

Pearson's *Grammar* became widely known[8] and he was elected to the Royal Society. **Newcomb**, as President of the forthcoming International Congress of Arts and Sciences (St. Louis, 1904), invited him to read a report on methodology of science[9]. Such scholars as **Boltzmann** and **Kapteyn** had participated there. Newcomb's attitude towards Pearson was also reflected in one of his pronouncements of 1903, see § 10.8-4.

S. L. Zabell somewhere noticed that **Mach** (1897, Introduction) had mentioned Pearson's *Grammar*…:

*The publication* [of the *Grammar*] *acquainted me with a researcher whose* erkenntnisskritischen [**Kantian**] *ideas on every important issue coincide with my own notions and who knows how to oppose, candidly and courageously, extra-scientific tendencies in science*.

I mention two more facts concerning Pearson. In 1921 – 1933 he had read a special course of lectures at University College and in 1978 his son, E. S. Pearson, published them making use of the extant notes and likely providing the title itself. After **Todhunter** (1865), this contribution was apparently the first considerable work in its field and on its first page the author expressed his regret that he did not study the history of statistics earlier (see one of the Epigraphs to this book).

E. S. Pearson supplied a Preface where he illustrated his father's interest in general history. But in my context it is more important to mention K. Pearson's fundamental biography of **Galton** (1914 – 1930), perhaps the most immense book from all works of such kind, wherever and whenever published.

Pearson also devoted several papers to the history of probability and statistics; I mentioned three of them (§§ 2.2.3, 3.2.3 and 7.1-5) and resolutely disagreed with the main conclusion of the second one. In two more articles Pearson (1920; 1928b) studied the history of correlation and maintained that some authors including **Gauss** could have applied the ideas and methods of correlation theory, but that it would be nevertheless wrong to attribute to them its beginnings. He often successfully attempted to introduce the statistical method and especially the theory of correlation into many branches of science. See however Bernstein's statement about that theory in § 14.1-4.

The work of Pearson and his followers [**Student** (real name, **Gosset**), **Yule** and others] is partly beyond the boundaries of my investigation and I shall only sketch the main directions of Pearson's subsequent (after about 1894) studies, of the person who (**Hald** 1998, p. 651)

*Between 1892 and 1911 […] created his own kingdom of mathematical statistics* [!] *and biometry in which he reigned supremely, defending its ever expanding frontiers against attacks*.

And, finally (Fisher 1937, p. 306), soon after Pearson's death:

K. Pearson's *plea of comparability* [between the methods of moments and maximum likelihood] *is […] only an excuse for falsifying the comparison*.

Egon Pearson kept silent.

Positive opinions had been certainly held as well (Mahalanobis 1936; Eisenhart 1974), see also Newcomb's letter to Pearson in § 10.8-4 and the statements of Bernstein and Kolmogorov below.

Pearson's main merits include the development of the principles of the correlation theory and contingency, the introduction of the "Pearsonian curves" for describing empirical distributions, rather than for replacing the normal law by another universal density, which was what **Newcomb** (§ 10.8.4) had attempted to accomplish, and the $^2$ test as well as the compilation of numerous statistical tables. Pearson (1896 with additions in 1901 and 1916) constructed the system of his curves in accordance with practical considerations but had not sufficiently justified it by appropriate stochastic patterns. That system was defined as the solution of the differential equation

$$y = \frac{x-k}{a+bx+cx^2} y \qquad (4)$$

with four parameters. The case $b = c = 0$ naturally led to the normal distribution; otherwise, 12 types of curves appeared of which at least some were practically useful. Pearson determined the parameters by the method of moments, – by four sample moments of the appropriate distribution[10]. Recall (§ 10.2-4) that the same term is used in the theory of probability in a quite another way. The "statistical" method of moments is opportune, but the estimates thus calculated often have an asymptotic efficiency much less than unity (**Cramér** 1946, § 33.1).

    **Bernstein** (1946, pp. 448 – 457) indicated a stochastic pattern (sampling with balls being added) that led to the most important Pearsonian curves. He referred to **Markov** and, on p. 337, to **Polya** (1931) as his predecessors. Markov (1917) had indeed considered the abovementioned pattern and mentioned the Pearsonian curves on his very first page; Bernstein, however, mistakenly indicated another of his papers.

    **Pearson** (**E. S. Pearson** 1936 – 1937, vol. 29, p. 208 with reference to his record of the lectures of his father) paid special attention to the notion of correlation and stated that

*The purpose of the mathematical theory of statistics is to deal with the relationship between 2 or more variable quantities, without assuming that one is a single-valued mathematical function of the rest.*

    **Abbe** (§ 9B-1) and then **Helmert** derived the $\chi^2$ distribution for revealing systematic influences in the theory of errors (§ 9B-1) whereas Pearson (1900) introduced the chi-squared test in the context of mathematical statistics. True, not at once, he began applying it for checking the goodness of fit, independence in contingency tables and homogeneity. In spite of its importance, the chi-squared test hardly "clears" empiricism of its dangers as **Fisher**, in 1922, claimed (**Hald** 1998, p. 714).

    I (§ 14.1-4) noted that the Continental statisticians were not recognizing **Pearson**, also see Ondar (1977/1981, p. 142) who quoted **Chuprov's** similar statement. Many of his colleagues, Chuprov wrote, "like **Markov**, shelve the English investigations without reading them". The cause of that attitude was the empiricism of the Biometric school (Chuprov 1918 – 1919, Bd. 2, pp. 132 – 133):

*The disinclination of English researchers for the concepts of mathematical probability and mathematical expectation caused much trouble. The refusal to use these basic notions obscured the stochastic statement of problems; on occasion, it even directed the attempts to solve them on a wrong track. If, however, this attire, so uninviting to the Continental eye, is shed and the discarded is picked up, then it will be distinctly seen that Pearson and Lexis often offer different in form but basically kindred methods for solving essentially similar problems.*

    *It seems that nowadays one of the most important problems of the flourishing stochastic theory of statistics is exactly to work out a synthesis eliminating the contradictions between the two currents which are directed to common goals. Not <u>Lexis against Pearson</u>, but <u>Pearson cleansed by Lexis and Lexis enriched by Pearson</u> should be the slogan of those who are not satisfied by the soulless empiricism of*

*the post-Queteletian statistics and strive for constructing its rational theory.*

**Fisher** (1922, pp. 311 and 329n) also indicated that Pearson had been confusing theoretical and empirical indicators. Similar pronouncements were due to **Anderson** (Sheynin 1990c/2011, p. 149), but I shall only quote **Bernstein** and **Kolmogorov**. Bernstein (1928/1964, p. 228), when discussing *a new cycle of problems in the theory of probability which comprises the theories of distribution and of the general non-normal correlation*, wrote:

*From the practical viewpoint the Pearsonian English school is occupying the most considerable place in this field. Pearson fulfilled an enormous work in managing statistics; he also has great theoretical merits, especially since he introduced a large number of new concepts and opened up practically important paths of scientific research. The justification and criticism of his ideas is one of the central problems of current mathematical statistics. Charlier and Chuprov, for example, achieved considerable success here whereas many other statisticians are continuing Pearson's practical work, definitely losing touch with probability theory ...*

*The modern period in the development of mathematical statistics began with the fundamental works of […] K. Pearson, Student, Fisher. ... Only in the contributions of the British school did the application of probability theory to statistics cease to be a collection of separate isolated problems and became a general theory of statistical testing of stochastic hypotheses …* Kolmogorov (1947, p. 63).

*The investigations made by Fisher, the founder of the modern British mathematical statistics, were not irreproachable from the standpoint of logic. The ensuing vagueness in his concepts was so considerable, that their just criticism led many scientists (in the Soviet Union, Bernstein) to deny entirely the very direction of his research.* Ibidem, p. 64.

Enumerating the "main weaknesses" of the **Pearsonian** school, Kolmogorov (1948/2002, pp. 68 – 69) indicated that

*Rigorous results concerning the proximity of empirical sample characteristics to theoretical related only to the case of independent trials. Notions held by the English statistical school about the logical structure of the theory of probability which underlies all the methods of mathematical statistics remained on the level of the eighteenth century. In spite of the immense […] work done […] the auxiliary tables used in statistical studies proved highly imperfect in respect of cases intermediate between "small" and "large" samples.*

On the statements of other scientists about Pearson see Sheynin (2010a). On Yule see M. G. Kendall (1951) and on Gosset (Student)

see E. S. Pearson (1990). His collected papers appeared in 1992. And just a few words about him from Irwin (1978, p. 409):

*He was one of the pioneers in the development of modern statistical method and its application to the design and analysis of experiments.*

Irwin (p. 410) also quoted Fisher: Gosset was *the Faraday of statistics*. Note that the German statistician Lüroth was the first to introduce the *t*-distribution (Pfanzagl & Sheynin 1996).

### 15.3. The Merging of the Continental Direction and the Biometric School?

So, did the two statistical streams merge, as **Chuprov** would have it? In 1923 he had become Honorary Fellow of the Royal Statistical Society and in 1926, after his death, the Society passed a resolution of condolence (Sheynin 1990c/2011, p. 156) which stated that his

*Contributions to science were admired by all* […]. *They did much to harmonise the methods of statistical research developed by continental and British workers.*

In § 14.1-4 I mentioned the unilateral and, for that matter, only partly successful attempts made by Chuprov, and the vain efforts of **Slutsky** to reconcile **Markov** with **Pearson's** works. And Bauer (1955, p. 26) reported that he had investigated, on Anderson's initiative, how both schools had been applying analysis of variance and concluded (p. 40) that their work was going on side by side but did not tend to unification. More details about Bauer`s study are contained in Heyde & Seneta (1977, pp. 57 – 58) where it also correctly indicated that, unlike the Biometric school, the Continental direction had concentrated on nonparametric statistics.

I myself (**Gnedenko** & Sheynin 1978/2001, p. 275) suggested that mathematical statistics properly originated as the coming together of the two streams; even now I think that that statement was not original (but am unable to mention anyone). However, now I correct myself. At least until the 1920s, say, British statisticians had continued to work all by themselves. **E. S. Pearson** (1936 – 1937), in his study of the work of his father, had not commented on Continental statisticians and the same is true about other such essays (**Mahalanobis** 1936; Eisenhart 1974). We only know that **K. Pearson** regretted his previous neglect of the history of statistics (see one of my Epigraphs).

I believe that English, and then American statisticians for the most part only accidentally discovered the findings already made by the Continental school. Furthermore, the same seems to happen nowadays as well. Even **Hald** (1998) called his book *History of Mathematical Statistics*, but barely studied the work of that school. In 2001, *Biometrika* (vol. 88) published five essays devoted to its centenary but not a word was said in any of them about the Continent, not once was **Chuprov** mentioned. It is opportune to add that **Cramér** (1946, Preface) aimed to unite, in his monograph, English and American statistical investigations (and, in the first place, the work of **Fisher**) with the new, purely mathematical theory of probability created "largely owing to the work of" French and Russian mathematicians.

Later he (1981, p. 315) stated that

*In the years before the war* [before WWI] *it had seemed to me that continental mathematicians and the Anglo-Saxon statisticians were working without sufficient mutual contact, and that it might be useful to try to join both these lines of research.*

In 1919 there appeared in *Biometrika* an editorial remarkably entitled *Peccavimus*! (we were guilty). Its author, undoubtedly **Pearson**, corrected his mathematical and methodological mistakes made during several years and revealed mostly by **Chuprov** (Sheynin 1990c/2011, p. 75) but he had not taken the occasion to come closer to the Continental statisticians. Continental statisticians had not been better, suffice it to mention Markov.

### Notes

**1. Bortkiewicz'** paper appeared in a rare Russian political periodical published abroad. I discovered that journal in the Rare books section of the (former) Lenin State Library in Moscow. A few other copies of the periodical's same issue, which I since found in Germany, do not, however, contain the paper in question; perhaps it was only included in a part of the edition.

**2. Oskar Nikolaevich Anderson** (1887 – 1960), a Russian German, emigrated in 1920. In 1924 – 1942 he lived and worked in Bulgaria, then in Germany (in West Germany), was the leading statistician in both these countries and a founder-member of the international *Econometric Society* (Anderson 1946). Also see Sheynin (1990c/2011, pp. 80 – 82), H. & R. Strecker (2001) and his collected works (Anderson 1963). Anderson was the last representative of the Continental direction.

**3.** In 1897 **Bortkiewicz** also unsuccessfully attempted to publish his work in Russian, in a periodical of the Petersburg Academy of Sciences. His request was refused since the contribution was to appear elsewhere (although only in German), see Sheynin (1990c/2011, p. 61).

**4.** This charge was not proved; furthermore, it contradicts our perception of his personality.

**5.** An author of a paper published in 1896, of a book devoted to biometry which appeared in 1899, and of two subsequent notes (**M. G. Kendall** & Doig 1968).

**6.** Already in **Darwin's** times, a theory was supposed to be quantitatively corroborated. Darwin, however, provided nothing of the sort and it would be more proper to say, as in § 10.8.2, "hypothesis of the origin of species".

**7.** In 1912, **Slutsky** had submitted two manuscripts to **Pearson** who rejected both. Three letters from Slutsky to Pearson (but no replies) are extant (Univ. College London, Pearson Papers 856/4 and 856/7; Sheynin 1999c). In this connection Slutsky had corresponded with **Chuprov** (Sheynin 1990c/2011, pp. 65 – 67) and soon published elsewhere (1914) one of his manuscripts, – the one, whose refusal by Pearson he called a misunderstanding. It was Chuprov who recommended him a proper outlet.

**8.** Here is **Neyman's** remarkable recollection (**E. S. Pearson** 1936 – 1937, vol. 28, p. 213): in 1916, he read the *Grammar of Science* on advice of his teacher at Kharkov University, **S. N. Bernstein**, and the book greatly impressed "us".

It is not difficult to imagine that **Pearson** was given a hostile reception in the Soviet Union. This issue is beyond my chronological boundaries and I only mention two episodes (Sheynin 1998c, pp. 536 and 538, note 16).

a) Maria Smit, the future Corresponding Member of the Academy of Sciences (!), 1930: the Pearsonian curves are based

*On a fetishism of numbers, their classification is only mathematical.*

*Although he does not want to subdue the real world as ferociously as it was attempted by […] Gaus [**Gauss**], his system nevertheless rests only on a mathematical foundation and the real world cannot be studied on this basis at all.*

b) A. Ya. Boiarsky, L. Zyrlin, 1947: they blasphemously charged **Pearson** with advocating racist ideas that "forestalled the Göbbels department".

Only somewhat more reserved was the anonymous author in the *Great Sov. Enc.*, 2$^{nd}$ ed., vol. 33, 1955, p. 85.

Soviet statistics and statisticians endured real suffering. The same Smit, in 1931 (Sheynin 1998c, p. 533, literal translation): "the crowds of arrested saboteurs are full of statisticians". She herself had probably helped to assure the success of that process.

And let us recall Schlözer (1804, § 15): *Statistics and despotism do not get along together.*

**9. Pearson** refused to come because of his financial problems and unwillingness to leave his Department at University College London under "less complete supervision" (Sheynin 2002b, pp. 143 and 163, Note 8).

**10.** Here is a passage from the extant part of an unsigned and undated letter certainly written by **Slutsky** to **Markov**, likely in 1912 (Sheynin 1999c, p. 132):

*…are not independent in magnitude from the sum of the already accumulated deviations or that the probabilities of equal deviations are not constant, we shall indeed arrive at the formula* [Slutsky wrote down formula (4) with $k = 0$ and $F(x)$ instead of the trinomial in the denominator]. […] *Much material* [already shows that the **Pearsonian** curves are useful but] […] *it seems desirable also for the asymmetric Pearson curves* […] *to provide a theoretical derivation which would put* [them] *in the same line as the **Gauss** curve on the basis of the theory of probability* (*hypergeometric series*).

## Supplement: Axiomatization

I present a bibliographic survey of some important points.

The main essays are Barone & Novikoff (1978) and Hochkirchen (1999) and among the lesser known authors is Bernstein (1917). After Hilbert (1901), Kolmogorov (1933) made the decisive step and Freudenthal & Steiner (1966, p. 190) commented: he *came with the Columbus' egg*. As the legend goes, Columbus cracked an egg which enabled it to stand firmly on his table. Among the new sources I list Hausdorff (2006) who left an important unpublished contribution, see Girlich (1996), Shafer & Vovk (2001) and Krengel (2011) who stressed the role of Bohlmann. Vovk & Shafer (2003, p. 27) characterized their book:

*We show how the classical core of probability theory can be based directly on game-theoretic martingales, with no appeal to measure theory. Probability again becomes* [a] *secondary concept but is now defined in terms of martingales*.

In concluding, I quote Boole (1854/1952, p. 288):

*The claim to rank among the pure sciences must rest upon the degree in which it* [the theory of probability] *satisfies the following conditions: 1° That the principles upon which its methods are founded should be of an axiomatic nature*.

Boole formulated two more conditions of a general scientific essence.

2012a (Gauss)
2016a (essence of statistics)

# 3. General Literature

**Daw R. H.** (1980), J. H. Lambert, 1727 – 1777. *J. Inst. Actuaries*, vol. 107, pp. 345 – 350.
**Dawid Ph.** (2005), Statistics on trial. *Significance*, vol. 2, No. 1, pp. 6 – 8.
**Daxecker F.** (2004), *The Physicist and Astronomer Christopher Scheiner*. Innsbruck.
--- (2006), *Die Hauptwerk des Astronomen P. Christoph Scheiner JJ*. Innsbruck.
**DeCandolle Alph.** (1855), *Géographie botanique raisonnée*, tt. 1 – 2. Paris.
**DeCandolle Aug. P.** (1832), *Physiologie végétale*, tt. 1 – 3. Paris.
**Delambre J. B. J.** (1819), Analyse des travaux de l'Académie … pendant l'année 1817, partie math. *Mém. Acad. Roy. Sci. Inst. de France*, t. 2 pour 1817, pp. I – LXXII of the *Histoire*.
*Demografichesky* (1985), *Demografichesky Enziklopedichesky Slovar* (Demographic Enc. Dict.). Moscow.
**De Moivre A.** (1712, in Latin). English transl.: De mensura sortis or the measurement of chance. ISR, vol. 52, 1984, pp. 236 – 262. Commentary (A. Hald): Ibidem, pp. 229 – 236.
--- (1718), *Doctrine of Chances*. Later editions: 1738, 1756. References in text to reprint of last edition: New York, 1967.
--- (1725), *Treatise on Annuities on lives*. London. Later edition of 1743 incorporated in the *Doctrine* (1756, pp. 261 – 328) in a somewhat improved form, as its anonymous editor stated on p. xi. German transl.: Wien, 1906.
--- (1730), *Miscellanea Analytica de Seriebus et Quadraturis*. London. French translation: Paris, 2009.
--- (1733, Latin), Transl. by author: A method of approximating the sum of the terms of the binomial $(a + b)^n$ expanded into a series from whence are deduced some practical rules to estimate the degree of assent which is to be given to experiments. Incorporated in subsequent editions of the *Doctrine* (in 1756, in an extended version, on pp. 243 – 254).
--- (1756), This being the last edition of the *Doctrine*.
**De Montessus R.** (1903), Un paradoxe du calcul des probabilités. *Nouv. Annales Math.*, sér. 4, t. 3, pp. 21 – 31.
**De Morgan A.** (1845), Theory of probabilities. *Enc. Metropolitana*, Pure sciences, vol. 2. London, pp. 393 – 490.
--- (1847), *Formal Logic, or the Calculus of Inference, Necessary and Probable*. Second ed.: Chicago, 1926.
--- (1864), On the theory of errors of observation. *Trans. Cambr. Phil. Soc.*, vol. 10, pp. 409 – 427.
**De Morgan Sophia Elizabeth** (1882), *Memoir of Augustus De Morgan*. London.
**Descartes R.** (1644), *Les principes de la philosophie. Œuvres,* t. 9, pt. 2 (the whole issue). Paris, 1978. Reprint of the edition of 1647.
**DeVries W. F. M.** (2001), Meaningful measures: indicators on progress, progress on indicators. ISR, vol. 69, pp. 313 – 331.
**De Witt J.** (1671, in vernacular), Value of life annuities in proportion to redeemable annuities. In Hendriks (1852, pp. 232 – 249).
**Dietz K.** (1988), The first epidemic model: historical note on P. D. Enko. *Austr. J. Stat.*, vol. 30A, pp. 56 – 65.
**Dietz K., Heesterbeek J. A. P.** (2000), D. Bernoulli was ahead of modern epidemiology. *Nature*, vol. 408, pp. 513 – 514.
**---** (2002), D. Bernoulli's epidemiological model revisited. *Math. Biosciences*, vol. 180, pp. 1 – 21.
**O'Donnell T.** (1936), *History of Life Insurance*. Chicago.
**Doob J. L.** (1989), Commentary on probability. In *Centenary of Math. in America*, pt. 2. Providence, Rhode Island, 1989, pp. 353 – 354. Editors P. Duren et al.
**Dorfman Ya. G.** (1974), *Vsemirnaia Istoria Fiziki* (Intern. History of Physics). Moscow.
**Dormoy E.** (1874), Théorie mathématique des assurances sur la vie. *J. des actuaries française*, t. 3, pp. 283 – 299, 432 – 461.
--- (1878), Same title, t. 1. Paris. Incorporates his paper of 1874.
**Dorsey N. E., Eisenhart C.** (1969), On absolute measurement. In Ku (1969, pp. 49 – 55).
275

--- (1823b, Latin), Theorie der den kleinsten Fehlern unterworfenen Combination der Beobachtungen, pts. 1 – 2. Ibidem, pp. 1 – 53.

--- (1823c), Anwendung der Wahrscheinlichkeitsrechnung auf eine Aufgabe der practischen Geometrie. W-9. Göttingen – Leipzig, 1903, pp. 231 – 237.

--- (1826, German), Preliminary author's report about Gauss (1828). Ibidem, pp. 200 – 204.

--- (1828a, Latin), Supplement to Gauss (1823b). German transl.: Ibidem, pp. 54 – 91.

--- (1828b), Bestimmung des Breitenunterschiedes zwischen den Sternwarten von Göttingen und Altona etc. W-9, pp. 5 – 64. **S, G,** 72.

--- (1845; Nachlass), Anwendung der Wahrscheinlichkeitsrechnung auf die Bestimmung der Bilanz für Witwenkassen. W-4, 1873, pp. 119 – 183.

--- (1855), *Méthode des moindres carrés.* Paris.

--- (1863 – 1930), *Werke*, Bde 1 – 12. Göttingen a. o. Reprint: Hildesheim, 1973 – 1981.

--- (1887), *Abhandlungen zur Methode der kleinsten Quadrate*. Hrsg, A. Börsch & P. Simon. Latest ed.: Vaduz, 1998.

**Gavarret J.** (1840), *Principes généraux de statistique médicale.* Paris.

**Gerardy T.** (1977), Die Anfänge von Gauss' geodätische Tätigkeit. *Z. f. Vermessungswesen*, Bd. 102, pp. 1 – 20.

**Gingerich O.** (1983), Ptolemy, Copernicus, and Kepler. In Adler M. G. & van Doren J., ditors (1983), *Great Ideas Today.* Chicago, pp. 137 – 180.

--- (2002), The trouble with Ptolemy. *Isis*, vol. 93, pp. 70 – 74.

**Gini C.** (1946), Gedanken von Theorem von Bernoulli. *Z. für Volkswirtschaft u. Statistik*, 82. Jg., pp. 401 – 413.

**Gnedenko B. V.** (1951, Russian), On Ostrogradsky's work in the theory of probability. IMI, vol. 4, pp. 99 – 123. **S, G,** 5.

--- (1954, Russian). *Theory of probability.* Moscow, 1973. [Providence, RI, 2005.] First Russian ed., 1950.

--- (1958), Main stages in the history of the theory of probability. *Actes VIIIe Congrès Hist. Sci. 1956.* N. p., 1958, vol. 1, pp. 128 – 131.

--- (1959, Russian), On Liapunov's work on the theory of probability. IMI, vol. 12, pp. 135 – 160.

**Gnedenko B. V., Gikhman I. I.** (1956, Russian), Development of the theory of probability in the Ukraine. IMI, vol. 9, pp. 477 – 536.

**Gnedenko B. V., Sheynin O. B.** (1978). See Sheynin (1978a).

**Goldstein B. R.** (1985), *The 'Astronomy' of Levi ben Gerson (1288 – 1344)*. New York.

**Golitzin Iv.** (1807), *Statisticheskie Tablitzy Rossiyskoi Imperii* (Stat. Tables of the Russian Empire). Moscow.

**Gordevsky D. Z.** (1955, Russian), *K. A. Andreev.* Kharkov.

**Gorrochurn P.** (2016), *Classic Topics on the History of Modern Mathematical Statistics*. Hoboken NJ.

**Gower B.** (1993), Boscovich on probabilistic reasoning and the combination of observations. In Bursill-Hall (1993, pp. 263 – 279).

**Gowing R.** (1983), *Roger Cotes – Natural Philosopher.* Cambridge.

**Graunt J.** (1662), *Natural and Political Observations Made upon the Bills of Mortality*. Baltimore, 1939. Editor, W. F. Willcox.

*Great Books* (1952), *Great Books of the Western World*, vols 1 – 54. Chicago.

**Greenwood M.** (1936), Louis and the numerical method. In author's *Medical Dictator.* London, pp. 123 – 142.

--- (1940), A statistical mare's nest? *J. Roy. Stat. Soc.*, vol. 103, pp. 246 – 248.

--- (1941 – 1943), Medical statistics from Graunt to Farr. *Biometrika*. Reprint: Pearson & Kendall (1970, pp. 47 – 120).

**Gridgeman N. T.** (1960), Geometric probability and the number . *Scripta Math.*, t. 25, pp. 183 – 195.

**Grodzensky S. Ya.** (1987, Russian), *A. A. Markov*. Moscow.

**Gusak A. A.** (1961, Russian), La préhistoire et les débuts de la théorie de la représentation approximative des fonctions. IMI, vol. 14, pp. 289 – 348.

**Guy W. A.** (1852), Statistics, medical. In *Cyclopaedia of Anatomy and Physiology*, vol. 4. London, pp. 801 – 814.
279

# Index of Names

This Index does not cover the General Literature and my own name is also absent.
The numbers refer to subsections rather than to pages.